\documentclass[fleqn,a4paper]{paper}

\newtheorem{theorem}{Theorem}
\newtheorem{lemma}{Lemma}
\newtheorem{conclusion}{Conclusion}
\newtheorem{auxillary_lemma}{Auxiliary Lemma}
\newtheorem{note}{Note}
\newtheorem{definition}{Definition}
\newtheorem{technical_lemma}{Technical Lemma}
\newtheorem{algorithm}{Algorithm}

\newenvironment{narrow}[2]{%
\begin{list}{}{%
\setlength{\topsep}{0pt}%
\setlength{\leftmargin}{#1}%
\setlength{\rightmargin}{#2}%
\setlength{\listparindent}{\parindent}%
\setlength{\itemindent}{\parindent}%
\setlength{\parsep}{\parskip}}%
\item[]}{\end{list}}

\usepackage[]{graphicx}
\usepackage[]{subfigure}
\usepackage[]{caption2}
\usepackage[section]{placeins}
\usepackage{color}
\usepackage{amssymb}

\begin{document}

\title{Smooth B\'ezier Surfaces over Unstructured Quadrilateral Meshes}

\author{Michel Bercovier$^1$\thanks{E-Mail :\texttt {berco@cs.huji.ac.il};Corresponding author}   and Tanya Matskewich$^2$ \\ $^1$ \small{The Rachel and Selim Benin School of Computer Science and Engineering} ,\\ \small{Hebrew University of Jerusalem,Israel.}   \\$^2$ \small{Microsoft Corp. Redmond, WA, USA,}}
%
%
%
\date{}

\newcommand{\pd}{\partial}
\newcommand{\DT}[2]{\frac{\partial^2{#1}}{\partial{#2}^2}}
\newcommand{\DD}[3]{\frac{\partial^2{#1}}{\partial{#2}\partial{#3}}}
\newcommand{\D}[2]{\frac{\partial{#1}}{\partial{#2}}}
\newcommand{\CNK}[2]{\left({\scriptsize\begin{array}{c}\!{#1}\cr{#2}\!\end{array}}\right)}
\newcommand{\Cnk}[2]{\scriptsize \left(\begin{tabular}{@{}c@{}} {#1}\cr {#2}\end{tabular}\right)}
\newcommand{\cnk}[2]{({\tiny\begin{tabular}{@{}c@{}} {#1}\cr {#2}\end{tabular}})}
\newcommand{\VV}[2]{\left(\begin{array}{c}{#1}\cr{#2}\end{array}\right)}
\newcommand{\MM}[4]{\left(\begin{array}{cc}{#1}&{#2}\cr{#3}&{#4}\end{array}\right)}
\newcommand{\VT}[3]{\left(\!\!\begin{array}{c}{#1}\cr{#2}\cr{#3}\end{array}\!\!\right)}
\newcommand{\imp}{$\clubsuit\clubsuit\clubsuit$}
\newcommand{\T}{\tilde}
\newcommand{\LL}{\tilde\lambda}
\newcommand{\RR}{\tilde\rho}
\newcommand{\GG}{\tilde\gamma}
\newcommand{\TTe}[1]{\tilde{e}^{({#1})}}
\newcommand{\TTL}{\tilde T_{\lambda}}
\newcommand{\TTR}{\tilde T_{\rho}}
\newcommand{\Dlt}{\Delta}
\newcommand{\dL}{\Delta L}
\newcommand{\dR}{\Delta R}
\newcommand{\dC}{\Delta C}
\newcommand{\FUN}[1]{\bar{\cal FUN}^{(#1)}}
\newcommand{\PAR}[1]{\T{\cal PAR}^{(#1)}}
\newcommand{\CP}[1] {\T{\cal CP}^{(#1)}}
\newcommand{\MDS}[1]{\T{\cal B}^{(#1)}} 
\newcommand{\GCP}[1] {\T{\cal CP}_G^{(#1)}}
\newcommand{\GMDS}[1]{\T{\cal B}_G^{(#1)}}
\newcommand{\FCP}[1] {\T{\cal CP}_F^{(#1)}}
\newcommand{\lmtT}[3]{\tiny\begin{array}{c}{#1}\cr{#2}\cr{#3}\end{array}}
\newcommand{\twolines}[2]{\scriptsize\begin{array}{c}\!\!\!{#1}\!\!\!\cr\!\!\!{#2}\!\!\!\end{array}}
\newcommand{\<}{\langle}
\renewcommand{\>}{\rangle}
\renewcommand{\kappa}{\chi}

\newcommand{\eop}{\noindent$\sqcup\!\!\!\!\sqcap$}

\parskip=.050in

%
\renewcommand{\thepage}{\roman{page}}
%
%
%
%

\maketitle
\addcontentsline{toc}{section}{Abstract}

\noindent

\begin{abstract}

We study the following problem: given a polynomial order of approximation $n$ and the corresponding B\'ezier tensor product patches  over an  unstructured quadrilateral mesh made of convex quadrilaterals with vertices of  any valence , is there a solution to the $G{^1}$  ( and as a consequence the $C{^1}$ ) approximation   (resp. interpolation  ) problem ? To illustrate the interpolation case , constraints  defining regularity conditions across patches have to be satisfied. The resulting number of free degrees of freedom must be such that the interpolation problem has a solution! This is  similar to studying the minimal determining set (MDS)  for a  $C{^1}$ continuity construction. 

Based on the equivalence of  $G{^1}$ and $C{^1}$  we introduce a sufficient  $G{^1}$  condition that is better adapted to the present problem. Boundary conditions are then analysed including normal derivative constraints ( common in FEM but not  in CAGD. )

The MDS are constructed for both polygonal meshes and meshes with $G{^1}$ -smooth piecewise B\'ezier cubic global boundary. 
The main results are that such MDS exists always for patches of order $\ge 5 $. For $n = 4$ criterions for  mesh structures avoiding under constrained situations are analysed . This leads to the construction of bases by solving a well defined linear system, which allows the solution of the problem for large families  structures of  planar meshes, \emph{without using macro-elements or subdivisions.}.

A complete solution for cubic $C{^1}$ boundaries is given , again by a constructive algorithm. We also show that one can mixes quartic and quintic patches. constraints. Explicit construction is provided for important types of interpolation/boundary Finally some numerical examples are given . As a conclusion from  a practical point of view, the present paper provides a way to solve  $C^1$ interpolations/approximations and fourth order  partial differential problems on arbitrary structures of quadrilateral meshes. Defining a global in-plane parametrisation a priory allows the introduction of a  "physical" energy functional, as is done in Isogeometric Analysis, energy that relates to the functional representation of the surface.   
\end{abstract}
\eject

\vspace{-0.3in}
\tableofcontents

\eject
\listoffigures
\addcontentsline{toc}{section}{List of Figures}

\eject

\setcounter{page}{1}
\renewcommand{\thepage}{\arabic{page}}

\part{Introduction}
\label{part:introduction}

\section{Problem definition , current sate and aim of the research }
\label{sect:problem_definition}

Definitions of the current Section are not always self-contained. The precise mathematical definitions of the notions used in this section  will be given in subsequent sections,(example: the precise definition of $G^1$ ($C^1$)-continuity ,is in Part ~\ref{part:review_3D_G1}).

\subsection{Description of the problem}
\label{subsect:problem_description}

\textbf {Input}
\emph{
Let $OXYZ$ a Cartesian coordinate system, $\T\Omega$  a planar simply-connected domain lying in  the  $XY$-plane, together with its regular quadrangulation into  convex elements ( by regular we mean that two elements have either an empty intersection, an edge  or a  vertex in common .) The mesh may have an \underline{arbitrary} \underline{topological and geometrical} structure (see Figure ~\ref{fig:fig1}) and may have either a polygonal or a  $G^1$-smooth piecewise-cubic B\'ezier parametric global boundary. }

\textbf{Output}
\emph{
Compute a $3D$ \underline{piecewise B\'ezier tensor-product} \underline{parametric} surface of order $n$, so that 
\begin{itemize}
\item[{\bf (1)}]
There is a \underline{one-to-one correspondence} between the planar mesh elements and the patches of the resulting surface.The orthogonal projection of every $3D$ patch onto the $XY$-plane defines a bijection between the patch and the corresponding element of the planar mesh (see Figure ~\ref{fig:fig2}).
\item[{\bf (2)}]
The union of every two adjacent patches is  \underline{$G^1$ ($C^1$)-regular}.
\item[{\bf (3)}]
The resulting surface is obtained as  the constrained minimum for a given energy functional, relative to the cartesian plane ,(  \underline{functional} representation of the surface.)
\item[{\bf (4)}]
One can impose  \underline{additional interpolation/boundary constraints} . For example, the resulting surface may be required to interpolate some $3D$ initial data at vertices (see Figure ~\ref{fig:fig3}).
\item[{\bf (5)}]
Determine the dimension of the underlying approximation space.
\end{itemize}
}
\label{def:general_problem_def}
An important result due to J. Peters ~\cite{jorg} establishes that  $G^1$  continuity over any unstructured mesh implies actually $C^1$ . Hence we are free to use 
 $G^1$ or $C^1$-continuity requirements , whichever is more adapted to the cases we consider! 
We always suppose that the additional constraints noted above  fit the smoothness requirements . The surface we thus define or compute can be seen as a \emph{functional surface} , that is a surface defined as $Z= f(X,Y)$ over $\T\Omega$

%
\subsection{Brief review of  related works}
\label{subsect:review_general}

The current work lies at the junction of Computer Aided Geometric Design and the Finite Element Methods. Both fields have extremely extensive bibliography that makes it impossible to present a full list of related works. 
 
Reviews and wide lists of references can be found, for example, in ~\cite{elber}, ~\cite{farin}, ~\cite{hoschek}, 
  ~\cite{tiller},  for CAGD methods) and in ~\cite{ciarlet_handbook}, ~\cite{zienkiewicz}, ~\cite{timoshenko_main} (for Finite Element methods).
The recent Isogeometric Analysis (IGA) breakthrough ~\cite{hughes} has brought together the fields of multi patch surface handling  , higher order smooth  order approximations  and Finite Element Methods (FEM)..

The related problems involve many "parameters"  
like : type of the considered meshes, choice of  the functional space, the energy functional, smoothness , hence an abundant literature  around the present subject.

A classic problem in CAGD  deals with the construction of piecewise parametric B\'ezier or B-spline surfaces or their extension to Non Uniform Rational B-Spline (NURBS) ,and  it often requires (at least) $G^1$($C^1$)-continuity of the resulting surface. The researches may be subdivided into two large categories: the first category concentrates on the continuity conditions; the second category studies the different types of the energy (fairness) functionals. 
The purpose of the present review is to outline briefly some fundamental concepts of the related approaches. We review mainly  the case of quadrilateral patches and smooth surfaces, and do not consider NURBS. 
 
\subsubsection{CAGD based constructions}

\paragraph{Continuity conditions}
One can onsider the construction of a smooth surfaces interpolating a given $3D$ mesh of curves (e.g. ~\cite{liu} ~\cite{peters_main}, ~\cite{piper}, ~\cite{sarraga_old}, ~\cite{sarraga_old_errata}) or start with the construction of the mesh of curves which interpolates the given data (e.g. ~\cite{nielson_functional}, ~\cite{pottmann}).
Usually, the initial triangular or quadrilateral mesh is not required to be regular. However, it appears (see ~\cite{peters_main}) that the piecewise parametric $G^1$-smooth interpolant does not necessarily exist for every mesh. Then either some restrictive assumption on the mesh of curves is introduced or some modification of the mesh is made.

Localizing the propagation of continuity constraints by refining surfaces is necessary for some cases.
Subdivision of some mesh elements (see~\cite{FVS}, ~\cite{clough_toucher}, ~\cite{doo_sabin_subdivision}, ) is commonly applied in order to improve the mesh quality (e.g. ~\cite{peters_main}, ~\cite{piper},~\cite{reif_spline}) and to make the mesh admissible for interpolation by a smooth surface. Another techniques ,~\cite{hahman}, is based on macro-elements to keep a low order order approximation .

Subdivision of an initial mesh element clearly implies that the resulting surface for the element is composed of several (polynomial/spline) pieces. The first step is to check that  the mesh of curves is  admissible,next one proceeds with filling the "faces". Both the weight functions and the inner control points in the under-constraint situations are defined by application of some (usually local) heuristics, such as the least-square or averaging techniques (see, for example, ~\cite{hahman}, ~\cite{nielson_functional}, ~\cite{peters_main}, ~\cite{reif_spline},  and the references herein. Application of the local heuristics allows to construct a resulting surface by the local methods and to avoid any complicated computations. 

To avoid macro elements or subdivision one needs higher order tensor product patches for $G{^1}$ construction of surfaces. The first candidate is the bi quartic patch , but as we will show existence and uniqueness of a solution depends on the underlying mesh structure .
$G^{1}$ construction of bi quintic B-spline surfaces over arbitrary topology is done in ~\cite{shi_wang_quint} and ~\cite{xshi_rec}. The aim is  to simplify surface representations by an approximation with such patches. To do that they derive many local $G^{1}$ properties similar to the one we will introduce. The authors do claim that the bi- quintic quadrilateral is of  the lowest order possible, but do not give any demonstration of this statement. They show that the problem has no solution over general meshes for bi-cubic patches  . Furthermore the actual dimension of the resulting basis is not studied and the functional used for approximating the given surfaces are not defined. Similarly,  ~\cite{adaquint}, construct a $G^{1}$ surface by  patch-by-patch scheme  smoothly stitching bi-quintic B\'ezier patches .  While the techniques described above are generally sufficient in order to construct nicely looking surfaces,by approximation or  interpolation of  the given data, they usually require some preprocessing  and  the nature of the local heuristics may not reflect any geometrical characteristic of the resulting surface.  ~\cite{demko} gives a higher order construction based on sixth order polynomials, we shall not consider this here.  

\paragraph{A study of the energy functionals}
The second category of techniques allows controlling the shape of the surface by minimisation of some (usually global) energy functional. The works, which study different forms of the energy functional, usually deal with intrinsically smooth functional bases (e.g. B-spline basis) and in any case assume the regular structure of the mesh. The user is required to enter only some essential interpolation data, the rest of the degrees of freedom are defined by the energy minimisation. 
The energies used in Computer Aided Geometrical Design commonly relate to the parametric representation of the surface and include the partial derivatives with respect to parameters. The spectrum of the energies is very wide; the most advanced techniques compute the energies using some initial approximation of the resulting surface which lead to a good approximation of the "natural" geometrical characteristics, such as total curvature (e.g. (~\cite{greiner_geometric_energy}, ~\cite{sarraga_greiner_extension}).

\subsubsection{Interaction between  FEM and CAGD}
\label{subsubsect:FEM}

CAGD and FEM are related domains, the  main link being the two way exchanges between geometry and meshes . 
Higher order approximations are often used in FEM, based on higher order polynomials 
( $  p  $ \emph{method})
 triangular, regular quadrilateral or so-called macro-elements ( splitting of convex quadrilaterals (see ~\cite{FVS}) or triangles (see ~\cite{clough_toucher})), for details see ~\cite{ciarlet_old}, ~\cite{zienkiewicz}. Moreover CAD representations have been used for the numerical solutions of partial differential equations (PDEs), ~\cite{reifweb} ,~\cite{anath},  and conversely some FEM methods have been used for the design of geometrical objects  ,~\cite{volp}.  The real convergence is given the Isogeometric Analysis ~\cite{hughes} , where PDEs are approximated by NURBS in the \emph {physical space} , using the geometric transformation that defines the domain  and not  the parametric( reference )  one.

\paragraph{The Bivariate Triangular Spline Finite Elements}

The construction of the Bivariate Triangular Spline Finite Elements (BSFE) is closely connected to the approach of the current work. The BSFE approach combines B\'ezier-Bernstein representation of the polynomials and the requirement of $C^r$ ($r>0$) smoothness. It leads to the problem of determining {\it minimal determining sets (MDS)}. This will be at the center of the present work  , so let us introduce this notion as it was for BSFE. 

Let a triangulation of a simply-connected planar polygonal domain $\T\Omega$ be given.  By definition, for integers $n$ and $0\le r\le n-1$, {\it space $S^r_n$} consists of $C^r(\T\Omega)$ smooth functions which are piecewise polynomials of total degree at most $n$ over each triangle with respect to the barycentric coordinates. 
 $S^r_n$ is called a bivariate spline finite element space. Note, that although for every triangle the polynomials are represented in their B\'ezier-Bernstein form with respect to the barycentric coordinates, they also can be considered as polynomials in the functional sense. 

For a triangle with vertices $A$, $B$, $C$, the  $XY$-coordinates of the B\'ezier control points are given by $\bar{P}_{i,j,k}=\frac{1}{n}(i A + j B + k C)$, $i+j+k=n$. Let $Z(\bar{P}_{i,j,k})$ denote the Z coordinate of the control point $ \T P_{i,j,k}$. Since at least $C^0$-continuity is assumed, the B\'ezier control points of shared edges are unique, determining the dimension of the space $S^0_n$.

 The dimension of $S^r_n$ is given by the number of control points in a {\it minimal determining set (MDS)}, i.e. a minimal set of points nodal points $D$ such that  ( see definition  ~\ref{def:CPNMDS}, below) :
 \[
\forall P_{i,j,k} \in D  , \quad Z(\bar{P}_{i,j,k})=0 ,\textrm{  and } Z  \in C{^r }  \Rightarrow  Z \equiv 0 ;
\]
The problem of  the dimension of $S^r_n$ was initiated with a conjecture of Strang in ~\cite{strang}. The dimensionality depends on both the topological and geometrical structure of the mesh; an arbitrary small perturbation of the mesh vertices may lead to changes in structure of the minimal determining set. The first important result was achieved by Morgan and Scott ~\cite{morgan_scott}, with the dimension formula and explicit basis for space $S^1_n$, $n\ge 5$. The minimal determining sets (and therefore bases of the underlying spline spaces) were explicitly constructed for $S^1_4$ for all triangulations (see Alfred et all.  ~\cite{alfeld_piper_schumaker}); $S^r_n$, $n\ge 3r+2$ for general triangulations (see ~\cite{hong}, ~\cite{ibrahim_schumaker}); $S^r_n$, $n=3r+1$ for almost all triangulations (see ~\cite{alfeld_schumaker}). We are not aware of results for the case $r=1$ and  $n \le 3 $ .  In the latter case ,  subdivisions of the initial triangulation lead to the construction of the cubic spline finite element space. For the convex quadrangulation, the space $S^r_n$ is defined by triangulation obtained by inserting the diagonals of each quadrilateral (see ~\cite{lai}). For the approximate solution of boundary-value problem, the spaces of type $S^r_{3r}$, $r\ge 1$ for convex quadrangulations appears to be of particular interest, since they possess full approximation power (in contrast to spaces based on general triangulations), have relatively low dimensions and may be locally refined (see ~\cite{lai_schumaker}).

Generalisation of the bivariate spline finite element space - parametric spline finite element space, composed of such functions that every one of $X$, $Y$, $Z$ coordinates belongs to $S^r_n$ - is given in ~\cite{schumaker_parametric}. There surfaces are build by interpolation and avoid the vertex enclosure problem (see below). However,  parametrisation in the $XY$-plane can not be fixed a priory, which makes the approach unsuitable for minimisation of the energy relating to the functional representation of the surface.

\subsection{The principal aim of the present work}
\label{subsect:principal_aim}
\subsubsection{Contribution}

\paragraph{Generalisation of BSFE approach}
The current work generalises the BSFE approach on unstructured non degenerate convex quadrilateral meshing of a given planar domain $\T\Omega$ with subparametric Bezi\'er tensor-product "Finite Elements"  (FE) on each quadrilateral ( were we define an in plane parametrisation by a bilinear mapping from a reference element.) This  in-plane parametrisation leads to the  linearisation of the  $C^1$-continuity conditions 
and  reduces the problem to a  linear constrained minimisation (see Part ~\ref{part:general_linearisation}). We also extend our results to the case where the edge of a quadrilateral on the boundary of $\T\Omega$ is given by a cubic parametrisation.

 In addition, it provides a natural set up for explicitly   computing  the minimal determining set of the control points $\MDS{n}$ (see Subsection ~\ref{subsubsect:math_formulation}). 
We compute the $\MDS{n}$ for the space of $C^1(\T\Omega)$-smooth, piecewise parametric polynomials of degree $n\ge 4$.\\

 In common with the BSFE approach described in subsection (~\ref{subsubsect:FEM} )a single patch of the resulting surface is given by  a polynomial subparametric FE. However,the  degrees of freedom which guarantee the $C^1$-smooth concatenation between adjacent patches  are not local and the MDS depends on the topology and geometry of the mesh. The resulting construction is done   for all possible unstructured mesh quadrangulations (both from a topological and geometrical points of view).
 
The principal \underline{\it differences} from the standard BSFE approach are the following.
\begin{itemize}
\item
The elements are defined over a square rather than a triangular reference element and has a tensor-product polynomial form.
\item
Mapping between the reference element and the corresponding element in $\T\Omega$ is of at least of bilinear order. The resulting surface is given by a functional minimisation or by interpolation. In  classic BSFE the "energies' are expressed in the parametric space, not in a "physical" one. In our case, like  for isoparametric finite elements , the space of functions does not coincide with the space of functional polynomials over quadrilaterals. 
\item
The MDS depends on the choice of mappings between the reference square and the mesh elements. For a polygonal quadrilateral element the  canonic bilinear  mapping is used but  a boundary mesh element with one curvilinear side requires a special analysis in order to choose the mapping in an optimal way. 
\end{itemize}

The main contributions of the current work are listed below.
\begin{itemize}
\item

The current approach works for quadrilateral meshes with any valence for the vertices (hence it can use standard FEM quad mesh generators .)
\item
The MDS are constructed for both polygonal meshes and meshes which at the boundary of  $\T\Omega$ consists of $G^1$-smooth piecewise B\'ezier cubics . In the last case,  mappings between reference and boundary mesh elements are of higher order. Handling curved boundaries lead to better  approximate solutions of  partial differential equations.
\item
While the dimension of the MDS is uniquely defined, the choice of control points, which participate in MDS ( i.e the basis functions), can be made in different ways. The current research analyses different MDS which are suitable for different "additional" conditions to cover the main types of interpolation and boundary conditions.
\end{itemize}

The current work is restricted to an analysis of the MDS; it does not analyse stability nor the approximation order of the solution, though the experimental results seem quite accurate.

\paragraph{A study of the continuity conditions}
As noted above by ~\cite{jorg} , on any quadrilateral mesh  imposing of $G^1$-continuity conditions is equivalent to the requirement of $C^1$-smoothness.
Hence, it is sufficient to analyse the $G^1$ continuity conditions for the inner edges in terms of control points of adjacent patches.

Moreover the study of $G^1$-continuity conditions for the control points adjacent to a mesh vertex leads to results which fit the general Vertex Enclosure theory (see Section ~\ref{sect:vertex_enclosure}). The results have elegant geometrical formulations, closely related to the structure of the planar mesh. 

In addition, a detailed analysis is made for the "middle" control points adjacent to an inner edge. All possible configurations of the adjacent mesh elements are studied and the nice dependencies between the geometry of the elements and the available degrees of freedom are defined.

\paragraph{Choice of the energy functional}

The shape of the resulting surface (in addition to $G^1$-continuity and "additional" constraints) is defined by minimisation of a "natural physical" global energy functional, no local heuristics are used. Fixation of in-plane parametrisation a priory allows to define a functional form of energy, which makes the solution applicable to a wide range of problems in PDEs.

\subsubsection{Definitions and Hypotheses}
\label{subsubsect:math_formulation}

\begin{definition}
Consider the tensor product Bernstein-B\'ezier polyomials  of degree $n$ over the unit square  $[u,v]\in [0,1]^2$. Let  \underline{${\cal POL}^{(n)}$}, \underline{$\T {\cal POL}^{(n)}$} and \underline{$\bar {\cal POL}^{(n)}$} be the corresponding spaces ,obtained by defining, as coefficients of these polynomials ,scalar, $2D$ and $3D$ control points respectively.
 \end{definition}
We deal with  B\'ezier  functions,(resp. plannar domains and surfaces) , following the CAGD conventions , ~\cite{farin}, such objects will be defined as being of \emph{order}  $n$, (resp. ($ n,m$) and ($n, m,r$) ),  where  $n$ (resp. $ n,m$, and  $n, m,r$) defines the maximum (formal) degree(s)  of a polynomial, the actual 
degree(s) may be less. For example by degree elevation of its first tensor term , a  B\'ezier  bilinear quadrilateral [order ($1,1$) ], can be considered as an order  ($2,1$ ) quadrilateral. When there is no confusion we will use indifferently \emph{order} or \emph{degree}.

For a planar simply-connected domain $\T\Omega$, lying in $XY$-plane, and its quadrangulation $\T{\cal Q}$ into non degenerate convex elements, the following definitions and notations will be used. 

\begin{definition}
Let $m$ be some integer, a piecewise-polynomial $2D$ function $\T\Pi$ is an \underline{order $m$ global regular in-plane parametrisation} of domain $\T\Omega$ if  
\begin{itemize}
\item[{\bf (1)}]
For every mesh element $\T q\in \T{\cal Q}$, the restriction    $\T Q=\T\Pi|_{\T q}$ belongs to  $\T{\cal POL}^{(m)}$ and defines a regular mapping (see Definition ~\ref{def:reg_param}) between the reference square and the planar element $\T q$.
\item[{\bf (2)}]
For two adjacent mesh elements $\T q$ and $\T q'$, the $2D$ control points of $\T Q = \T\Pi|_{\T q}$ and $\T Q' = \T\Pi|_{\T q'}$ coincide along the common edge of the elements.
\end{itemize} 
The space of all order $m$ regular in-plane parametrisation will be denoted by $\PAR{m}$.
\label{def:reg_param_all}
\end{definition}

\begin{definition}
A piecewise-parametric $3D$ function $\bar\Psi$ \underline{agrees} with a given global in-plane parametrisation $\T\Pi$ if for every mesh element $\T q\in \T{\cal Q}$ the restriction of the function $\bar Q = \bar\Psi|_{\T q}$ defines a mapping from the unit square into the $3D$ space $\bar Q: [0,1]^2\rightarrow{\bf R^3}$ and the  $(X,Y)$ coordinates of $\bar Q$ coincide with the restriction of the global in-plane parametrisation $\T\Pi|_{\T q}$.
\label{def:agrees_param}
\end{definition}

\begin{definition}
Let $n$ and $m<n$ be some integers and $\T\Pi\in \T{\cal PAR}^{(m)}$ be some fixed degree $m$ global regular in-plane parametrisation of domain $\T\Omega$. Then \underline{space $\FUN{n}(\T\Pi)$} is by definition composed of piecewise-parametric $3D$ functions $\bar\Psi$ so that
\begin{itemize}
\item[{\bf (1)}]
$\bar\Psi$ agrees with the in-plane parametrisation $\T\Pi$.
\item[{\bf (2)}]
For every mesh element $\T q\in\T{\cal Q}$, the $Z$-coordinate of the restriction $\bar Q =\bar\Psi|_{\T q}$ belongs to ${\cal POL}^{(n)}$ (and hence $\bar\Psi|_{\T q}$ is a subparametric FE).
\item[{\bf (3)}]
$\bar\Psi$ is a $C^1$-smooth function in the functional sense over  $\T\Omega$: $\bar\Psi\in C^1(\T\Omega)$.
\end{itemize}
\label{def:space_fun}
\end{definition}
It is important to note that space $\FUN{n}(\T\Pi)$ depends on the chosen in-plane parametrisation $\T\Pi$, although in what follows it will be usually clear which underlying in-plane parametrisation is considered and the space will be usually denoted by $\FUN{n}$.

\begin{definition}
Let \underline{$\CP{n}(\T\Pi)$} be a set of in-plane B\'ezier control points of all patches which result from degree elevation of a global regular in-plane parametrisation $\T\Pi$ up to degree $n$ for every patch. Since the B\'ezier control points of the in-plane parametrisation always coincide along the shared edges, they are unambiguously well defined. 

A \underline{determining set} D is a subset of $\CP{n}$ so that :
 \[
\forall P \in D  , \quad Z(P)=0   \Rightarrow  \bar\Psi\ \equiv 0 ;
\]
A determining set is called \underline{minimal determining set (MDS) $\MDS{n}$} if there is no  determining  which size is smaller. 
\label{def:CPNMDS}
\end{definition}
The subset $\MDS{n}$ depends on the chosen in-plane parametrisation $\T\Pi$ and is not necessarily uniquely defined for a fixed in-plane parametrisation, but all instances  have the same size, equal to the dimension of the vector space generated by the corresponding basis functions.

The purpose of the current work is to choose an in-plane parametrisation in some optimum way and describe the MDS $\MDS{n}$ for all $n\ge 4$ and for all possible mesh configurations. Moreover, the control points which belong to the MDS may be chosen in different ways, defining different {\it instances} of $\MDS{n}$.
The principal goal of the work is to analyse the different instances of the MDS according to different "additional" constraints.

More details regarding different instances of the MDS, relations between MDS and the "additional" constraints and other important definitions relevant for the current approach are given in Section ~\ref{sect:general_flow}.

\subsection{Domains of application}

\subsubsection{Solution of an interpolation/approximation problem using the functional form of the shape functional}
\label{subsubsect:interpolation_problem_definition}

Let a set of $3D$ points be given and the topological structure of a $3D$ quadrilateral mesh be defined. The $3D$ mesh itself is "virtual" in the meaning that the quadrilateral faces of the mesh are defined in a symbolic manner, boundary curves of $3D$ patches are not specified. 

The goal is to construct such a {\it $G^1$-smooth}, {\it piecewise B\'ezier parametric} surface that {\it interpolates/approximates} the given $3D$ points and satisfies some optional additional conditions  (for example normals of the tangent planes at the mesh vertices or the boundary curve of the whole mesh may be specified \footnote{More detailed description of the possible kinds of additional constraints is given in Section ~\ref{subsect:linear_interp}}); the shape of the surface is defined (in addition to the interpolation/approximation conditions) by some {\it energy functional} which relates to the {\it functional} representation of the surface $Z=Z(X,Y)$.

Define  $3D$ quadrilateral elements  by connecting by straight segments  the endpoints of every edge of the "virtual" mesh (see Figure ~\ref{fig:fig3}). If the construction is such that  the orthogonal {\it projection} of these $3D$ elements on the $XY$-plane defines a {\it bijection} and forms a planar mesh of  {\it convex quadrilaterals} then our algorithm can be applied.

We consider in details the main kind of interpolation problems and provides a general approach, so that any interpolation/approximation problem can be handled in the same manner. The solution does not use the auxiliary construction of a $3D$ mesh of curves; the vertex enclosure constraints (see Section ~\ref{sect:vertex_enclosure}
) are intrinsically satisfied by the construction of the MDS, hence we can solve the problem for any structure of the mesh.

A comparison of the current approach and the standard  techniques based on interpolation o f $3D$ mesh of curves is presented in Appendix, Section ~\ref{sect:compare_interpolation}.

\subsubsection{Approximate solution of a partial differential equation over an arbitrary quadrilateral mesh}
\label{subsubsect:pde_problem_definition}


Consider 4th order  partial differential equations (PDEs) (for example, the Thin Plate Problem or a biharmonc operator ): to have a conforming FEM one needs  a $C^1$ basis.  An approximate solution  can then be found by constrained minimisation of a corresponding  energy functional (see  ~\cite{ciarlet_old} ~\cite{zienkiewicz},  ~\cite{timoshenko_main}). Constraints (fixed explicitly or implicitly) are  used for the imposition of  {\it boundary conditions}, the {\it energy} functional is derived from  the PDE and the computed solution can be represented  as  the {\it functional representation} of the resulting surface (an example of the Thin Plate functional is given in Section ~\ref{sect:application_example}).
 
\begin{figure}[!ph]
\centering
\includegraphics[clip,height=2.5in]{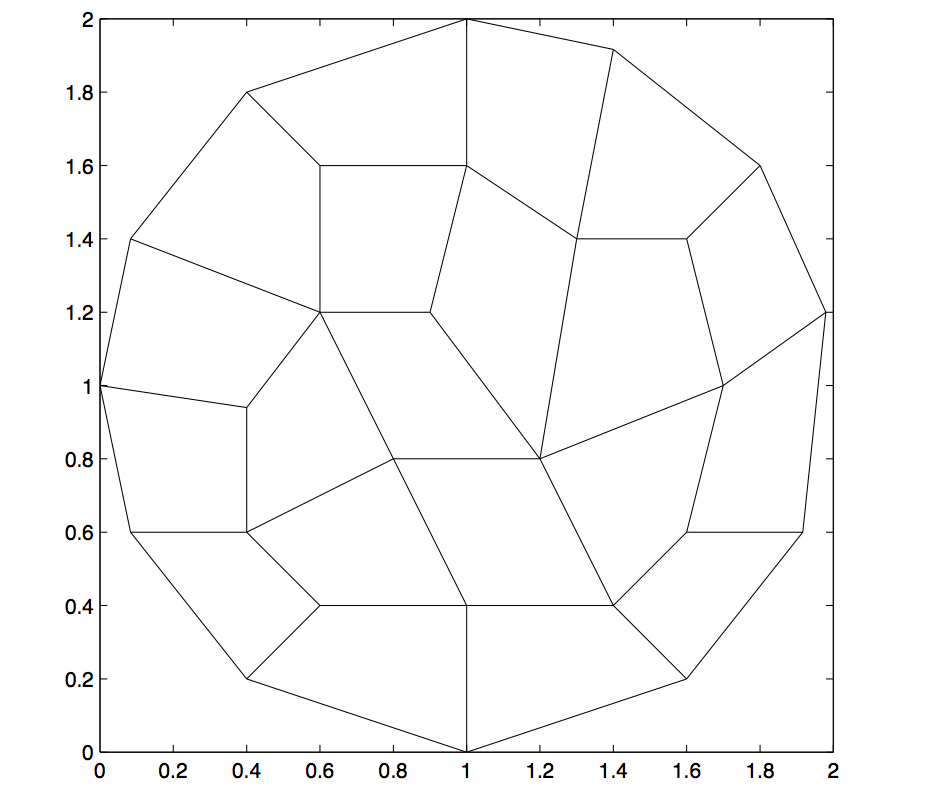}
\caption{Irregular quadrilateral mesh for a circular domain.}
\label{fig:fig1}
\end{figure}

Usually subdivision of the domain into elements is done by classical 2D meshing techniques ( see for instance ~\cite{blaker} ) and is included in the input of the problem together with the domain itself. The current research provides a possibility to construct a {\it $C^1$-smooth piecewise-polynomial approximate solution of a 4rth order partial differential equation for quadrilateral meshes with an arbitrary geometrical and topological structure}, like the mesh shown in Figure ~\ref{fig:fig1}. (Requirement of $C^1$-smoothness leads to a conforming approximate solution for  fourth-order partial differential equations.) The solution is constructed and the different boundary conditions are analysed for a planar mesh with a polygonal global boundary and for a planar mesh with piecewise-cubic $G^1$-smooth B\'ezier parametric global boundary. Although the error estimation lies beyond the scope of the current work, practical results (see Section ~\ref{sect:application_example}) show that the approach leads to an approximate solution of a high quality.  In the  spirit of  Isogeometric Analysis one could also approximate 2nd order PDEs by mean of these $C^1$ bases.

\vspace{0.3in}
\subsection{Structure of the work}
\label{sect:structure_of_the_work}

\vspace{0.15in}
\paragraph{Contents}

 Section ~\ref{sect:problem_definition} (Part ~\ref{part:introduction}) describes the problem and the general approach to its solution, introduces essential concepts and formulates the principal goals of the research. We  compare the current research with related works, highlights our contributions  and describe the domains of application.

Part ~\ref{part:review_3D_G1} presents fundamental results and definitions from the common theory of smooth surfaces, closely connected to the subject of the present work.
The general method of solution, adopted in the current research, is described in Part ~\ref{part:general_linearisation}.
Two central Parts, ~\ref{part:linear_boundary} 
and ~\ref{part:smooth_boundary} , apply the general method for two different types of planar mesh configurations. These parts contain the most important theoretical results, related to the existence and to the explicit construction of the solution, both in regular and in all possible degenerated cases. Full proofs of theoretical results as well as implementation algorithms are provided.
Part ~\ref{part:mds_mixed_4_5} shows that definitions based on the common fundamental concepts can be naturally generalised. The Generalisation leads to the composite solution with a wide range of applications. 

Part ~\ref{part:illustrative_examples} presents the computational examples, which allow to illustrate the correctness of the approach, and the results of application of the general solution to the Thin Plate problem. Part ~\ref{part:further_research} discusses possible topics for further research.

\vspace{0.15in}

In order to free the main text from the technical details and long computations as much as possible, all proofs, less important statements or statements which require a complicated algebraic formulation (Technical  and Auxiliary  Lemmas  ) and some Algorithms are given in an Appendix (Sections ~\ref{sect:proofs}, ~\ref{sect:technical_lemmas} and ~\ref{sect:algorithms} respectively). In addition, the Appendix contains examples of the solution of the partial differential equation (Section 
~\ref{sect:figures_tp}) and several  Sections  that includes some auxiliary material (Sections ~\ref{sect:energy_example_bilinear}, ~\ref{sect:mds_to_algebraic_solution} and ~\ref{sect:compare_interpolation}).

{\it Figures.}
Generally, the Figures are placed in the end of every Section. Figures which present results of the practical application of the approach are given in Appendix, Section ~\ref{sect:figures_tp}. The full list of figures is given at the beginning of the work.

{\it Notations.}
A special Section (Section ~\ref{sect:Notations}) presents a list of the most common and useful formal notations and definitions.  In addition, Section ~\ref{sect:Notations} contains an index list for some essential definitions used in the current work.

{\it Fonts and underlines.} 

Tems and notions ,with a precise definition, are usually written in italic font and/or appear in quotes. (For example, the {\it "Middle"} system of equations or the {\it middle} control points).

\eject

\FloatBarrier

\begin{figure}[!pt]
\vspace{-1.3in}
\begin{narrow}{-0.2in}{0.0in}
\centering
\includegraphics[clip,height=3.7in]{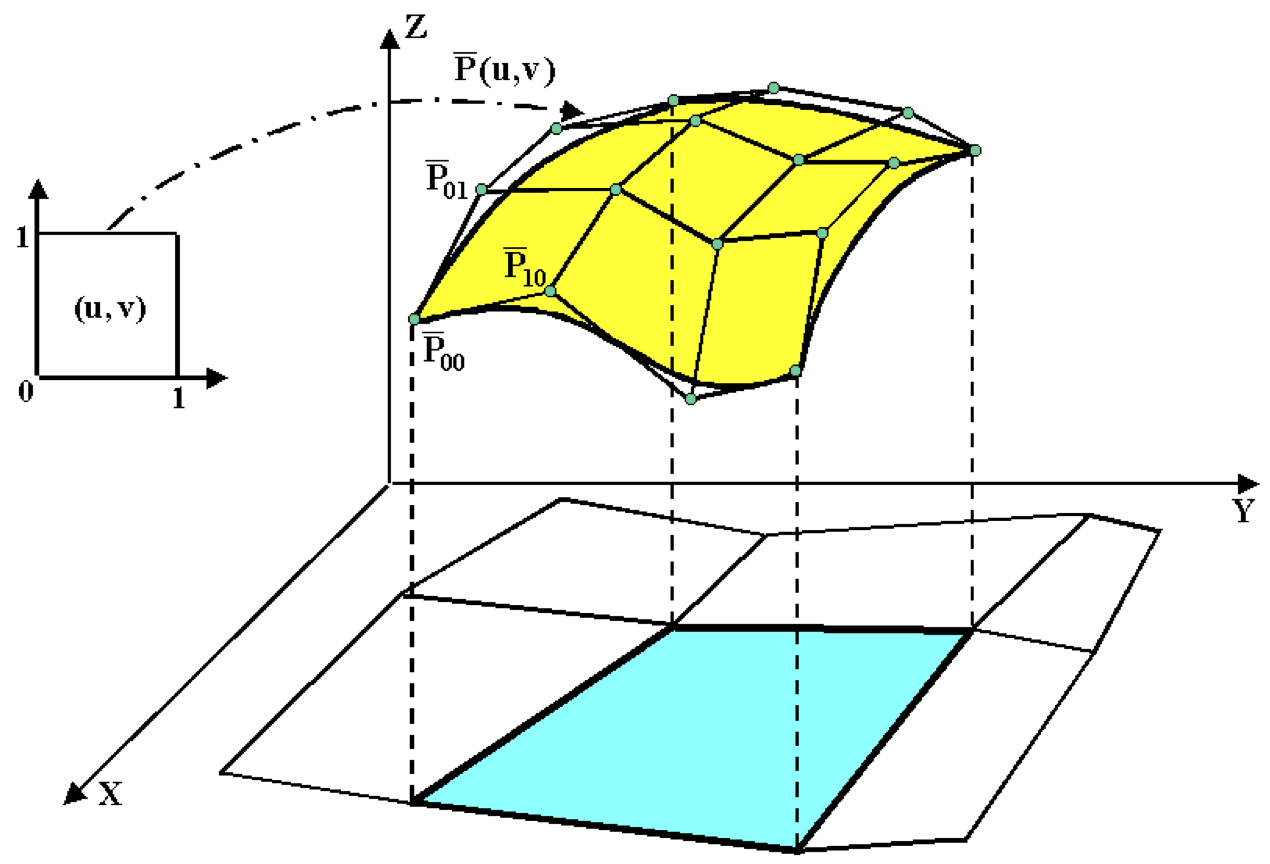}
\end{narrow}
\vspace{-0.1in}
\caption{$3D$ B\'ezier parametric patch and the corresponding planar mesh element.}
\label{fig:fig2}
\end{figure}

\begin{figure}[!pb]
\vspace{-0.1in}
\centering
\includegraphics[clip,height=3.6in]{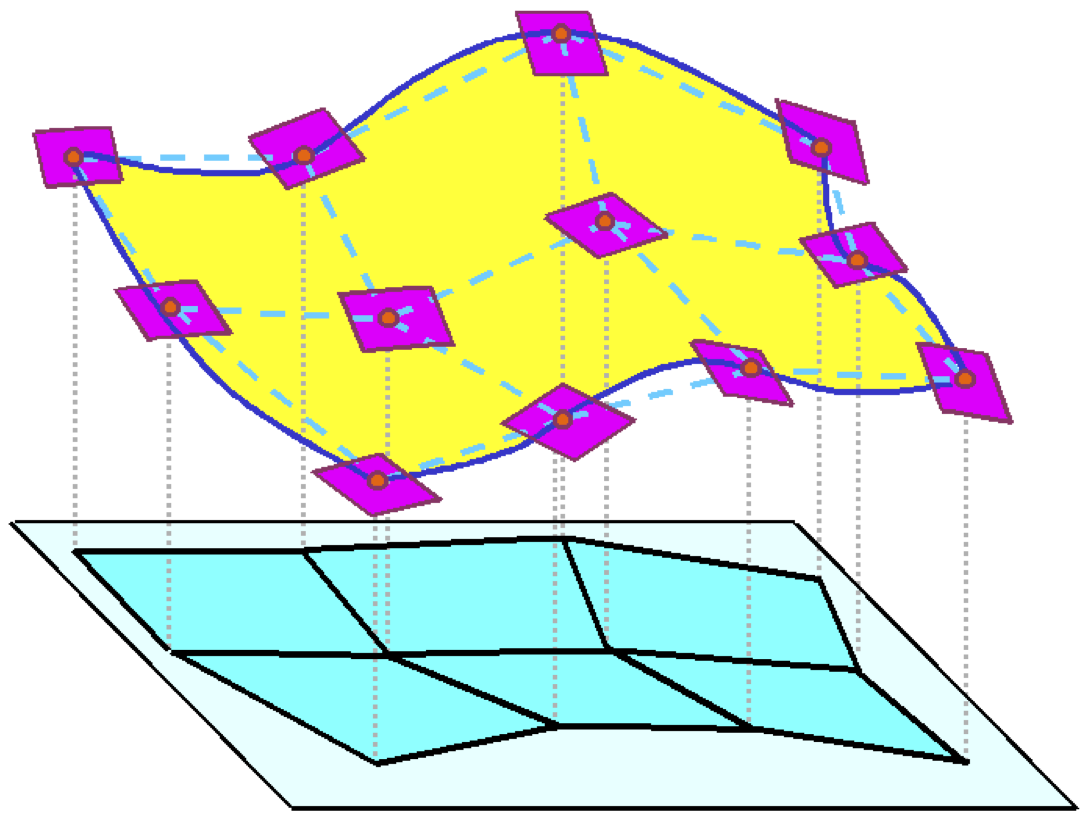}
\vspace{-0.1in}
\caption{Planar mesh and $3D$ surface interpolating $3D$ points and tangent planes at the mesh vertices.}
\label{fig:fig3}
\end{figure}

\FloatBarrier

\eject
\section{Notations}
\label{sect:Notations}
Here is  presents a list of the most common and useful formal notations and definitions. 

\subsection{Points and vectors}

\begin{description}

\item{$\bar A$ = $(\tilde A, A)$ = $(A_{X}, A_{Y}, A)$}
 - a $3D$ point or vector, where $\tilde A=(A_{X}, A_{Y})$ is its $2D$ component in
$XY$-plane and $A$ is its $Z$-coordinate. In general, $\ {\bar{}}\ $ will be used for $3D$ objects and $\ {\tilde{}}\ $ for $2D$ objects in $XY$-plane.

\item{$||\bar A||$ or $||\T A||$} - Euclidean norm of $3D$ or $2D$ vector.

\item{$<\bar A, \bar B>$} - cross product of two $3D$ vectors.

\item{$<\tilde A, \tilde B>$} - $Z$-coordinate of the cross product of two vectors lying in 
$XY$-plane (signed length of the cross product of two $2D$ vectors).

\item{$(\bar A, \bar B)$} - scalar product of two $3D$ vectors.

\item
$
\begin{array}{l}
mix\left(\!\!\begin{array}{ccc}\bar A\cr \bar B\cr \bar C\end{array}\!\!\right)\!=\!(\bar A, <\bar B, \bar C>) =
det\left(\!\!
\begin{array}{ccc}
A_{X},\! & A_{Y}\! & A \cr
B_{X}\! & B_{Y}\! & B \cr
C_{X}\! & C_{Y}\! & C \cr
\end{array}
\!\!\right)\cr
\end{array}
$ - "mixed" product of three $3D$ vectors. 
\end{description}

\subsection{Polynomials}
\label{subsect:def_polynomials}
\begin{description}
\item
$B^n_i(u)=\Cnk{n}{i} u^i(1-u)^{n-i}$, $u\in[0,1]$, $i=0,\ldots,n$  - degree $n$ Bernstein polynomial of one variable.
\item
$B^n_{ij}(u,v)=\Cnk{n}{i}\Cnk{n}{j} u^i(1-u)^{n-i} v^j(1-v)^{n-j}$, $(u,v)\in[0,1]^2$, $i,j=0,\ldots,n$ -
degree $n$ tensor-product Bernstein polynomial of two variables.
\item 
$P(u)=\sum_{i=0}^n P_i B_i^n(u)$ - B\'ezier polynomial of degree $n$, 
\item
$P_i$ , $\T P_i$, $\bar P_i$ - B\'ezier control points ( in $1D,2D,3D$).
\item

\item
$deg(P)$ - {\it actual} degree of a polynomial; the lowest integer such that $P(u)$ can be represented in the form $P(u)=\sum_{i=0}^{deg(p)}{\alpha}_i u^i$, with ${\alpha}_{deg(p)}  \neq 0$  .
\item 
$P(u,v)=\sum_{i,j=0}^n P_{ij} B_{ij}^n(u,v)$ -
tensor-product B\'ezier polynomial of order $n$, 

$P_{ij}$ (or $\T P_{i,j}$, $\bar P_{i,j}$) - B\'ezier control points (see Figure ~\ref{fig:fig2}). 
\end{description}

\subsection{Planar mesh data}
\subsubsection{Vertices, edges and twist characteristics}
\label{subsect:def_vertices_edges_twists_planar}

\paragraph{Vertices, edges and twist characteristics of a single mesh element}

\begin{description}
\item
$\T{A},\T{B},\T{C},\T{D}$ - four vertices of a convex quadrilateral planar mesh element (see Figure ~\ref{fig:fig8}).
\item
$\T{t}(\T{A},\T{B},\T{C},\T{D})\!=\!
\T{A}\!-\!\T{B}\!+\!\T{C}\!-\!\T{D}$ - twist characteristic of the element, double difference between the midpoints of the diagonals of the quadrilateral (see Figure ~\ref{fig:fig8}).
\end{description}

\paragraph{Vertices, edges and twist characteristics of two adjacent mesh elements}

\begin{description}
\item
$\LL,\LL',\GG, \GG', \RR, \RR' $ -  vertices of two adjacent planar mesh elements (see Figure ~\ref{fig:fig11}). 
\item
$\T e^{(R)}=\RR-\GG$,\ \ \  $\T e^{(C)}=\GG'-\GG$,\ \ \ $\T e^{(L)}=\LL-\GG$ - directed planar mesh edges adjacent to vertex $\GG$ (see Figure ~\ref{fig:fig11}).
\item
$\T t^{(R)}=\T t(\GG,\RR,\RR',\GG')$, $\T t^{(L)}=\T t(\GG,\GG',\LL',\LL)$, where $\T t$ is the 
twist characteristic of the planar mesh elements (see Figure ~\ref{fig:fig15}).
\end{description}

\paragraph{Edges and twist characteristics of elements adjacent to a given vertex}

Let planar mesh elements share the common vertex $\T V$ of degree $val(V)$ (see Figure ~\ref{fig:fig10})
\begin{description}
\item
$\T e^{(j)}=\T V^{(j)}-\T V$ ($j=1,\ldots,val(V)$) - directed planar mesh edges adjacent to vertex $\T V$, ordered counter clockwise.
\item
$\T t^{(j)} = \T V-\T V^{(j)}+\T F^{(j)}-\T V^{(j+1)}$ 
($j=1,\ldots,val(V)$ for an inner vertex and $j=1,\ldots,val(V)-1$ for a boundary vertex) - twist characteristics of planar elements adjacent to vertex $\T V$. 
\end{description}

\subsection{Partial derivatives of in-plane parametrisations}
\label{subsect:pd_planar_points_cubic}
Let two adjacent planar mesh elements be parametrized by $\T L(u,v)$ and $\T R(u,v)$ respectively and let $\T L(1,v)\equiv\T R(0,v)$ (see Figure ~\ref{fig:fig7}). 
\begin{description}
\item{$\T L_v = \T R_v = \D{\T L}{v}(1,v) = \D{\T R}{v}(0,v)$} - partial derivatives of the in-plane parametrisations in the direction along the common edge.

\item{$\T L_u=\D{\T L}{u}(1,v),\ \ \ \T R_u=\D{\T R}{u}(0,v)$} - partial derivatives of the in-plane parametrisations along the common edge in the cross direction.

\item{$\T\lambda_i,\T\rho_i$} - coefficients of the polynomials $\T L_u$ and $\T R_u$ with respect to the B\'ezier basis.

\item{$\T\lambda^{(power)}_i,\T\rho^{(power)}_i$} - coefficients of the polynomials $\T L_u$ and $\T R_u$ with respect to the power basis.

\end{description}

\subsection{Weight functions}
\label{subsect:def_weight_functions}
\begin{description}
\item{$c(v),l(v),r(v)$} - scalar weight functions from the definition of $G^1$-continuity (see Definition ~\ref{def:G1}). Note that these are B\'ezier \emph{functions.}

\item{$c_j, l_j, r_j$} - coefficients of the weight functions with respect to the B\'ezier basis.

\item{$c^{(power)}_j, l^{(power)}_j, r^{(power)}_j$} - coefficients of the weight functions with respect to the power basis.

\item{$(order(c),order(l),order(r))$ (or $(deg(c),deg(l),deg(r))$}) - match which is defined by formal (or actual) degrees of the weight functions.

\item{$max\_deg(l,r)$} - maximum among $deg(r)$ and $deg(l)$.

\end{description}

\subsection{$3D$ data of the resulting surface}

\subsubsection{B\'ezier control points of two adjacent patches}
\label{subsect:def_Bezier_cp_two_adjacent_patches}


\begin{description}
\item{$\bar C_j$, $\bar L_j$, $\bar R_j$} ($j=0,...,n$) - for two adjacent patches, B\'ezier control points along the common edge and of the rows adjacent to the common edge in the left and the right patch respectively. Control points $\bar C_j$  are called "central" control points and control points $\bar L_j$, $\bar R_j$ are called "side" control points  (see Figure ~\ref{fig:fig5}).

\item{$\Delta\bar C_j,\Delta\bar L_j,\Delta\bar R_j$} -
first-order differences of the control points  (see Figure ~\ref{fig:fig5})
\vspace{-0.1in}
\begin{equation}
\begin{array}{lll}
\Delta\bar C_j = \bar C_{j+1} - \bar C_{j},& j=0,\ldots,n-1 \cr
\Delta\bar L_j = \bar C_j - \bar L_j,& j=0,\ldots,n \cr
\Delta\bar R_j = \bar R_j - \bar C_j,& j=0,\ldots,n \cr
\end{array}
\label{eq:deltas}
\end{equation}
\vspace{-0.15in}
\end{description}

\subsubsection{B\'ezier control points adjacent to some mesh vertex}
\label{subsect:notation_vertex_adjacent_cp}
Let $\T V$ be a planar mesh vertex of valence $val(V)$ and let $\T V$ have at least one adjacent inner edge. The following notations are used for the control points adjacent to the vertex and participating in $G^1$-continuity conditions for at least one inner edge (see Figure ~\ref{fig:fig6})
\begin{description}
\item
$\bar V$ - control point corresponding to the mesh vertex, the control point (or its components) will be called  $V$-type control point.
\item
$\bar E^{(1)},\ldots,\bar E^{(val(V))}$ for an inner vertex or
$\bar E^{(2)},\ldots,\bar E^{(val(V)-1)}$ for a boundary vertex
- the first control points adjacent to $\bar V$ along the inner edges emanating from the vertex; these control points (or their components) will be called tangent or $E$-type control points.
\item
$\bar D^{(1)},\ldots,\bar D^{(val(V))}$ for an inner vertex 
$\bar D^{(2)},\ldots,\bar D^{(val(V)-1)}$ for a boundary vertex
- the second control points adjacent to $\bar V$ along the inner edges emanating from the mesh vertex; these control points (or their components) will be called $D$-type control points.
\item
$\bar T^{(1)},\ldots,\bar T^{(val(V))}$ for an inner vertex or $\bar T^{(1)},\ldots,\bar T^{(val(V)-1)}$ for a boundary vertex, which are adjacent to $\bar V$ and do not lie at any edge; these control points (or their components) will be called twist or $T$-type control points.
\end{description}

\subsubsection{Partial derivatives of patches at the common vertex}
\label{subsect:def_deriv_patches_at_common_vertex}

\paragraph{Partial derivatives of two adjacent patches at the common vertex}

Let two adjacent patches be parametrized as shown in Figure ~\ref{fig:fig4} and let the parametrisations agree along the common edge. 
\begin{description}
\item
$\bar\epsilon^{(R)}$, $\bar\epsilon^{(C)}$, $\bar\epsilon^{(L)}$ -
the first-order derivatives along right, central and left edges computed at vertex $\bar V$.
\vspace{-0.05in}
\begin{equation}
\begin{array}{l}
\bar\epsilon^{(R)} = \D{\bar R}{u}(0,0) = n\Delta\bar R_0\cr
\bar\epsilon^{(C)} = \D{\bar R}{v}(0,0) = \D{\bar L}{v}(1,0) = n\Delta\bar C_0\cr
\bar\epsilon^{(L)} = -\D{\bar L}{u}(1,0)= -n\Delta\bar L_0
\end{array} 
\end{equation}
\vspace{-0.05in}
\item
$\bar\tau^{(R)}$ - the second-order mixed partial derivatives of the left and the right patches computed at vertex $\bar V$.
\vspace{-0.05in}
\begin{equation}
\begin{array}{l}
\bar\tau^{(R)}=\DD{\bar R}{u}{v}(0,0)=
n^2(\Delta\bar R_1-\Delta\bar R_0)\cr
\bar\tau^{(L)}=-\DD{\bar L}{u}{v}(1,0)=
-n^2(\Delta\bar L_1-\Delta\bar L_0) 
\end{array}
\end{equation}
\vspace{-0.05in}
\item
$\delta^{(C)}=\DT{R}{v}(0,0)=\DT{L}{v}(1,0)=n(n-1)(\Delta C_1-\Delta C_0)=n(n-1)(C_2-2 C_1+C_0)$ 
- $Z$-component of the second-order partial derivative along the central edge computed at vertex $\bar V$.
\end{description}

\paragraph{Partial derivatives of all patches sharing a common vertex}

Let patches sharing a common vertex $\bar V$ be parametrized as shown in Figure ~\ref{fig:fig6} and let parametrisations of adjacent patches agree along the common edges.


\begin{description}
\item
$\bar\epsilon^{(j)} = n(\bar E^{(j)}-\bar V) = \D{\bar P^{(j-1)}}{v}(0,0) = 
\D{\bar P^{(j)}}{u}(0,0)$ 
- first-order partial derivative in  the direction  of edge $\T e^{(j)}$, computed at vertex $\bar V$.
\item
$\bar\tau^{(j)} = n^2 (\bar V-\bar E^{(j)}+\bar T^{(j)}-\bar E^{(j+1)}) = 
\DD{\bar P^{(j)}}{u}{v}(0,0)$ 
- second-order mixed partial derivative of patch 
$\bar P^{(j)}$ computed at vertex $\bar V$.
\item
$\delta^{(j)} = n(n-1) (D^{(j)}-2E^{(j)}+V) =
\DT{P^{(j-1)}}{v}(0,0) = \DT{P^{(j)}}{u}(0,0)$-\\ 
$Z$-component of the second-order partial derivative in the direction of edge $\T e^{(j)}$, computed at vertex $\bar V$.
\end{description}

\subsection{Definitions of special sets, spaces and equations}
\label{subsect:def_sets_spaces_equations}

\paragraph{Definitions of some sets of control points}

\begin{description}
\item
$\MDS{n}$, $\GMDS{n}$, $\MDS{4,5}$ - see Definitions ~\ref{def:CPNMDS}, ~\ref{def:cp_subsets} and ~\ref{def:mds_mixed_4_5} respectively.
\item
$\CP{n}$, $\GCP{n}$, $\FCP{n}$, $\CP{4,5}$ - see Definitions 
~\ref{def:CPNMDS}, ~\ref{def:cp_subsets}, ~\ref{def:cp_subsets} and ~\ref{def:cp_mixed_4_5} respectively.
\item
{\it middle} control points - see Definition ~\ref{def:middle_bilinear} for a mesh with a polygonal global boundary and Lemma ~\ref{lemma:unchanged_middle_bicubic} and Definition ~\ref{def:middle_edge_one_boundary_bicubic} for a mesh with a smooth global boundary.
\end{description}

\paragraph{Definitions of some functional spaces}

\begin{description}
\item
$\FUN{n}$, $\FUN{4,5}$ - see Definition ~\ref{def:space_fun} and Part ~\ref{part:mds_mixed_4_5} respectively.
\item
$\PAR{m}$ - see Definition ~\ref{def:reg_param_all}.
\item
$S^r_n$ - see Subsection ~\ref{subsubsect:FEM}.
\end{description}

\paragraph{Definitions of some special equations and systems of equations}
\begin{description}
\item
{\it "Eq(s)"} indexed equation - see general Definition 
~\ref{def:indexed_equation}, Lemma ~\ref{lemma:weights_equations_bilinear} for a mesh with a polygonal global boundary and Lemma ~\ref{lemma:Eq_sumC_linear_system} for a mesh with a smooth global boundary.
\item
{\it "Eq(s)-type"} equations - see Definition ~\ref{def:indexed_equation}.
\item
{\it "Middle"} system of equations - see Definition ~\ref{def:middle_bilinear} for a mesh with a polygonal global boundary and Lemma ~\ref{lemma:unchanged_middle_bicubic} and Definition ~\ref{def:middle_edge_one_boundary_bicubic} for a mesh with a smooth global boundary.
\item
{\it "Restricted Middle"} system of equations - see Definition 
~\ref{def:middle_edge_one_boundary_bicubic}.
\item
{\it "sumC-equation","C-equation"} - see Definition  ~\ref{def:sumLR_sumC} .
\end{description}

\vspace*{0.2in}
\eject
\begin{figure}[!pb]
\centering
\includegraphics[clip,height=2.5in]{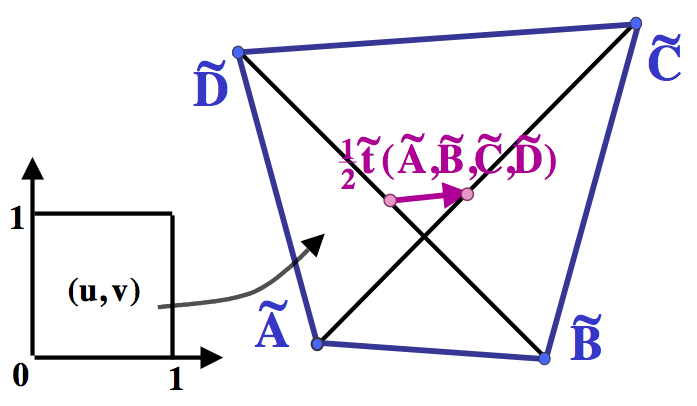}
\vspace{-0.1in}
\caption{A planar mesh element.}     
\label{fig:fig8}
\end{figure}


\begin{figure}[!ph]
\centering
\includegraphics[clip,height=2.8in]{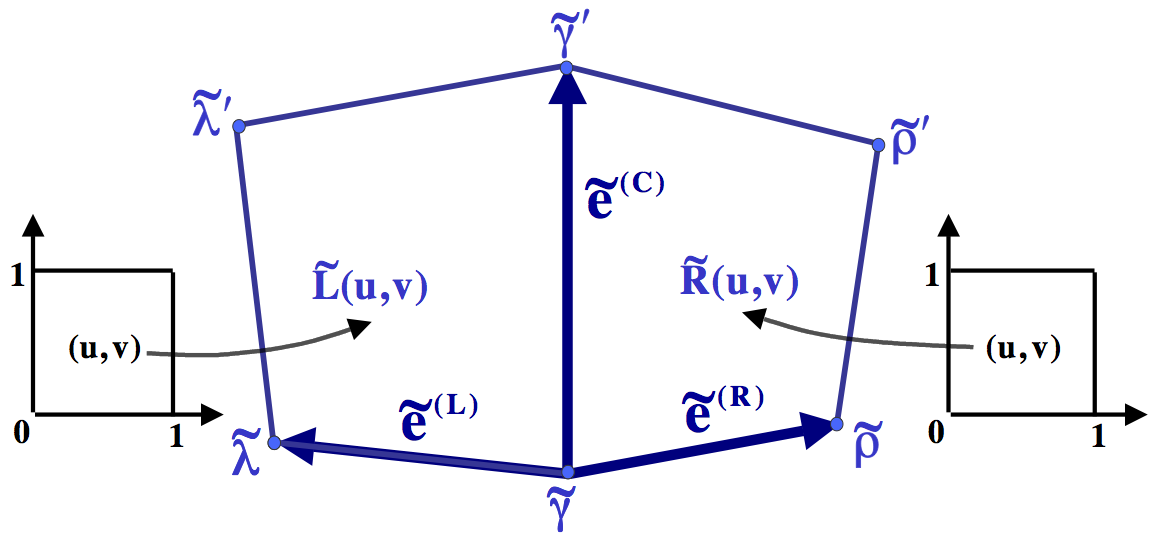}
\vspace{-0.1in}
\caption{Two adjacent planar mesh elements.}     
\label{fig:fig11}
\end{figure}

\FloatBarrier

\begin{figure}[!pt]
\vspace{-1.0in}
\centering
\includegraphics[clip,totalheight=2.8in]{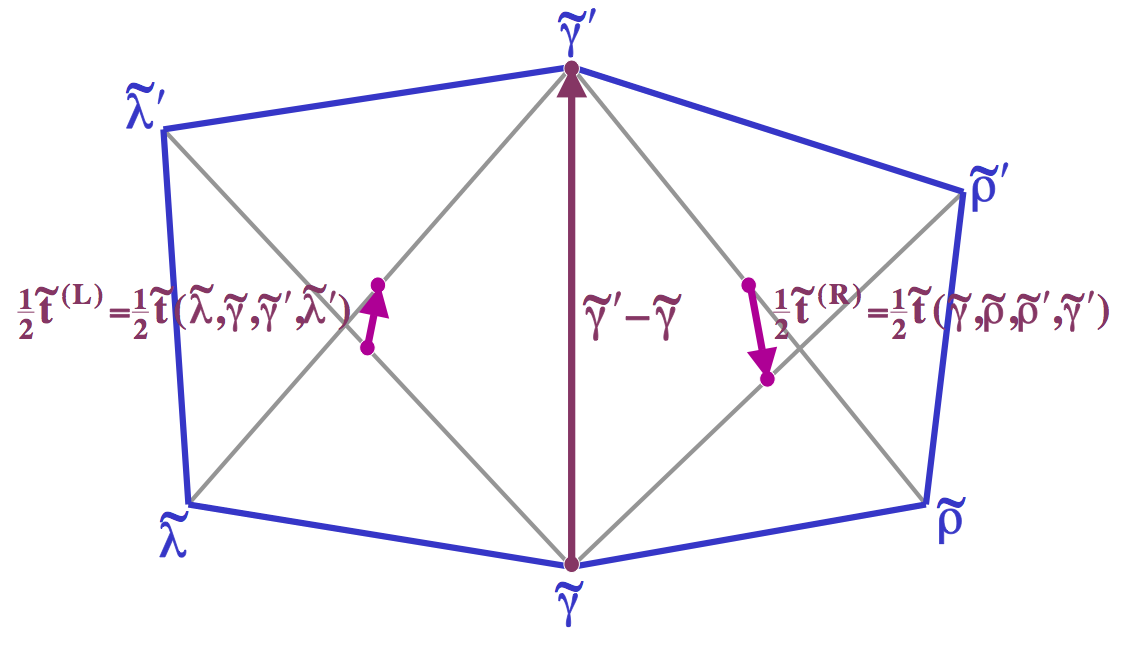}
\vspace{-0.15in}
\setcaptionwidth{5.2in}
\caption{Vectors important for computation of the actual degrees of the weight functions in the case of the bilinear parametrisation of two adjacent mesh elements.}
\label{fig:fig15}
\end{figure}

\begin{figure}[!ph]
\vspace{-0.05in}
\centering
\includegraphics[clip,height=2.3in]{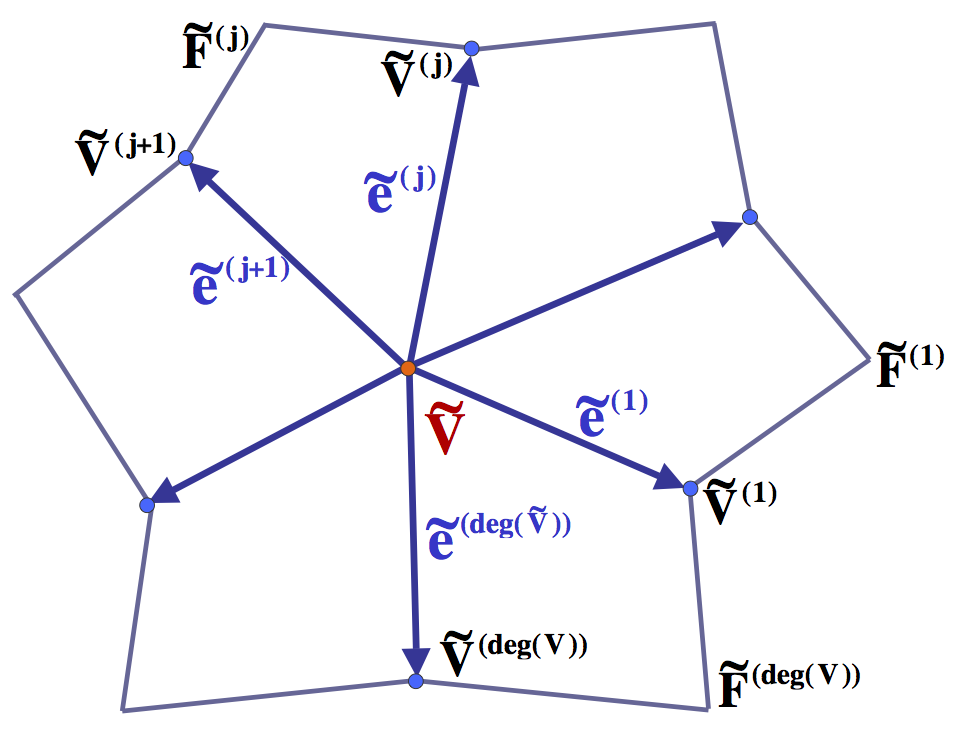}
\caption{Planar mesh elements adjacent to the common vertex.}
\label{fig:fig10}
\end{figure}


\begin{figure}[!pt]
\vspace{0.1in}
\begin{narrow}{-0.1in}{0.0in}
\centering
\includegraphics[clip,totalheight=2.0in]{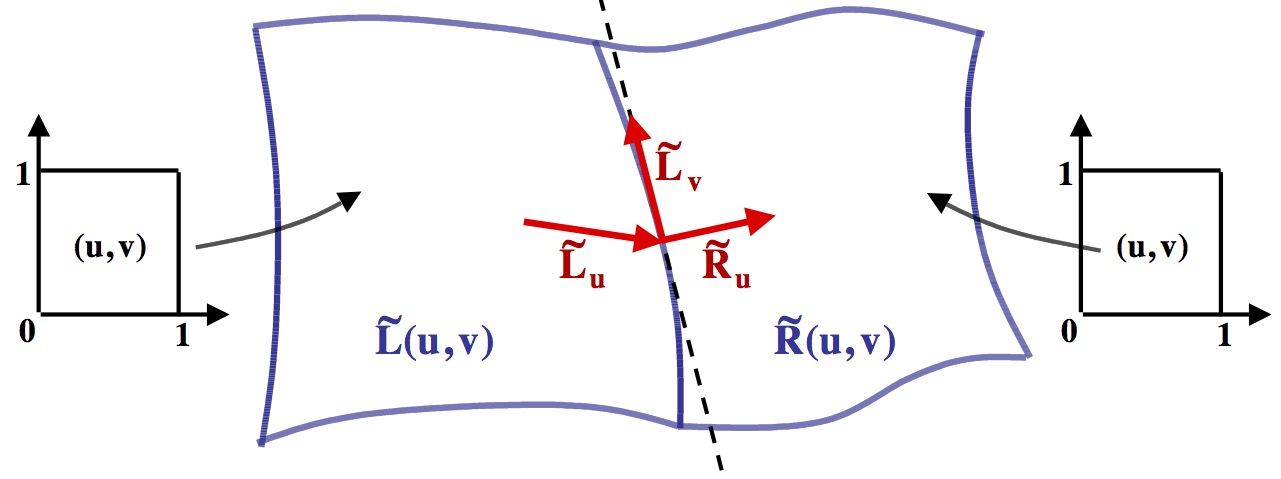}
\end{narrow}
\vspace{-0.15in}
\caption{Partial derivatives of in-plane parametrisations for two adjacent patches.}
\label{fig:fig7}
\end{figure}

\begin{figure}[!ph]
\vspace{0.2in}
\begin{narrow}{-0.1in}{0.0in}
\centering
\includegraphics[clip,height=3.3in]{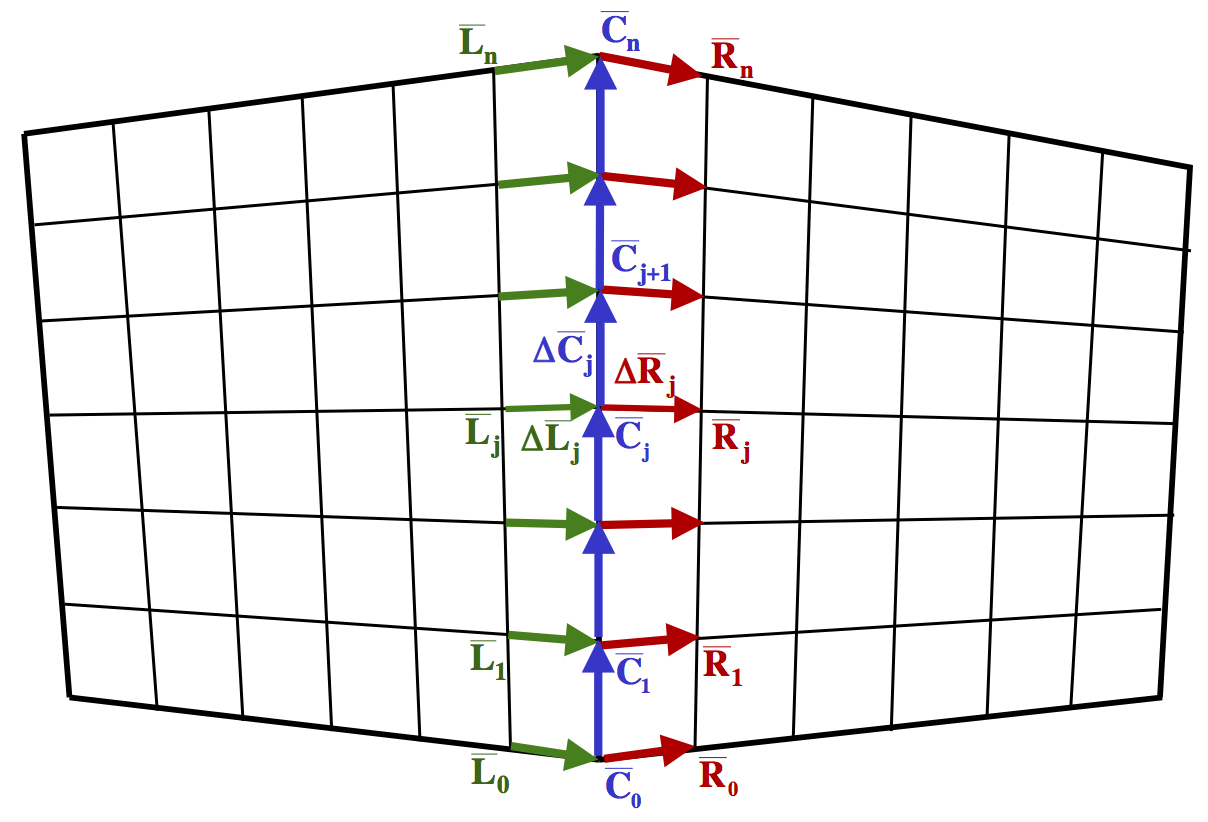}
\end{narrow}
\caption{B\'ezier control points adjacent to the common edge of two patches.}
\label{fig:fig5}
\end{figure}

\FloatBarrier

\begin{figure}[!pt]
\includegraphics[clip,height=3.9in]{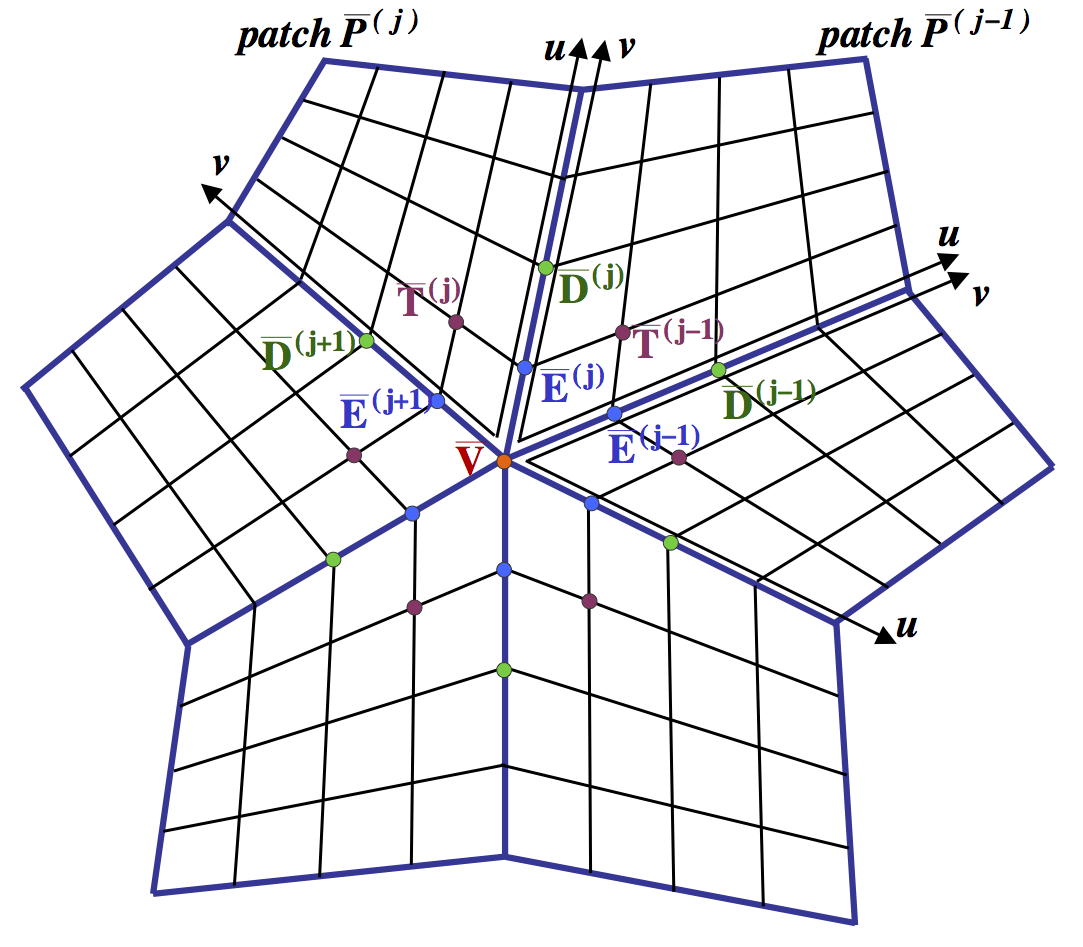}
\caption{B\'ezier control points adjacent to some mesh vertex.}
\label{fig:fig6}
\end{figure}

\vspace{0.1in}
                            
\begin{figure}[!pb]
\begin{narrow}{-0.3in}{0.0in}
\includegraphics[clip,height=2.0in]{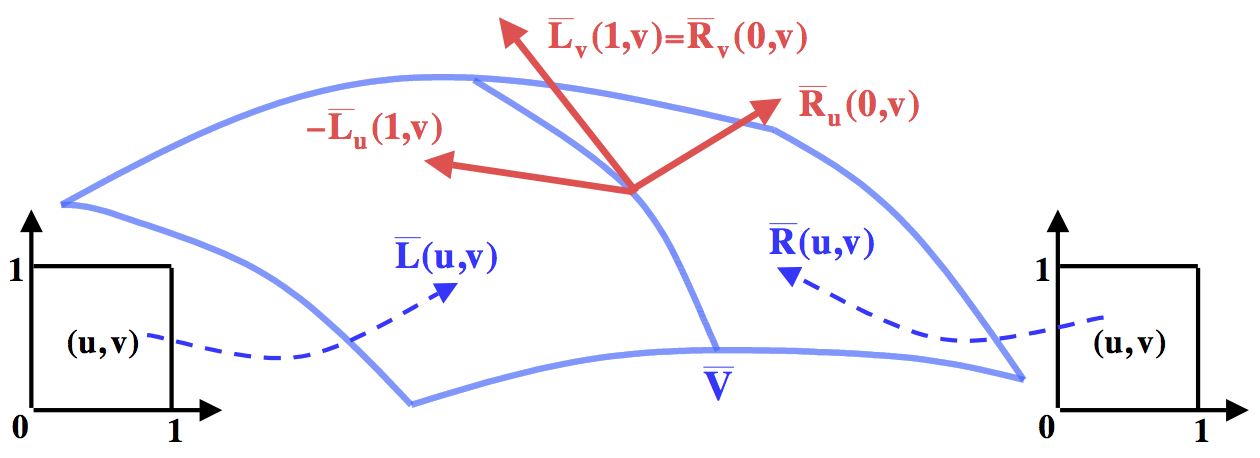}
\end{narrow}
\caption{ Two adjacent parametric patches with $G^1$-smooth concatenation along the common boundary.}
\label{fig:fig4}
\end{figure}

\FloatBarrier

\eject
\part{Some fundamental results regarding $G^1$-smooth surfaces}
\label{part:review_3D_G1}

Construction of the MDS (defined in Subsection ~\ref{subsubsect:math_formulation}) is based on the analysis of the smoothness conditions between adjacent patches. The current Section contains the formal definitions of the different kinds of smoothness and presents the general theoretical results of the vertex enclosure problem, which are closely connected to the analysis of the local structure of the MDS.

\section{Basic definitions related to smoothness of the surface}
\label{sect:smoothness_definitions}
%
 Two kinds of smoothness - functional and parametric ones - will be involved. The following standard definitions are used.

\begin{definition}
[$C^1$-smoothness of a functional surface]
A functional surface $Z(X,Y)$ is \underline{$C^1$-smooth} over domain $\tilde\Omega$ if for every point
$(X,Y)\in\tilde\Omega$ the first-order partial derivatives $\frac{\partial Z}{\partial X}(X,Y)$,
$\frac{\partial Z}{\partial Y}(X,Y)$ are well defined and 
 continuous over $\tilde\Omega$.
\label{def:C1_functional}
\end{definition}

In order to define $G^1$ parametric smoothness, let us consider two adjacent quadrilateral patches $\bar L(u,v)$ and $\bar R(u,v)$. Let every one of the patches be parametrized by the unit square (see Figure ~\ref{fig:fig4}), such that their  parametrisations agree along the common edge and the concatenation  between the patches is $G^0$-smooth (continuous)
\begin{equation}
\bar L(1,v)=\bar R(0,v) {\rm \ for\ every\ } v\in[0,1]
\label{eq:C0_agree}
\end{equation}  
In addition, every patch is supposed to be sufficiently smooth, with at least a continuous first order   partial derivatives along the common edge . Equation ~\ref{eq:C0_agree} implies that 
\begin{equation}
\bar L_v(1,v)=\bar R_v(0,v) {\rm \ for\ every\ } v\in[0,1]
\label{eq:deriv_along_bound_equal}
\end{equation} 

\begin{definition}
[$G^1$-smooth concatenation  between two parametric patches]
Patches $\bar L(u,v)$ and $\bar R(u,v)$ join  \underline{$G^1$-smoothly} along the common edge if and only if there exist a scalar-valued weight functions $l(v)$, $c(v)$, $r(v)$ such that for every 
$v\in[0,1]$ 
\begin{equation}
\bar L_u(1,v) l(v) + \bar R_u(0,v) r(v) + \bar L_v(1,v) c(v) = 0 
\label{eq:G1_eq_general_1}
\end{equation}
\begin{equation}
l(v)r(v) < 0 
\label{eq:G1_eq_general_2}
\end{equation}
\begin{equation}
<\bar L_u(1,v),\bar L_v(1,v)> \neq 0 
\label{eq:G1_eq_general_3}
\end{equation}
(see Definition given in ~\cite{peters_main}).
\label{def:G1}
\end{definition}

Geometrically $<\bar L_u,\bar L_v>$ and $<\bar R_u,\bar R_v>$ define the tangent plane normals for the left and the right patches respectively. Equation ~\ref{eq:G1_eq_general_1} means that the tangent planes of the adjacent patches are co-planar along the common edge. Equation  ~\ref{eq:G1_eq_general_3} means  normal to the  the tangent plane does not vanish and Equation ~\ref{eq:G1_eq_general_2} controls the orientation of the patches
in order to avoid cusps.
We have the following Lemma :
\begin{lemma}
Let two patches with a linear common edged join $G^1$ smoothly   then their respective parametrisation join $C^1$ continuously. 
\label{lemma:jorg}
\end{lemma}
{\bf Proof}
The common face of the two patches being a linear segment, one can trivially reparametrize one as image of the second  ,this  is a special case Peters'    fundamental  Lemma ~\cite{jorg}\noindent. 
\nopagebreak 
\eop${}_{{\bf Lemma ~\ref{lemma:jorg}}}$


It is then possible to combine the functional representation of the surface defining the energy functional and the parametric representation of the surface in order to impose the $G^1$smoothness constraints in parametric form.

\section{The vertex enclosure problem}
\label{sect:vertex_enclosure}

The current Section mainly relates to the work ~\cite{peters_main}, devoted to smooth interpolation of a given $3D$ mesh of curves. The work is chosen as the main reference, because it formulates and analyses in details the general vertex enclosure constraint. The satisfaction of the vertex enclosure constraints determines the existence of a $G^1$-smooth interpolant for a given $3D$ mesh of curves. (Generalisation of the vertex enclosure problem to the case of concatenation of a few patches around a common vertex with a definite degree of smoothness can be found in ~\cite{peters_complicated}.)

\subsection{General formulation of the vertex enclosure constraint}

Let a $3D$ mesh of polynomial curves be given, ( we only study  meshes  which faces are  $4$-sided,) construct  a $G^1$-smooth piecewise B\'ezier tensor-product interpolanting these curves. ( ~\cite{peters_main} contains a full analysis for the mixed triangular/quadrilateral meshes and shows that from a theoretical point of view the quadrilateral or triangular form of a patch does not lead to essential differences for the vertex enclosure constraint).
In the problem formulated above, the boundary curves of every patch (mesh curves) are given and the inner B\'ezier control points of every patch play the role of unknowns. These unknowns should satisfy $G^1$-continuity constraints, which means that the weight functions from Definition ~\ref{def:G1} should exist for every inner edge of the mesh. In particular these functions should exist for every one of the edges that share a common inner vertex. Consideration of $G^1$-continuity constraints together for all edges adjacent to the same vertex leads to a so called vertex enclosure problem. 
The vertex enclosure constraint is met at a mesh vertex $\bar V$ 
if weight functions could be simultaneously defined for each mesh edge emanating from $\bar V$.

For an inner mesh vertex Equations ~\ref{eq:G1_eq_general_1} applied to all edges emanating from the vertex have a circulant nature (the "left" patch of the "first" edge is also the "right" patch of the "last" edge) and lead to a linear system of equations such that the matrix of the system has a circulant structure. Independently of the order and geometry of the mesh curves, the matrix is always invertible at the odd vertices and rank deficient at the even vertices. At the even vertices the rank of the matrix is equal to its size minus one, which generally means that one additional constraint for every even mesh vertex should be satisfied in order to allow a $G^1$-smooth interpolation. A mesh is called {\it admissible} if a $G^1$ smooth interpolant can be constructed
(or in other words, if weight functions for all inner edges can be defined without contradictions).

In Peters, ~\cite{peters_main}, the vertex enclosure constraint is considered in its most general form (for example, the mesh curves sharing the common vertex may have different polynomial degrees), which leads to quite complicated equations. The constraint is not written explicitly, sufficient conditions that allow concluding that a given mesh is {\it admissible} are supplied.

We will show that in our  case, the explicit form of the vertex enclosure constraint becomes very simple and elegant. The Subsection presents  formulas from Peters ~\cite{peters_main} in order to verify later that results of the current work fit the general theory. The general results are formulated in notations of the current work (see Section ~\ref{sect:Notations}). Although it makes the presentation quite different from its original, the conversion between different forms of presentation is purely technical and straightforward. In order to make the formulas more compact and clear, some minor simplifying assumptions, which are always satisfied in the current work, will be used.

According to the problem definition, a quadrilateral mesh of curves of degree $m$ should be $G^1$-smoothly interpolated by piecewise tensor-product B\'ezier patches of degree $n$.  $G^1$-continuity between a pair of adjacent patches implies that the following two equations should be satisfied.

{\it "Tangent Constraint"}
\begin{equation}
c_0 \Delta\bar C_0 + l_0 \Delta\bar L_0 + r_0 \Delta\bar R_0 = 0
\label{eq:Peters_tangent}
\end{equation}

{\it "Twist Constraint"}
\begin{equation} 
\begin{array}{ll}
(n-1)c_0\Delta\bar C_1+deg(c)c_1\Delta\bar C_0+\cr
n\ l_0\Delta\bar L_1+deg(l)l_1\Delta\bar L_0+
n\ r_0\Delta\bar R_1+deg(r)r_1\Delta\bar R_0=0
\end{array}
\label{eq:Peters_twist_1}
\end{equation}
Here notations from Section ~\ref{sect:Notations} are used. In particular $\bar C_j, \bar L_j, \bar R_j$ are the control points of B\'ezier patches (see Figures ~\ref{fig:fig4}) and $c_j,l_j,r_j$ are coefficients of the weight functions.

In the interpolation problem formulated in work ~\cite{peters_main},  tangents $\Dlt\bar C_0$, $\Dlt\bar L_0$, $\Dlt\bar R_0$ and boundary curve control points $\bar C_j$ ($j=0,\ldots,n$) are given, and twist control points $\bar L_1$ and $\bar R_1$ as well as coefficients of the weight functions serve as unknowns.

For a vertex $\bar V$ with $val(V)$ emanating curves, superscript $j$ will be used when tangent or twist constraint is considered for the curve with order number $j=1,\ldots,val(V)$. 
Control points $\bar C_i^{(j)}$,$\bar L_i^{(j)}$,$\bar R_i^{(j)}$ for $i=0,1$ participate in different roles in equations for more than one curve. In order to avoid an ambiguity, notations $\bar V$, $\bar E^{(j)}$ and $\bar T^{(j)}$ will be used respectively  for the vertex, tangent and twist control points (see Subsection ~\ref{subsect:notation_vertex_adjacent_cp}).

The {\it "Tangent Constraint"} defines (up to a scale factor) the zero-indexed coefficients of the weight functions. These coefficients depend on the geometry of the tangent vectors of curves emanating from the vertex.

The {\it "Twist Constraint"} for curve with order number $j$ can be rewritten in the form
\begin{equation}
m\ l_0^{(j)} \bar T^{(j)} - m\ r_0^{(j)} \bar T^{(j-1)} = \bar A^{(j)}
\label{eq:Peters_twist_2}
\end{equation}
where
\begin{equation}
\begin{array}{lll}
\bar A^{(j)} = &(n-1)c_0^{(j)}\Delta\bar C_1^{(j)} + deg(c^{(j)}) c_1^{(j)}\Delta\bar C_0^{(j)}+\cr
&l_0^{(j)}\left(m \bar E^{(j)}-(m-n)\Delta\bar L_0^{(j)}\right)+ deg(l^{(j)}) l_1^{(j)}\Delta\bar L_0^{(j)}-\cr
&r_0^{(j)}\left(m \bar E^{(j)}+(m-n)\Delta\bar R_0^{(j)}\right)+ deg(r^{(j)}) r_1^{(j)}\Delta\bar R_0^{(j)}\cr
\end{array}
\label{eq:Peters_twist_right_part}
\end{equation}
For an {\it inner} vertex $\bar V$, the {\it "Twist Constraint"} applied simultaneously to all  the $val(V)$ edges emanating from the vertex, leads to a circulant linear system of equations. 
\begin{equation}
MT=A
\label{eq:Peters_twist_matrix_1}
\end{equation}
Here matrix $M$ has a circulant structure
\begin{equation}
M= \left(\begin{array}{ccccccc}
 l_0^{(1)} & 0 & 0 & \ldots & 0 & 0 & -r_0^{(1)}\cr
-r_0^{(2)} & l_0^{(2)} & 0 & \ldots & 0 & 0 & 0 \cr
0 & -r_0^{(3)} & l_0^{(3)} & \ldots & 0 & 0 & 0 \cr
\vdots & \vdots & \vdots & \ddots & \vdots & \vdots & \vdots \cr
0 & 0 & 0 & \ldots &  l_0^{(val(V)-2)} & 0 & 0 \cr 
0 & 0 & 0 & \ldots & -r_0^{(val(V)-1)} & l_0^{(val(V)-1)} & 0    \cr
0 & 0 & 0 & \ldots & 0 & -r_0^{(val(V))} & l_0^{(val(V))} \cr
\end{array}
\right)
\label{eq:Peters_twist_matrix_2}
\end{equation}
and
\begin{equation}
T = \left(\begin{array}{llllll}
\bar T^{(1)}\cr
\vdots\cr 
\bar T^{(val(V))}\cr
\end{array}
\right)
\ \ \ \ 
A = \left(\begin{array}{llllll}
\bar A^{(1)}\cr
\vdots\cr
\bar A^{(val(V))}\cr
\end{array}
\right)
\label{eq:Peters_twist_vectors_2}
\end{equation}
Equation ~\ref{eq:Peters_twist_matrix_1} together with the dependency defined by the {\it "Tangent Constraint'} lead to the general Parity Phenomenon. The following two theorems are due to Peters ~\cite{peters_main}.

\begin{theorem} 
[The General Parity Phenomenon]
For an inner vertex $\bar V$, matrix $M$ (see Equation ~\ref{eq:Peters_twist_matrix_1}) is of full rank if and only if $\bar V$ is an odd vertex. Otherwise its rank is equal to $val(V)-1$.
\label{theorem:parity_phenomenon}
\end{theorem}
From the algebraic point of view, the Parity Phenomenon means that in case of an even vertex some linear combination of the right sides of equations should be equal to zero. In the mesh interpolation problem the only free variables of the right-side expressions are coefficients of the weight functions (more precisely - coefficients with index $1$ and scaling factor for the coefficients with index $0$). If these coefficients are fixed in advance, then the Parity Phenomenon implies that the given mesh should satisfy one additional constraint for every even vertex.
In general, the existence of $G^1$-smooth interpolant depends on the solvability of a quite complicated equation in terms of coefficients of the weight functions. It appears that in order to conclude that a given mesh of curves is {\it admissible}, the complicated equation should not be solved explicitly. Instead of it, one can verify whether the mesh satisfies the sufficient conditions formulated below.

\begin{theorem}
[Sufficient conditions for the vertex enclosure constraint]

If at every inner even vertex $\bar V$ of a given mesh of curves either of the following holds
\begin{itemize}
\item
(Colinearity Condition)
The vertex is a $4$-vertex and all odd-numbered and all even-numbered curves emanating from $V$ have colinear tangent vectors (that is the tangent vectors form an 'X').
\item
(Sufficiency of $C^2$ Data) 
The mesh curves emanating from $\bar V$ are compatible with a second fundamental form at $\bar V$.
\end{itemize}
then the mesh is admissible.
\label{th:suff_vert_enclosure_P}
\end{theorem}

\eject
\part{General linearisation method}
\label{part:general_linearisation}
\vspace{-0.1in}

This Part shows that  the functional space $\FUN{n}$ (see Definition ~\ref{def:space_fun}) reduces  the problem formulated in paragraph ~\ref{def:general_problem_def} to the solution of some {\it linear} constrained minimisation problem.
Later the general method for linearisation will be applied to cases when a planar mesh has a polygonal (Part ~\ref{part:linear_boundary}) or piecewise-cubic $G^1$-smooth global boundary (Part ~\ref{part:smooth_boundary}).

In addition this Part provides some important definitions and notations related to the general flow of the solution and to the analysis of the MDS for the different kinds of "additional" constraints.

\section{Linearisation of the minimisation problem}
\label{sect:Linearisation}

\subsection{Linearisation of the smoothness condition}
\label{subsect:linearisation_G1_general}

Properties of a global regular in-plane parametrisation (see Definition ~\ref{def:reg_param_all}) play a principal role in the following discussion. We now complete this definitionn by :

\begin{definition}
{\bf (Regular parametrisation of a mesh element)}\\
In-plane parametrisation $\T P(u,v)=(P_{X}(u,v),P_{Y}(u,v))$, $(u,v)\in[0,1]^2$ of a single mesh element $\T p\in\T{\cal Q}$ is called \underline{regular} if and only if
\begin{itemize}
\item[{\bf (1)}]
$\T P(u,v)$ is a bijective mapping between the unit square and $\T p$ 
\item[{\bf (2)}]
$\T P(u,v)$ is at least $C^1$-smooth
\item[{\bf (3)}]
The Jacobian $J^{(\T P)}(u,v)$ of the mapping has no singular points: for every $(u,v)\in[0,1]^2$
\vspace{-0.05in}
\begin{equation}
det(J^{(\T P)}(u,v)) = 
det\left(\begin{array}{cc} 
\D{P_{X}}{u} & \D{P_{Y}}{u} \cr
\D{P_{X}}{v} & \D{P_{Y}}{v} \cr
\end{array}\right)
\neq 0 
\label{eq:reg_param}
\end{equation}
\vspace{-0.05in}
\end{itemize}
\label{def:reg_param}
\end{definition}

\vspace{-0.3in}
\begin{theorem}
Let $\T\Pi$ be a fixed global regular in-plane parametrisation and $\bar\Psi$ a piecewise parametric $3D$ function 
agreeing with $\T\Pi$ (see Definition ~\ref{def:agrees_param}).
Let two adjacent mesh elements be parametrized as shown in Figure ~\ref{fig:fig4}, and let $\bar L=(\T L,L)$ and $\bar R=(\T R,R)$ denote the restriction of $\bar\Psi$ on the elements. Then the  $G^1$-continuity of $\T\Pi$ along the common edge is equivalent to the satisfaction of the equation
\begin{equation}
L_u(v)l(v)+R_u(v)r(v)+L_v(v)c(v)=0
\label{eq:G1_Z_equation}
\end{equation}
where 
\vspace{-0.05in}
\begin{equation}
\begin{array}{lll}
c(v)=<\tilde L_u, \tilde R_u>, &
l(v)=<\tilde R_u, \tilde L_v>, &
r(v)=-<\tilde L_u, \tilde L_v> 
\end{array}
\label{eq:weights_general}
\end{equation}
\vspace{-0.05in}
are fixed scalar-valued functions and $L_u(v), R_u(v), L_v(v)$ are the $Z$-components of the partial derivatives along the common edge.

Functions $c(v)$, $l(v)$, $r(v)$ computed according to Equation 
~\ref{eq:weights_general} will be called \underline{conventional weight functions}.
\label{theorem:G1_equiv_Z_equation}
\end{theorem}
{\bf Proof} See Appendix, Section ~\ref{sect:proofs}.

%

Lemma ~\ref{lemma:jorg} and Theorem ~\ref{theorem:G1_equiv_Z_equation} imply that the requirement of $C^1$-smoothness in definition of space $\FUN{n}$ (Definition ~\ref{def:space_fun}, Item {\bf (3)}) may be substituted by the weaker requirement of $G^1$-smoothness, and, furthermore, by the requirement that Equation ~\ref{eq:G1_Z_equation} is satisfied for every inner edge. Thus the smoothness conditions can be studied in terms of the $Z$-components of the B\'ezier control points of the adjacent patches. Moreover, we have  the following Lemma .
 
\begin{lemma}
Let $\T\Pi\in\PAR{m}$ and $\bar\Psi\in\FUN{n}(\T\Pi)$ for $n>m$.
Then the sum in Equation ~\ref{eq:G1_Z_equation} is a B\'ezier polynomial of (formal) degree $n+2m-1$; its coefficients are linear functions in terms of $Z$-components of the control points. In order to satisfy the $G^1$-continuity condition, it is sufficient to impose $NumEqFormal = n+2m$ linear equations of the form :

\begin{equation}
\hspace{-0.4in}
\sum\limits_{\lmtT{j+k=s}{0\le j\le n}{0\le k\le 2m-1}}
\!\!\!\!\!\!\!\!\!
\Cnk{n}{j}\Cnk{2m-1}{k} (l_k\Delta L_j+r_k\Delta R_j)+
\!\!\!\!\!\!\!\!\!
\sum\limits_{\lmtT{j+k=s}{0\le j\le n-1}{0\le k\le 2m}}
\!\!\!\!\!\!\!\!\!
\Cnk{n-1}{j}\Cnk{2m}{k} c_k\Delta C_j=0
\label{eq:G1_Z_general}
\end{equation}

\vspace{-0.1in}

where $s=0,\ldots,n+2m-1$. Here $\Dlt L_j$, $\Dlt R_j$ ($j=0,\ldots,n$) and $\Dlt C_j$ ($j=0,\ldots,n-1$) are first order differences between $Z$-components of the control points, defined in Subsection ~\ref{subsect:def_Bezier_cp_two_adjacent_patches}.

The system may contain redundant equations. The number of necessary and sufficient equations is given by
\begin{equation} 
NumEqActual = n+max\{max\_deg(l,r)+1,\ \ deg(c)\}
\label{eq:num_eq_sufficient}
\end{equation} 
\label{lemma:linear_1}
\end{lemma}
{\bf Proof} See Appendix, Section ~\ref{sect:proofs}.

\begin{definition}
The following notations will be used
\begin{itemize}
\item
\underline{Indexed equation "Eq(s)"} - the equation which follows from equality to zero of the B\'ezier coefficient with index $s$
\item
\underline{NumInd} - the number of the indexed equations (for one edge). Although in most situations $NumInd$ is equal to the total number of equations, the equality should not be necessarily satisfied, because one may like to separate some equations with a special meaning from the homogeneous system of the indexed equations. 
\item
\underline{{\it "Eq(s)"}-type equations} - pair of equations {\it "Eq(s)"} and {\it "Eq(NumInd-1-s)"} ($s=0,\ldots,\lceil(NumInd-1)/2\rceil$).
\end{itemize}
\label{def:indexed_equation}
\end{definition}
Definition of {\it "Eq(s)"}-type equations clearly makes sense in cases when equations {\it "Eq(s)"} and {\it "Eq(NumInd-s)"} are symmetric (for example, in the case when in-plane parametrisations of the elements adjacent to an edge have a symmetric form and both vertices of the edge are inner). The formal definition of {\it "Eq(s)"}-type equations in non-symmetric cases will also be used.

Of course, the linear systems should be considered for all inner edges. The global system may have a sufficiently complicated structure: some control points participate in $G^1$-continuity equations for more than one edge. In addition, a relatively high degree of in-plane parametrisation results in high degrees of the weight functions and increases both the total number of equations and the complexity of each equation. 

A study of the global system is equivalent to a study of the MDS, since it is composed of the control points, which correspond to free variables of the linear system. In addition to an analysis of the dimensionality, one should clearly be careful of a "uniform distribution" of the basic control points. Dimensionality and structure of the MDS for concrete choices of in-plane parametrisation will be analysed in detail in Parts ~\ref{part:linear_boundary} and ~\ref{part:linear_boundary}. Special attention will be paid to a study of the geometrical meaning of equations involved in the linear system.

\subsection{Linear form of "additional" constraints}
\label{subsect:linear_additional}

The current Subsection describes the commonly used types of interpolation and boundary constraints and shows which control points become fixed as a result of application of the constraint. The relation between MDS and the chosen type of "additional" constraints will be explained in greater detail in Subsection ~\ref{subsect:relation_mds_additional}.

\subsubsection{Interpolation constraints}
\label{subsect:linear_interp}
The following constraints (the first one and optionally the second or/and the third ones) are usually applied in the case of an interpolation problem.

\noindent
{\bf (Vertex)-interpolation.} The resulting surface should pass through the given $3D$ point at every mesh vertex. 
The constraint involves $V$-type control points (see Subsection ~\ref{subsect:notation_vertex_adjacent_cp} for definition) of the mesh vertices.

\noindent
{\bf (Tangent plane)-interpolation.} The normal of the tangent plane at every $3D$ vertex should have a specified direction. At every mesh vertex, in addition to $V$-type control point, the constraint involves $E$-type control points. The constraint automatically fits the requirement of $G^1$-smoothness. The assignment values for two $E$-type control points of two non-colinear edges at every vertex is sufficient. 

\noindent
{\bf (Boundary curve)-interpolation.} The resulting surface should interpolate a given $3D$ curve along the global boundary. The constraints result in assignment values for all control points lying on the global boundary of the mesh.

The following notations will be used in order to specify the kind of interpolation problem: round brackets $()$ mean that the type of interpolation is applied, square brackets $[]$ mean that the type of interpolation is optional. For example, the (vertex)[tangent plane]-interpolation problem means that the resulting surface should pass through given $3D$ points at vertices and in addition the normals of the tangent planes at vertices might be specified.

\subsubsection{Boundary conditions}
\label{subsect:linearisation_boundary_conditions}

The following standard boundary conditions are imposed when an approximate solution of some partial differential equation should be found. Here $\bar P=(\T P,P)$ is the restriction of the resulting  function to some boundary mesh element $\T p\in \T{\cal Q}$

\begin{figure}[!pb]
\vspace{-0.5in}
\centering
\includegraphics[clip,totalheight=1.5in]{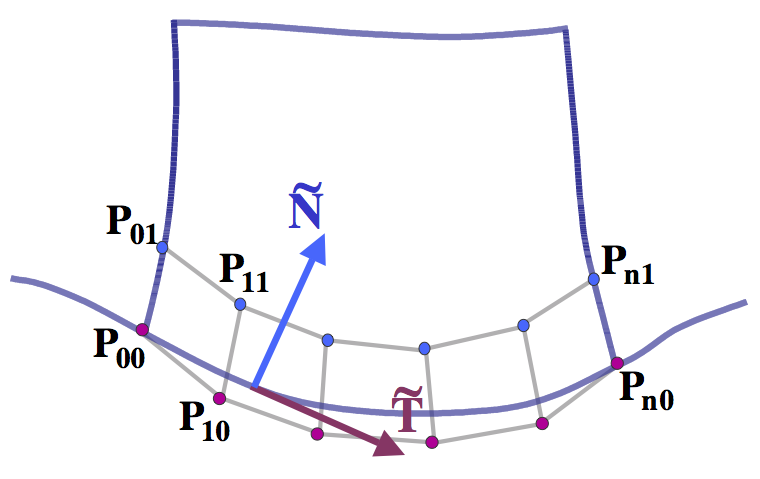}
\caption{Control points involved into the boundary conditions.}
\label{fig:fig9}
\end{figure}

\vspace{-0.1in}
\paragraph{Simply-supported boundary condition.} The standard simply-supported boundary constraint implies that
\vspace{-0.05in}
\begin{equation}
P(boundary)=0
\label{eq:simply_supported_Z}
\end{equation}
\vspace{-0.05in}
should be explicitly fixed. Let $\bar P(u,v)$ be  a patch, such  that $\T P(u,0)$ lies on the global boundary of domain $\Omega$ (see Figure ~\ref{fig:fig9}),then the simply-supported boundary condition means that 
\vspace{-0.05in}
\begin{equation}
P_{i0}=0\ \ {\rm for\ every\ } i=0,\ldots n
\label{eq:simply_supported_cp}
\end{equation}

\vspace{-0.1in}
\paragraph{Clamped boundary condition.} The standard clamped boundary condition means that  
\vspace{-0.1in}
\begin{equation}
P(boundary)=0\ \ {\rm and}\ \ \D{P}{\T{N}}(boundary)=0
\label{eq:clamped_Z}
\end{equation}
\vspace{-0.1in}
where $\T{N}$ is the unit planar normal to the boundary of the domain. Let patch $\bar P(u,v)$ have a regular in-plane parametrisation $\T P(u,v)$. Then
\vspace{-0.05in}
\begin{equation}
\begin{array}{ll}
\D{P}{\T{N}}= 
\D{P}{v}\left(\D{v}{P_{X}} N_{X}+ \D{v}{P_{Y}} N_{Y}\right)=
\D{P}{v}\frac{||\T P_u||}{det(J^{\T P})}\cr
\end{array}
\end{equation}
\vspace{-0.05in}
where $||\T P_u(boundary)||\neq 0$, because otherwise $\T N(boundary)$ is not correctly defined. Condition $\D{P}{v}(boundary)=0$ together with condition $P(boundary)=0$ imply equality to zero of the $Z$-components of two rows of the boundary control points 
(see Figure ~\ref{fig:fig9})
\vspace{-0.05in}
\begin{equation}
P_{ij}=0\ \ {\rm for\ every\ } i=0,\ldots n,\ \ j=0,1 
\end{equation} 

Non homogeneous boundary conditions can also be treated given functions $F$ and $G$, consider  $P(boundary)=F$ and $\D{P}{\T{N}}(boundary)=G$ .

If functions $F$ and $G$ are represented or approximated by
B\'ezier parametric polynomials of some degrees (for example degree of $F$ should be less or equal to the chosen degree of polynomial for $Z$-component of the patch), then this general condition does not lead to additional complications. For the current approach it is important that the control points which are defining the boundary conditions be free from dependencies which follow from $G^1$-smoothness conditions. A simply-supported boundary condition affects the control points along the global boundary of domain; a clamped boundary condition affects the control points along the global boundary of the domain and the control points adjacent to the boundary. 

\subsection{Quadratic form of the energy functional}

Let $\T\Pi$ be a fixed global regular parametrisation, $\bar\Psi\in\FUN{n}(\T\Pi)$ and $\bar P=(\T P,P)=\bar\Psi|_{\T p}$ be the restriction of $\bar\Psi$ on some mesh element $\T p\in\T{\cal Q}$.  The in-plane parametrisation $\T P=\T\Pi|_{\T p}$ is fixed, therefore all partial derivatives of $Z$ with respect to $X$ and $Y$ become linear in terms of $Z$-components of B\'ezier control points. Any energy functional defined as the integral of some quadratic expression of the partial derivatives has a quadratic form in terms of $Z$-components of the control points, hence we have a quadratic  minimisation problem. 

The expression for the energy functional depends on the chosen in-plane parametrisation of the mesh element and may be different for different elements. Although the basic kind of parametrisation considered in the present work (the bilinear parametrisation) leads to very simple formulas, a separate computation is generally required for every mesh element.  

Section ~\ref{sect:energy_example_bilinear} (see Appendix) presents an example of computation of energy functional in the case of bilinear in-plane parametrisation.

\FloatBarrier

\eject 
\section{Principles of construction of MDS}
\label{sect:general_flow}

\subsection{Special subsets of control points and their dimensionality}

\begin{figure}[!ph]
\vspace{-0.1in}
\begin{narrow}{-0.5in}{-0.5in}
\begin{minipage}[]{0.33\linewidth}
\subfigure[]
{
\includegraphics[clip,totalheight=2in]{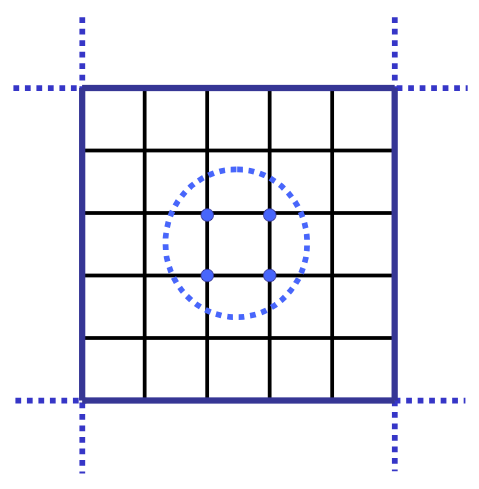}
\label{fig:fig12a}
}
\end{minipage}
\begin{minipage}[]{0.33\linewidth}
\subfigure[]
{
\includegraphics[clip,totalheight=2in]{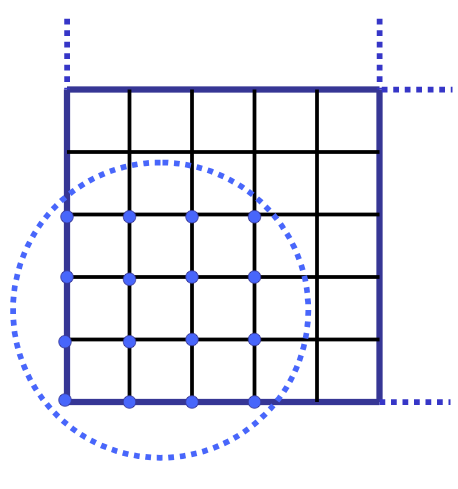}
\label{fig:fig12b}
}
\end{minipage}
\begin{minipage}[]{0.33\linewidth}
\subfigure[]
{
\includegraphics[clip,totalheight=2in]{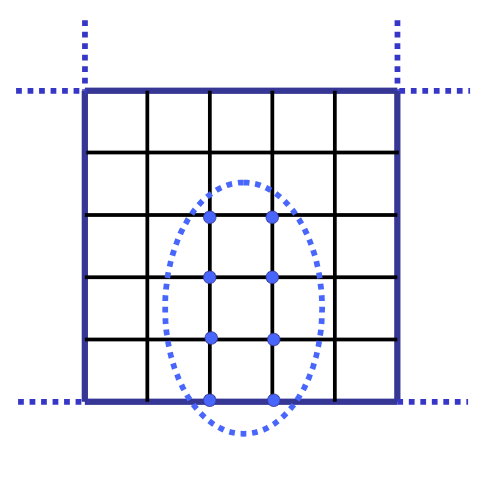}
\label{fig:fig12c}
}
\end{minipage}
\end{narrow}
\vspace{-0.2in}
\caption{Control points, which do not participate in $G^1$-continuity conditions (a) Inner element (b) Corner element (c) Boundary non-corner element.}
\label{fig:fig12}
\end{figure}

\begin{figure}[!pb]
\begin{narrow}{-0.4in}{-0.4in}
\centering
\includegraphics[clip,height=3in]{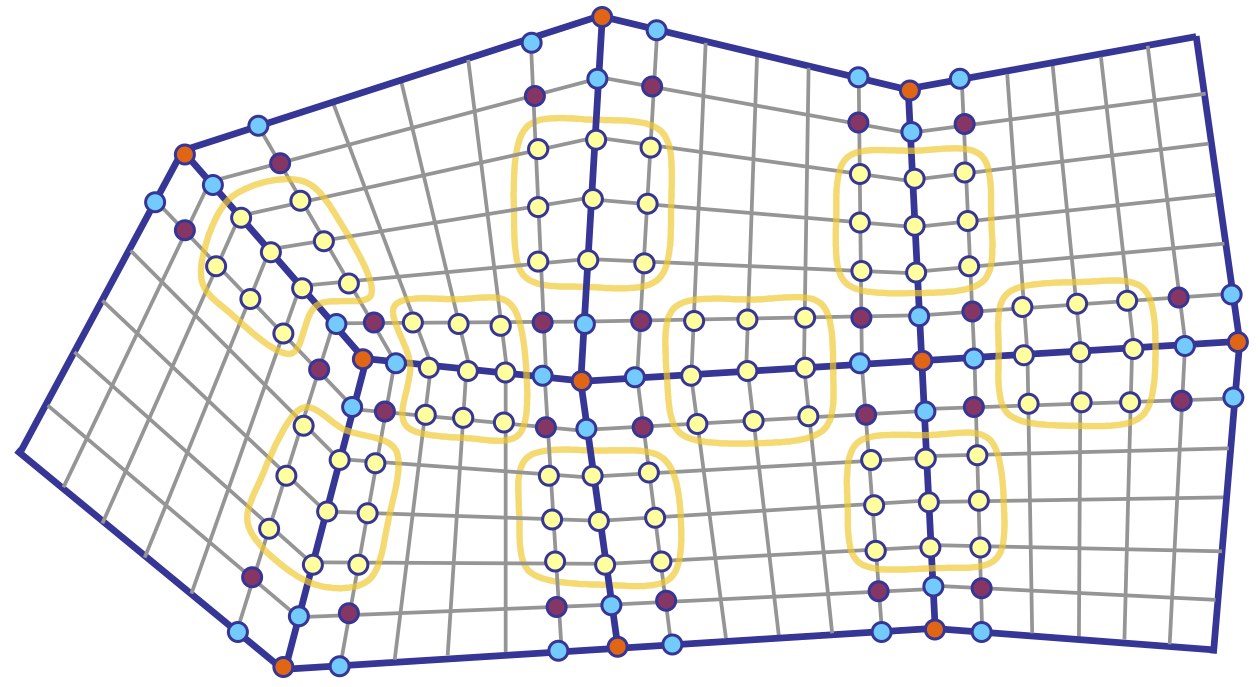}
\end{narrow}
\caption{Control points of a global in-plane parametrisation.}
\label{fig:fig38}
\end{figure}

To clarify the discussion, and  restrict their number  that should be analysed we introduce some subsets of control points.
\begin{definition}
Let (see Figures ~\ref{fig:fig12} and  ~\ref{fig:fig38})
\begin{description}
\item
\underline{$\FCP{n}$} : subset of in-plane control points $\CP{n}$ which are not involved in the $G^1$-continuity conditions; 
\item
\underline{$\GCP{n}$ ($\GMDS{n}$)} : subset of in-plane control points $\CP{n}$ (\emph{resp}. minimal determining set $\MDS{n}$) which participate in $G^1$-continuity conditions.
\end{description}
\label{def:cp_subsets}
\end{definition}
The set $\FCP{n}$  of control points (in other words, all control points which do not lie at some inner edge or adjacent to it) clearly belong to any determining set. 

\begin{lemma}
Dimensions of the subsets $\GCP{n}$ and $\FCP{n}$ are given by the following formulas (here $|\ \ |$ denotes dimension of a set)
\begin{equation}
\hspace{-0.35in}
\begin{array}{ll}
|\FCP{n}|=&
(n\!-\!3)^2|Face_{inner}|+
(n\!-\!3)(n\!-\!1)|Face_{\twolines{boundary}{non\!-\!corner}}|+\cr
&(n\!-\!1)^2|Face_{corner}|\cr
|\GCP{n}|=&
|Vert_{non-corner}|+
3|Vert_{\twolines{boundary}{non\!-\!corner}}|+
(3n\!-\!5)|Edge_{inner}|\cr
\end{array}
\end{equation}
The following relations take place
\begin{equation}
\hspace{-0.24in}
|\CP{n}|=|\GCP{n}|+|\FCP{n}|,\ \ \ \ \ \ \ \ \ \ \ \ \ \ \ \  
|\MDS{n}|=|\GMDS{n}|+|\FCP{n}|
\end{equation}
\label{lemma:GMDS_importance}
\end{lemma}
The important conclusion from Lemma ~\ref{lemma:GMDS_importance} is that we only need to study the structure and dimensionality of $\GMDS{n}$ - {\it subset of the minimal determining set which participates in $G^1$-continuity conditions}.

\subsection{Relation between MDS and the "additional" constraints}
\label{subsect:relation_mds_additional}

The definition of  paragraph ~\ref{def:general_problem_def}, implies that any "additional" constraints is assumed to be consistent and to fit the $G^1$-continuity requirements.

An "additional" constraint result in some definite control points being fixed. These control points get their values according to the "additional" constraints and can not influence the satisfaction of $G^1$ continuity conditions.

\begin{definition}
The minimal determining set $\MDS{n}$ is said to \underline{fit} a given  "additional" constraint if any control point which should be fixed according to this "additional" constraint either
\begin{itemize}
\item
Belongs to $\MDS{n}$ or
\item
Does not belong to $\MDS{n}$ but depends only on the control points which belong to $\MDS{n}$
\end{itemize}
An MDS is called \underline{"pure"} if it is constructed according to $G^1$-conditions alone and is not required to fit any specific  "additional" constraint. 
\label{def:mds_fits_additional_constraint}
\end{definition}

\subsection{Principle of locality in construction of MDS}
\label{subsect:locality_principle}
As it was mentioned above, MDS is not uniquely defined. According to the current approach, construction of the MDS will be built up gradually and will follow two (closely connected) kinds of locality concepts.

At every step of MDS construction, some subset of the linear equation will be considered. The first principle of locality requires that the subset includes indexed equations (see Definition ~\ref{def:indexed_equation}) with successive indices. In addition, the analysis starts from the application of the equations locally, for example, to edges sharing some common vertex or to control points participating in $G^1$-continuity equation for a given edge. 
Control points, which get their status (basic or dependent) during the construction step, clearly obey the principle of  geometrical locality. The local set of the control points which are classified according to the local application of some set of equations, is called \underline{\it "local template"} of MDS (see Figure ~\ref{fig:fig39}).

As soon as the local analysis is completed, one should define the order in which the local templates should be constructed and take care to put together the local templates without contradiction (different local templates may intersect!).

\subsection{Aim of the classification process}
Construction of the MDS implies assignment of the definite status to every one of the control points. Control points, which belong to the MDS, are \underline{\it basic} control points and the rest are \underline{\it dependent} control points. For a given "additional" constraint, the basic control points which get their values according to the "additional" constraints are called \underline{\it basic fixed}, the remaining basic control points are \underline{\it basic free}.

Classification process by definition includes
\begin{itemize}
\item[-]
Construction of the minimal determining set MDS (or several instances of the MDS).
\item[-]
Description of the dependency of every one of the dependent control points on the basic ones. (More precisely, dependency of $Z$-component corresponding to the dependent control point on $Z$-components corresponding to the basic control points).
\item[-]
For a given "additional" constraint, the choice of the instance of MDS that fits the constraint. 
\end{itemize}
It is important to note, that although usually several different configurations of MDS are considered, construction of MDS follows some definite principles (see Subsection ~\ref{subsect:locality_principle}) and of course does not cover all possible configurations.
According to the current approach, in case none of the constructed instances of $\MDS{n}$ fits some "additional" constraint, MDS of a higher degree will be considered. However, failure to choose a suitable instance of MDS does not necessarily imply that a "pure" algebraic solution of the constrained linear system does not exist in space $\FUN{n}$. For example, it will always be assumed that any $V$-type control point belongs to MDS, while an algebraic solution may use such a control point as a dependent one in non-interpolating problems.

In order to make the discussion precise, the following definition of the {\it stages} of the classification process is introduced.
\begin{definition}
A \underline{Stage} is usually a large part of the classification process, which is defined by some set of equations 
and so that at the end of the stage:
\begin{itemize}
\item[{\bf (1)}]
All control points which participate (or may participate under definite geometrical conditions) in these equations 
are classified (as basic or dependent ones) and the status of every one of the control points is final, it can not be changed during the next stages of the classification.  
\item[{\bf (2)}]
Any dependent control point depends only on the control points with the final basic classification status. 
\item[{\bf (3)}]
All these equations 
 are satisfied by classification of the control points. 
\end{itemize}
\label{def:stage_step}
\end{definition}

\section{From MDS to solution of the linear minimisation problem}

As soon as for a given "additional" constraint, a suitable MDS is constructed, dependencies of the dependent control points are defined and energy for every mesh element is computed, construction of the solution of the linear minimisation problem is made in straightforward algebraic manner (see Appendix, Section
~\ref{sect:mds_to_algebraic_solution} for more details).

\FloatBarrier

\eject
\part{MDS for a quadrilateral mesh with a polygonal global boundary}
\label{part:linear_boundary}

\section{Mesh limitations}
\label{sect:mesh_limitations_bilinear}
\begin{figure}[!pb]
\vspace{-0.7in}
\begin{narrow}{-0.5in}{-0.5in}
\begin{minipage}[]{0.25\linewidth}
\subfigure[]
{
\includegraphics[clip,totalheight=1.5in]{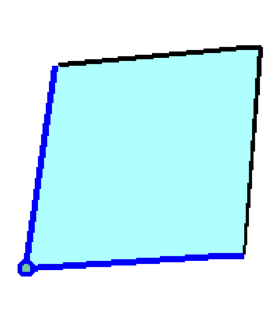}
\label{fig:fig14a}
}
\end{minipage}
\begin{minipage}[]{0.35\linewidth}
\subfigure[]
{
\includegraphics[clip,totalheight=1.5in]{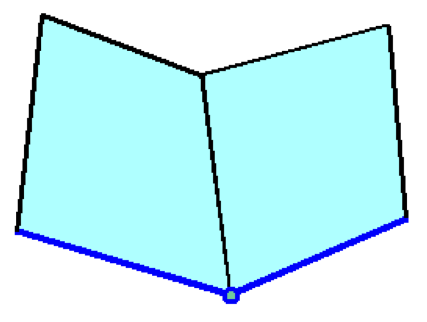}
\label{fig:fig14b}
}
\end{minipage}
\begin{minipage}[]{0.35\linewidth}
\subfigure[]
{
\includegraphics[clip,totalheight=1.5in]{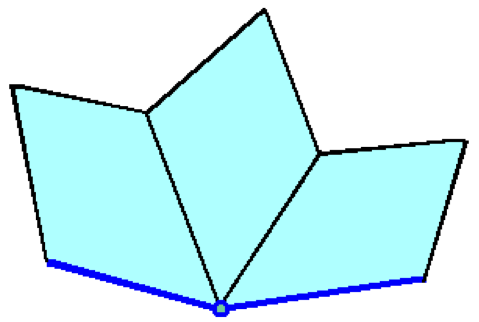}
\label{fig:fig14c}
}
\end{minipage}
\end{narrow}
\vspace{-0.2in}
\caption{An illustration for the mesh limitations.}
\label{fig:fig14}
\end{figure}

The following minor mesh limitations are always supposed to be satisfied
\begin{itemize}
\item[-]
The mesh consists of {\it strictly convex} quadrilaterals. Every mesh element is a convex quadrilateral and angle between any two sequential edges is strictly less
than $\pi$.
\item[-]
Boundary vertices have valence $2$ (a corner vertex) or $3$ (see Figures ~\ref{fig:fig14a},
~\ref{fig:fig14b}). The situation shown in Figure ~\ref{fig:fig14c} is not allowed.
\item[-]
Any inner edge has at most one boundary vertex.
\end{itemize}
The limitations are naturally satisfied in most of the practical situations. In addition, a standard technique of necklacing (see ~\cite{necklacing}) may be applied to a mesh in order to achieve the second and the third requirement.

\begin{figure}[!pt]
\centering
\includegraphics[clip,totalheight=1.4in]{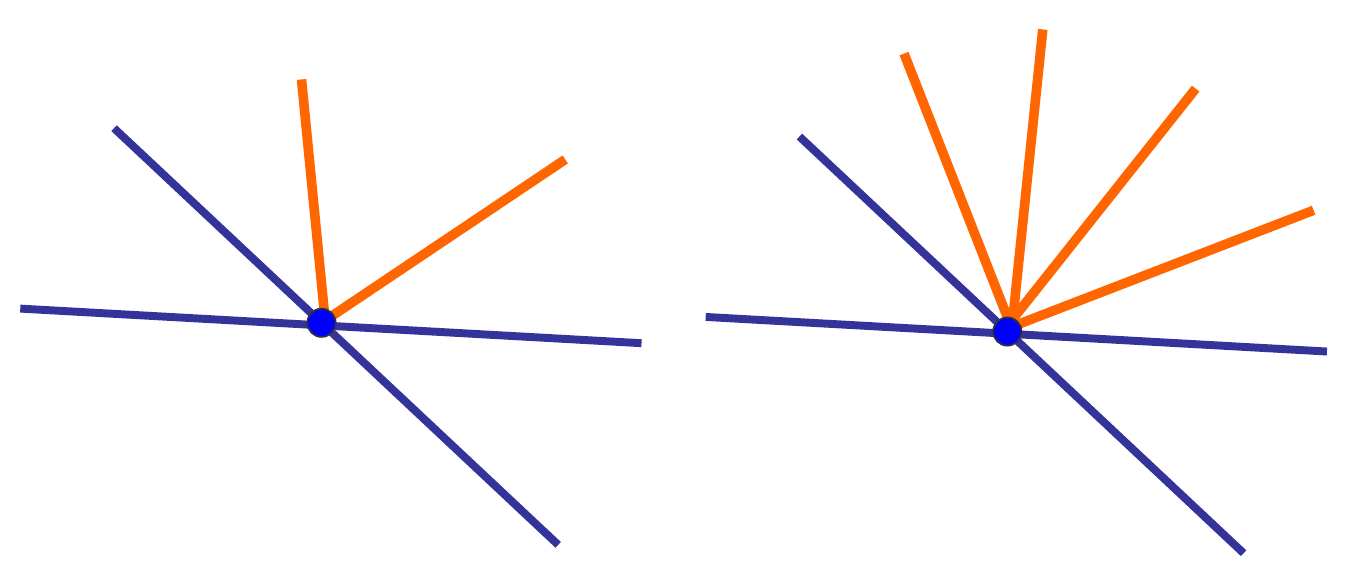}
\vspace{-0.1in}
\caption{Examples of meshes, which do not satisfy the {\it "Uniform Edge Distribution Condition"}.}
\label{fig:fig23}
\end{figure}

It is convenient to introduce an additional minor mesh limitation, which is required to be satisfied only when MDS of degree $n=4$ is considered. In this case a planar mesh should satisfy the "Uniform Edge Distribution Condition", defined as follows
\begin{definition}
The mesh is said to satisfy the {\it "Uniform Edge Distribution Condition"} if for any even vertex
of degree $\ge 6$ which has two pairs of colinear edges, the remaining edges
($2$ edges in case of a $6$-vertex, $4$ edges in case of a $8$-vertex and so on) do
not all belong to the same quadrant formed by lines containing the colinear edges
(see Figure ~\ref{fig:fig23}).
\label{def:uniform_edge_distribution}
\end{definition}

\FloatBarrier

\section{In-plane parametrisation}
\label{sect:parametrisation_bilinear}

For a quadrilateral planar mesh element, the bilinear in-plane parametrisation will be considered. It is a natural choice because for a general quadrilateral (without curvilinear sides) it clearly provides a parametrisation of the minimal possible degree, which finally leads to the minimal possible number of $G^1$-continuity equations.

For a convex quadrilateral planar mesh element with vertices $\T{A},\T{B},\T{C},\T{D}$ (see Figure ~\ref{fig:fig8}) the parametrisation is given by explicit formula 
\begin{equation}
\T{P}(u,v) = 
\T{A}(1-u)(1-v)+\T{B}u(1-v)+\T{C}uv+\T{D}(1-u)v
\label{eq:parametrisation_bilinear}
\end{equation}
\vspace{-0.05in}
and   $det(J^{(\T P)} (u,v)>0$ for every $(u,v)\in[0,1]^2$ , hence the following Lemma holds.

\begin{lemma}
For a strictly convex planar quadrilateral element, 
the bilinear in-plane parametrisation of the element is regular.
\label{lemma:reg_bilinear}
\end{lemma}
\vspace{-0.1in}

The bilinear parametrisations for all mesh elements clearly satisfies the requirements of Definition ~\ref{def:reg_param_all} and define a degree $1$ global regular in-plane parametrisation $\T\Pi^{(bilinear)}\in\PAR{1}$.

Explicit formulas for in-plane control points which belong to $\GCP{n}(\T\Pi^{(bilinear)})$ are given in Technical Lemma ~\ref{tl:degree_elevation_bilinear} (see Appendix, Section ~\ref{sect:technical_lemmas}).

\section{Conventional weight functions and linear form of $G^1$-continuity conditions}
\label{sect:weight_functions_and_linear_system_bilinear}

Application of the general linearisation method (Theorem ~\ref{theorem:G1_equiv_Z_equation} and Lemma ~\ref{lemma:linear_1}) to a particular case of bilinear in-plane parametrisation leads to the next Lemma.
\begin{lemma}
Let two adjacent mesh elements with vertices $\LL,\LL',\GG,\GG',\RR,\RR'$ (see Figure ~\ref{fig:fig11}) each with a bilinear in-plane parametrisation, then
\begin{itemize}
\item[{\bf (1)}]
The conventional weight functions $l(v),r(v)$ and $c(v)$ along the common edge are B\'ezier polynomials of (formal) degrees $1$,$1$ and $2$ respectively and their coefficients with respect to the B\'ezier basis depend on the geometry of the planar elements in the following way 
\begin{equation}
\hspace{-0.3in}
\begin{array}{ll}
l_0=\<\RR-\GG,\GG'-\GG\>     
& c_0=\<\RR-\GG,\LL-\GG\>\cr
l_1=\<\GG-\GG',\RR'-\GG'\>   
& c_1=\frac{1}{2}
\left[\<\RR-\GG,\LL'-\GG'\>-\<\LL-\GG,\RR'-\GG'\>\right]\cr
r_0=-\<\GG'-\GG,\LL-\GG\> & 
c_2=-\<\LL'-\GG',\RR'-\GG'\> \cr
r_1=-\<\LL'-\GG',\GG-\GG'\>
\end{array}
\label{eq:weights_bilinear}
\end{equation}
\item[{\bf (2)}]
The system of $n+2$ linear equations {\it "Eq(s)"} for $s=0,\ldots,n+1$ is sufficient in order to satisfy the $G^1$- continuity condition along the common edge 
\begin{equation}
\hspace{-0.35in}
{\it "Eq(s)"}\ \ 
\begin{array}{l}
(n+1-s)(l_0\Delta L_s+r_0\Delta R_s)+
s(l_1\Delta L_{s-1}+r_1\Delta R_{s-1})+\cr
\frac{(n-s)(n+1-s)}{n}c_0\Delta C_s\!+\!
\frac{2s(n+1-s)}{n}c_1 \Delta C_{s-1}\!+\!
\frac{s(s-1)}{n}c_2 \Delta C_{s-2}\!=\! 0
\end{array}
\label{eq:G1_general_bilinear}
\end{equation}
Here $\Delta L_j$, $\Delta R_j$ for $j<0$ or $j>n$ and 
$\Delta C_j$ for $j<0$ or $j>n-1$ are assumed to be equal to zero.
\end{itemize}
\label{lemma:weights_equations_bilinear}
\end{lemma}
Weight functions $l(v)$, $r(v)$ and $c(v)$ may have lower actual degrees than $1$, $1$ and $2$.
Weight function $l(v)$ becomes constant if $\T t^{(R)}$ is parallel to $\GG'-\GG$; $r(v)$ becomes a constant if $\T t^{(L)}$ is parallel to $\GG'-\GG$  and the actual degree of $c(v)$ is at most $1$ if $\T t^{(R)}$ and $\T t^{(L)}$ are parallel (see Figure ~\ref{fig:fig15} and Subsection ~\ref{subsect:def_vertices_edges_twists_planar} for definition of $\T t^{(R)}$ and $\T t^{(L)}$). For example, in the case of two adjacent square elements $deg(l)=deg(r)=0$, $deg(c)\le 1$. It implies that
$n+1$ linear equations are sufficient in order to guarantee  $G^1$-smooth concatenation and therefore an additional degree of freedom is available.

\section{Local MDS} 
\label{sect:local_mds_bilinear}

As stated in Subsection ~\ref{subsect:locality_principle},
construction of the MDS follows the principle of locality. All control points are subdivided into several types: $V$,$E$,$D$ and $T$-type control points adjacent to some mesh vertex (see Subsection ~\ref{subsect:notation_vertex_adjacent_cp}) and the {\it middle}  control points adjacent to some mesh edge (see Definition ~\ref{def:middle_bilinear}). Every type of the control points is responsible for the satisfaction of some definite subset of the linear equations.

The current Subsection is devoted to an analysis of equations applied to a separate mesh vertex or edge, a possible influence of the other equations is ignored. The analysis results in construction of local MDS, templates, which locally define which control points belong to MDS and describe the dependencies of the dependent control points. The same set of equations may define several structures of the MDS, suitable for the different mesh geometry and different types of "additional" constraints.

\begin{definition}
We  say that two local templates are \underline{different}, if they contain a different number of basic control points or if there is a difference in the types or a principal difference in the location of the control points. 
\end{definition}

\begin{note}
Sometimes, the templates do not uniquely specify which control points should be classified as basic. In case of ambiguity, classification of the control points is made arbitrarily. The local geometric characteristics, such as edge lengths or angles, plays an important role in stabilizing the solution and can be a matter of additional research.
\end{note}

\subsection{Local classification of $E$,$V$-type control points for a separate vertex based on $"Eq(0)"$-type equations}
\label{subsect:local_mds_VE_bilinear}

\subsubsection{Formal equation and geometrical formulation}
\label{subsub: forlema}

Formal substitution of $s=0$ in Equation ~\ref{eq:G1_Z_general} gives
\vspace{-0.07in}
\begin{equation}
"Eq(0)"\ \ \ \ l_0\Delta L_0 + r_0\Delta R_0+c_0\Delta C_0 = 0
\label{eq:s_0_bilinear}
\end{equation}
\vspace{-0.07in}
It is precisely the general {\it "Tangent Constraint"} (Equation ~\ref{eq:Peters_tangent}) applied to $Z$-components of the control points. The difference is that in the curve mesh interpolation problem $3D$ tangents $\Delta \bar L_0$, $\Delta\bar R_0$, $\Delta\bar C_0$ are given and coefficients of the weight functions are unknown. In the current case on the contrary, coefficients of the weight functions are fixed a priory and $Z$-components of the control points play the role of unknowns. 

{\it "Eq(0)"}-type equations have a very simple geometrical meaning.
Let $\T{V}$ be a planar mesh vertex of degree $val(V)$ and $\T e^{(j)}$ ($j=1,\ldots, val(V)$) be a directed planar mesh edges emanating from $\T{V}$ (see Figure ~\ref{fig:fig10}). Then for the edge $\T e^{(j)}$, zero-indexed coefficients of the weight functions can be rewritten as follows
\begin{equation}
\hspace{-0.3in}
l_0^{(j)}=\<\T e^{(j-1)}, \T e^{(j)}\>,\ \ \ 
r_0^{(j)}=-\<\T e^{(j)}, \T e^{(j+1)}\>,\ \ \ 
c_0^{(j)}=\<\T e^{(j-1)}, \T e^{(j+1)}\>
\label{eq:coeff_0_geometrical}
\end{equation}
Note that  $\bar E^{(j-1)}-\bar V$, $\bar E^{(j)}-\bar V$, and $\bar E^{(j+1)}-\bar V$ are colinear
if and only if 
\begin{equation}
\begin{array}{ll}
0=&mix\left(\begin{array}{c}
\bar E^{(j-1)}-\bar V\cr 
\bar E^{(j)}-\bar V\cr 
\bar E^{(j+1)}-\bar V\cr 
\end{array}\right)= 
\frac{1}{n^2}
mix\left(\begin{array}{cc}
\T{e}^{j-1} & E^{(j-1)}-V\cr 
\T{e}^{j} & E^{(j)}-V\cr 
\T{e}^{j+1} & E^{(j+1)}-V\cr 
\end{array}\right)=\cr
&l_0^{(j)} (E^{(j+1)}-V)+c_0^{(j)} (E^{(j)}-V)+r_0^{(j)} (E^{(j-1)}-V)
\end{array}
\end{equation}
which exactly means that the {\it "Eq(0)"}-type equation for the edge with order number $j$ is satisfied,  thus $\bar E^{(j)}- \bar V $ , ($ j=1,\ldots,val(V)$ ) are coplanar ; 
this can be summed up in the lemma:
\begin{lemma}
Let $ \bar V$ be a $3D$ vertex control point and let $\bar E^{(j)}$ ,($j=1,\ldots, val(V)$) be the edge control points adjacent to the vertex (see Subsection ~\ref{subsect:notation_vertex_adjacent_cp}). 
Then for {\it"Eq(0)"}-type equations applied simultaneously to all edges sharing vertex $\bar V$, the tangent vectors $\bar E^{(j)}-\bar V$ ($j=1,\ldots,val(V)$) should be coplanar.   
\label{lemma:tangent_coplanarity}
\end{lemma}
\vspace{-0.1in}

\subsubsection{Degrees of freedom and dependencies}

At every vertex, the tangent plane is defined by any three noncolinear control points lying in it. Let $V$-type control point and such a pair of $E$-type control points $\T E^{(i)}$, $\T E^{(j)}$ ($1\le i,j\le val(V)$), that $\T{e}^{(i)}$ and $\T{e}^{(j)}$ are not colinear, be classified as \emph{basic}.
Any other $E$-type control point $\T E^{(k)}$ is classified as dependent. Its dependency (dependency of the corresponding $Z$-component) is defined by system of {\it "Eq(0)"}-type equations and has the following explicit form
\begin{equation}
\begin{array}{ll}
E^{(k)}= \frac{1}{\<\TTe{i},\TTe{j}\>}
&\left\{
-E^{(i)}\<\TTe{j},\TTe{k}\>
-E^{(j)}\<\TTe{k},\TTe{i}\>+\right.\cr
&\left.V\left(\<\TTe{i},\TTe{j}\>+\<\TTe{j},\TTe{k}\>+\<\TTe{k},\TTe{i}\>\right)
\right\}
\end{array}
\label{eq:tangent_coplanarity_cp}
\end{equation}

\subsubsection{Local templates for a separate inner vertex}
\label{subsect:local_VE_inner_vertex_bilinear}

For any inner vertex, we classify as \emph{basic} any $V$ type control point that has the properties above. 

 The remaining $E$-type control points depend on the basic control points according to Equation ~\ref{eq:tangent_coplanarity_cp}.
The correspondent local template is shown in Figure ~\ref{fig:fig42Aa}. 

This local MDS clearly fits all types of considered "additional" constraints, including the (vertex)(tangent plane)-interpolation condition. The basic control points can be easily classified into free and fixed, depending on the kind of the "additional" constraints. 

\subsubsection{Local templates for a separate boundary vertex}
\label{subsect:local_VE_boundary_vertex_bilinear}

According to the mesh limitations, any non-corner boundary vertex $\T V$ has exactly one adjacent inner edge $\TTe{2}$. 

\begin{figure}[!ph]
\vspace{-0.1in}
\begin{narrow}{0.0in}{0.0in}
\begin{minipage}[]{0.32\linewidth}
\subfigure[]
{
\includegraphics[clip,totalheight=1.3in]{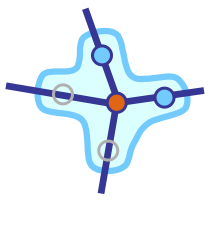}
\label{fig:fig42Aa}
}
\end{minipage}
\begin{minipage}[]{0.32\linewidth}
\subfigure[]
{
\includegraphics[clip,totalheight=1.3in]{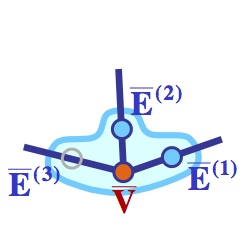}
\label{fig:fig42Ab}
}
\end{minipage}
\begin{minipage}[]{0.32\linewidth}
\subfigure[]
{
\includegraphics[clip,totalheight=1.3in]{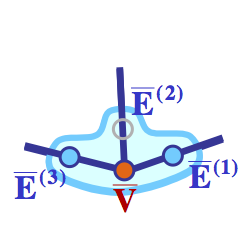}
\label{fig:fig42Ac}
}
\end{minipage}
\end{narrow}
\vspace{-0.2in}
\caption{Local templates for the classification of $V$,$E$-type control points in case of global bilinear in-plane parametrisation $\T\Pi^{(bilinear)}$.}
\label{fig:fig42A}
\end{figure}

The following two local templates are defined:  
\begin{description}
\item[]
\underline{$TB0^{(V,E)}$} (Figure ~\ref{fig:fig42Ab}). The local MDS contains $\T V$, the  boundary control point $\T E^{(1)}$ and the inner control point $\T E^{(2)}$. This template is always used when the boundary edges are colinear.
\item[]
\underline{$TB1^{(V,E)}$} (Figure ~\ref{fig:fig42Ac}). The local MDS contains $\T V$ and two boundary control points $\T E^{(1)}$, $\T E^{(3)}$. This template is always used in the case of  boundary curve interpolation and simply supported boundary conditions, provided the boundary edges are not colinear.
\end{description}
The following two examples show that for the considered "additional" constraint at least one of $TB0^{(V,E)}$, $TB1^{(V,E)}$ provides the local MDS which fits the constraint and
present classification of the basic control points into free and fixed.

{\it (Vertex)(Boundary curve)-interpolation condition}.
If the boundary edges are not colinear, then $TB1^{(V,E)}$ is used. $\T V$, $\T E^{(1)}$ and $\T E^{(3)}$ are basic fixed control points; $\T E^{(2)}$ is dependent, the corresponding $Z$-component is computed according to Equation 
~\ref{eq:tangent_coplanarity_cp}. If the boundary edges are colinear, then $TB0^{(V,E)}$ is used. $\T V$ and $\T E^{(1)}$ are basic fixed control points , $\T E^{(2)}$ is a basic free one. In this case, one should verify that the data of the boundary curve fits the {\it "Tangent Constraint"}: the given value of $E^{(3)}$ should be equal to the value computed according to Equation ~\ref{eq:tangent_coplanarity_cp}, using the given values of $V$ and $E^{(1)}$.

{\it Clamped boundary condition}. It is always possible to make use of $TB0^{(V,E)}$. All basic control points are fixed. The standard clamped boundary constraint clearly satisfies the {\it "Tangent Constraint"}. In case of a more complicated clamped boundary condition, classification of the control points remains unchanged. One should verify that the boundary condition and the {\it "Tangent Constraint"} fit together. Value of $E^{(3)}$ computed according to Equation ~\ref{eq:tangent_coplanarity_cp}, should be equal to the value given by the boundary condition.

\subsection{Local classification of $D$,$T$-type control points for a separate vertex based on $"Eq(1)"$-type equations}
\label{subsect:local_mds_DT_bilinear}

In the current Subsection it will always be assumed that $V$,$E$-type control points are classified and {\it "Eq(0)"}-type equations are satisfied by choice of an appropriate template. 

\subsubsection{Formal equation and geometrical formulation}

Substitution of $s=1$ in Equation ~\ref{eq:G1_general_bilinear} leads to the formula 
\begin{equation}
\begin{array}{ll}
"Eq(1)"\ \ \ \ &
n(l_0\Delta L_1+r_0\Delta R_1)+(l_1\Delta L_0+r_1\Delta R_0)+\cr
&(n-1)c_0\Delta C_1+2 c_1\Delta C_0 = 0
\end{array}
\label{eq:s_1_bilinear}
\end{equation}
This is a particular case of the general {\it "Twist Constraint"} (Equation ~\ref{eq:Peters_twist_1}) applied to $Z$-components of  the control points. An advantage of the current particular case is that the coefficients of the weight functions have a clear geometrical meaning, closely connected to the structure of the initial planar mesh. It allows rewriting {\it "Eq(1)"} in a more meaningful form.


Let $\T{V}$ be a planar mesh vertex of degree $val(V)$ and $\T e^{(j)}$ ($j=1,\ldots, val(V)$) be directed planar mesh edges emanating from $\T{V}$ (see Figure ~\ref{fig:fig10}). 

Let $\bar\Psi\in\FUN{n}(\T\Pi^{(bilinear)})$ and patch $\bar P^{(j)}$ denote the restriction of $\bar\Psi$ on the mesh element adjacent to $\T V$ ,  containing edges $\T e^{(j-1)}$ and $\T e^{(j)}$ (see Figure ~\ref{fig:fig6}) as part of its boundary.
The following important relations between $XY$-components of the first and second-order partial derivatives of the patches and the initial mesh data hold

\begin{equation}
\hspace{-0.3in}
\begin{array}{l}
\bar\epsilon^{(j)}=(\T e^{(j)}, n(E^{(j)}-V))\cr
\bar\tau^{(j)}=(\T t^{(j)}, n^2(D^{(j)}-2E^{(j)}+V))=
(\T t^{(j)}, \frac{n}{n-1}\delta^{(j)})
\end{array}
\label{eq:deriv_planar_connection}
\end{equation}
Here $\bar\epsilon^{(j)}$ and $\bar\tau^{(j)}$, $\delta^{(j)}$ are the first and second order partial derivatives (see Subsection 
~\ref{subsect:def_deriv_patches_at_common_vertex}) and $\T e^{(j)}$, $\T t^{(j)}$ are directed planar edges and twist characteristics of the planar mesh elements (see Subsection ~\ref{subsect:def_vertices_edges_twists_planar}). Relations given in Equation ~\ref{eq:deriv_planar_connection} allow to conclude that the following Lemma holds.


\begin{lemma}
Let all {\it "Eq(0)"}-type equations for all inner edges adjacent to vertex $\T V$ be satisfied by classification of  $V$ and $E$-type control points. 
Then for inner edge $\T e^{(j)}$, {\it "Eq(1)"}-type equation applied to the control points adjacent to vertex $\T V$, has the following geometrical form 
\begin{equation}
tw^{(j)}+tw^{(j-1)} = coeff^{(j)}\delta^{(j)}
\label{eq:eq1_for_one_edge} 
\end{equation}
Where 
\vspace{-0.1in}
\begin{equation}
\hspace{-0.35in}
tw^{(j)}\!=\!\frac{1}{\<\TTe{j},\TTe{j+1}\>^2}
mix\VT{\bar\tau^{(j)}}{\bar\epsilon^{(j)}}{\bar\epsilon^{(j+1)}},\ \ 
coeff^{(j)}\!=\!\frac{\<\TTe{j-1},\TTe{j+1}\>}{\<\TTe{j-1},\TTe{j}\>\<\TTe{j},\TTe{j+1}\>}
\label{eq:tw_delta_def}
\end{equation}
\vspace{-0.2in}
\label{lemma:eq1_geom_edge_j}
\end{lemma}
\vspace{-0.05in}
{\bf Proof} See Appendix, Section ~\ref{sect:proofs}.

It is important to note, that $tw^{(j)}$ and $\delta^{(j)}$ are linear expressions in terms of $V$-type, $E$-type and $T$-type or $D$-type control points; $tw^{(j)}$ contains a {\it single} non-classified $T$-type point $T^{(j)}$ and $\delta^{(j)}$ contains a {\it single} non-classified $D$-type point $D^{(j)}$.

\subsubsection{Theoretical results for an inner vertex}
\label{subsect:theory_vertex_enclosure_bilinear}

\paragraph{The Parity Phenomenon.}

Application of Lemma ~\ref{lemma:eq1_geom_edge_j} to all edges emanating from a common inner vertex leads to the following Theorem.

\begin{theorem}
Let $\T V$ be an inner vertex of degree $val(V)$ and let {\it "Eq(0)"}-type equations for all edges adjacent to the vertex be satisfied. Then 
\begin{itemize}
\item[{\bf (1)}]
The system of {\it "Eq(1)"}-type equations applied simultaneously to all edges adjacent to $\T V$ has the following form
\begin{equation}
M
\left(\begin{array}{c}
tw^{(1)}\cr \vdots\cr tw^{(val(V))}\cr
\end{array}\right) =
\left(\begin{array}{c}
coeff^{(1)}\delta^{(1)}\cr  
\vdots\cr 
coeff^{(val(V))}\delta^{(val(V))} 
\end{array}\right)
\label{eq:parity_phenomenon}
\end{equation}
\vspace{-0.1in}
where $M$ is the matrix with a simple circulant structure
\begin{equation}
M=\left(\begin{array}{ccccccc}
1      & 0      & \ldots & 0      & 0      & 1      \cr
1      & 1      & \ldots & 0      & 0      & 0      \cr
\vdots & \vdots & \ddots & \vdots & \vdots & \vdots \cr
0      & 0      & \ldots & 1      & 1      & 0      \cr
0      & 0      & \ldots & 0      & 1      & 1      \cr
\end{array}\right)
\end{equation}
\item[{\bf (2)}] 
In the case of an even \textcolor{red} {odd} 
vertex matrix $M$ is of full rank.
In the case of an odd  \textcolor{red} {even} vertex $rank(M)=val(V)-1$ and the system has a solution if and only if the following additional condition is satisfied
\vspace{-0.05in}
\begin{equation}
\hspace{-0.35in}
"Circular\ Constraint"\ \ \ \ 
\sum_{j=1}^{val(V)} (-1)^j coeff^{(j)}\delta^{(j)}=0
\label{eq:eq1_condition}
\end{equation}
\vspace{-0.15in}
\end{itemize}
\label{theorem:eq1_all_edges_of_vertex}
\end{theorem}
\vspace{-0.1in}

\begin{note}
The {\it "Circular Constraint"} does not involve $T$-type control points. It establishes some dependency between $D$-type control points adjacent to a given vertex. (Under the assumption that all $V$-type and $E$-type control points are already classified according to the first stage of the classification process).
\end{note}

Results of Theorem ~\ref{theorem:eq1_all_edges_of_vertex} clearly fit the general Parity Phenomenon (Theorem ~\ref{theorem:parity_phenomenon}). The {\it "Circular Constraint"} corresponds to the necessary condition which should be satisfied for the right sides of the general {\it "Twist Constraint"} (see Equation ~\ref{eq:Peters_twist_right_part}) for an even vertex. The main advantage of the current particular case is a very elegant and geometrically meaningful form of the {\it "Circular Constraint"}.

The way in which the {\it "Circular Constraint"} is applied presents the second important difference between the current approach and the standard techniques for interpolation by $3D$ smooth piecewise parametric surface. Usually some initial data ($3D$ mesh of curves in work ~\cite{peters_main}) is tested to satisfy the necessary condition. In case of negative answer, a $G^1$-smooth surface cannot be constructed. In the current approach one may take advantage of the fact that even in the case of (vertex)(tangent plane)-interpolation, a boundary curve of two adjacent patches (which has at least degree $4$) is not totally fixed. At least one control point in the middle of every curve remains non-fixed. The purpose is to construct the MDS in such a manner, that every vertex at which the {\it "Circular Constraint"} should be satisfied, has at least one "own" basic $D$-type control point.

\begin{note}
Coefficient $coeff^{(j)}$ of $\delta^{(j)}$ ($D^{(j)}$) in the {\it "Circular Constraint"} may be equal to zero; it happens if the planar mesh edges $\TTe{j-1}$ and $\TTe{j+1}$ are colinear. In this case $D^{(j)}$ does not contribute to the {\it "Circular Constraint"}.
\end{note}

\begin{definition}[Regular $4$-vertex]
Vertex of valence $4$ is called \underline{4-regular} if $4$ planar edges emanating from the vertex form two colinear pairs: $\TTe{1}$ is colinear to $\TTe{3}$ and 
$\TTe{2}$ is colinear to $\TTe{4}$.
\end{definition}

\begin{lemma}
Regular $4$-vertex is the only possible configuration of the edges adjacent to some inner even vertex when all coefficients $coeff^{(j)}$ $(j=1,\ldots,val(V))$ are equal to zero and the {\it  "Circular Constraint"} is satisfied automatically.
\label{lemma:regular_vertex}
\end{lemma}

\vspace{0.32in}
\noindent
{\bf Proof of the  Lemma}
The strict convexity of the mesh elements implies that any inner even vertex $\T V$ has degree $val(V)=4$ at least.

Let $coeff^{(j)}=0$ for every $j=1,\ldots,val(V)$.
In particular, $coeff^{(2)}=0$ and so $\TTe{1}$ and $\TTe{3}$ are colinear and lie on some straight line $\T l^{(1,3)}$;  $coeff^{(3)}=0$ and so $\TTe{2}$ and $\TTe{4}$ are colinear and lie on some straight line $\T l^{(2,4)}$ (see Figure ~\ref{fig:fig16}).
Therefore for $val(V)=4$ the vertex is proven to be regular.

It remains to show that $val(V)$ could not be greater than $4$. Indeed, let $val(V)>4$. Then, $\TTe{2}$ should be colinear to both $\TTe{4}$ and $\TTe{val(V)}$ (because both $coeff^{(3)}$ and $coeff^{(1)}$ are equal to zero). But $\TTe{4}$ and $\TTe{val(V)}$ can not be colinear because $\TTe{val(V)}$ lies strictly between 
$\TTe{4}$ and $\TTe{deg(1)}$ which span an angle less than $\pi$ due to the strict convexity of the mesh elements.
\eop${}_{{\bf Lemma ~\ref{lemma:regular_vertex}}}$.

\paragraph{Some necessary and sufficient conditions for the satisfaction of the {\it "Circular Constraint"} at a separate inner even vertex}

Results of the current paragraph correspond to the sufficient vertex enclosure conditions formulated in Theorem
~\ref{th:suff_vert_enclosure_P}.
Although the results do not contribute directly to the construction of the MDS, they provide an additional confirmation that the present approach fits the general theory of $G^1$-smooth piecewise parametric surfaces. 
 
Lemma ~\ref{lemma:regular_vertex} from the previous paragraph shows that for an inner vertex of degree $4$ colinearity of two pairs of emanating edges is a sufficient condition for the satisfaction of the {\it "Circular Constraint"}. Necessary and sufficient conditions are presented in the following Lemma.

\begin{lemma}
Let $\T V$ be an inner even vertex
\begin{itemize}
\item[{\bf (1)}]
If $\delta^{(j)}$, $j=1,\ldots,val(V)$
($Z$-components of the second-order derivatives in the directions of the planar edges) are chosen in such a manner that they are compatible at $\T V$ with second-order partial derivatives of some functional surface, then the {\it "Circular Constraint"} is satisfied. 
\item[{\bf (2)}]
For a non-regular $4$-vertex $\T V$, compatibility of $\delta^{(j)}$,  $j=1,\ldots,4$ with second-order partial derivatives of some functional surface is not only a sufficient but also a necessary condition for the satisfaction of the {\it "Circular Constraint"}.
\end{itemize}
\label{lemma:C2_sufficient}
\end{lemma}
\vspace{-0.1in}
{\bf Proof} See Appendix, Section ~\ref{sect:proofs}.

\begin{note}
Lemma ~\ref{lemma:C2_sufficient} does not mean that the resulting surface is necessarily $C^2$-smooth at the vertex.  
Besides values of $\delta^{(j)}$ ($j=1,\ldots,val(V)$) there is always at least one additional degree of freedom ($T$-type control point) which implies that the second-order partial derivatives in the functional sense are not necessarily well defined at the vertex.
\end{note}

\subsubsection{Local templates for a separate inner vertex}
\label{subsect:local_DT_inner_vertex_bilinear}

\paragraph{Odd vertex}
A local template for an inner odd vertex is shown in Figure ~\ref{fig:fig42Ba}. All $D$-type control points are classified as basic and all $T$-type control points are dependent. There are  $val(V)$ basic control points in all.

The correctness of the classification and dependencies of $T$-type control points are explained below.

As stated in Theorem ~\ref{theorem:eq1_all_edges_of_vertex}, matrix $M$ is invertible for an odd inner vertex. Therefore all $D$-type control points can be classified as basic and $T$-type control points depend on them (and on $V$-type and $E$-type basic control points) according to equation
\begin{equation}
\left(\begin{array}{c}
tw^{(1)}\cr \vdots\cr tw^{(val(V))}\cr
\end{array}\right) = M^{-1}
\left(\begin{array}{c}
coeff^{(1)}\delta^{(1)}\cr  
\vdots\cr  
coeff^{(val(V))}\delta^{(val(V))} 
\end{array}\right)
\label{eq:eq_T_from_D_odd}
\end{equation}
In greater detail, $D$-type control points together with $V$-type and $E$-type basic control points fully define values of $\delta^{(j)}$ ($j=1,\ldots,val(V)$). Equation ~\ref{eq:eq_T_from_D_odd} defines dependency of $tw^{(j)}$ ($j=1,\ldots,val(V)$) on
$\delta^{(j)}$ ($j=1,\ldots,val(V)$), and finally the values of $T$-type control points are given by
\begin{equation}
\begin{array}{ll}
T^{(j)}=&\frac{1}{n^2}\<\TTe{j},\TTe{j+1}\>tw^{(j)}-
V\left(1+\frac{1}{n}\frac{\<\T t^{(j)},\TTe{j+1}-\TTe{j}\>}{<\TTe{j},\TTe{j+1}\>}\right)+\cr
&E^{(j)}\left(1+\frac{1}{n}\frac{\<\T t^{(j)},\TTe{j+1}\>}{<\TTe{j},\TTe{j+1}\>}\right)+
E^{(j+1)}\left(1-\frac{1}{n}\frac{\<\T t^{(j)},\TTe{j}\>}{<\TTe{j},\TTe{j+1}\>}\right)
\end{array}
\label{eq:T_cp_dependency_inner_even}
\end{equation}
The classification may be formally subdivided into {\it two steps}: at the first step all $D$-type control points are classified as basic, at the second step all $T$-type control points are classified as dependent and their dependencies are established. 

\paragraph{Even vertex excluding the regular $4$-vertices}

The local template for an inner odd vertex, excluding regular $4$-vertices, is shown in Figure ~\ref{fig:fig42Bb}. $val(V)-1$ among $D$-type control points and one $T$-type control point are classified as basic; $val(V)$ basic control points in all. $D$-type control of some edge may be chosen to be dependent only if two neighboring edges of the edge are not colinear.

The correctness of the classification and dependencies of the dependent control points are explained below.

According to Theorem ~\ref{theorem:eq1_all_edges_of_vertex}, the  {\it "Circular Constraint"} should be imposed on $D$-type control points adjacent to an inner even vertex, excluding regular $4$-vertices. Classification of $D$-type and $T$-type control points  can be made as follows. 

At the {\it first step} $D$-type control points are classified. One $D$-type control point with a non-zero coefficient, say $\T D^{(k)}$, is chosen. The remaining $D$-type control points are classified as basic and $D^{(k)}$ depends on them (and $V$, $E$-type basic control points) according to the {\it "Circular Constraint"}
\vspace{-0.05in}
\begin{equation}
\begin{array}{l}
D^{(k)}=2 E^{(k)}-V+\frac{1}{n(n-1)coeff^{(k)}}\sum_{\tiny \begin{array}{c}j=1\cr j\neq k\end{array}}^{val(V)}
(-1)^{j+k+1} coeff^{(j)}\delta^{(j)}
\end{array}
\label{eq:eq1_delta_dependency}
\end{equation}
\vspace{-0.05in}

At the {\it second step} $T$-type control points are classified.
Rank-deficiency of matrix $M$ means that one of $T$-type control points, for example $\T T^{(val(V))}$, can be classified as a basic control point. The remaining $T$-type control points depend on this control point and $D$-type control points (which are classified during the previous step) according to the following equation
\begin{equation}
\hspace{-0.4in}
\left(\!\!
\begin{array}{c}
tw^{(1)}\cr
tw^{(2)}\cr
\vdots\cr 
tw^{(val(V)-1)}\end{array}
\!\!\right)=
\left(\!
M^{\tiny\begin{array}{c}1,\ldots,val(V)-1\cr 1,\ldots,val(V)-1\end{array}}
\!\!\right)
^{-1}
\left(\!\!
\begin{array}{c}
coeff^{(1)}\delta^{(1)}-tw^{(val(V))}\cr
coeff^{(2)}\delta^{(2)}\cr
\vdots\cr
coeff^{(val(V)-1)}\delta^{(val(V)-1)}\end{array}
\!\!\right)
\label{eq:eq1_tw_dependency}
\end{equation}
where $M^{\tiny\begin{array}{c}1,\ldots,val(V)-1\cr 1,\ldots,val(V)-1\end{array}}$ is 
a $(val(V)-1)\times(val(V)-1)$ square matrix which contains $val(V)-1$ first rows and lines
of the matrix $M$.

\paragraph{Regular $4$-vertex}

Local template for a regular $4$-vertex is shown in Figure ~\ref{fig:fig42Bc}. All $D$-type control points and one $T$-type control point are classified as basic, the remaining $T$-type control points are dependent; there are $val(V)+1$ basic control points in all.

In case of a regular $4$-vertex the {\it "Circular Constraint"} should not be explicitly imposed. At the {\it first step} all $D$-type control points are classified as basic. At the {\it second step} one of the  $T$-type control points is classified as basic and others as dependent ones. $T$-type control points do not depend on $D$-type control points; dependency between $T$-type control points is defined by the relation
\vspace{-0.07in}
\begin{equation}
tw^{(1)} = -tw^{(2)}= tw^{(3)}= -tw^{(4)}
\label{eq:T_cp_dependency_4_regular}
\end{equation}

\subsubsection{Local templates for a separate boundary vertex}
\label{subsect:local_mds_DT_boundary_bilinear}

In case of a boundary non-corner vertex there is exactly one adjacent inner edge $\T e^{(2)}$ and so it is sufficient to consider a single constraint in the form of Equation ~\ref{eq:eq1_for_one_edge}. The following two local templates are defined:  
\begin{description}
\item[]
\underline{$TB0^{(D,T)}$} (Figure ~\ref{fig:fig42Ca}). The local MDS contains a $D$-type control point $\T D^{(2)}$ and one $T$-type control points ($\T T^{(1)}$). This template can always be used , excluding the case when the boundary edges are not co-linear and the clamped boundary condition is imposed.
\item[]
\underline{$TB1^{(D,T)}$} (Figure ~\ref{fig:fig42Cb}). The local MDS contains two $T$-type control points $\T T^{(1)}$ and $\T T^{(2)}$. This template is used only in the case when the boundary edges are not co-linear and the clamped boundary condition is imposed.
\end{description}

Below the different "additional" conditions are considered and $TB0^{(D,T)}$ is shown to fit almost all situations.

{\it Any kind of Interpolation/Simply-supported boundary condition.}

These boundary conditions do not involve neither $D$-type nor $T$-type control points adjacent to the boundary vertex. While $D$-type control point do not contribute to Equation ~\ref{eq:eq1_for_one_edge} in case of colinear boundary edges, the coefficient of $T^{(1)}$ is equal to $-n\<\T e^{(2)}, \T e^{(3)}\>$ and never vanishes. Therefore it is very natural to choose a $T$-type control point as a dependent variable (template $TB0^{(D,T)}$), its dependency is defined by Equation ~\ref{eq:eq1_for_one_edge}. Both basic control points are classified as free control points.\

{\it Clamped boundary condition.} 

The clamped boundary condition involves both $T$-type control points and do not involves $D$-type control point. Therefore a MDS which fits the clamped boundary condition can not define $T$-type control point as dependent on $D$-type control point (see ~\ref{def:mds_fits_additional_constraint}).

If the boundary edges are colinear, $D$-type control point do not participate in Equation ~\ref{eq:eq1_for_one_edge}. Equation ~\ref{eq:eq1_for_one_edge} defines dependency between $T^{(1)}$ and $T^{(2)}$
\vspace{-0.07in}
\begin{equation}
tw^{(1)}+tw^{(2)}=0
\end{equation} 
\vspace{-0.07in}
In case of the standard clamped boundary condition (Equation ~\ref{eq:clamped_Z}) $tw^{(1)}=tw^{(2)}=0$ and the dependency is automatically satisfied. Any other clamped boundary condition should be checked to be compatible with  the $G^1$-smoothness requirement. Template $TB0^{(D,T)}$ can be used: the basic $D$-type control point is classified as free and the basic $T$-type control point is classified as fixed.

If the boundary edges are not colinear, one should use template $TB1^{(D,T)}$ (template $TB0^{(D,T)}$ does not longer fit the "additional" condition). Both $T$-type basic control points are fixed. Dependency of $D^{(2)}$-type is defined according to Equation ~\ref{eq:eq1_for_one_edge}
\begin{equation}
\begin{array}{l}
D^{(2)}=2E^{(2)}-V+\frac{1}{n(n-1)coeff^{(2)}}(tw^{(1)}+tw^{(2)})
\end{array}
\label{eq:D_dependency_boundary_bilinear}
\end{equation}

\subsection{Local classification of the {\it middle} control points for a separate edge}
\label{subsect:local_mds_middle_bilinear}

\begin{definition}
Let the global bilinear in-plane parametrisation be considered.

For an inner edge, the set of $n-2$ indexed equations {\it "Eq(s)"} for $s=2,\ldots,s=n-1$ is said to compose the \underline{\it "Middle"} system of equations.

Control points $\T L_2,\ldots,\T L_{n-2}$, $\T R_2,\ldots,\T R_{n-2}$ and $\T C_3,\ldots,\T C_{n-3}$ are respectively called the \underline{"side" middle} and the \underline{"central" middle} control points.
\label{def:middle_bilinear}
\end{definition}

\subsubsection{Existence and types of the local templates}
\label{subsect:templates_middle_bilinear}

Let all $V$, $E$, $T$ and $D$-type control points be already classified and all {\it "Eq(0)"}-type and {\it "Eq(1)"}-type equations be satisfied. The only non-classified control points which participate in $G^1$-continuity equations for an inner edge are the {\it middle} control points and the {\it "Middle"} system is not yet satisfied. We shall show that the {\it middle} control points are sufficient in order to satisfy the {\it "Middle"} system of equations for an inner edge. It provides the possibility to define the local MDS separately for every inner edge.
The following geometric condition  plays a principal role in the analysis of the {\it "Middle"} system of equations.
\begin{definition}
For two adjacent patches with vertices $\LL$,$\LL'$,$\GG$,$\GG'$,$\RR$,$\RR'$ (see Figure ~\ref{fig:fig27}), let ${(\LL-\GG)}^{(proj)}$, ${(\LL'-\GG')}^{(proj)}$, ${(\RR-\GG)}^{(proj)}$, ${(\RR'-\GG')}^{(proj)}$ denote the lengths of projections of the corresponding planar vectors onto direction perpendicular to $\GG'-\GG$. (In other words, distances from vertices $\LL$,$\LL'$,$\RR$,$\RR'$ to the line $(\GG,\GG')$)
\begin{itemize}
\item
It will be said that the \underline{"Projections Relation"} holds if
\vspace{-0.1in}
\begin{equation}
\begin{array}{l}
\frac{{(\RR-\GG)}^{(proj)}}{{(\RR'-\GG')}^{(proj)}}=
\frac{{(\LL-\GG)}^{(proj)}}{{(\LL'-\GG')}^{(proj)}}
\end{array}
\end{equation} 
\item
For the bilinear global in-plane parametrisation $\T\Pi^{(bilinear)}$, the {\it "Projections Relation"} means that the coefficients of the conventional weight functions $l(v)$, $r(v)$ satisfy the following equation:
\vspace{-0.1in}
\begin{equation}
l_0 r_1 - r_0 l_1 = 0
\end{equation}
\end{itemize} 
\label{def:projection_relation}
\end{definition}

The structure of the local MDS depends on the geometrical configuration of two adjacent mesh elements. The following two local templates for classification of the {\it middle} control points are defined.
\begin{description}
\item[]
\underline{$TM0$} (Figure ~\ref{fig:fig42Da}). The local MDS contains all "central" middle control points (if any) and $n-4$ "side" middle control points. The template is applied if the "Projections Relation" does not hold.
\item[]
\underline{$TM1$} (Figure ~\ref{fig:fig42Db}). The local MDS contains all "central" middle control points (if any) and $n-3$ "side" middle control points, one in every pair $(\T L_j,\T R_j)$,  
$j=2,\ldots,n-2$. The template is applied if the "Projections Relation" holds.
\end{description}
Theorem ~\ref{theorem:middle_equations_bilinear} justifies the choice of the local MDS and proves the correctness of the classification of the control points.

\begin{theorem}
Let $\T L$ and $\T R$ be the restrictions of the global in-plane parametrisation $\T\Pi^{(bilinear)})$ on two adjacent mesh elements. 
Let {\it "Eq(0)"}-type and {\it "Eq(1)"}-type equations be satisfied for the common edge and all $E$,$V$,$D$,$T$-type control points be classified. Then for any $n\ge 4$ \begin{itemize}
\item[{\bf (1)}] {\bf Consistency.}
The {\it "Middle"} system of equations has a  solution in terms of the {\it middle} control points.
\item[{\bf (2)}] {\bf Classification of the {\it middle} control points.} The following classification of the {\it middle} control points guarantees that the {\it "Middle"} system of equations for an inner edge is satisfied.
\begin{itemize}
\item
All "central" {\it middle} control points (if there are any) are classified as the basic (free) ones. 
\item
If the {\it "Projections Relation"} does not hold then there are $n-4$ basic (free) and $n-2$ dependent control points among $2(n-3)$ "side" {\it middle} control points.
\item
If the {\it "Projections Relation"} holds then there are $n-3$ basic (free) and $n-3$ dependent control points among $2(n-3)$ "side" {\it middle} control points; in every pair $(\T L_j,\T R_j)$ $(j=2,\ldots,n-2)$ one control point is basic and another one is dependent. 

\end{itemize}
\end{itemize}
\label{theorem:middle_equations_bilinear}
\end{theorem} 
Dependencies of the dependent middle control points are described in the theorem's proof.\\
\vspace{-0.1in}

{\bf Proof} See Appendix, Section ~\ref{sect:proofs}.

\subsubsection{Example of the local MDS for $n=4$ and $n=5$}
\label{subsect:example_4_5_middle_bilinear}

\paragraph{}
For \underline{$n=4$} the {\it "Middle"} system consists of two indexed equations {\it "Eq(2)"} and {\it "Eq(3)"} and there are two {\it middle} control points $\T L_2$ and $\T L_3$.

Theorem ~\ref{theorem:middle_equations_bilinear} implies that if the {\it "Projections Relation"} is not satisfied, then the local MDS is empty. Both {\it middle} control points are dependent and their dependencies on the basic $E$,$V$,$D$,$T$-type control points are defined by the {\it "Middle"} system which has a $2\times2$ matrix of  full rank.

If the {\it "Projections Relation"} holds then equations {\it "Eq(2)"} and {\it "Eq(3)"} are no longer independent and an additional degree of freedom is available; one of the control points $\T L_2$, $\T R_2$ becomes basic (free). 

\paragraph{}
For \underline{$n=5$}, if the {\it "Projections Relation"} does not hold then one of the control points $\T L_2, \T L_3, \T R_2, \T R_3$ is basic (free) and the others are dependent. If the {\it "Projections Relation"} is satisfied then there are two basic (free) control points, one in every pair $(\T L_2,\T R_2)$ and $(\T L_3,\T R_3)$.

\begin{figure}[!pt]
\centering
\includegraphics[clip,totalheight=1.5in]{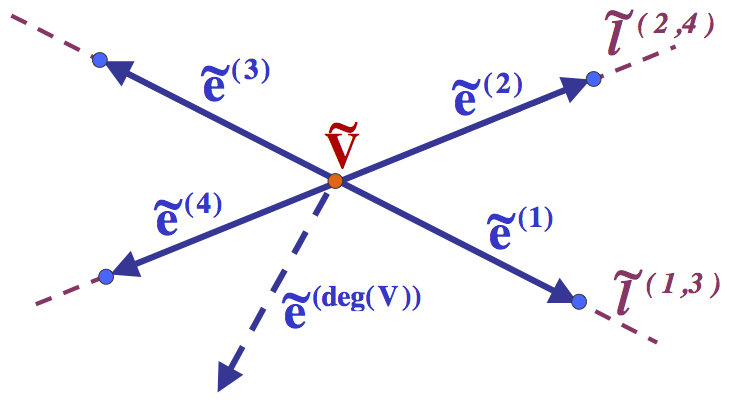}
\caption{Possible mesh configuration, which automatically satisfies the {\it "Circular Constraint"}.}
\label{fig:fig16}
\end{figure}

\begin{figure}[!ph]
\vspace{-0.1in}
\begin{narrow}{0.0in}{0.0in}
\begin{minipage}[]{0.31\linewidth}
\subfigure[]
{
\includegraphics[clip,totalheight=1.3in]{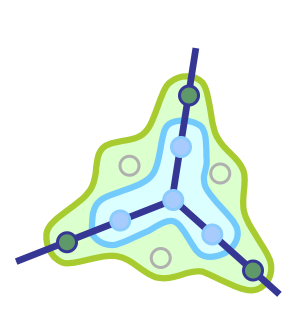}
\label{fig:fig42Ba}
}
\end{minipage}
\begin{minipage}[]{0.31\linewidth}
\subfigure[]
{
\includegraphics[clip,totalheight=1.3in]{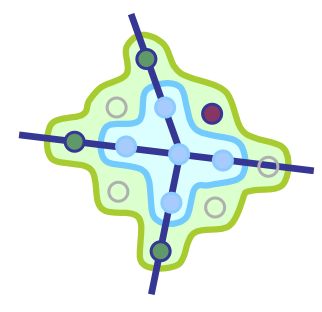}
\label{fig:fig42Bb}
}
\end{minipage}
\begin{minipage}[]{0.31\linewidth}
\subfigure[]
{
\includegraphics[clip,totalheight=1.3in]{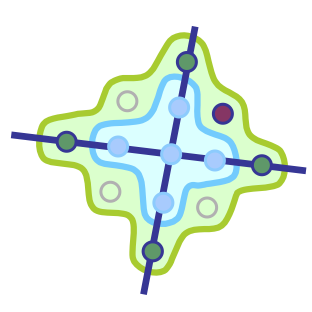}
\label{fig:fig42Bc}
}
\end{minipage}
\end{narrow}
\vspace{-0.2in}
\caption{Local templates for the classification of $D$,$T$-type control points adjacent to an inner vertex in case of global bilinear in-plane parametrisation $\T\Pi^{(bilinear)}$.}
\label{fig:fig42B}
\end{figure}

\FloatBarrier

\begin{figure}[!ph]
\vspace{-0.1in}
\begin{narrow}{0.1in}{0.0in}
\begin{minipage}[]{0.5\linewidth}
\subfigure[]
{
\includegraphics[clip,totalheight=1.1in]{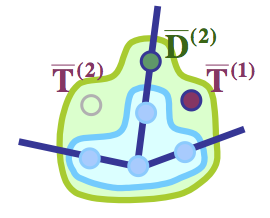}
\label{fig:fig42Ca}
}
\end{minipage}
\begin{minipage}[]{0.5\linewidth}
\subfigure[]
{
\includegraphics[clip,totalheight=1.1in]{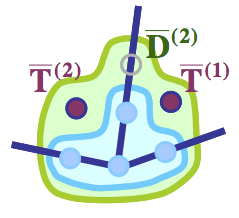}
\label{fig:fig42Cb}
}
\end{minipage}
\end{narrow}
\vspace{-0.2in}
\caption{Local templates for the classification  of $D$,$T$-type control points adjacent to a boundary vertex in case of global bilinear in-plane parametrisation $\T\Pi^{(bilinear)}$.}
\label{fig:fig42C}
\end{figure}

\begin{figure}[!ph]
\vspace{-0.1in}
\centering
\includegraphics[clip,totalheight=2.3in]{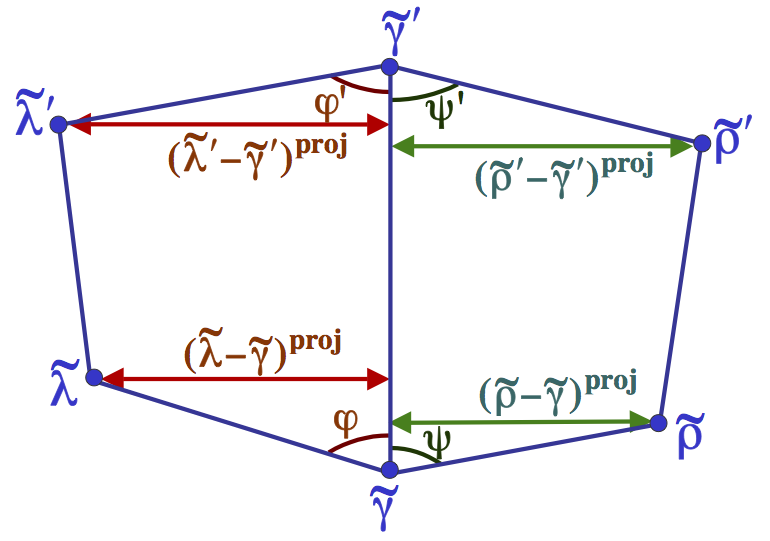}
\caption{An illustration for the {\it "Projections Relation"}.}
\label{fig:fig27}
\end{figure}
\begin{figure}[!ph]
\vspace{-0.1in}
\begin{narrow}{0.1in}{0.0in}
\begin{minipage}[]{0.5\linewidth}
\subfigure[]
{
\includegraphics[clip,totalheight=1.3in]{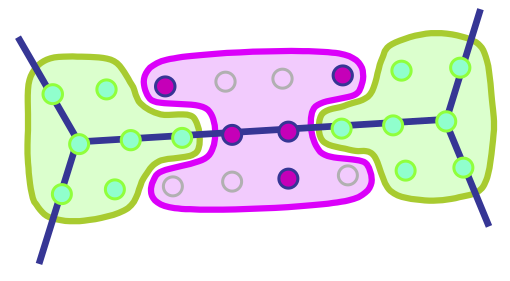}
\label{fig:fig42Da}
}
\end{minipage}
\begin{minipage}[]{0.5\linewidth}
\subfigure[]
{
\includegraphics[clip,totalheight=1.3in]{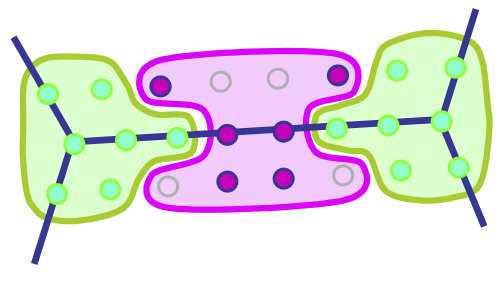}
\label{fig:fig42Db}
}
\end{minipage}
\end{narrow}
\vspace{-0.2in}
\caption{Local templates for the classification of the {\it middle} control points in case of global bilinear in-plane parametrisation $\T\Pi^{(bilinear)}$.}
\label{fig:fig42xD}
\end{figure}

\FloatBarrier

\eject
\section{Global MDS}
\label{sect:global_mds_bilinear}

The main purpose of the current Section is to define the global MDS based on the local analysis given in Section ~\ref{sect:local_mds_bilinear}.
One should try to "put together" local templates without contradictions. More precisely , the {\it order} of construction of the local MDS should be defined. It is important to remember, that even if local templates do not intersect geometrically,  the construction of local MDS is usually based on the assumption that some control points are already classified and it is possible to define dependencies on these control points. Therefore, the order of construction of the local MDS plays the principal role in the construction of the global MDS.

Although a {\it "pure"} global MDS (MDS based on $G^1$-conditions only) always exists, application of an "additional" constraint requires construction of some {\it definite} suitable local templates, which do not always fit together. 
If the global classification succeeds then it will be said that an instance of global MDS, which fits the "additional" constraints, is constructed. Otherwise, according to the current approach, MDS of a higher degree should be considered; an attempt to rebuild the local templates defined in Section ~\ref{sect:local_mds_bilinear} is never done.

\subsection{MDS of degree $n\ge 5$}
\label{subsect:global_mds_n_ge_5_bilinear}

For $n\ge 5$ the local templates for different vertices and/or edges never involve the same control points (see Figure ~\ref{fig:fig39}). For any "additional" constraint (any choice of the local templates), construction of the global MDS can be made according to the following Algorithm.

\begin{algorithm}
\label{algorithm:mds_global_deg_ge_5_bilinear}
\end{algorithm}
\centerline{\underline{\it Algorithm for construction of global MDS}}
\centerline{\underline{\it in case of a global bilinear in-plane parametrisation $\T\Pi^{(bilinear)}$}}

\begin{description}
\item[]
{\bf "Stage 1"} For every non-corner mesh vertex, construct a local MDS for the  classification of $V$,$E$-type control points (see Subsection ~\ref{subsect:local_mds_VE_bilinear}). The choice of the local template for every vertex is made according to the local mesh structure and should fit the given "additional" constraint. At the end of the stage all $V$,$E$-type control points are classified and all {\it "Eq(0)}-type equations are satisfied.

\item[]
{\bf "Stage 2"}
For every non-corner mesh vertex, construct a local MDS for classification of $D$,$T$-type control points (see Subsection ~\ref{subsect:local_mds_DT_bilinear}). The choice of the local template for every vertex should fit the given "additional" constraint. At the end of the stage all $V$,$E$,$D$,$T$-type control points are classified and all {\it "Eq(0)} and {\it "Eq(1)"}-type equations are satisfied.

\item[]
{\bf "Stage 3"}
For every inner edge, construct a local MDS for the classification of the {\it middle} control points (see Subsection ~\ref{subsect:local_mds_middle_bilinear}). The choice of the local template for every edge depends on the local geometrical structure of the mesh. At the end of the stage all control points are classified and all $G^1$-continuity equations are satisfied.
\end{description}

\begin{lemma}
For any $n\ge 5$ and for any type of "additional" constraints, there exists  an instance $\MDS{n}(\T\Pi^{(bilinear)})$ of the global MDS, that fits the "additional" constraints.
\label{lemma:deg_n_ge_5_mds_existence_bilinear}
\end{lemma}
Dimensionality of the global MDS is studied in Subsection ~\ref{subsect:mds_dimensionality_bilinear}.

\subsection{MDS of degree $n=4$}

\subsubsection{Principal role of classification of $D$-type control points}

For $n=4$, construction of the global MDS becomes more complicated since the local templates defined in Section ~\ref{sect:local_mds_bilinear} may intersect. The intersection occurs between local templates for classification of $D$,$T$-type control points adjacent to end-vertices of the same inner edge (see Figure ~\ref{fig:fig41}). The intersection contains the single control point - $D$-type control point of the edge; the problem arises if both templates define the control point as dependent. 

In general, Algorithm ~\ref{algorithm:mds_global_deg_ge_5_bilinear} remains valid for $n=4$. 
Moreover, the only stage which should be made more accurate is {\it "Stage 2"}. At this stage for $n\ge 5$, the local templates for vertices may be defined in any order and for every local template the choice of dependent $D$-type control point (if needed) may be made arbitrarily in case of ambiguity. For $n=4$ the {\it order} in which the vertices are traversed and the choice of the dependent $D$-type control points (if needed) play the principal role.

According to the current approach, the global classification of $D$,$T$-type is said to exist if it is possible to "put together" the local templates without contradiction. In order to do it, one should specify the traversal order of the vertices and the choice of dependent $D$-type control point (if any) for every template. The classification involves a global analysis and may fail for some kinds of the "additional" constraints. If the classification succeeds, the global MDS exists and the rest of the control points are classified precisely as in the case of $n\ge 5$.

It is important to note, that for all local templates, excluding $TB1^{(D,T)}$ (see Subsection ~\ref{subsect:local_mds_DT_boundary_bilinear}), $D$-type control points can be classified prior to classification of $T$-type control points. Furthermore, as soon as a global classification of $D$-type control points is made, classification of $T$-type control points can be made locally, separately for every vertex and can not lead to any conflicts between different vertices.
It implies that global classification of $D$-type control points is the most difficult and the important step in the construction of the global MDS. The following notations will be used.

\vspace*{0.3in}
\begin{definition}
\label{definition:D_relevant}
\nopagebreak
\hspace{0.5in}
\nopagebreak
\begin{itemize}
\item[]
A vertex  $\T V$ is called \underline{$D$-relevant} if the local template for classification of $D$,$T$-type control points contains  a dependent $D$-type control point.
\item[]
It will be said, that the {\it $D$-relevant} vertex $\T V$ \underline{uses} an adjacent $D$-type control point $\T D^{(j)}$ if $D^{(j)}$ enters Equation ~\ref{eq:eq1_for_one_edge} with a non-zero coefficient.
\item[]
A $D$-type control point is said to be \underline{assigned} to an adjacent $D$-relevant vertex if it is classified as dependent according to the local template at the vertex. The vertex is called "owner" of the assigned $D$-type control point.
\end{itemize}
\end{definition}
There are two types of {\it $D$-relevant} vertices
\begin{itemize}
\item[-]
Inner even vertices, excluding regular $4$-vertices.
\item[-]
Boundary vertices which use template $TB1^{(D,V)}$, or in other words, boundary vertices with non-colinear adjacent boundary edges in the case of a clamped boundary condition. 
\end{itemize} 

\subsubsection{Examples of possible difficulties}
The following small examples show what kind of difficulties one may encounter during construction of global MDS of degree $4$ which fits a given "additional" constraint.

In the case of a clamped boundary condition, the
three meshes presented in Figure ~\ref{fig:fig18} have no sufficient $D$-type control points in order to {\it assign} an "owned" 
$D$-type control point to every {\it $D$-relevant} mesh vertex (arrows in the Figure show the "owner" vertices for $D$-type control points). 
Figure ~\ref{fig:fig18a} - the inner $4$-vertex has not its own $D$-type control point.
Figure ~\ref{fig:fig18b} - none of the inner $4$-vertices has its own $D$-type control point because $D$-points of the dashed edges do not contribute to the {\it "Circular Constraint"} for the inner vertices.
Figure ~\ref{fig:fig18c} - there are no sufficient $D$-type control points for either one of two inner $4$-vertices.

An additional example of Figure ~\ref{fig:fig19} shows the situation when an inner $3$-vertex $\T V$ (which does not use the adjacent $D$-type control points itself) also "does not help" to get free $D$-type control points for the neighboring $4$-vertices. This is because the  $D$-type control points of the dashed edges do not
contribute to the {\it "Circular Constraint"} for vertices $\T V^{(1)}$, $\T V^{(2)}$ and $\T V^{(3)}$.


\subsubsection{Sufficient conditions and algorithms for the global classification of $D$,$T$-type control points}
\label{subsect:global_classification_DT_bilinear}

Classification of $D$,$T$-type control points will explicitly be made for any kind of "additional" constraints, excluding the clamped boundary condition. For most of the mesh configurations, the general algorithms work as well in the case of a clamped boundary condition. Different techniques, which may help in case of failure of the general algorithm, are provided. The classification process heavily uses different graph-like structures and graph-theory algorithms (see ~\cite{cormen}). The next Definition introduces some special "graph-related" notations, which will be used in the current discussion.

\begin{definition}
\hspace{0.1in}
\begin{itemize}
\item
\underline{Mesh vertex} or \underline{primary vertex} - vertex of the initial mesh.
\item
\underline{Secondary vertex} - auxiliary (symbolic) vertex in the middle of the edge.
\item
\underline{Mesh edge} or \underline{full edge}  - edge of the initial mesh.
\item
\underline{Half-edge} - edge connecting a primary vertex and an adjacent secondary vertex.
\end{itemize}
\end{definition}
First, an undirected $D$-dependency graph will be constructed.
Then the spanning tree algorithm will be applied to every connected component of the dependency graph and a directed $D$-dependency forest will be built. In most cases, the dependency 
forest allows to define the traversal order for the vertices and to assign $D$-type control point to every $D$-relevant vertex. 

\paragraph{Construction of an undirected $D$-dependency graph}

An undirected $D$-dependency graph is constructed according to the following Algorithm.

\begin{algorithm}
\label{algorithm:dependency_graph_construction}
\end{algorithm}
\centerline{\underline{\it Algorithm for construction of the $D$-dependency graph}}

\begin{itemize}
\item[{\bf (1)}]
Add {\it secondary vertex} $\T S^{(ij)}$ in the middle of every {\it mesh edge} $(\T V^{(i)},\T V^{(j)})$ (see  Figure ~\ref{fig:fig20a}).
\item[{\bf (2)}]
Delete {\it half-edge} $(\T V^{(i)},\T S^{(ij)})$ if vertex $\T V^{(i)}$ does not use the $D$-type control point of the {\it mesh edge} $(\T V^{(i)},\T V^{(j)})$ (see Definition ~\ref{definition:D_relevant} and Figure ~\ref{fig:fig20b}).
\item[{\bf (3)}]
Delete a {\it secondary vertex} if after elimination of {\it half-edges} (step {\bf(2)} of the Algorithm) no {\it half-edge} connected to it remained.
\item[{\bf (4)}]
Delete a {\it primary} vertex if it does not use any $D$-type control point (if it is not {\it $D$-relevant}).
\end{itemize}

Now the $D$-dependency graph consists of {\it $D$-relevant} {\it primary vertices} and such {\it secondary vertices} that the correspondent $D$-type control point is used by some {\it primary vertex}. A {\it primary} and a {\it secondary vertex} are connected by a {\it half-edge} if and only if the {\it primary vertex} uses the $D$-type control point corresponding to the {\it secondary vertex}. It will be said that two {\it primary vertices} $\T V^{(i)}$ and $\T V^{(j)}$ are connected by the {\it full edge} if both of them use $D$-type control point of the {\it mesh edge}  $(\T V^{(i)}, \T V^{(j)})$ (both {\it half-edges} $(\T V^{(i)}, \T S^{(ij)})$ and  $(\T S^{(ij)}, \T V^{(j)})$ belong to the $D$-dependency graph). If $\T S^{ij}$ is connected to only one of 
the {\it primary vertices}, say $\T V^{(i)}$, then the {\it half-edge} $(\T V^{(i)}, \T S^{(ij)})$ will be called a \underline{\it dangling half-edge}.

From the topological point of view the graph consists of a few (possible zero or one) connected components. Each connected component may contain pairs of the {\it primary vertices} connected by the {\it full edges} and the {\it secondary vertices} connected to the {\it primary} ones by the {\it dangling half-edges}.

\paragraph{Construction of directed a $D$-dependency forest when any connected component of $D$-dependency graph has a dangling half-edge}

\begin{algorithm}
\label{algorithm:dependency_forest_construction}
\end{algorithm}
\centerline{\underline{\it Algorithm for the construction of the $D$-dependency  forest}}
\centerline{\underline{\it when every connected component of the $D$-dependency graph}}
\centerline{\underline{\it has at least one dangling half-edge}}

\noindent
For every connected component of $D$-dependency graph, define a $D$-dependency tree in the following way.

\begin{itemize}
\item[{\bf (1)}]
Choose a {\it primary vertex} which has a \underline{\it dangling half-edge} (at least one such vertex exists for every connected component), denote this vertex by $\T R$.
\item[{\bf (2)}]
Build a \underline{spanning tree} of the {\it primary vertices} of the connected component using the {\it full edges} only (this is possible because the component obviously remains connected with respect to the {\it primary vertices} after elimination of all {\it dangling half-edges} (see Figure ~\ref{fig:fig21b}).
\item[{\bf (3)}]
Define vertex $\T R$ to be the \underline{root} of the spanning tree.
\item[{\bf (4)}]
\underline{Direct} every {\it mesh edge} of the tree from the upper vertex to a lower one (according to the hierarchy of the spanning tree, see Figure ~\ref{fig:fig21c}).
\item[{\bf (5)}]
\underline{Direct} the {\it mesh edge} corresponding to a {\it dangling  half-edge} at $\T R$
towards $\T R$. From this moment this {\it mesh edge} is defined as belonging to the $D$-dependency tree of $\T R$ (the correctness of this definition is shown in Lemma
~\ref{lemma:dependency_forest_properties}).
\end{itemize}

\begin{note}
The $D$-dependency forest consists of the {\it primary} $D$-relevant vertices and directed {\it mesh edges}; at this step one may forget about auxiliary secondary vertices and half-edges.
\end{note}

\begin{note}
Assigning directions to the mesh edges which correspond to the dangling half-edges of the root vertices may lead to the situation when $D$-dependency trees of different components of the $D$-dependency graph become connected by these directed edges (see Figure ~\ref{fig:fig21A}). 
Such trees still will be referred to as the \underline{different} trees of the $D$-dependency forest. 
The mesh vertex and its adjacent directed mesh edge may belong to  different $D$-dependency trees; a directed \underline{mesh edge} always belongs to the \underline{$D$-dependency tree of the vertex it points to}.
\end{note}

\begin{lemma}
The directions assigned to {\it mesh edges} as a result of the construction of $D$-dependency forest, obey the following rules
\begin{itemize}
\item[{\bf (1)}]
Every mesh edge is either undirected or its direction is correctly defined. In particular, it implies that every mesh edge belongs to at most one $D$-dependency tree.
\item[{\bf (2)}]
Every {\it $D$-relevant} mesh vertex uses $D$-type control points of either undirected mesh edges or of the directed mesh edges that belong to the same $D$-dependency tree as the vertex itself.
\item[{\bf (3)}]
Every {\it $D$-relevant} mesh vertex $\T V$ uses exactly one $D$-type control point which belongs to the directed edge pointing towards $\T V$. This $D$-type control point is denoted $\T D(\T V)$ (see Figure ~\ref{fig:fig21c}).
\end{itemize}
\label{lemma:dependency_forest_properties}
\end{lemma}
\vspace{-0.07in}
{\bf Proof} See Appendix, Section ~\ref{sect:proofs}.

\paragraph{Explicit classification of $D$,$T$-type control points when any connected component of $D$-dependency graph has a dangling half-edge}

As soon as the $D$-dependency forest is constructed, the global classification of $D$,$T$-type control points can be made according to Algorithm ~\ref{algorithm:D_T_classification_degree_4}. The correctness of classification is justified by Lemma ~\ref{lemma:dependency_forest_properties}. In particular, Lemma  ~\ref{lemma:dependency_forest_properties} implies that classification of $D$-type control points for the different trees of the $D$-dependency forest can be made independently, without any conflicts.

\begin{algorithm}
\label{algorithm:D_T_classification_degree_4}
\end{algorithm}
\centerline{\underline{\it Algorithm for the classification of $D$,$T$-type control points}}
\centerline{\underline{\it based on $D$-dependency forest}}

\vspace{0.05in}
\noindent
{\bf "Step 1"}
\vspace{-0.05in}
\begin{itemize}
\item[{\bf (1)}] 
For every $D$-relevant boundary vertex, classify $T$-type control points as basic according to template $TB1^{(D,V)}$.
\item[{\bf (2)}] 
Classify as basic $D$-type control points of all mesh edges which remained undirected after construction of the $D$-dependency forest.
\item[{\bf (3)}] 
For every tree in the $D$-dependency forest do the following

\noindent
$\bullet$\ \ \
{\it Assign} the $D$-type control point of every directed {\it mesh edge} to the {\it mesh vertex} this edge points to. 

\noindent
$\bullet$\ \ \
Traverse the $D$-dependency tree down-up, level by level. 

\noindent
$\bullet$\ \ \ 
At every level for every $D$-relevant vertex $\T V$ (the order of vertices belonging to the same level is not important) choose the assigned $D$-type control point $\T D(\T V)$ as dependent in the local template of the vertex. Define the dependency of the corresponding $Z$-component $D(\T V)$ according to the local template (Equation ~\ref{eq:eq1_delta_dependency} for an inner vertex and Equation ~\ref{eq:D_dependency_boundary_bilinear} for a boundary vertex). $D(\T V)$ may depend on
\begin{itemize}
\item[-]
$V$-type and $E$-type control points which are classified as basic during the first stage of the classification process.
\item[-]
$D$-type control points which are classified as basic during traversing the lower levels of the tree. 
\item[-]
$T$-type basic fixed control points adjacent to $D$-relevant boundary vertices. 
\end{itemize}
\end{itemize}

\noindent
{\bf "Step 2"}

\vspace{-0.07in}
\begin{itemize}
\item[]
For every non-corner vertex (excluding $D$-relevant boundary vertices), classify $T$-type control points locally according to the local template (see Subsections ~\ref{subsect:local_DT_inner_vertex_bilinear} and ~\ref{subsect:local_mds_DT_boundary_bilinear}). 
\end{itemize}
For example, in the mesh fragment shown in Figure ~\ref{fig:fig22}, $\T V$ is a leaf-vertex of some $D$-dependency tree, $\T V^{(4)}$ is the "father" of $\T V$ in the $D$-dependency tree and $\T V^{(1)}$ is the root vertex of some other $D$-dependency tree. Then 
$\T D^{(2)}$,$\T D^{(3)}$,$\T D^{(5)}$,$\T D^{(6)}$ are basic $D$-type control points. $\T D^{(1)}$ is a dependent control point, however it does not contribute to the {\it "Circular Constraint"} for $\T V$ since edges  numbered  $2$ and $6$ are colinear.This {\it "Circular Constraint"} defines the  dependency of $D^{(4)}=D(\T V)$,
 $D^{(2)}$,$D^{(3)}$,$D^{(5)}$ (and $E$ and $V$-type basic control points) according to the equation
\vspace{-0.05in}
\begin{equation}
\begin{array}{l}
\sum_{j=2}^6\frac{\<\TTe{j-1},\TTe{j+1}\>}{\<\TTe{j-1},\TTe{j}\>\<\TTe{j},\TTe{j+1}\>}
(D^{(j)}-2E^{(j)}+V) = 0
\end{array}
\end{equation}

\paragraph{Absence of $D$-relevant boundary vertices as a sufficient condition for global classification of $D$,$T$-type control points}

Theorem ~\ref{theorem:sufficient_DT_classification_bilinear}
presents a sufficient condition for the existence of global classification of $D$,$T$-type control points.

It is important to pay attention to the fact that the presence of inner $D$-relevant vertices depends on the mesh structure, while the presence of boundary $D$-relevant vertices depends both on mesh structure and on the type of applied "additional" constraints. 

\begin{theorem}
Absence of $D$-relevant boundary vertices implies that the $D$-dependency graph is either empty or every connected component of $D$-dependency graph contains at least one dangling half-edge.
Therefore the global classification of $D$,$T$-type control points exists and can be made by successive application of Algorithms ~\ref{algorithm:dependency_graph_construction}, ~\ref{algorithm:dependency_forest_construction} and ~\ref{algorithm:D_T_classification_degree_4}.
\label{theorem:sufficient_DT_classification_bilinear}
\end{theorem}
{\bf Proof} See Appendix, Section ~\ref{sect:proofs}.

\paragraph{Case when there exist $D$-relevant boundary vertices}

In most of practical situations, even in the presence of $D$-relevant boundary vertices, every connected component of the $D$-dependency graph has at least one dangling half-edge due to the neighborhood of inner odd vertices or regular $4$-vertices. 
In this case the classification of $D$,$T$-type control points exists and is made according to Algorithms ~\ref{algorithm:dependency_graph_construction}, ~\ref{algorithm:dependency_forest_construction} and ~\ref{algorithm:D_T_classification_degree_4}. If some connected component has no dangling half-edges, then Algorithm ~\ref{algorithm:dependency_forest_construction} fails to construct the $D$-dependency tree for this component. 

Algorithm ~\ref{algorithm:dependency_forest_round} (see Appendix, Section ~\ref{sect:algorithms}) presents a simple modification of Algorithms ~\ref{algorithm:dependency_forest_construction} and ~\ref{algorithm:D_T_classification_degree_4} which may help to classify $D$,$T$-type control points for those components of the $D$-dependency graph, which have no dangling half-edges. The construction is made separately for every such component and does not affect the construction of the classification solution for components with dangling half-edges.

The principal reason of failure of Algorithm ~\ref{algorithm:dependency_forest_construction} for a connected  component $\T{\cal C}$ which has no dangling half-edges is an impossibility to {\it assign} $D$-type control point to the root vertex of the spanning tree. Algorithm ~\ref{algorithm:dependency_forest_round} tries to overcome this problem by splitting a full edge of a component so that the edge does not participate in the spanning tree of $\T{\cal C}$. Algorithm ~\ref{algorithm:dependency_forest_round} does not necessarily succeed. If the Algorithm fails (for one of the reasons which are explained below), any one of the following possibilities may be chosen.
\begin{itemize}
\item
Pass to MDS of degree $5$ which is well defined for any kind of "additional" constraints. 
\item
Pass to the solution which combines patches of degree $5$ along the boundary and patches of degree $4$ in the inner part of the mesh (see Part ~\ref{part:mds_mixed_4_5}).
\item
Try to create a dangling half-edge by a local change of the initial mesh, for example by the regularisation of one of the inner $4$-vertices, as shown in Figure ~\ref{fig:fig25}.
\end{itemize}

\begin{algorithm}
See Appendix, Section ~\ref{sect:algorithms}.
\label{algorithm:dependency_forest_round}
\end{algorithm}

\subsubsection{The existence of MDS. Analysis of different "additional" constraints.}

The results of Subsection ~\ref{subsect:global_classification_DT_bilinear} lead to the following conclusions.

\begin{lemma}
\hspace{0.4in}
\begin{itemize}
\item[{\bf (1)}]
Global {\it "pure"} MDS\ \ $\MDS{4}(\T\Pi^{(bilinear)})$ (MDS which relates to $G^1$ continuity constraints alone) is well defined for any mesh configuration.

\item[{\bf (2)}]
For any "additional" constraint which does not involve $D$-relevant boundary mesh vertices, an instance of $\MDS{4}(\T\Pi^{(bilinear)}),$ which fits the "additional" constraint, exists for any mesh configuration. In particular, a suitable $\MDS{4}(\T\Pi^{(bilinear)})$ always exists for any kind of interpolation and for simply-supported boundary conditions.
\end{itemize}
\label{lemma:deg_4_mds_existence_bilinear}
\end{lemma}
The MDS is constructed according to Algorithm ~\ref{algorithm:mds_global_deg_ge_5_bilinear}, where {\it "Stage 2"} is accomplished by successive application of Algorithms ~\ref{algorithm:dependency_graph_construction}, ~\ref{algorithm:dependency_forest_construction} and ~\ref{algorithm:D_T_classification_degree_4}. The dimensionality of $\MDS{4}(\T\Pi^{(bilinear)})$ is analysed in Subsection ~\ref{subsect:mds_dimensionality_bilinear}.

For any "additional" constraint which involve $D$-relevant boundary vertices (ex. clamped boundary  condition for vertices with non-colinear adjacent boundary edges), the current approach may fail to construct an instance of $\MDS{4}(\T\Pi^{(bilinear)})$ which fits this "additional" constraint. Nevertheless even in the presence of $D$-relevant boundary vertices, the general classification algorithm or its modification (Algorithm ~\ref{algorithm:dependency_forest_round}) will work for most mesh configurations. A solution which combines patches of degrees $4$ and $5$ is presented in Part ~\ref{part:mds_mixed_4_5}.

\subsection{Dimensionality of MDS}
\label{subsect:mds_dimensionality_bilinear}

\begin{theorem}
For a global bilinear in-plane parametrisation $\T\Pi^{(bilinear)}$ and for any $n\ge 4$ , the dimension of $\GMDS{n}$ (subset of MDS which participates in the  $G^1$-continuity condition) and the dimension of $\MDS{n}$ (full dimension of MDS) are given by the following formulas (see Definition ~\ref{def:cp_subsets} and Lemma  ~\ref{lemma:GMDS_importance})
\begin{equation}
\hspace{-0.4in}
\begin{array}{ll}
|\GMDS{n}| = &
3|Vert_{non-corner}|+
|Vert_{\twolines{boundary}{non\!-\!corner}}|+
(2n-7) |Edge_{inner}|+\cr
&|Vert_{\twolines{inner}{4\!-\!regular}}|+
|Edge_{\twolines{inner,}{"Projections\ Relation"\ holds}}|
\end{array}
\label{eq:dimension_gmds_bilinear}
\end{equation}

\begin{equation}
\hspace{-0.4in}
\begin{array}{lll}
|\MDS{n}| = &|\GMDS{n}|+&|\FCP{n}|=\cr
&|\GMDS{n}|+&(n\!-\!3)^2|Face_{inner}|+
(n\!-\!3)(n\!-\!1)|Face_{\twolines{boundary}{non\!-\!corner}}|+\cr
&&(n\!-\!1)^2|Face_{corner}|
\end{array}
\label{eq:mds_dimension_bilinear}
\end{equation}
\label{theorem:mds_dimension_bilinear}
\end{theorem}
\vspace{-0.1in}
{\bf Proof} See Appendix, Section ~\ref{sect:proofs}.

\begin{figure}[!ph]
\vspace{-0.1in}
\begin{narrow}{-0.4in}{-0.4in}
\centering
\includegraphics[clip,height=2.8in]{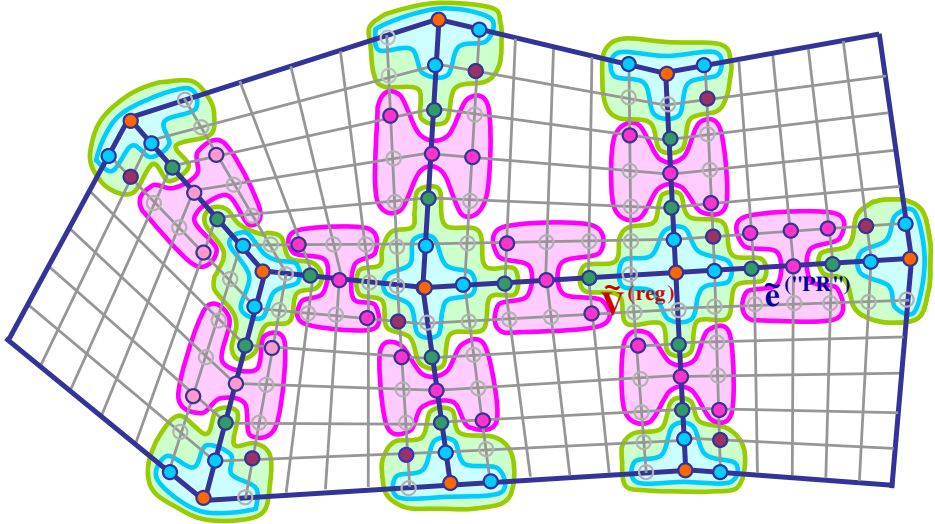}
\end{narrow}
\caption{An example of a "pure" global MDS $\MDS{5}(\T\Pi^{(bilinear)}$.}
\label{fig:fig39}
\end{figure}

\begin{figure}[!pb]
\vspace{-0.25in}
\begin{narrow}{-0.4in}{-0.4in}
\centering
\includegraphics[clip,height=2.8in]{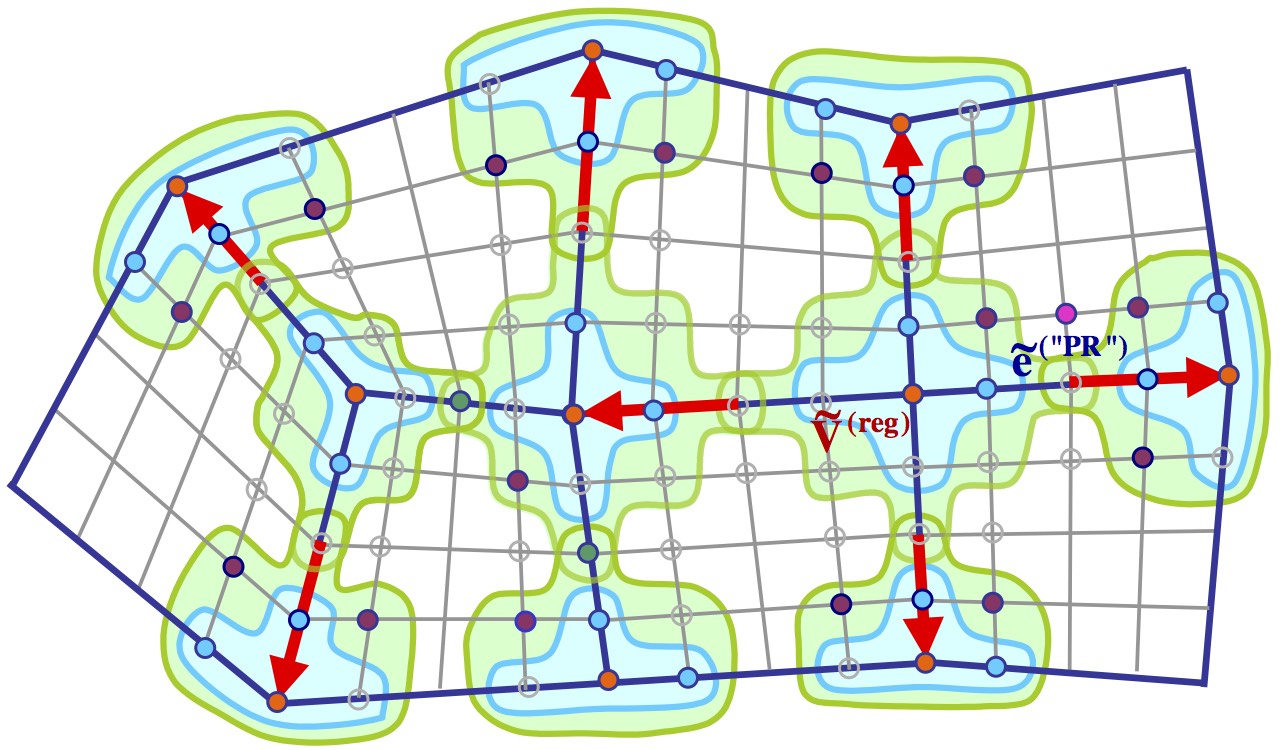}
\end{narrow}
\caption{An example of a global MDS $\MDS{4}(\T\Pi^{(bilinear)}$, which fits the clamped boundary conditions.}
\label{fig:fig41}
\end{figure}

\begin{figure}[!pt]
\vspace{-0.4in}
\begin{narrow}{-0.5in}{-0.5in}
\begin{minipage}[]{0.23\linewidth}
\subfigure[]
{
\includegraphics[clip,totalheight=1.4in]{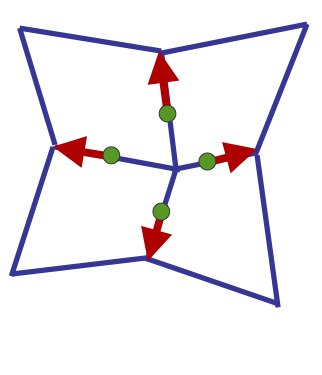}
\label{fig:fig18a}
}
\end{minipage}
\begin{minipage}[]{0.46\linewidth}
\subfigure[]
{
\includegraphics[clip,totalheight=1.4in]{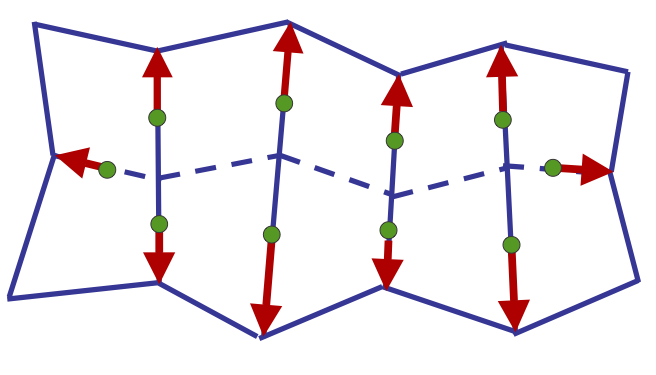}
\label{fig:fig18b}
}
\end{minipage}
\begin{minipage}[]{0.28\linewidth}
\subfigure[]
{
\includegraphics[clip,totalheight=1.4in]{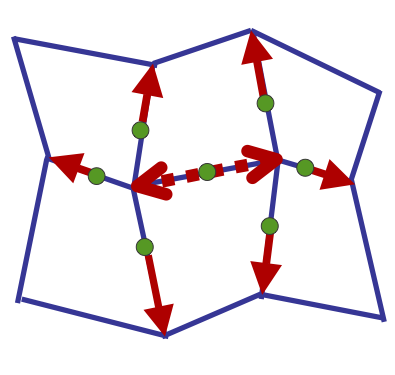}
\label{fig:fig18c}
}
\end{minipage}
\end{narrow}
\vspace{-0.2in}
\caption{Examples of such meshes that not every {\it $D$-relevant} mesh vertex has its own $D$-type control point.}
\label{fig:fig18}
\end{figure}

\begin{figure}[!ph]
\centering
\includegraphics[clip,totalheight=2.0in]{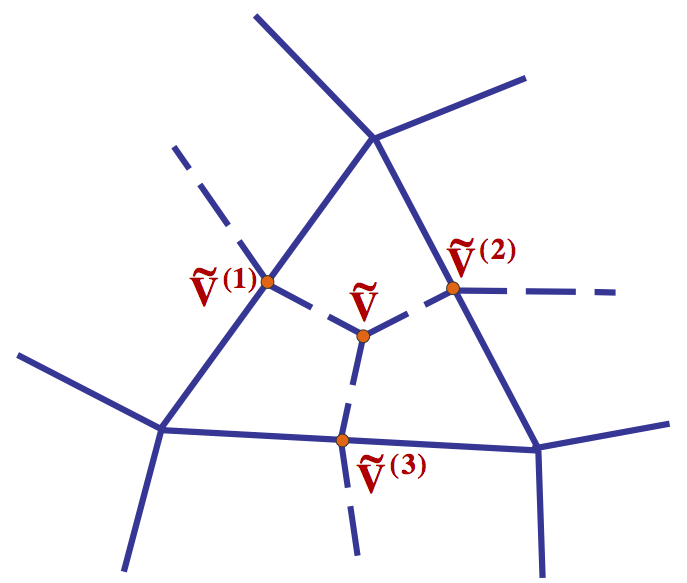}
\caption{An example of mesh configuration when $3$-vertex does not contribute to the {\it "Circular Constraint"} for adjacent $4$-vertices.}
\label{fig:fig19}
\end{figure}

\begin{figure}[!ph]
\begin{narrow}{-0.5in}{-0.5in}
\begin{minipage}[]{0.3\linewidth}
\subfigure[]
{
\includegraphics[clip,totalheight=1.0in]{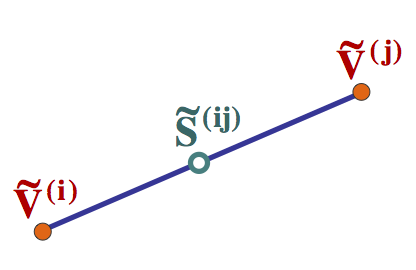}
\label{fig:fig20a}
}
\end{minipage}
\begin{minipage}[]{0.65\linewidth}
\subfigure[]
{
\includegraphics[clip,totalheight=1.2in]{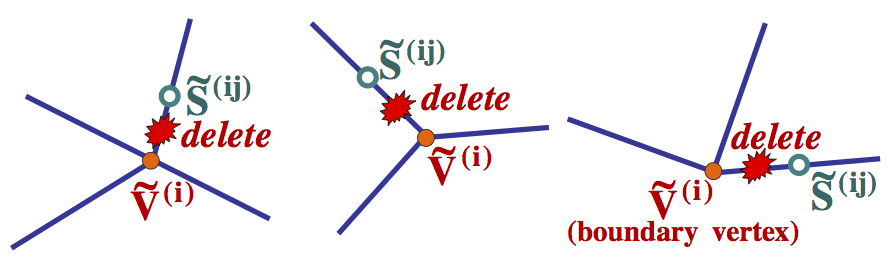}
\label{fig:fig20b}
}

\end{minipage}
\end{narrow}
\caption{An illustrations for the construction of the $D$-dependency graph.}
\label{fig:fig20}
\end{figure}

\begin{figure}[!pt]
\vspace{-0.7in}
\begin{narrow}{-0.5in}{-0.5in}
\begin{minipage}[]{0.33\linewidth}
\subfigure[]
{
\includegraphics[clip,totalheight=2.5in]{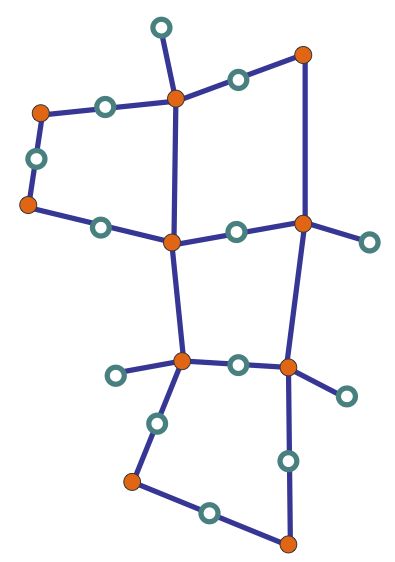}
\label{fig:fig21a}
}
\end{minipage}
\begin{minipage}[]{0.33\linewidth}
\subfigure[]
{
\includegraphics[clip,totalheight=2.5in]{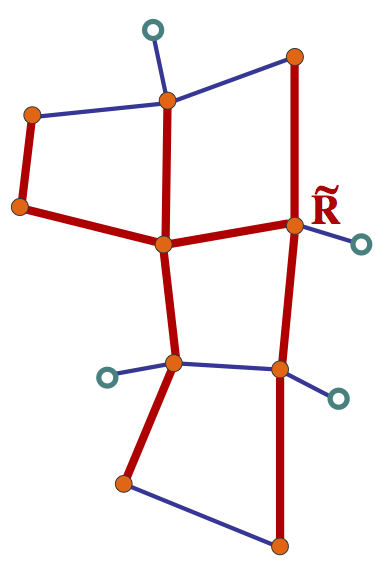}
\label{fig:fig21b}
}
\end{minipage}
\begin{minipage}[]{0.33\linewidth}
\subfigure[]
{
\includegraphics[clip,totalheight=2.5in]{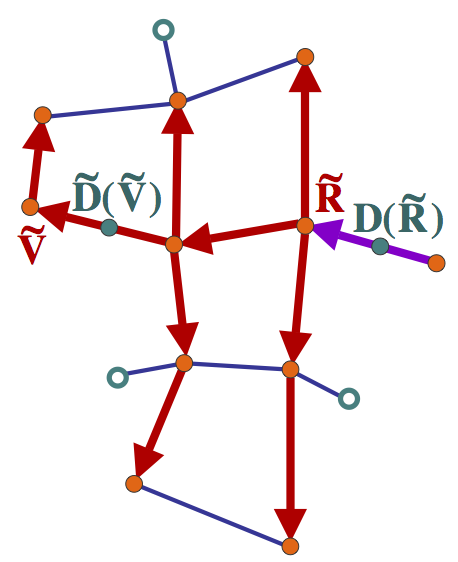}
\label{fig:fig21c}
}
\end{minipage}
\end{narrow}
\vspace{-0.2in}
\setcaptionwidth{5.4in}
\caption{Construction of $D$-dependency tree for a connected component of $D$-dependency graph.
\ \ (a) Connected component of $D$-dependency graph (here the structure of the component is not correct in the meaning that it does not correspond to any planar mesh, the Figure serves only 
as an illustration for Algorithm ~\ref{algorithm:dependency_forest_construction}).
\ \ (b) Spanning tree of the connected component.
\ \ (c) $D$-dependency tree of the connected component.
}
\label{fig:fig21}
\end{figure}

\begin{figure}[!pb]
\vspace{-0.7in}
\begin{narrow}{-0.5in}{-0.5in}
\begin{minipage}[]{0.33\linewidth}
\subfigure[]
{
\includegraphics[clip,totalheight=2.5in]{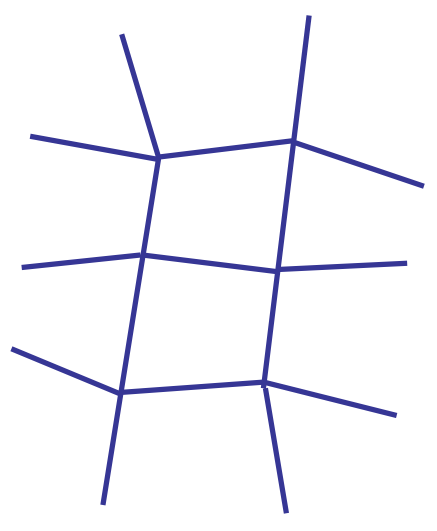}
\label{fig:fig21Aa}
}
\end{minipage}
\begin{minipage}[]{0.33\linewidth}
\subfigure[]
{
\includegraphics[clip,totalheight=2.5in]{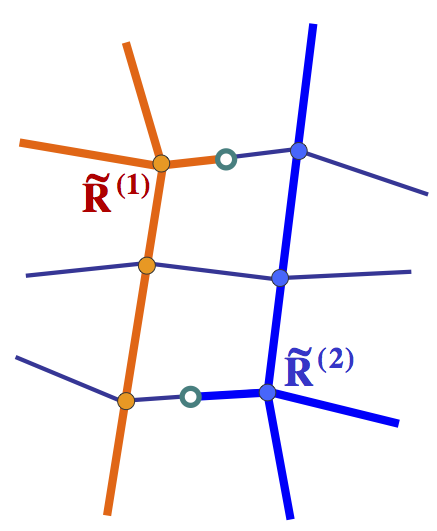}
\label{fig:fig21Ab}
}
\end{minipage}
\begin{minipage}[]{0.33\linewidth}
\subfigure[]
{
\includegraphics[clip,totalheight=2.5in]{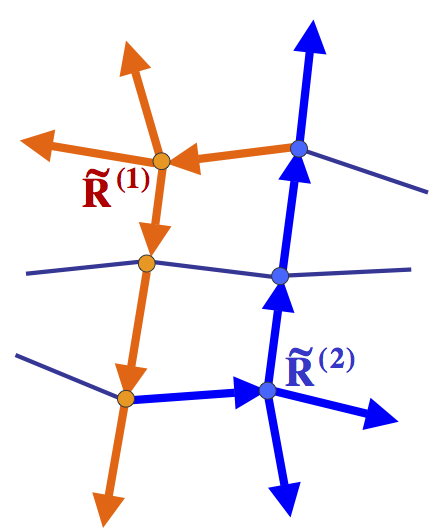}
\label{fig:fig21Ac}
}
\end{minipage}
\end{narrow}
\vspace{-0.2in}
\setcaptionwidth{5.4in}
\caption{An example of two different $D$-dependency trees connected by the directed edges at the root vertices.
\ \ (a) Planar mesh vertices and edges.
\ \ (b) $D$-dependency graph consisting of two connected components.
\ \ (c) $D$-dependency trees for two components of $D$-dependency graph.
}
\label{fig:fig21A}
\end{figure}

\begin{figure}[!pt]
\vspace{-0.7in}
\centering
\includegraphics[clip,totalheight=1.9in]{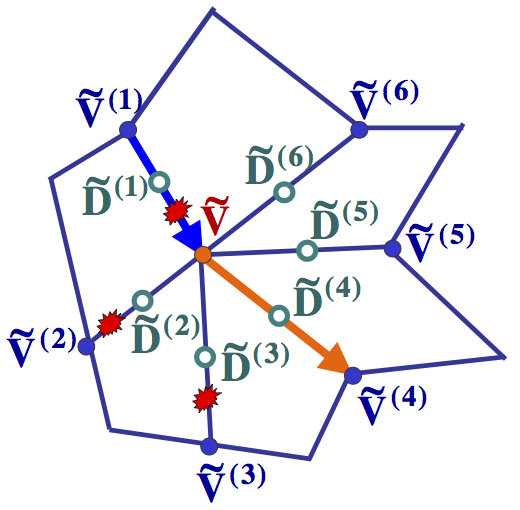}
\vspace{-0.1in}
\caption{An example of the classification of $D$-type control points.}
\label{fig:fig22}
\end{figure}

The examples given in Figures ~\ref{fig:fig39}, ~\ref{fig:fig40}, ~\ref{fig:fig41} provide illustrations to the formulas given in Equation ~\ref{eq:dimension_gmds_bilinear}. The Figures present instances of $\GMDS{5}$ (Figure ~\ref{fig:fig39}) and $\GMDS{4}$ (Figures ~\ref{fig:fig40} and ~\ref{fig:fig41}) for the same mesh.
The mesh contains one regular $4$-vertex $\T V^{(reg)}$ and the {\it "Projections Relation"} holds for edge $\T e^{("PR")}$. Control points which belong to $\GMDS{n}$ are marked by filled circles of different colors (orange, blue, green, violet and pink for $V$,$E$,$D$,$T$-type and {\it middle} control points respectively). Arrows in Figures ~\ref{fig:fig40} and ~\ref{fig:fig41} show which $D$-type control points are {\it assigned} to $D$-relevant vertices.

Figure ~\ref{fig:fig39} shows a {\it "pure"} instance of $\GMDS{5}$ (no "additional" constraints are applied). $\GMDS{5}$ contains $67$ control points (among $117$ control points of $\GCP{5}$), which precisely fits Equation ~\ref{eq:dimension_gmds_bilinear}. Figures ~\ref{fig:fig40} and ~\ref{fig:fig41} respectively show a "pure" instance of $\GMDS{4}$ and an instance which fits the clamped boundary condition. Both instances have the same dimension and contain $47$ control points.

An additional example, which allows to verify the correctness of  Equation ~\ref{eq:mds_dimension_bilinear}, is given in Section ~\ref{sect:MDS_example}.

\FloatBarrier

\FloatBarrier

\FloatBarrier

\begin{figure}[!ph]
\vspace{-0.2in}
\centering
\includegraphics[clip,totalheight=2.0in]{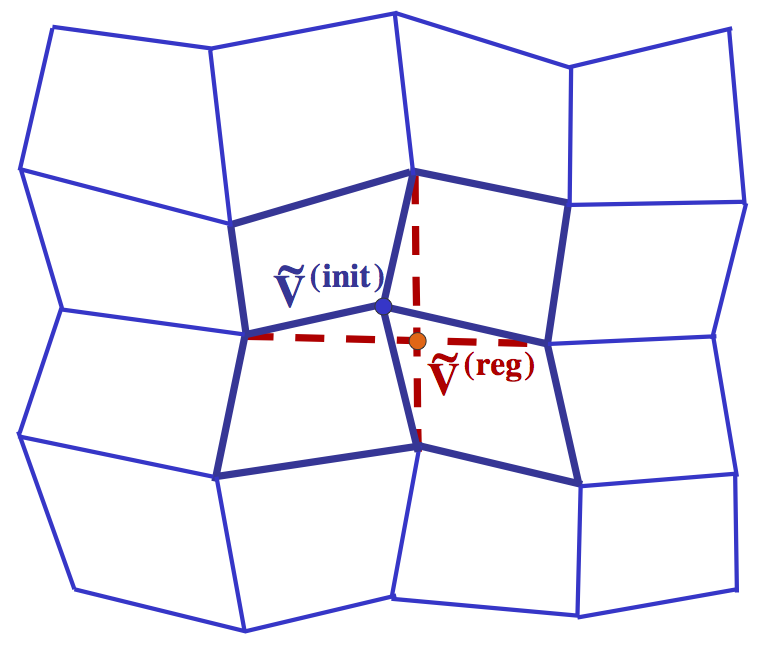}
\vspace{-0.15in}
\setcaptionwidth{5.3in}
\caption{Regularization of an inner $4$-vertex. Here $\T{V}^{(init)}$ is a vertex of the initial mesh and $\T{V}^{(reg)}$ is the corresponding regularized vertex. 
In $D$-dependency graph, $\T{V}^{(reg)}$ helps to obtain dangling half-edges for connected components of the adjacent vertices.}
\label{fig:fig25}
\end{figure}

\begin{figure}[!pb]
\vspace{-0.3in}
\begin{narrow}{-0.4in}{-0.4in}
\centering
\includegraphics[clip,height=3in]{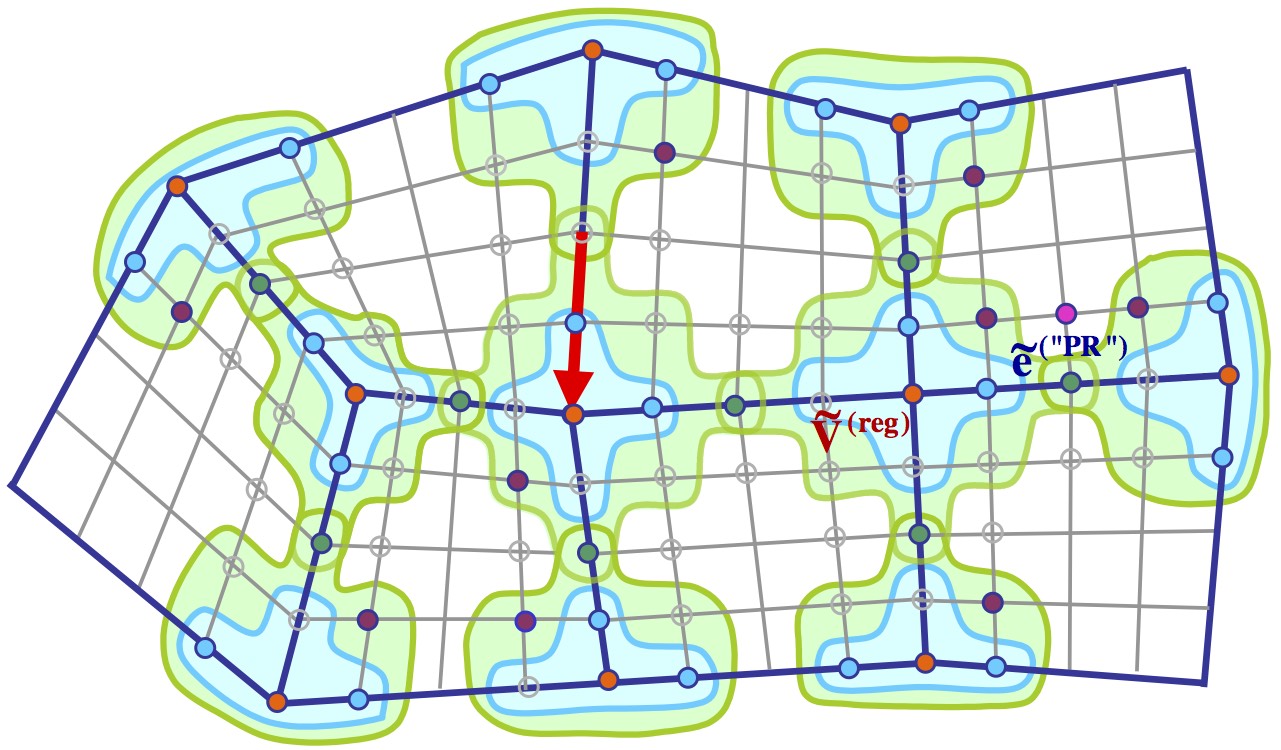}
\end{narrow}
\vspace{-0.1in}
\caption{An example of a "pure" global MDS $\MDS{4}(\T\Pi^{(bilinear)}$.}
\label{fig:fig40}
\end{figure}

\FloatBarrier

\eject
\part{MDS for a quadrilateral mesh with a smooth global boundary}
\label{part:smooth_boundary}

\begin{figure}[!pb]
\centering
\includegraphics[clip,totalheight=2.9in]{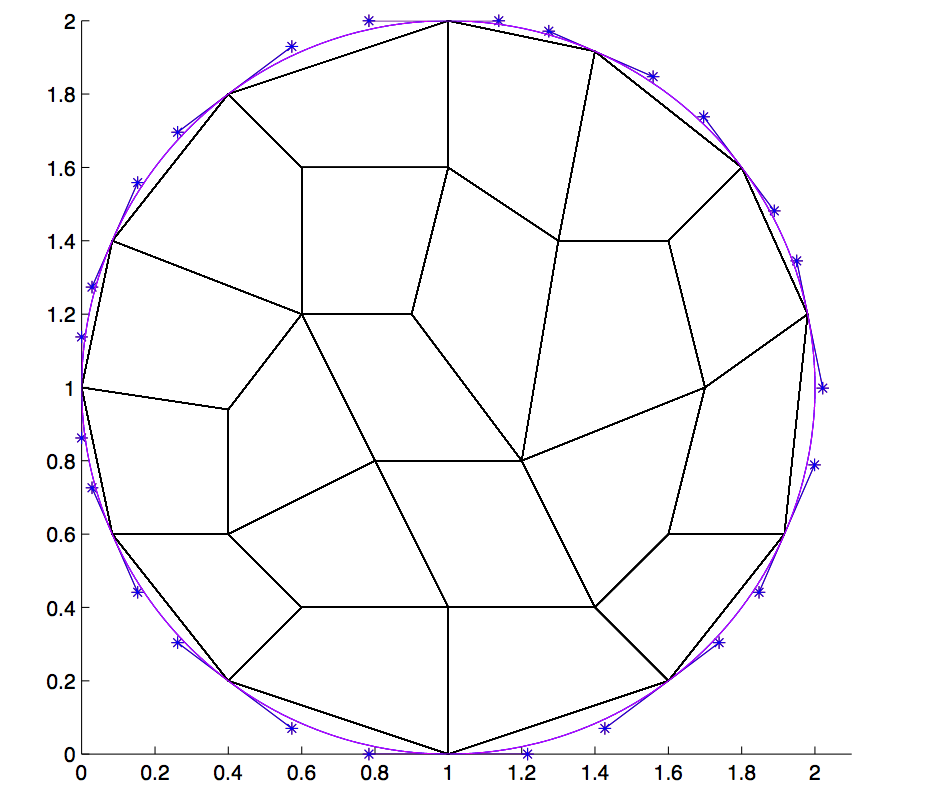}
\caption{An example  of a planar domain with a smooth global boundary.}
\label{fig:fig28}
\end{figure}

To construct higher order approximations/interpolations,  one needs also to handle  smooth boundaries without reducing them to polygonal lines!  This chapter deals with such constructions :
 planar meshes with a smooth global boundary (like the mesh shown in Figure ~\ref{fig:fig28}).
 
The bilinear in-plane parametrisation is no longer sufficient at the boundary. However, the global in-plane parametrisation
is constructed in such a manner, that the local templates should be changed only for the boundary vertices and for inner edges adjacent to the boundary. Therefore, a study of the continuity constraints for the edges adjacent to the boundary plays the principal role in construction of an MDS.

Like in the case of a polygonal global boundary, a {\it "pure"} $\MDS{n}$ is constructed for any $n\ge 4$ and for any mesh configuration.
Although it will be shown, that in the current case a boundary vertex is never $D$-relevant, a relatively high degree of the weight functions for the inner edges adjacent to the boundary may result in the failure of $\MDS{4}$ construction for some "additional" constraints.

\eject

\section{ Definitions,Mesh limitations and In-Plane parametrisation}
\label{sect:mesh_limitations_bicubic}
Let us recall some definitions from ...
\subsection{Definitions and In-plane parametrisation}
\label{sect:parametarization_bicubic}

\subsubsection{Points important for in-plane parametrisation of the boundary mesh elements in the case of a mesh with a smooth global boundary}
\label{subsect:notations_planar_points_cubic}
\paragraph{A single boundary element}
\begin{description}
\item
$\T{A},\T{B},\T{C},\T{D}$ - four vertices of a boundary mesh element, where $\T D$ and $\T C$ are the boundary vertices 
(see Figure ~\ref{fig:fig29a}).
\item
$\T E$, $\T F$ - two inner control points of the cubic boundary curve (see Figure ~\ref{fig:fig29a}).
\end{description}

\paragraph{Two adjacent boundary elements}
\begin{description}
\item
$\LL,\LL',\GG, \GG', \RR, \RR' $ -  vertices of two adjacent boundary mesh elements (see Figure ~\ref{fig:fig31})
\item
$\TTL',\TTR'$ - control points of the boundary curves adjacent to 
the vertex $\GG'$ (see Figure ~\ref{fig:fig31}).
\end{description}

The notion of  planar quadrilateral mesh is slightly generalised and it is assumed that 
\begin{itemize}
\item
The geometry of every inner edge is described by a straight segment.
\item
The edges along the global boundary of the planar domain have  a cubic parametric B\'ezier representation and the concatenation between any pair of adjacent edges is (non-degenerated) $G^1$-smooth. 
\end{itemize}
All basic assumptions listed in Section ~\ref{sect:mesh_limitations_bilinear} remain valid in the current 
case. Moreover, for simplicity, the mesh is assumed to have no corner elements, or in other words, every boundary vertex is supposed to have exactly one adjacent inner edge.

\subsection{In-plane parametrisation}
\label{sect:parametarization_bicubic}

\subsubsection{In-plane parametrisation for a boundary mesh element}
\label{subsect:cubic_param_definition}
There is no reason to change the bilinear type of parametrisation for inner mesh elements. For a boundary mesh element, we introduce 
a bicubic parametrisation $\T P(u,v)=\sum_{i,j=0}^3 \T P_{ij} B^3_{ij}$  in order to fit the boundary curve data. The following principles define the choice of parametrisation for a boundary mesh element.

\begin{figure}[!pb]
\vspace{-0.4in}
\begin{narrow}{-0.4in}{-0.4in}
\subfigure[]
{\includegraphics[clip,height=2.0in]{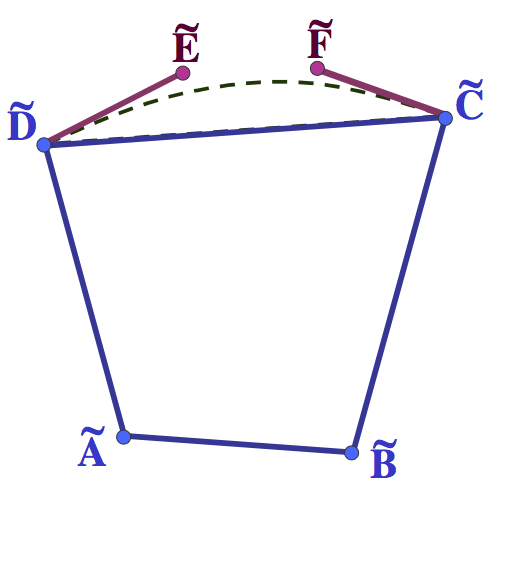}
\label{fig:fig29a}}
\subfigure[]
{\includegraphics[clip,height=2.0in]{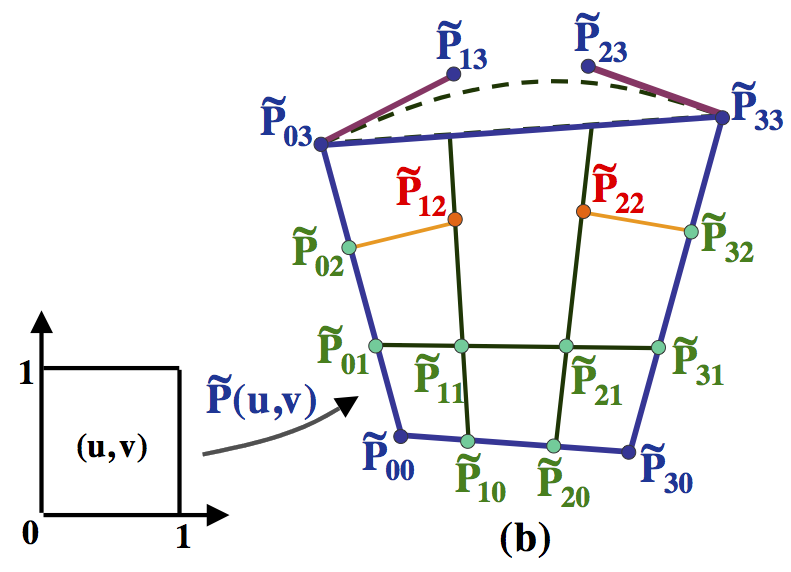}
\label{fig:fig29b}}
\end{narrow}
\vspace{-0.15in}
\setcaptionwidth{5.2in}
\caption{A planar mesh element, control points of a boundary curve and control points of in-plane parametrisation in the case of a mesh with a smooth global boundary.}
\label{fig:fig29}
\end{figure}

\begin{narrow}{-0.1in}{0.0in}
\noindent\vspace{-0.2in}
\begin{itemize}
\item
In order to make the definition of a global in-plane parametrisation possible,  it should be linear along the edge with two inner vertices and the parametrisations of two boundary elements should agree along edges with a common boundary vertex.
\item
The parametrisation should have a minimal possible influence on $G^1$-continuity equations for edges with two inner vertices.
It makes it possible to reuse the local templates constructed in the case of global bilinear in-plane parametrisation.
\item
Partial derivatives of the parametrisation along the edges with one boundary vertex should have the minimal possible degrees. It allows decreasing the number of  linear equations which are sufficient in order to guarantee that the $G^1$-continuity condition is satisfied.
\end{itemize}
\end{narrow}
Let the upper edge of a boundary mesh element lie on the global boundary of the planar domain. The following choice of the control points for in-plane parametrisation of the element is adopted (see Figure ~\ref{fig:fig29b}).

\begin{narrow}{-0.1in}{0.0in}
\noindent\vspace{-0.2in}
\begin{itemize}
\item
$\T P_{00},\T P_{30},\T P_{03}, \T P_{33}$ and $\T P_{13}, \T P_{23}$ are given respectively by the vertices of the planar element and by the control points of the cubic boundary curve. 
\item
$\T P_{10}, \T P_{20}, \T P_{01}, \T P_{11}, \T P_{21}, \T P_{31}$ are given by the degree elevation of the bilinear parametrisation up to degree $3$. This allows the  linearisation of the  $G^1$-continuity for the lower edge exactly as it was done in the case of a global bilinear in-plane parametrisation.
\item
$\T P_{02}, \T P_{32}$ are also given by the degree elevation of the bilinear parametrisation. It leads to the minimal possible degrees of the partial derivatives of in-plane parametrisation in directions along the left and the right edges.
\item
$\T P_{12}, \T P_{23}$ are chosen in such a manner, that the partial derivatives of in-plane parametrisation in the cross direction for the left and the right edges have the minimal possible formal degrees (degree $2$). Here it is assumed that prior to the choice of $\T P_{12}$ and $\T P_{23}$ all other control points of in-plane parametrisation are already fixed according to the three first items.
\end{itemize}
\end{narrow}
The explicit formulas for the control points of in-plane bicubic parametrisation of a boundary mesh element are given in Technical Lemma ~\ref{tl:control_points_param_bicubic} (see Appendix, Section ~\ref{sect:technical_lemmas}).

\subsubsection{Regularity of the bicubic in-plane parametrisation}
\label{subsect:parametrisation_regularity_bicubic}

Technical Lemma ~\ref{tl:regularity_cubic} (see Appendix, Section ~\ref{sect:technical_lemmas}) describes a sufficient geometrical condition for the regularity of the chosen in-plane parametrisation for the boundary mesh elements. A study of the regularity in the current work is restricted to this simple condition which is s quite natural , meaning that it is satisfied for the "non-degenerated" structures of the boundary elements. Of course, a more detailed analysis of the regularity conditions can be made.

\subsection{Global in-plane parametrisation $\T\Pi^{(bicubic)}$}
For every boundary element, let a bicubic in-plane parametrisation be constructed according to Technical Lemma ~\ref{tl:control_points_param_bicubic} (see Appendix, Section ~\ref{sect:technical_lemmas}) and let it be regular (the requirements of Technical Lemma ~\ref{tl:regularity_cubic} (see Appendix, Section ~\ref{sect:technical_lemmas}) are satisfied). Then according to Definition ~\ref{def:reg_param_all}, Paragraph ~\ref{subsubsect:math_formulation}, the collection of  bicubic parametrisations for boundary mesh elements and bilinear parametrisations for inner mesh elements defines a global regular in-plane parametrisation $\T\Pi^{(bicubic)}\in\PAR{3}$. Clearly, degree $3$ is just a formal degree of the parametrisation;  all inner mesh elements are of degree $1$ and parametrisations of the boundary mesh element have actual degree $3$ only along the global boundary.
\\
\\
\\
\eject
\section{Conventional weight functions}
\label{sect:weight_functions_bicubic}

\subsection{Weight functions for an edge with two inner vertices}

Bilinear parametrisation of the inner mesh elements and a special choice of the control points for in-plane parametrisation of the boundary mesh elements lead to the following Lemma.
 
\begin{lemma}
Let a mesh edge have two inner vertices. Then the conventional weight functions defined according to the global in-plane parametrisations $\T\Pi^{(bicubic)}$ and $\T\Pi^{(bilinear)}$, coincide. 
\label{lemma:wf_for_edge_two_inner_bicubic}
\end{lemma}

\subsection{Weight functions for an edge with one boundary vertex}

\subsubsection{Partial derivatives of in-plane parametrisations along the edge}

Let $\T P$ be the restriction of a global in-plane parametrisation $\T\Pi^{(bicubic)}$ on a boundary mesh element. For the edges with one boundary vertex, the first-order partial derivatives of $\T P$ have degree $0$ in the direction along the edges and formal degree $2$ in the cross direction. The explicit formulas of the partial derivatives are given in Technical Lemma ~\ref{tl:partial_derivatives_in_plane_bicubic} (see Appendix, Section ~\ref{sect:technical_lemmas}).

Let further $\T L(u,v)$ and $\T R(u,v)$ be the restrictions of the global in-plane parametrisation $\T\Pi^{(bicubic)}$ on the adjacent boundary elements. Although the formal degree of $\T L(u,v)$ and $\T R(u,v)$ is $m=3$, none of the partial derivatives $\T L_u$, $\T R_u$, $\T L_v$ has the full actual degree. Therefore, the linear system of $NumEqFormal=n+2m=n+6$ equations, described in Lemma ~\ref{lemma:linear_1},paragraph ~\ref{subsect:linearisation_G1_general} ,contains redundant equations. An analysis of the conventional weight functions allows to decrease the number of sufficient linear equations and finally to compute precisely the rank of the corresponding linear system. The polynomial representations of the partial derivatives $\T L_u$, $\T R_u$, $\T L_v$ in both B\'ezier and power basis play a very important role in the study of the weight functions.

\begin{lemma}
Let $\T L(u,v)$ and $\T R(u,v)$ be the restrictions of global in-plane
 parametrisation $\T\Pi^{(bicubic)}$ on two adjacent boundary mesh elements. Then the first-order partial derivatives of the in-plane parametrisations along the common edge have the  polynomial representations that follows in terms of the initial mesh data.
Here $\LL,\LL',\GG,\GG',\RR,\RR'$ and $\TTL'$ and $\TTR'$ denote respectively vertices of the elements and control points of the boundary curves (see Subsection ~\ref{subsect:notations_planar_points_cubic}  and  
~\ref{subsect:pd_planar_points_cubic} and 
Figure ~\ref{fig:fig31}  for the definitions).

\begin{narrow}{-0.15in}{0.0in}
\noindent\vspace{-0.2in}
\begin{itemize}
\item
$\T L_v=\GG'-\GG$ is a constant (polynomial of degree $0$).
\item
$\T L_u, \T R_u$ are polynomials of formal degree $2$. The coefficients of $\T L_u$ and $\T R_u$ with respect to B\'ezier and power basis are given below  (cf section~\ref{subsect:pd_planar_points_cubic}):

\begin{equation}
\hspace{-0.4in}
\begin{array}{lcl}
\T\lambda_0\!\!=\!\GG\!-\!\LL & \ \ &
\T\lambda^{(power)}_0\!\!=\!\T\lambda_0\!=\!\GG\!-\!\LL\cr
\T\lambda_1\!\!=\!\frac{1}{2}[(\GG\!-\!\LL)\!+\!(\GG'\!-\!\LL')]&& 
\T\lambda^{(power)}_1\!\!=\!2(\T\lambda_1\!-\!\T\lambda_2)\!=\!
(\GG'\!-\!\GG)\!-\!(\LL'\!-\!\LL)\cr
\T\lambda_2\!\!=\!3(\GG'\!-\!\TTL') &&
\T\lambda^{(power)}_2\!\!=\!\T\lambda_0\!\!-\!2\T\lambda_1\!\!+\!\T\lambda_2\!=\!
3(\GG'\!\!-\!\TTL')\!-\!(\GG'\!\!-\!\LL')\cr
\T\rho_0\!\!=\!\RR\!-\!\GG &&
\T\rho^{(power)}_0\!\!=\!\T\rho_0\!=\!\RR\!-\!\GG\cr
\T\rho_1\!\!=\!\frac{1}{2}[(\RR\!-\!\GG)\!+\!(\RR'\!-\!\GG')]&&
\T\rho^{(power)}_1\!\!=\!2(\T\rho_1\!-\!\T\rho_2)\!=\!
(\RR'\!-\!\RR)\!-\!(\GG'\!-\!\GG)\cr
\T\rho_2\!\!=\!3(\TTR'\!-\!\GG') &&
\T\rho^{(power)}_2\!\!=\!\T\rho_0\!-\!2\T\rho_1\!+\!\T\rho_2\!=\!
3(\TTR'\!\!-\!\GG')\!-\!(\RR'\!\!-\!\GG')
\end{array}
\end{equation}
\end{itemize}
\end{narrow}
\label{lemma:coeff_partial_bicubic}
\end{lemma}
Vectors $\T\lambda_i,\T\rho_i$ $(i=0,\ldots,2)$ are shown in Figure ~\ref{fig:fig32a}. 
Figure ~\ref{fig:fig32b} shows vectors $\T L_v$,\ $-\T\lambda^{(power)}_0$,\ 
$-\T\lambda^{(power)}_1+\T L_v$,\ $\T\lambda^{(power)}_2$,\ 
$\T\rho^{(power)}_0$,\ $\T\rho^{(power)}_1+\T L_v$,\ $\T\rho^{(power)}_2$. These vectors define the actual degrees of the weight functions and minimal number of  linear equations which are sufficient in order to satisfy the $G^1$-continuity condition. 

\subsubsection{Weight functions}
\label{subsect:def_weight_functions_bicubic}
Lemma ~\ref{lemma:coeff_weight_bicubic}
describes relations between coefficients of the weight functions and coefficients of the partial derivatives $\T L_u,\T R_u,
\T L_v$. Explicit formulas for coefficients of the weight functions in terms of the initial mesh data with respect to the  B\'ezier and the power basis immediately follow from Lemmas 
~\ref{lemma:coeff_partial_bicubic} and ~\ref{lemma:coeff_weight_bicubic}. Although the explicit formulas are useful in the following analysis, they become relatively complicated in the case of bicubic parametrisation and are given in Technical Lemmas ~\ref{tl:weight_coeff_cubic_explicit} and ~\ref{tl:weight_coeff_cubic_Bezier_power} (see Appendix, Section ~\ref{sect:technical_lemmas}). 

\begin{lemma}
Let $L(u,v)$, $R(u,v)$ be the restrictions of the global in-plane parametrisation $\T\Pi^{(bicubic)}$ on two adjacent boundary mesh elements. Then for the common edge of the elements, conventional weight functions $c(v)$, $l(v)$ and $ r(v)$ have formal degrees $4$, $2$ and $4$ respectively. Relations between the coefficients of the weight functions with respect to the B\'ezier and power bases and coefficients of the partial derivatives $\T L_u,\T R_u,\T L_v$ are given below.
\begin{equation}
\hspace{-0.4in}
\begin{array}{lcl}
l_i=\<\T \rho_i,\GG'-\GG\> & 
\ \ \ \ \ \ \ \ \ \ \ \ \ \ \ \ \ \ \ \ \ &
l^{(power)}_i=\<\T \rho^{(power)}_i,\GG'-\GG\> \cr
r_i=-\<\T \lambda_i,\GG'-\GG\>& &
r^{(power)}_i=-\<\T \lambda^{(power)}_i,\GG'-\GG\>
\end{array}
\label{eq:lr_lambda_rho}
\end{equation}
\begin{equation}
\hspace{-0.4in}
\begin{array}{lcl}
c_k=\frac{1}{\Cnk{4}{k}}
\!\!\!
\sum\limits_{\tiny\begin{array}{c}
{i\!+\!j\!=\!k}\cr{0\!\le\! i,j\!\le\! 2}\end{array}}
\!\!\!\!\!\!
\Cnk{2}{i}\Cnk{2}{j}\<\T\lambda_i,\T\rho_j\> & \ \ \ \ &
c^{(power)}_k=
\!\!\!\!\!
\sum\limits_{\tiny\begin{array}{c}
{i\!+\!j\!=\!k}\cr{0\!\le\! i,j\!\le\! 2}\end{array}}
\!\!\!\!\!\!
\<\T\lambda^{(power)}_i,\T\rho^{(power)}_j\>
\end{array}
\label{eq:c_lambda_rho}
\end{equation}
Here $i=0,\ldots,2$ in formulas for coefficient of $l(v)$ and $r(v)$ and $k=0,\ldots,4$ in formulas for coefficients of $c(v)$.
\label{lemma:coeff_weight_bicubic}
\end{lemma}

\begin{note}
From Lemma ~\ref{lemma:coeff_weight_bicubic} and according to the geometrical meaning of $\T\lambda_0$, $\T\lambda_2$, $\T\rho_0$, $\T\rho_2$ (see Lemma ~\ref{lemma:coeff_partial_bicubic}) it follows that
\begin{equation}
\begin{array}{lll}
l_0=l^{(power)}_0\neq 0,&\ \ \ \ \ &l_2\neq 0\cr
r_0=r^{(power)}_0\neq 0,&          &r_2\neq 0
\end{array}
\end{equation}
\vspace{-0.05in}
\begin{equation}
\ c_4=0
\end{equation}
\label{note:coeff_neq_zero}
\end{note}
A very important relation between the actual degrees of the weight functions is given in Lemma ~\ref{lemma:c_lr_degree_relation} (see Subsection ~\ref{subsect:def_weight_functions} for definitions).

\begin{lemma}
Let the conventional weight functions $c(v)$, $l(v)$, $r(v)$ for the common edge of two adjacent boundary mesh elements be defined according to the global in-plane parametrisation $\T\Pi^{(bicubic)}$. Then the actual degrees of the weight functions are connected by the following relation:
\begin{equation}
deg(c)\le max\_deg(l,r)+2
\end{equation}
\label{lemma:c_lr_degree_relation}
\end{lemma}
\vspace{-0.1in}
{\bf Proof} See Appendix, Section ~\ref{sect:proofs}.

\subsubsection{Geometrical meaning of the actual degrees of the weight functions}
\label{subsect:geometrical_meaning_weight_bicubic}
The different geometrical configurations of two adjacent boundary elements result in the  actual different degrees of the weight functions. Examples shown in Figure ~\ref{fig:fig33} provide some intuition regarding the geometrical meaning of the actual degrees of the weight functions.

\begin{itemize}
\item[(a)]
Example of degrees  $(4,2,2)$ for the weight functions $c,l,r$.
Vectors $\T\lambda^{(power)}_2$ and $\T\rho^{(power)}_2$
are not parallel to each other and both of them are not parallel to the vector $\GG'-\GG$. 
\item[(b)]
Example of degrees  $(3,2,2)$ for the weight functions $c,l,r$.
Vectors $\T\lambda^{(power)}_2$ and $\T\rho^{(power)}_2$
are parallel to each other but they are not parallel to the vector $\GG'-\GG$. $deg(c)\le 3$ {\it if and only if} $\T\lambda^{(power)}_2$ and $\T\rho^{(power)}_2$ are parallel.
\item[(c)]
Example of degrees  $(3,1,1)$ for the weight functions $c,l,r$.
Three vectors $\T\lambda^{(power)}_2$, $\T\rho^{(power)}_2$ and $\GG'-\GG$ are parallel. $max\_deg(l,r)\le 1$ {\it if and only if} vectors $\T\lambda^{(power)}_2$ and $\T\rho^{(power)}_2$ are parallel to $\GG'-\GG$. 
\item[(d)]
Example of  degrees $(2,1,1)$ for the weight functions $c,l,r$.
Vectors $\T\lambda^{(power)}_2$, $\T\rho^{(power)}_2$ and $\GG'-\GG$ are parallel. Equality $\T\lambda^{(power)}_1 = \T\rho^{(power)}_1$ implies that $deg(c)\le 2$. 
\item[(e)]
Example of degrees $(2,0,0)$ for the weight functions $c,l,r$.
Vectors $\T\lambda^{(power)}_2$, $\T\rho^{(power)}_2$, $\T\lambda^{(power)}_1$, $\T\rho^{(power)}_1$ and $\GG'-\GG$ are parallel. $l^{(power)}_1=r^{(power)}_1=0$ {\it if and only if} vectors $\T\lambda^{(power)}_2$, $\T\rho^{(power)}_2$, $\T\lambda^{(power)}_1$, $\T\rho^{(power)}_1$ are parallel to $\GG'-\GG$.
\end{itemize}
The examples show that different geometrical configurations lead to a variety of the triples $(deg(c),deg(l),deg(r))$. The actual degree of the weight function $c$ is not uniquely defined by the  degrees of the weight functions $l$ and $r$.

\begin{figure}[!ph]
\begin{narrow}{-0.2in}{-0.1in}
\centering
\includegraphics[clip,totalheight=2.0in]{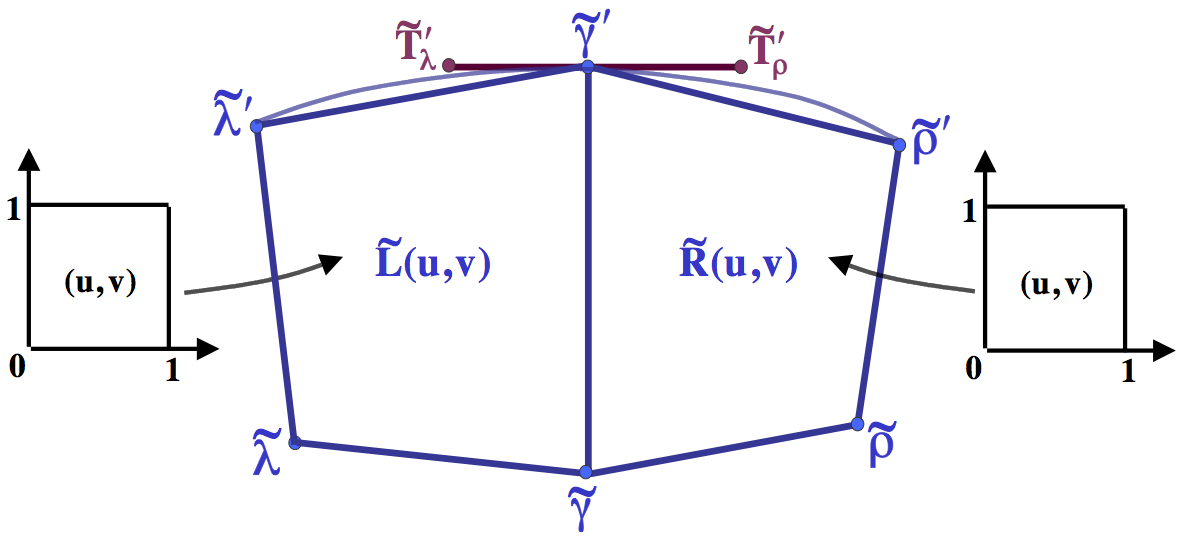}
\end{narrow}
\vspace{-0.2in}
\caption{Two adjacent boundary mesh elements in the case of a mesh with a smooth global boundary.}
\label{fig:fig31}
\end{figure}

\begin{figure}[!pb]
\vspace{-0.3in}
\begin{narrow}{-0.7in}{-0.7in}
\begin{minipage}[]{0.5\linewidth}
\subfigure[]
{
\includegraphics[clip,totalheight=2.1in]{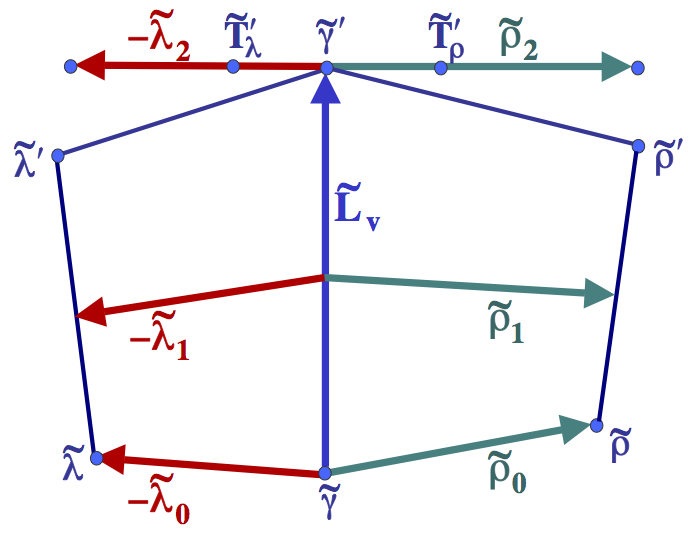}
\label{fig:fig32a}
}
\end{minipage}
\begin{minipage}[]{0.5\linewidth}
\subfigure[]
{
\includegraphics[clip,totalheight=2.1in]{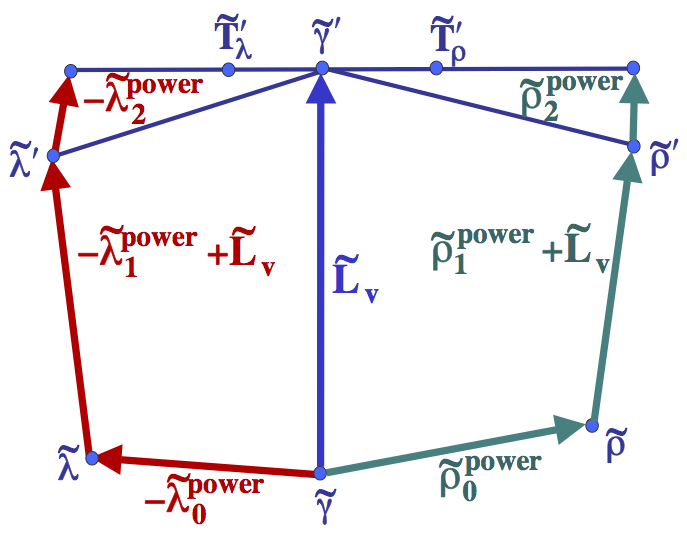}
\label{fig:fig32b}
}
\end{minipage}
\end{narrow}
\vspace{-0.2in}
\begin{narrow}{-0.25in}{0.0in}
\caption{Coefficients of $\T L_u,\T R_u,\T L_v$ with respect to B\'ezier and to power bases for two adjacent mesh elements in the case of global in-plane parametrisation $\T\Pi^{(bicubic)}$.}
\end{narrow}
\label{fig:fig32}
\end{figure}

\FloatBarrier

\begin{figure}[!pt]
\vspace{-0.5in}
\centering
\begin{narrow}{-0.9in}{-0.9in}
\begin{minipage}[]{0.33\linewidth}
\subfigure[(4,2,2)]
{
\includegraphics[clip,totalheight=1.8in]{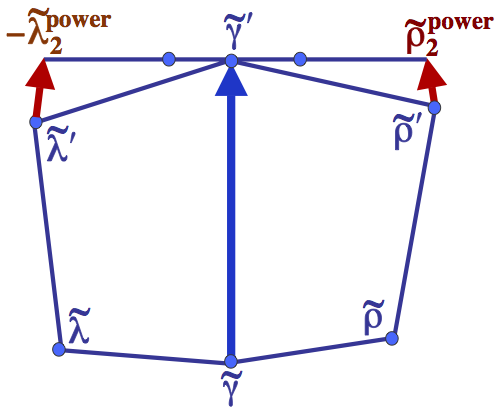}
\label{fig:fig33a}
}
\end{minipage}
\begin{minipage}[]{0.33\linewidth}
\subfigure[(3,2,2)]
{
\includegraphics[clip,totalheight=1.8in]{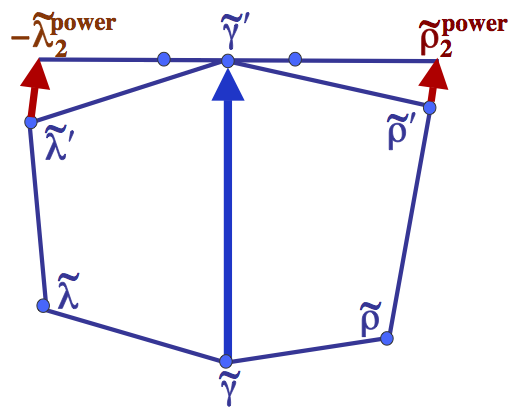}
\label{fig:fig33b}
}
\end{minipage}
\begin{minipage}[]{0.33\linewidth}
\subfigure[(3,1,1)]
{
\includegraphics[clip,totalheight=1.8in]{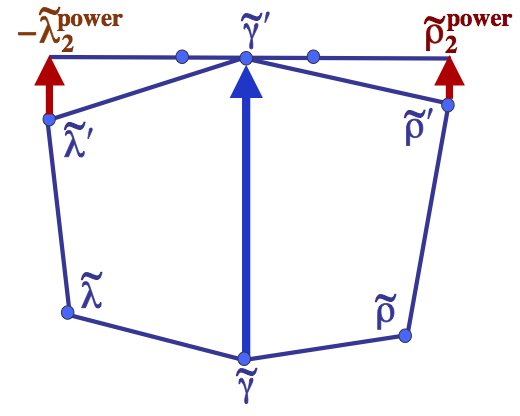}
\label{fig:fig33c}
}
\end{minipage}\\
\vspace{-0.2in}
\centering
\begin{minipage}[]{0.4\linewidth}
\subfigure[(2,1,1)]
{
\includegraphics[clip,totalheight=1.8in]{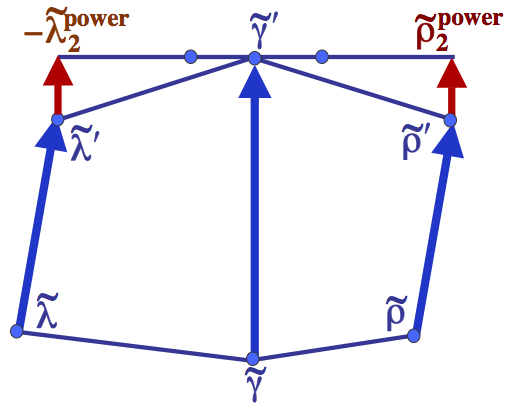}
\label{fig:fig33d}
}
\end{minipage}
\begin{minipage}[]{0.4\linewidth}
\subfigure[(2,0,0)]
{
\includegraphics[clip,totalheight=1.8in]{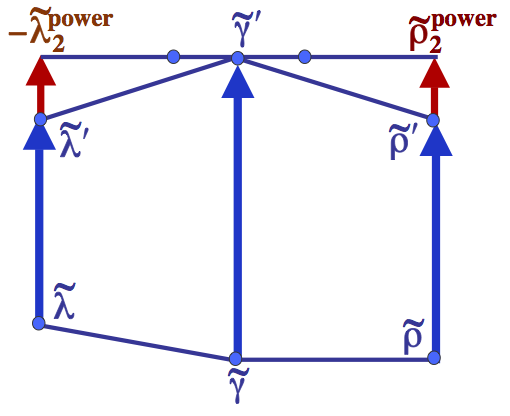}
\label{fig:fig33e}
}
\end{minipage}
\end{narrow}
\vspace{-0.2in}
\caption{Different geometrical configurations of two adjacent 
boundary mesh elements lead to different actual degrees of the conventional weight functions.}
\label{fig:fig33}
\vspace{0.8in}
\end{figure}

\FloatBarrier

\eject
\section{Linear form of $G^1$-continuity conditions}
\label{sect:linearisation_G1_bicubic}

\subsection{$G^1$-continuity conditions for an edge with two inner vertices}

Lemma ~\ref{lemma:wf_for_edge_two_inner_bicubic} clearly implies that the following Lemma holds.

\begin{lemma}
For any edge with two inner vertices, the linear system of $G^1$-continuity equations in the case of global in-plane parametrisation $\T\Pi^{(bicubic)}$ remains unchanged with respect to the case of bilinear global in-plane parametrisation $\T\Pi^{(bilinear)}$.
\label{lemma:G1_edge_with_two_inner_vertices_bicubic}
\end{lemma}

\subsection{$G^1$-continuity conditions for an edge with one boundary vertex}

\subsubsection{Formal construction of the linear system of equations}

In the case of global in-plane parametrisation $\T\Pi^{(bicubic)}$, the conventional weight functions $c(v)$, $l(v)$, $r(v)$ have formal degrees $4$, $2$, $2$ for an edge with one boundary vertex (see Lemma ~\ref{lemma:coeff_weight_bicubic}). Therefore the linearized $G^1$-continuity equation (Theorem 
~\ref{theorem:G1_equiv_Z_equation}, Equation ~\ref{eq:G1_Z_equation}) has the following form
\begin{equation}
\hspace{-0.4in}
\begin{array}{l}
\sum\limits_{i=0}^n \dL_i B^n_i \sum\limits_{j=0}^2 l_j B^2_j+
\sum\limits_{i=0}^n \dR_i B^n_i \sum\limits_{j=0}^2 r_j B^2_j+
\sum\limits_{i=0}^{n-1} \dC_i B^n_i 
\sum\limits_{j=0}^4 c_j B^4_j=0\cr
\end{array}
\label{eq:G1_cubic_diff_degrees}
\end{equation}
Unlike the case of global bilinear in-plane parametrisation $\T\Pi^{(bilinear)}$, not all summands of the last equation have the same formal degrees: the degree of the two first summands is $n+2$ while the degree of the last one is $n+3$. Of course, this difficulty can be easily overcome by application of the standard degree elevation to the first two summands. On the other hand it is important to remember, that $n+2$ and $n+3$ are just the {\it formal} degrees and that the {\it actual} degrees of the summands may coincide. For example in the case when the actual degrees of the weight functions $c,l,r$ are equal to $(3,2,2)$ the degree elevation is unnecessary, its application will lead to a redundant equation in the resulting linear system.

Generally, there are two different approaches to the analysis of Equation ~\ref{eq:G1_cubic_diff_degrees}. In the first one, all possible triples of the actual degrees of the weight functions should be considered separately. This approach "starts from geometry" and leads to the different algebraic systems, which correspond to the different geometric configurations. 
The second approach is a much more formal algebraic one. It starts by degree elevation, and then analyses the unique system of equations with the aid of the algebraic tools. This analysis eventually leads to different sub-cases, which of course correspond to different geometrical configurations.

Like in the case of global bilinear in-plane parametrisation, the second approach is adopted in the current work. However, while performing an algebraic analysis, one always should be aware of the alternative way which allows to verify the correctness of the results. For example in the case of  degrees $(3,2,2)$ for the weight functions, one may verify that the resulting system of the linear equations does not have the full row rank.

Application of the degree elevation to the first two summands of Equation ~\ref{eq:G1_cubic_diff_degrees} and writing down the coefficients of the resulting B\'ezier polynomial of degree $n+3$ lead to the following Lemma.
\begin{lemma}
In the case of global in-plane parametrisation $\T\Pi^{(bicubic)}$, the system of the following $n+4$ linear equations is sufficient in order to guarantee a $G^1$-smooth concatenation between two boundary patches 
\vspace{-0.1in}
\begin{equation}
\hspace{-0.45in}
\begin{array}{l}
\sum_{\lmtT{0\!\le\!i\!\le\! n}{0\!\le\!j\!\le\!2}{i\!+\!j\!=\!s}}
\!\!\!
(\dL_i l_j+\dR_i r_j)\Cnk{n}{i}\Cnk{2}{j}\!+\!
\sum_{\lmtT{0\!\le\!i\!\le\!n}{0\!\le\!j\!\le\!2}{i\!+\!j\!=\!s\!-\!1}}
\!\!\!\!\!
(\dL_i l_j+\dR_i r_j)\Cnk{n}{i}\Cnk{2}{j}+\cr
\sum_{\lmtT{0\le i\le n-1}{0\le j\le 4}{i+j=s}}
\dC_i c_j\Cnk{n-1}{i}\Cnk{4}{j} = 0
\end{array}
\label{eq:G1_cubic_formal}
\end{equation}
\vspace{-0.1in}
where $s=0,\ldots,n+3$.
\label{lemma:G1_cubic_formal}
\end{lemma}

\subsubsection{An equivalent system of the linear equations}
\label{subsect:equivalent_system_bicubic}

A linear system of equations, which is equivalent to the system given in Lemma ~\ref{lemma:G1_cubic_formal} and has a clearer and more intuitive structure, will be constructed in order to simplify the analysis of MDS. The following notations lead to a more compact form of the linear equations.
 
\begin{definition}
Let
\begin{equation}
\hspace{-0.4in}
\begin{array}{ll}
\underline{sumLR_s}=\sum_{\lmtT{0\le i\le n}{0\le j\le 2}{i+j=s}}
(\dL_i l_j+\dR_i r_j)\Cnk{n}{i}\Cnk{2}{j} & 
{\rm for\ \ } s=0,\ldots,n+2\cr
sumLR_{-1}=sumLR_{n+3}=0, & \cr
\underline{sumC_s}=\sum_{\lmtT{0\le i\le n-1}{0\le j\le 4}{i+j=s}}
\dC_i c_j\Cnk{n-1}{i}\Cnk{4}{j} & 
{\rm for\ \ } s=0,\ldots,n+3
\end{array}
\end{equation}
\label{def:sumLR_sumC}
\end{definition}
\vspace{-0.15in}
and
\vspace{-0.1in}
\begin{equation}
\hspace{-0.4in}
\begin{array}{lcl}
\underline{"Eq(s)"} & \ \ \ \ \ \ \ \ &
sumLR_s = (-1)^{s+1}\sum\limits_{k=0}^s (-1)^k sumC_k\cr
\underline{"sumC\!-\!equation"} &  & \sum\limits_{k=0}^{n+3}(-1)^k sumC_k=0
\end{array}
\label{eq:def_Eq_sumC_bicubic}
\end{equation}
(An expression in the left side of the last equation will also be called {\it "sumC-equation"}).

\begin{lemma}
The following linear system of $n+4$ equations is sufficient in order to guarantee a $G^1$-smooth concatenation  between two boundary patches in the case of global in plane parametrisation $\T\Pi^{(bicubic)}$ \begin{equation}
\left\{
\begin{array}{ll}
{\it"Eq(s)"}         & s=0,\ldots,n+2 \cr
{\it"sumC\!-\!equation"} & 
\end{array}
\right.
\end{equation}

\label{lemma:Eq_sumC_linear_system}
\end{lemma}

\vspace{-0.05in}

{\bf Proof} See Appendix, Section ~\ref{sect:proofs}.

Note, that {\it "sumC-equation"} might be formally written as {\it "Eq(n+3)"}, which would lead to a homogeneous system of equations. It was decided to separate this equation because of its role in the construction of MDS and because of its special properties, some of which are listed in the next Subsection.

\subsubsection{Some important properties of {\it "sumC-equation"}}
\label{subsect:sumC_equation_properties}

Unlike the indexed equations, {\it "sumC-equation"} deals only with the "central" control points $C_j$ ($j=0,\ldots,n$), none of the "side" control points $(L_j,R_j)$ ($j=0,\ldots,n$) contributes to it. 

An another important property of {\it "sumC-equation"} is given in Lemma ~\ref{lemma:sumC_expansion}.
\begin{lemma}
Let equation {\it "C-equation"} be defined as follows
\begin{equation}
\hspace{-0.4in}
\begin{array}{lcl}
\underline{\it "C\!\!-\!\!equation"} & \ \ \ \ \ \ \ &  
\sum_{i=0}^{n-1}(-1)^i\Cnk{n-1}{i}\dC_i= 
-\sum_{i=0}^n (-1)^i \Cnk{n}{i}C_i = 0
\end{array}
\label{eq:C_equation_def}
\end{equation}
(An expression in the left side of the equation will also be called {\it "C-equation"}).
Then {\it "sumC-equation"} can be represented as
\begin{equation}
\hspace{-0.4in}
\begin{array}{lcl}
{\it "C\!\!-\!\!equation"}\ c^{(power)}_4 = 0
\end{array}
\label{eq:sumC_expansion}
\end{equation}
\label{lemma:sumC_expansion}
\end{lemma}
Lemma ~\ref{lemma:sumC_expansion} in particular means that {\it "sumC-equation"} is automatically satisfied when $deg(c)\le 3$. Moreover, it shows that if $deg(c)=4$ then values of $C_j$ $(j=0,\ldots,n)$ should necessarily satisfy the {\it "C-equation"}.

In addition, Equation ~\ref{eq:C_equation_def} implies that 
{\it "sumC-equation"} and {\it "C-equation"}  have a "global nature" in the meaning that every one of $C_j$ control points $(j=0,\ldots,n)$ participates in these equations with a non-zero coefficient.

\section{Local MDS}
\label{sect:local_mds_bicubic}

In the following discussion, it always will be assumed that edge with one boundary vertex is parametrized in such a manner, that $\bar C_0$ is the inner vertex of the edge and $\bar C_n$ is the boundary vertex. It implies, for example, that equations {\it "Eq(0)","Eq(1)"} relate to the inner vertex and equations {\it "Eq(n+1)", "Eq(n+2)"} to the boundary vertex of the edge (see Figure ~\ref{fig:fig34b}).

\subsection{Local templates for a separate vertex}

\subsubsection{Local templates for an inner vertex}
\label{subsect:local_mds_inner_vertex_bicubic}

\begin{lemma}
Consider global   $\T\Pi^{(bicubic)}$ parametrisation. Then for any inner vertex
\begin{itemize}
\item[\bf (1)]
{\it "Eq(0)"}-type equation remains unchanged with respect to the case of global bilinear in-plane parametrisation $\T\Pi^{(bilinear)}$.
\item[\bf (2)]
The couple of {\it "Eq(0)"}-type and {\it "Eq(1)"}-type equations is  equivalent to the couple of the corresponding equations in the case of global bilinear in-plane parametrisation $\T\Pi^{(bilinear)}$.
\end{itemize}
\label{lemma:Eq0_Eq1_bicubic}
\end{lemma}
\vspace{-0.1in}
{\bf Proof} See Appendix, Section ~\ref{sect:proofs}.

\noindent 
Lemma ~\ref{lemma:Eq0_Eq1_bicubic} leads to the following Conclusions

\begin{conclusion}
Consider a global in-plane parametrisation $\T\Pi^{(bicubic)}$. Then for any inner vertex, the local templates for classification of $V$,$E$,$D$,$T$-type control points remain unchanged with respect to the case of global bilinear in-plane parametrisation $\T\Pi^{(bilinear)}$. 

The templates are constructed precisely as described in Subsections ~\ref{subsect:local_VE_inner_vertex_bilinear} and ~\ref{subsect:local_DT_inner_vertex_bilinear} and are responsible for the satisfaction of {\it "Eq(0)"}-type and {\it "Eq(1)"}-type equations at the inner vertex.
\end{conclusion}

\begin{conclusion}
All theoretical results related to the Vertex Enclosure problem
(see Section ~\ref{subsect:theory_vertex_enclosure_bilinear}) remain valid in the case of global in-plane parametrisation $\T\Pi^{(bicubic)}$.
\end{conclusion}

\subsubsection{Local templates for a boundary vertex}
\label{subsect:local_mds_boundary_vertex_bicubic}

\paragraph{Formal representation of equations {\it "Eq(n+2)"} and {\it "Eq(n+1)"} for an edge with one boundary vertex}

In the case of a global in-plane parametrisation $\T\Pi^{(bicubic)}$, equations {\it "Eq(n+2)"} and {\it "Eq(n+1)"} have the following representations for an edge with one boundary vertex.

\noindent
\underline{{\it "Eq(n+2)"} equation.}
\begin{equation}
\dL_n l_2 +\dR_n r_2=0
\label{eq:Eq_n_plus_2_bicubic}
\end{equation}
Geometrically {\it "Eq(n+2)"} means that  the  control points $\bar L_n,\bar C_n,\bar R_n$ are colinear. 
Equation {\it "Eq(n+2)"} {\it never involves $C_{n-1}$} ($E$-type control point of the edge adjacent to the boundary vertex).

\noindent
\underline{{\it "Eq(n+1)"} equation.} 
\begin{equation}
n(\dL_{n-1} l_2 + \dR_{n-1} r_2)+ 
2(\dL_n l_1 + \dR_n r_1) +
4\dC_{n-1} c_3 = 0
\label{eq:Eq_n_plus_1_bicubic}
\end{equation}
This equation involves only control points lying on the mesh boundary ($L_n,C_n,R_n$) or adjacent to it 
($L_{n-1},C_{n-1},R_{n-1}$). The control point $C_{n-2}$ ($D$-type control point of the edge adjacent to the boundary vertex) {\it does not participate} in this equation.

\paragraph{Local templates and different types of "additional" constraints}

The choice of the local templates for a boundary vertex $\T V$ in the case of global in-plane parametrisation $\T\Pi^{(bicubic)}$ is based on the following two principles.
\begin{itemize}
\item
The template should include all control points adjacent to $\T V$ and participating in specified "additional" constraints.
\item
The template is responsible for the satisfaction of some sub-system of the indexed equations (the sub-system contains either the last indexed equation or the last pair of the indexed equation).
\item
The template should contain a minimal possible number of the control points in order to provide additional degrees of freedom to the "Middle" system, which becomes relatively complicated in the case of $\T\Pi^{(bicubic)}$ for edges with one boundary vertex.
\end{itemize}

\begin{figure}[!pt]
\vspace{-0.1in}
\begin{narrow}{0.0in}{0.0in}
\begin{minipage}[]{0.5\linewidth}
\subfigure[]
{
\includegraphics[clip,totalheight=1.0in]{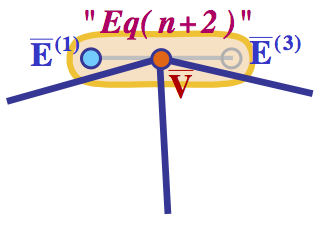}
\label{fig:fig34Aa}
}
\end{minipage}
\begin{minipage}[]{0.5\linewidth}
\subfigure[]
{
\includegraphics[clip,totalheight=1.0in]{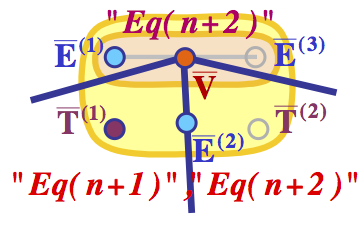}
\label{fig:fig34Ab}
}
\end{minipage}
\end{narrow}
\vspace{-0.15in}
\caption{Local templates for classification of $V$,$E$-type control points adjacent to a boundary vertex in the case of global in-plane parametrisation $\T\Pi^{(bicubic)}$.}
\label{fig:fig34A}
\end{figure}

\noindent
The following two local templates are defined:  
\begin{description}
\item[]
\underline{$TB0^{(bicubic)}$} (Figure ~\ref{fig:fig34Aa}). The template includes control points lying at the global boundary ($V$-type control point and two boundary $E$-type control points).
The template is responsible for the satisfaction of equation  {\it "Eq(n+2)"} and "additional" constraints, provided the "additional" constraints involve only control points lying at the global boundary.

A local MDS contains $\T V$ and one of $E$-type control points ($\T E^{(1)}$). Control point $\T E^{(3)}$ is dependent and dependency of $E^{(3)}$ is defined according to Equation ~\ref{eq:Eq_n_plus_2_bicubic}. 
\item[]
\underline{$TB1^{(bicubic)}$} (Figure ~\ref{fig:fig34Ab}). 
The template includes the control points lying at the global boundary and adjacent to it ($V$,$E$ and $T$-type control points). The template is responsible for the satisfaction of equations  {\it "Eq(n+2)"}, {\it "Eq(n+1)"} and for any type of considered "additional" constraints.

A local MDS contains $\T V$, one of $E$-type control points ($\T E^{(1)}$) and one of $T$-type control points ($\T T^{(1)}$). Control points $\T E^{(3)}$ and $\T T^{(2)}$ are dependent. First, dependency of $E^{(3)}$ is defined according to Equation ~\ref{eq:Eq_n_plus_2_bicubic} and then dependency of $T^{(2)}$ is defined according to Equation ~\ref{eq:Eq_n_plus_1_bicubic}.
\end{description}

\noindent
The template $TB0^{(bicubic)}$ is used when the "additional" constraints involve only control points lying at the global boundary {\it and} degree $4$ MDS is considered. In all other cases template $TB0^{(bicubic)}$ is used. 

\subsection{Local classification of the {\it middle} control points for a separate edge}
\label{subsect:local_mds_middle_bicubic}

\subsubsection{Local templates for an edge with two inner vertices}

Lemma ~\ref{lemma:G1_edge_with_two_inner_vertices_bicubic} implies that the following Lemma holds.

\begin{lemma}
Consider a global in-plane parametrisation $\T\Pi^{(bicubic)}$ and let  for any edge with two inner vertices the {\it "Middle"} system of equations and the {\it middle} set of the control points be defined according to Definition ~\ref{def:middle_bilinear}. 

Then for an edge with two boundary vertices, the local templates responsible for the classification of the {\it middle} control points are the same as for the global bilinear in-plane parametrisation $\T\Pi^{(bilinear)}$ (see Subsection ~\ref{subsect:templates_middle_bilinear}).
\label{lemma:unchanged_middle_bicubic}
\end{lemma}

\subsubsection{Local templates for an edge with one boundary vertex}
\label{subsect:local_mds_middle_edge_one_boundary_bicubic}

\paragraph{Definition of the {\it "Middle"} system of equations and {\it middle} control points for an edge with one boundary vertex}

The following notations are introduced in order to unify the description of the local MDS for $\FUN{4}$ and $\FUN{n}$, $n\ge 5$.

\begin{figure}[!pb]
\begin{narrow}{-0.5in}{-0.7in}
\begin{minipage}[]{0.31\linewidth}
\subfigure[]
{
\includegraphics[clip,totalheight=3.1in]{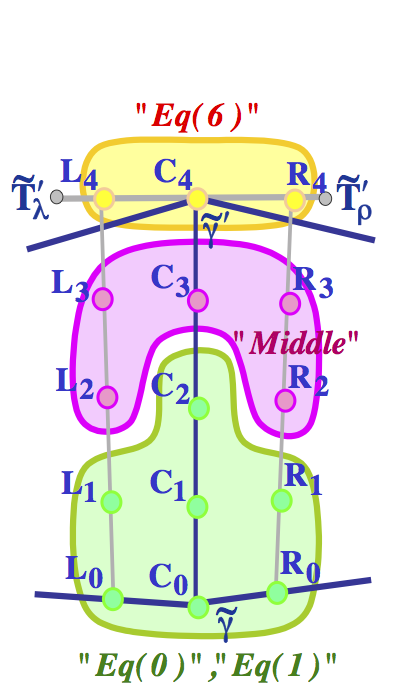}
\label{fig:fig34a}
}
\end{minipage}
\begin{minipage}[]{0.31\linewidth}
\subfigure[]
{
\includegraphics[clip,totalheight=3.1in]{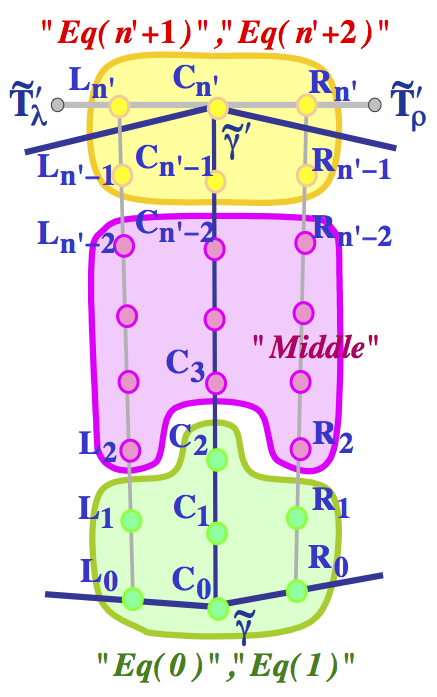}
\label{fig:fig34b}
}
\end{minipage}
\begin{minipage}[]{0.31\linewidth}
\subfigure[]
{
\includegraphics[clip,totalheight=3.1in]{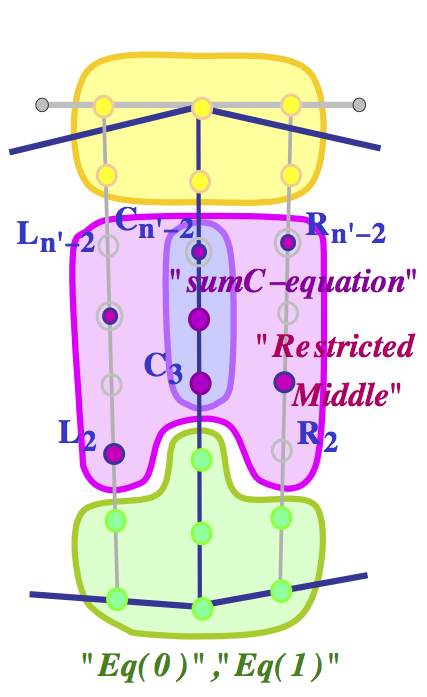}
\label{fig:fig34c}
}
\end{minipage}
\end{narrow}
\vspace{-0.15in}
\caption{Local templates for an edge with one boundary vertex in the case of global in-plane parametrisation $\T\Pi^{(bicubic)}$.}
\label{fig:fig34}
\vspace{-0.2in}
\end{figure}

\begin{definition}
Consider the functional space $\FUN{n}(\T\Pi^{(bicubic)})$, $n\ge 4$. For an edge with one boundary vertex, let
\vspace{-0.1in}
\begin{equation}
\underline{n'}=\left\{
\begin{array}{ll}
5 & if\ n=4\cr
n & if\ n\ge 5
\label{eq:deg_n_tag}
\end{array}
\right.
\vspace{-0.1in}
\end{equation}
\vspace{0.03in}
\noindent
The set of the \underline{middle} control points contains $2(n'-3)$ "side" {\it middle} control points $\T L_2,\ldots,\T L_{n'-2}$, $\T R_2,\ldots,\T R_{n'-2}$ and $n'-4$ "central" {\it middle} control points $\T C_3,\ldots,\T C_{n'-2}$ (see Figure ~\ref{fig:fig34a}, ~\ref{fig:fig34b}).

\vspace{0.03in}
\noindent
The \underline{"Middle"} system of equations consists of $n'$ equations: $n'-1$ indexed equations {\it "Eq(s)"} for $s=2,\ldots,n'$ and {\it "sumC-equation"} (see Equation ~\ref{eq:def_Eq_sumC_bicubic}).

\vspace{0.03in}
\noindent
The \underline{"Restricted Middle"} system of equations consists of $n'-1$ indexed equations {\it "Eq(s)"} for $s=2,\ldots,n'$.

\label{def:middle_edge_one_boundary_bicubic}
\end{definition}

Theorems ~\ref{theorem:C_equation_sufficiency} and ~\ref{theorem:LR_classification} show that the "central" {\it middle} control points are responsible for the consistency of the {\it "Middle"} system of equations and for the satisfaction of {\it "sumC-equation"}. Classification of the "side" {\it middle} control points is made after classification of the "central" control points and depends on the rank of the {\it "Restricted Middle"} system. 

The classification is based on the consistency and the rank analysis of the {\it "Middle"} system. The relatively complicated algebraic proofs of the main results deal with nice algebraic and geometric dependencies, which define the {\it actual} degrees of the weight functions and finally define the structure and the rank of the {\it "Middle"} system.

\paragraph{A necessary and sufficient   condition for the consistency of the {\it "Middle"} system. Classification of the "central" {\it middle} control points.}

\begin{theorem}
Given a global in-plane parametrisation $\T\Pi^{(bicubic)}$ ,for an edge with one boundary vertex, let all non-middle control points be classified and equations {\it "Eq(0)"}, {\it "Eq(1)"}, {\it "Eq(n'+1)"} for $n\ge 4$ and equation {\it "Eq(n'+2)"} for $n\ge 5$ be satisfied. Then for any $n\ge 4$
the satisfaction of {\it "C-equation"} in a case when 
\vspace{-0.1in}
\begin{equation}
deg(c)-max\_deg(l,r)=2
\end{equation}
\vspace{-0.1in}
\begin{itemize}
\item[{\bf(1)}]
Is a sufficient condition for the satisfaction of the {\it "sumC-equation"}.
\item[{\bf(2)}]
Is a necessary and sufficient condition for the consistency of the {\it "Middle"} system of equations.
\end{itemize}
\label{theorem:C_equation_sufficiency}
\end{theorem}
\vspace{-0.15in}
{\bf Proof} See Appendix, Section ~\ref{sect:proofs}.

\begin{lemma}
Consider the global in-plane parametrisation $\T\Pi^{(bicubic)}$.
For an edge with one boundary vertex, the "central" {\it middle} control points are classified prior to and independently of the classification of the "side" middle control points according to the following rules.
\begin{itemize}
\item
$\T C_t$ for $t=3,\ldots,n'-3$ (if any) are basic (free) control points.
\item
$\T C_{n'-2}$ is a dependent control point with the dependency defined by {\it "C-equation"} if
$deg(c)-max\_deg(l,r)=2$ and basic (free) otherwise.
Note, that $C_{n'-2}$ may depend only on the basic "central" {\it middle} control points
\end{itemize}
According to Theorem ~\ref{theorem:C_equation_sufficiency}, the classification guarantees that the {\it "sumC-equation"} is satisfied and that the {\it "Middle"} system of equations is consistent.
\label{lemma:C_classification}
\end{lemma}

\paragraph{Classification of the "side" {\it middle} control points}

Lemma ~\ref{lemma:C_classification} implies that after the classification of the "central" {\it middle} control points, {\it "sumC-equation"} is satisfied and it is sufficient to study the {\it "Restricted Middle"} system of equations.
The system is known to be consistent as a sub-system of the consistent {\it "Middle"} system. The {\it "Restricted Middle"} system contains $2(n'-3)$ non-classified "side" middle control points $(\T L_t,\T R_t)$ for $t=2,\ldots,n'-2$. It remains to study the rank of the system and to classify the "side" middle control points accordingly.

The next Definition generalises the {\it "Projections Relation"} (see Definition ~\ref{def:projection_relation}) that plays an important role in the rank analysis of the {\it "Middle"} system in case of global bilinear in-plane parametrisation $\T\Pi^{(bilinear)}$.

\begin{definition}
For an edge with one boundary vertex, let $l(v)$, $r(v)$ be conventional weight functions corresponding to a global in-plane parametrisation $\T\Pi^{(bicubic)}$. Scalars $g^{(ij)}$ ($i,j\in\{0,1,2\}$) are defined by the next formula in terms of coefficients of the weight functions
\begin{equation}
\underline{g^{(ij)}} = l_i r_j - l_j r_i
\end{equation}
\label{def:g_ij}
\end{definition}

The following three principal kinds of relations between 
$g^{(01)}$, $g^{(02)}$ and $g^{(12)}$ will be considered. These relations  correspond to the different possible values of the rank of the {\it "Restricted Middle"} system of equations in terms of non-classified "side" {\it middle} control points.

\begin{theorem}
Given the  global in-plane parametrisation $\T\Pi^{(bicubic)}$, for an edge with one boundary vertex, let
\begin{itemize}
\item[-]
All non-middle control points be classified and equations {\it "Eq(0)"}, {\it "Eq(1)"}, {\it "Eq(n'+1)"} for $n\ge 4$ and equation {\it "Eq(n'+2)"} for $n\ge 5$ be satisfied.
\item[-]
The "central" middle control points be classified according to Lemma ~\ref{lemma:C_classification}. 
\end{itemize}
Then the rank of the {\it "Restricted Middle"} system in terms of the "side" middle control points and the classification of the "side" middle control points depend on the relations between the $g^{(ij)}$
\begin{description}
\item[\bf "CASE 1"] 
If 
\vspace{-0.1in}
\begin{equation}
g^{(01)} g^{(12)} g^{(02)}\neq 0\ \ \  and\ \ \ 
\{g^{(02)}\}^2=4g^{(01)}g^{(12)}
\end{equation}
\vspace{-0.1in}
then $rank=n'-2$ and there are $n'-4$ basic (free) and $n'-2$ dependent control points among $2(n'-3)$ "side" middle control points.
\item[\bf "CASE 2"]
If
\vspace{-0.1in}
\begin{equation}
g^{(01)}=g^{(12)}=g^{(02)}=0
\end{equation}
\vspace{-0.1in}
then $rank=n'-3$ and there are $n'-3$ basic (free) and $n'-3$ dependent "side" middle control points; in every pair $(\T L_t,\T R_t)$ $(t=2,\ldots,n'-2)$ one control point is basic and one is dependent. 
\item[\bf "CASE 0"]
If none of previous conditions on $g^{(01)}, g^{(12)}, g^{(02)}$ are satisfied, then $rank=n'-1$ and there are $n'-5$ basic (free) and $n'-1$ dependent control points among $2(n'-3)$ "side" middle control points.
\end{description}
\label{theorem:LR_classification}
\end{theorem}
\vspace{-0.15in}
Like in the case of global bilinear in-plane parametrisation $\T\Pi^{(bilinear)}$, explicit dependencies between the control points are described in the proof of  the theorem.

\paragraph{The different kinds of  local templates}

A local template for an edge with one boundary vertex includes the {\it middle} control points and is responsible for the satisfaction of the {\it "Middle"} system of equations. Application of a local template assumes that all non-middle control points for the edge are classified and equations {\it "Eq(0)"}, {\it "Eq(1)"}, {\it "Eq(n'+1)"} for $n\ge 4$ and equation {\it "Eq(n'+2)"} for $n\ge 5$ are satisfied.

The choice of a local template depends entirely on the geometrical structure of the boundary mesh elements adjacent to the edge; a unique template corresponds to every geometrical configuration.
Classification of the {\it middle} control points is made as follows (see Figure ~\ref{fig:fig34c}).
\begin{description}
\item
{\it "Central" middle control points.} These control points are classified prior to and independently of the classification of "side" middle control points.
According to Lemma ~\ref{lemma:C_classification}, the local MDS contains all $n'-4$ "central" middle control points if $deg(c)-max\_deg(l,r)\neq 2$ and contains $n'-5$ control points otherwise.
\item
{\it "Side" middle control points.}
Classification of these control points depends on the rank of the {\it "Restricted Middle"} system. The number of basic control points varies from $n'-5$ (in {\bf "CASE 0"}) to $n'-3$ (in {\bf "CASE 2"}) (see Theorem ~\ref{theorem:LR_classification}).
\end{description}

\paragraph{Example of  local MDS for $n=4$ and $n=5$}
\label{subsect:classification_middle_particular_bicubic}

In both cases ($n=4$ or $n=5$) a local template contains the same control points: $\T C_3$, $\T L_2$, $\T R_2$, $\T L_3$, $\T R_3$.

Lemma ~\ref{lemma:C_classification} implies that the control point $\T C_3$ is dependent (and dependency of $C_3$ is defined by {\it "C-equation"}) if $deg(c)-max\_deg(l,r)=2$ and belongs to local MDS otherwise.

In {\bf "CASE 0"} all control points $\T L_2,\T R_2,\T L_3,\T R_3$ are dependent. In {\bf "CASE 1"} one of these four control points belongs to local MDS and others are dependent. In {\bf "CASE 2"} in every pair $(\T L_2,\T R_2)$ and $(\T L_3,\T R_3)$ one of the control points belongs to local MDS and the second one is dependent.

\eject

\FloatBarrier

\eject
\section{Global MDS}
\label{sect:global_mds_bicubic}
\subsection{Algorithm for construction}

Algorithm ~\ref{algorithm:mds_global_bicubic} defines order of application of the local templates described in Section ~\ref{sect:local_mds_bicubic}. The Algorithm allows to "put together" local templates without contradiction and to define global MDS.
The Algorithm reuses Algorithms ~\ref{algorithm:mds_global_deg_ge_5_bilinear} and ~\ref{algorithm:D_T_classification_degree_4} for the inner part of the mesh; local modifications for edges with one boundary vertex are made at the last stage of the classification process.

It was decided to present the common algorithm for construction of global MDS of degree $n=4$ and $n\ge 5$, because strategies for $n=4$ and $n\ge 5$ are based on the same principles and differ only in small details.

\begin{note}
MDS of degree $4$ can be constructed in situations, when "additional" constraints applied at a boundary vertex, involve only control points lying on the global boundary (local template $TB0^{(bicubic)}$, see Subsection ~\ref{subsect:local_mds_boundary_vertex_bicubic}, can be used for every boundary vertex). MDS of degree $n\ge 5$ can be constructed in any case.
\end{note}

\begin{algorithm}
\label{algorithm:mds_global_bicubic}
\end{algorithm}
\centerline{\underline{\it Algorithm for construction of global MDS}}
\centerline{\underline{\it in the case of global in-plane parametrisation $\T\Pi^{(bicubic)}$}}

\begin{description}
\item
[{\bf "Stage 1"\ \ }]
Classify all {\it $E$,$V$,$D$,$T$-type} control points adjacent to {\it inner vertices} exactly as it was done in the case of global bilinear in-plane parametrisation $\T\Pi^{(bilinear)}$ (see Subsections ~\ref{subsect:local_VE_inner_vertex_bilinear}, ~\ref{subsect:local_DT_inner_vertex_bilinear}, and Section ~\ref{sect:global_mds_bilinear}).

At the end of this stage, all {\it {\it "Eq(0)"}-type,{\it "Eq(1)"}-type} equations for all {\it inner vertices} are satisfied.

For an edge with one boundary vertex, control points $(\T L_0,\T C_0,\T R_0)$, $(\T L_1,\T C_1,\T R_1)$ and $\T C_2$ are classified and equations {\it "Eq(0)"} and {\it "Eq(1)"} are satisfied (see Figure ~\ref{fig:fig34a}, ~\ref{fig:fig34b}).

\item
[{\bf "Stage 2"\ \ }]
For every edge with two inner vertices, classify the {\it middle} control points exactly as it was done in the case of global bilinear in-plane parametrisation $\T\Pi^{(bilinear)}$ (see Subsection ~\ref{subsect:local_mds_middle_bilinear})

At the end of this stage {\it all} control points, lying on or adjacent to some {\it edge with two inner vertices}, and {\it $D$-type} control points adjacent to the {\it inner vertex} of an edge with one boundary vertex, are classified. All $G^1$-continuity equations excluding the {\it "Middle"} system of equations and indexed equations {\it "Eq(n'+1)"} for $n\ge 4$ and {\it "Eq(n'+2)"} for $n\ge 5$ are satisfied.

\item
[{\bf "Stage 3"}]
For every edge with one inner vertex, perform the following two steps of classification.
\begin{narrow}{0.4in}{0.0in}
\begin{itemize}
\item[{\bf "Step 1"}]
At the boundary vertex apply local template (see Subsection ~\ref{subsect:local_mds_boundary_vertex_bicubic})
\begin{itemize}
\item[-]
$TB0^{(bicubic)}$ if a given "additional" constraints, applied at the vertex, involve only control points lying at the global boundary {\it and} MDS of degree $4$ is considered.
\item[-]
$TB1^{(bicubic)}$ if MDS of degree $n\ge 5$ is considered.
\end{itemize}

At the end of this step, for every edge with one boundary vertex, all control points excluding {\it middle} control points are classified and indexed equations {\it "Eq(0)", "Eq(1)", "Eq(n'+1)"} for $n\ge 4$ and equation {\it "Eq(n'+2)"} for $n\ge 5$ are satisfied.

\item[{\bf "Step 2"}]
Apply local template for classification of the {\it middle} control points.

At the end of this step all control points are classified and all 
$G^1$-continuity equations are satisfied.
\end{itemize}
\end{narrow}
\end{description}

\noindent
{\bf Proof of the correctness} See Appendix, Section ~\ref{sect:proofs}.

\subsection{Existence of global MDS of degree $n\ge 5$}

Algorithm ~\ref{algorithm:mds_global_bicubic} allows to conclude that the following Lemma holds.

\begin{lemma}
For any $n\ge 5$ and for any type of "additional" constraints, there exists such an instance $\MDS{n}(\T\Pi^{(bicubic)})$ of the global MDS, that the instance fits the "additional" constraints.
\label{lemma:deg_n_ge_5_mds_existence_bicubic}
\end{lemma}
Dimensionality of $\MDS{n}(\T\Pi^{(bicubic)})$ is given in Subsection ~\ref{subsect:mds_dimensionality_bicubic}.

\subsection{Existence of global MDS of degree $n=4$. Analysis of different "additional" constraints.}

\begin{lemma}
\hspace{0.4in}
\begin{itemize}
\item[{\bf (1)}]
Global {\it "pure"} MDS\ \ $\MDS{4}(\T\Pi^{(bicubic)})$ (MDS which relates to $G^1$ continuity constraints alone) is well defined for any mesh configuration.

\item[{\bf (2)}]
Let an "additional" constraint, applied at any boundary vertex, involve only the control points lying at the global boundary. Then an instance of $\MDS{4}(\T\Pi^{(bicubic)}),$ which fits the "additional" constraint, exists for any mesh configuration. In particular, a suitable $\MDS{4}(\T\Pi^{(bicubic)})$ always exists for a (vertex)[boundary curve]-interpolation problem under the assumption that the (tangent plane)- interpolation is not required and for a partial differential equation with a simply-supported boundary condition. 
\end{itemize}
\label{lemma:deg_4_mds_existence_bicubic}
\end{lemma}
Dimensionality of $\MDS{4}(\T\Pi^{(bicubic)})$ is presented in Subsection ~\ref{subsect:mds_dimensionality_bicubic}.

There are two important differences with respect to the case of global bilinear in-plane parametrisation $\T\Pi^{(bilinear)}$.
The first evident difference is that an instance of global MDS of degree $4$, which fits the (tangent plane)-interpolation condition, is no longer known to exist. The second difference is that one never tries to construct MDS of degree $4$ in cases when its existence is not guaranteed. In a case of the (tangent plane)-interpolation and a case of the clamped boundary condition either MDS of degree $n\ge 5$ of mixed MDS of degrees $4$ and $5$ (see Part ~\ref{part:mds_mixed_4_5}) should be constructed.

Although in the case of global in-plane parametrisation $\T\Pi^{(bicubic)}$, no $D$-relevant boundary vertices are considered (and therefore the main reason for failure of MDS construction in the case of global bilinear parametrisation no longer exists), one encounters another problem. linearisation of the $G^1$-continuity condition for an edge with one boundary vertex leads to a more complicated {\it "Middle"} system of equations. Now construction of global MDS of degree $4$ may fail because of insufficient number of the control points which are necessary in order to satisfy the {\it "Middle"} system of equations for an edge with one boundary vertex. The existence of a suitable MDS of degree $4$ in the current case depends even more strongly on the kind of "additional" constraints than in the case of global bilinear in-plane parametrisation. 

\subsection{Dimensionality of MDS}
\label{subsect:mds_dimensionality_bicubic}

\begin{theorem}
For global in-plane parametrisation $\T\Pi^{(bicubic)}$ and for any $n\ge 4$ dimension of $\GMDS{n}$ (subset of MDS which participate in the $G^1$-continuity condition) and dimension of $\MDS{n}$ (full dimension of MDS) are given by the following formulas (see Definition ~\ref{def:cp_subsets} and Lemma  ~\ref{lemma:GMDS_importance})
\begin{equation}
\hspace{-0.4in}
\begin{array}{ll}
|\GMDS{n}| = &
3|Vert_{non-corner}|+
(2n-7) |Edge_{inner}|+\cr
&
|Vert_{\twolines{inner}{4\!-\!regular}}|+
|Edge_{\twolines{two\ inner\ vertices,}
                {"Projections\ Relation"\ holds}}|-\cr
&
|Edge_{\twolines{one\ boundary\ vertex,}
                {deg(c)-max\_deg(l,r)=2}}|+\cr
&
|Edge_{\twolines{one\ boundary\ vertex,}{{\bf "CASE 1"}}}|+
2|Edge_{\twolines{one\ boundary\ vertex,}{{\bf "CASE 2"}}}|
\end{array}
\label{eq:dimension_gmds_bicubic}
\end{equation}
\begin{equation}
\hspace{-0.4in}
\begin{array}{lll}
|\MDS{n}| = &|\GMDS{n}|+&|\FCP{n}|=\cr
&|\GMDS{n}|+&(n\!-\!3)^2|Face_{inner}|+
(n\!-\!3)(n\!-\!1)|Face_{\twolines{boundary}{non\!-\!corner}}|+\cr
&&(n\!-\!1)^2|Face_{corner}|
\end{array}
\label{eq:mds_dimension_bicubic}
\end{equation}
\label{theorem:mds_dimension_bicubic}
\end{theorem}

\eject
\part{Mixed MDS of degrees $4$ and $5$}
\label{part:mds_mixed_4_5}

\begin{figure}[!ph]
\vspace{-0.2in}
\centering
\subfigure[]
{
\includegraphics[clip,totalheight=2.3in]{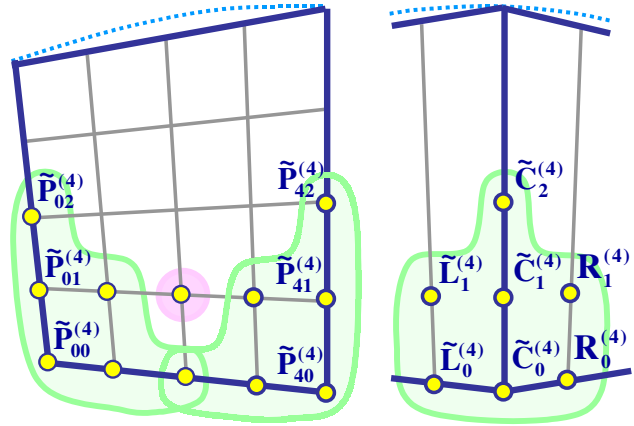}
\label{fig:fig36a}
}\\
\subfigure[]
{
\includegraphics[clip,totalheight=2.3in]{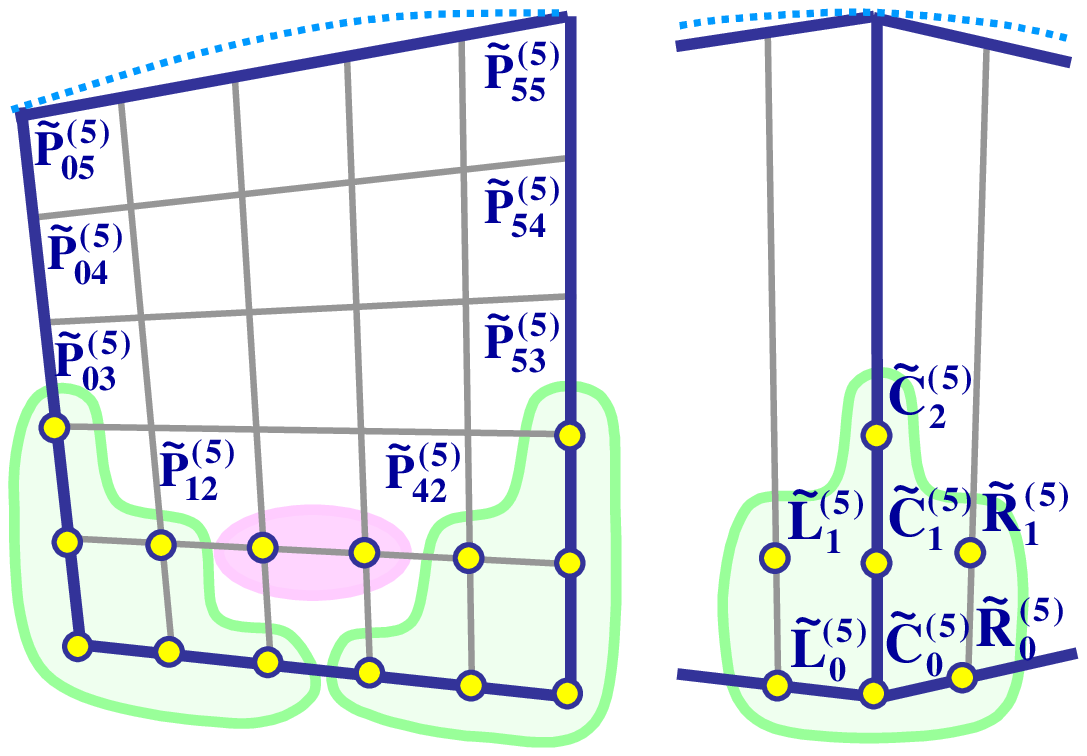}
\label{fig:fig36b}
}
\vspace{-0.15in}
\caption{An illustration for degree elevation for a boundary mesh element.}
\label{fig:fig36}
\end{figure}

\paragraph{Definition of mixed MDS}
For a global regular in-plane parametrisation $\T\Pi\in\PAR{m}$, $m<4$, mixed functional space \underline{$\FUN{4,5}(\T\Pi)$} is defined
according to Definition ~\ref{def:space_fun} where item {\bf (2)} is substituted with the following item
\begin{itemize}
\item[{\bf (2')}]
{\it For a mesh element $\T q\in\T{\cal Q}$, $Z$-coordinate of the restriction $\bar Q =\bar\Psi|_{\T q}$ belongs to ${\cal POL}^{(4)}$ if $\T q$ is an inner mesh element and belongs to ${\cal POL}^{(5)}$ if $\T q$ is a boundary mesh element.}
\end{itemize}
Definition of MDS (see Definition ~\ref{def:CPMDS}) also requires just a minor modification. 
\begin{definition}
\underline{$\CP{4,5}$} is composed of
\begin{itemize}
\item[-]
Subset of $\CP{4}$ for inner mesh elements.
\item[-]
"Degree 4" control points $\T P^{(4)}_{i,j}$, $i=0,\ldots,4$, $j=0,1$ and $\T P^{(4)}_{0,2}$, $\T P^{(4)}_{0,2}$ and
 
"degree 5" control points $\T P^{(5)}_{i,2}$, $i=1,\ldots,4$ and $\T P^{(5)}_{i,j}$, $i=0,\ldots,5$, $j=3,4,5$ for boundary mesh elements (see Figure ~\ref{fig:fig36}). 
\end{itemize}
\label{def:cp_mixed_4_5}
\end{definition}
\begin{definition}
MDS \underline{$\MDS{4,5}$} is such a minimal subset of $\CP{4,5}$ that for any function $\bar\Psi\in\FUN{4,5}$, equality to zero of $Z$-coordinates corresponding to all control points from the subset implies that $\bar\Psi\equiv 0$.
\label{def:mds_mixed_4_5}
\end{definition}

\paragraph{Algorithm for the construction of mixed MDS\ $\MDS{4,5}$}

For both mesh with a polygonal global boundary (Part ~\ref{part:linear_boundary}) and mesh with a smooth global boundary
(Part ~\ref{part:smooth_boundary}) a possible failure to construct MDS of degree $4$, which fits a given "additional" constraints, is always related to the boundary mesh elements. 
Algorithm ~\ref{algorithm:mixed_mds} (see Appendix, Section ~\ref{sect:algorithms}) shows that in a case of failure it is sufficient to elevate the degree up to $5$ for the boundary patches only, leaving the degree $4$ for all inner patches. In other words, it is sufficient to pass over to consideration of mixed MDS\ $\MDS{4,5}$.

Clearly, Algorithm ~\ref{algorithm:mixed_mds}
results in lower than $|\MDS{5}|$ dimension of MDS. However, Algorithm ~\ref{algorithm:mixed_mds} requires the global classification of $D$-type control points, like Algorithms for construction of MDS of degree $4$. Although a formal proof of the correctness of the Algorithm is not given, one can easily verify it by her/himself. 

\begin{algorithm}
\label{algorithm:mixed_mds}
See Appendix, Section ~\ref{sect:algorithms}.
\end{algorithm}

\paragraph{The existence of a suitable instance of mixed MDS for any type of "additional" constraints}

Algorithm ~\ref{algorithm:mixed_mds} allows to conclude that the following Lemma takes place.
\begin{lemma}
Let a planar mesh have either a polygonal or a smooth global boundary and global in-plane parametrisation $\T\Pi=\T\Pi^{(bilinear)}$ or $\T\Pi=\T\Pi^{(bicubic)}$ be considered.
Then for any type of "additional" constraints, there exists such an instance of mixed MDS\ $\MDS{4,5}(\T\Pi)$, that the instance fits the "additional" constraints.
\label{lemma:mixed_mds_existance}
\end{lemma}

\eject
\part{Computational examples}
\label{part:illustrative_examples}

The current Part illustrates theoretical analysis given in the previous Parts by a few computational examples. The first example presents MDS for some irregular $4$-element mesh. 

The rest of the examples study more complicated irregular meshes; these examples present approximate solutions of the Thin Plate problem under different boundary conditions. The main purpose of the examples is to verify the existence of MDS that fits the boundary constraints and to demonstrate quality of the resulting surface. Although some statistics connected to numerical precision of the solution are given, an error estimate of the approximate solution remains beyond the scope of the current research. 

\section{Examples of MDS}
\label{sect:MDS_example}

\begin{figure}[!pb]
\centering
\begin{narrow}{-0.8in}{-0.5in}
\begin{minipage}[]{0.5\linewidth}
\subfigure[]
{
\includegraphics[clip,totalheight=3in]{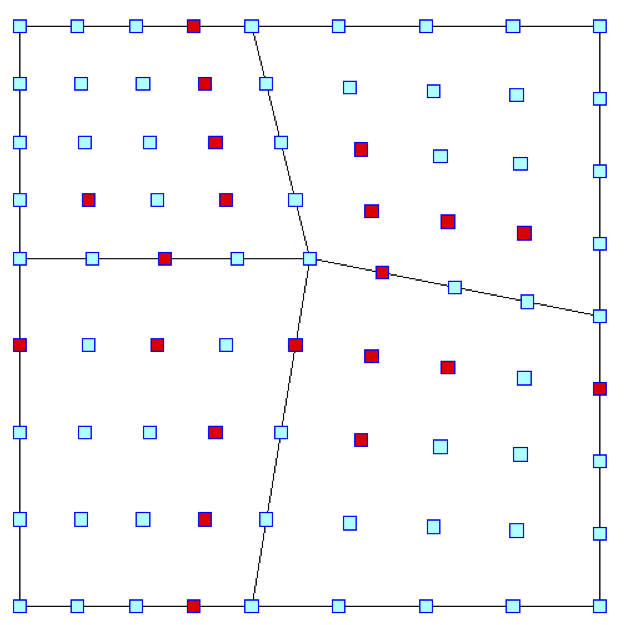}
\label{fig:fig_tp1a}
}
\end{minipage}
\hspace{0.1in}
\begin{minipage}[]{0.5\linewidth}
\subfigure[]
{
\includegraphics[clip,totalheight=3in]{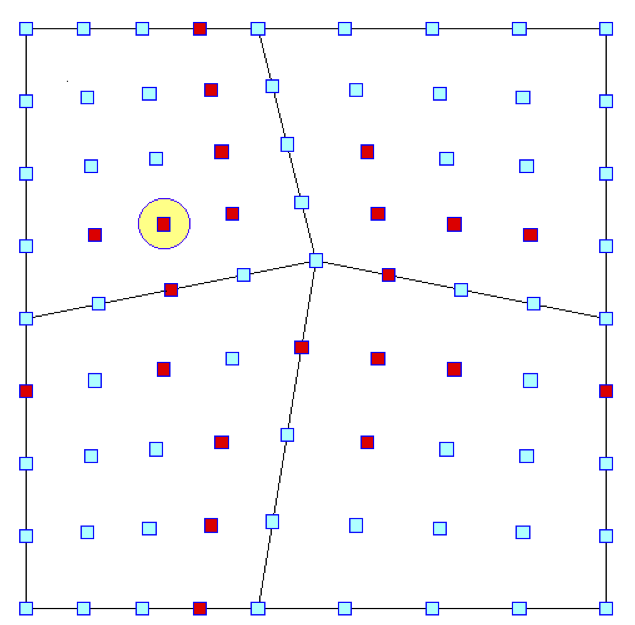}
\label{fig:fig_tp1b}
}
\end{minipage}
\end{narrow}
\vspace{-0.2in}
\setcaptionwidth{5.3in}
\caption{Examples of $\MDS{4}$ for irregular $4$-element mesh.}
\label{fig:fig_tp1}
\end{figure}

Figure ~\ref{fig:fig_tp1} presents an example of {\it "pure"} MDS $\ \MDS{4}$ for an irregular $4$-element mesh over a square domain. 
The basic control points are colored light blue and the dependent control points are colored red. The following control points belong to $\MDS{4}$.
\begin{narrow}{-0.15in}{0.0in}
\noindent\vspace{-0.2in}
\begin{itemize}
\item
\underline{$\FCP{4}$} contains $9$ corner control points for every mesh element. 
\item 
\underline{$\GMDS{4}$} contains the following control points
\begin{itemize}
\item[-]
{\it $V$-type control points.}
$V$-type control points for $5$ non-corner vertices.
\item[-]
{\it $E$-type control points.}
$2$ $E$-type control points for every one of the $5$ non-corner vertices. (For non-corner boundary vertices, the inner and one of the boundary $E$-type control points belong to MDS). 
\item[-]
{\it $D$-type control points.}
$3$ $D$-type control points adjacent to the inner vertex, the only $D$-relevant vertex for the considered mesh configuration. 
\item[-]
{\it $T$-type control points.}
One $T$-type control point for every one of $5$ non-corner vertices.
\item[-]
{\it Middle control points.}
One "side" {\it middle} control point for the left inner edge, since the {\it "Projections Relation"} is satisfied for mesh elements adjacent to the edge.
\end{itemize}
\end{itemize}
\end{narrow}
According to Theorem ~\ref{theorem:mds_dimension_bilinear},
\begin{equation}
\begin{array}{l}
|\GMDS{4}|= 3\times 5+ 4 + 1\times 4 + 0 + 1 = 24\cr
|\MDS{4}| = |\GMDS{4}|+|\FCP{4}| = 24+36 = 60
\end{array}
\end{equation}
\noindent
Correctness of the construction of the MDS and of the formula given in Theorem ~\ref{theorem:mds_dimension_bilinear} was verified by the straightforward construction of such a linear system of constraints that
\begin{itemize}
\item[]
All the control points of all the patches ($100$ control points altogether) are the unknowns.
\item[]
The system of constraints is formed by $G^1$-continuity equations for all inner edges.
\end{itemize}
The formal rank analysis of the system shows that its rank is equal to $40$. It means that there are $60$ free unknowns, which is precisely equal to the number of basic control points, obtained as a result of the construction of "pure" MDS.

The example shown in Figure ~\ref{fig:fig_tp1b} provides another illustration to the correctness of MDS construction. More precisely, it verifies the correctness of the classification rule for the {\it middle} control points (see Subsection ~\ref{subsect:local_mds_middle_bilinear}).
In Figure ~\ref{fig:fig_tp1b} the boundary vertex of the left inner edge is moved downwards with respect to the example shown in Figure ~\ref{fig:fig_tp1a}. Now the {\it "Projections Relation"} is no longer satisfied for the left edge and both {\it middle} control points adjacent to the edge are classified as dependent ones. It means that $|\MDS{4}|=59$, which again is precisely equal to the number of free unknowns resulting from the linear system of constraints for the current mesh configuration.

\eject
\vspace{-1.0in}
\section{Examples of an approximate solution of the Thin Plate problem for irregular quadrilateral meshes}
\label{sect:application_example}
The example given on a circular domain can be compared to the recent  survey ~\cite{nguyen} where several approches , including IGA with multi-patches  are evaluated. One will note the complete different approach in the underlying mesh structure.

\vspace{-0.1in}
\subsection{The Thin Plate problem}
\label{subsect:tp_definition}

Let the middle plane of a plate with uniform thickness lie in $(XY)$-plane. Then deflection of the plate under load $f$ - the solution of the Thin Plate problem - can be found by constrained minimisation of the energy functional
\vspace{-0.1in}
\begin{equation}
\begin{array}{l}
\hspace{-0.4in}
{\cal E}\!=\! 
\int\!\int\limits_{\!\!\!\!\!\!{\tilde\Omega}}\left[   
D\left\{\!
\left(\DT{Z}{X}\!+\!\DT{Z}{Y}\right)^2\!\!-\!
2(1\!-\!\nu)
\left(\DT{Z}{X}\DT{Z}{Y}\!-\!
\left(\DD{Z}{X}{Y}\right)^2\right)\!
\right\}- 
fZ\right] dX dY
\cr
\end{array}
\label{eq::thin_plate_functional_1}
\end{equation} 
\vspace{-0.1in}
where the constraints are defined by simply-supported of clamped boundary conditions.
\noindent
Here 
\vspace{-0.1in}
\begin{equation}
D=E h^3/\left[12(1-\nu^2)\right]
\label{eq:def_D_thin_plate}
\end{equation}
\vspace{-0.1in}
where 
$h$ - thickness of the plate, $\nu$ - Poisson's ratio and $E$ - modulus of elasticity.

References ~\cite{timoshenko_main}, ~\cite{zienkiewicz},
~\cite{ciarlet_old}, ~\cite{bernadou} provide a detailed analysis of the Thin Plate problem and discuss different physical assumptions and limitations for application of Equation ~\ref{eq::thin_plate_functional_1}. 

Next two Subsections contain examples of approximate solutions of the Thin Plate problem, which are based on the technique presented in the current work. In both examples, a highly irregular meshing of a simple planar domain is considered. It allows, on the one hand, to verify the correctness of the current approach for complicated meshes and, on the other hand, to compare the approximate solution with the exact one, which is known for  simple domains. In ~\cite{nguyen} , one will find an evaluation of many higher order methods on a similar problem, noting that none are on a completely unstructured mesh.

\vspace{-0.1in}
\subsection{Approximate solution over a circular domain}

The precise deflection equation for a uniformly loaded circular thin plate with the clamped boundary condition has the following form (see ~\cite{roarks}, ~\cite{timoshenko_main}) 
\vspace{-0.05in}
\begin{equation}
\begin{array}{l}
Z(X,Y)=f(R^2-r^2)^2/(64 D)
\end{array}
\label{eq:solution_circle_exact}
\end{equation}
\vspace{-0.05in}
where $R$-radius of the circle, $\T C=(CX,CY)$ - center of the circle, $r^2=(X-CX)^2+(Y-CY)^2$ - square distance from a point to the center of the circle. 

\begin{figure}[!pt]
\vspace{-0.7in}
\centering
\begin{narrow}{-0.8in}{-0.5in}
\begin{minipage}[]{0.5\linewidth}
\subfigure[]
{
\includegraphics[clip,totalheight=2.9in]{FIGURES_PS/fig_tp2a.png}
\label{fig:fig_tp2a}
}
\end{minipage}
\hspace{0.1in}
\begin{minipage}[]{0.5\linewidth}
\subfigure[]
{
\includegraphics[clip,totalheight=2.9in]{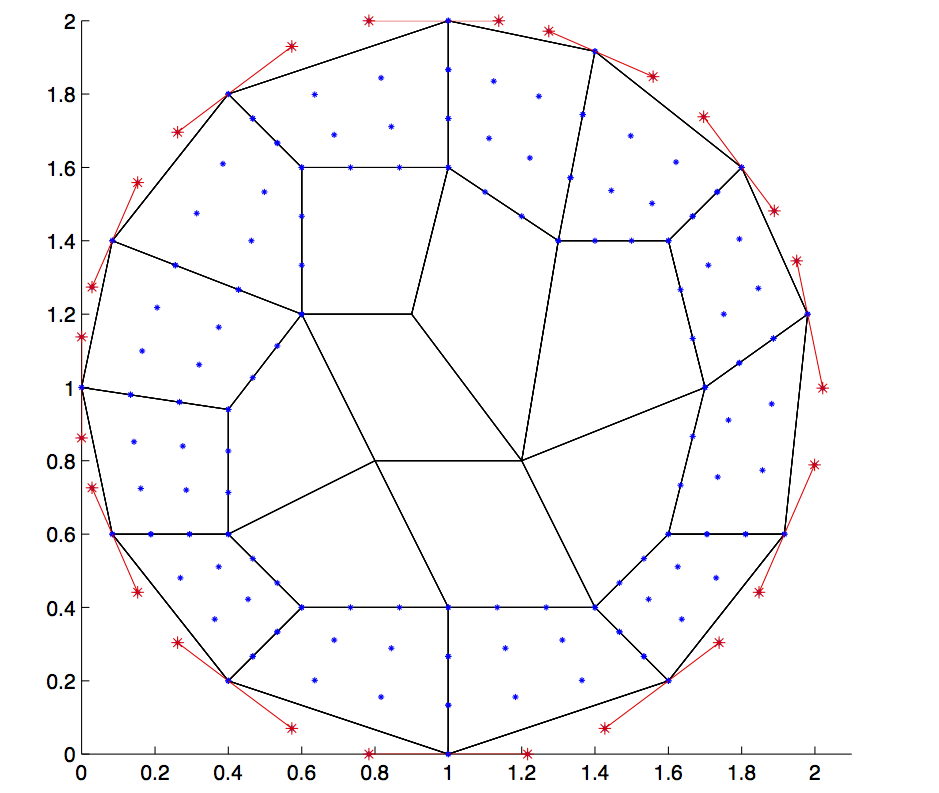}
\label{fig:fig_tp2b}
}
\end{minipage}
\end{narrow}
\vspace{-0.2in}
\caption{An irregular quadrilateral mesh for a circular domain.}
\label{fig:fig_tp2}
\end{figure}

An approximate solution of degree $5$ is constructed for the irregular mesh presented in Figure ~\ref{fig:fig_tp2a}. 
The mesh contains $20$ quadrilateral elements; the global boundary of the domain is approximated by a piecewise-cubic parametric B\'ezier curve (control points of the boundary curves are shown in Figure ~\ref{fig:fig_tp2a}). According to the current approach, global in-plane parametrisation $\T\Pi^{(bicubic)}$ is considered. Control points of bicubic in-plane parametrisation for the boundary mesh elements are shown in Figure ~\ref{fig:fig_tp2b}.

In the example $R=1 in$, $\T C=(1 in,1 in)$ and let $h=0.04 in$, $\nu=0.3$, $E=40\times 10^6 lb/in^2$. From the formal point of view, $f$ is just a coefficient of proportionality in the deflection equation. An unreasonably large value of $f=20\times 10^3 lb/in^2$ is taken in order to visualize the resulting surface, while all other computations and statistics correspond to $f=20 lb/in^2$.

\begin{figure}[!pb]
\vspace{-0.1in}
\centering
\includegraphics[clip,totalheight=3.7in]{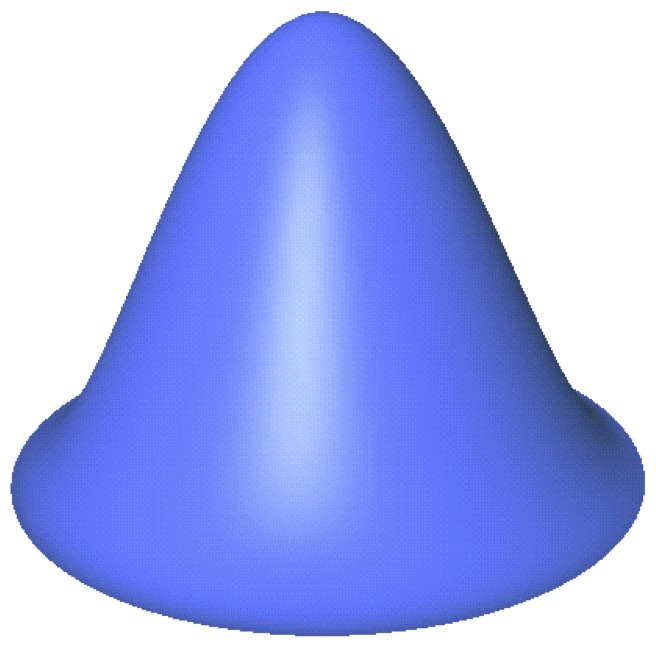}
\caption{The resulting smooth surface (case of the circular domain, irregular mesh).}
\label{fig:fig_tp3}
\end{figure}

The resulting smooth surface is shown in Figure ~\ref{fig:fig_tp3}. The total number of the basic free control points is equal to $208$.
Level lines for the resulting surface and for its first-order derivatives (see Figure ~\ref{fig:fig_tp4},  Appendix, Section ~\ref{sect:figures_tp}) confirm $C^1$-smoothness of the surface.
A comparison between the approximate and the exact solutions along segment $X=1 in$ is shown in Figures ~\ref{fig:fig_tp5} and ~\ref{fig:fig_tp6}(see Appendix, Section ~\ref{sect:figures_tp}).
One can see, that the approximate solution fits the exact one with a high precision. Errors in the center of the circle for the solution itself and for the bending moments do not exceed one hundredth percent. 

\vspace{-0.1in}
\subsection{Approximate solution over a square domain}
Analysis of an approximate solution of degree $4$ for a square domain with the simply-supported boundary condition is made according to the same scenario. The exact solution for the uniform load is (see ~\cite{roarks}, ~\cite{timoshenko_main})
\vspace{-0.05in}
\begin{equation}
Z(X,Y)=\frac{16 f a^4}{\pi^6 D}
\sum_{\lmtT{m,n=1}{m,n-odd}{}}^\infty
\frac{sin\frac{m\pi X}{a}sin\frac{n\pi Y}{b}}{mn(m^2+n^2)^2}
\label{eq:solution_quadratic_exact}
\end{equation}
\vspace{-0.15in}
where the square domain is axes-aligned, the origin of the coordinate system for the $XY$-plane coincides with the lower-left corner of the square and $a$ is the length of the side.

\begin{figure}[!pt]
\vspace{-0.7in}
\centering
\includegraphics[clip,totalheight=2.9in]{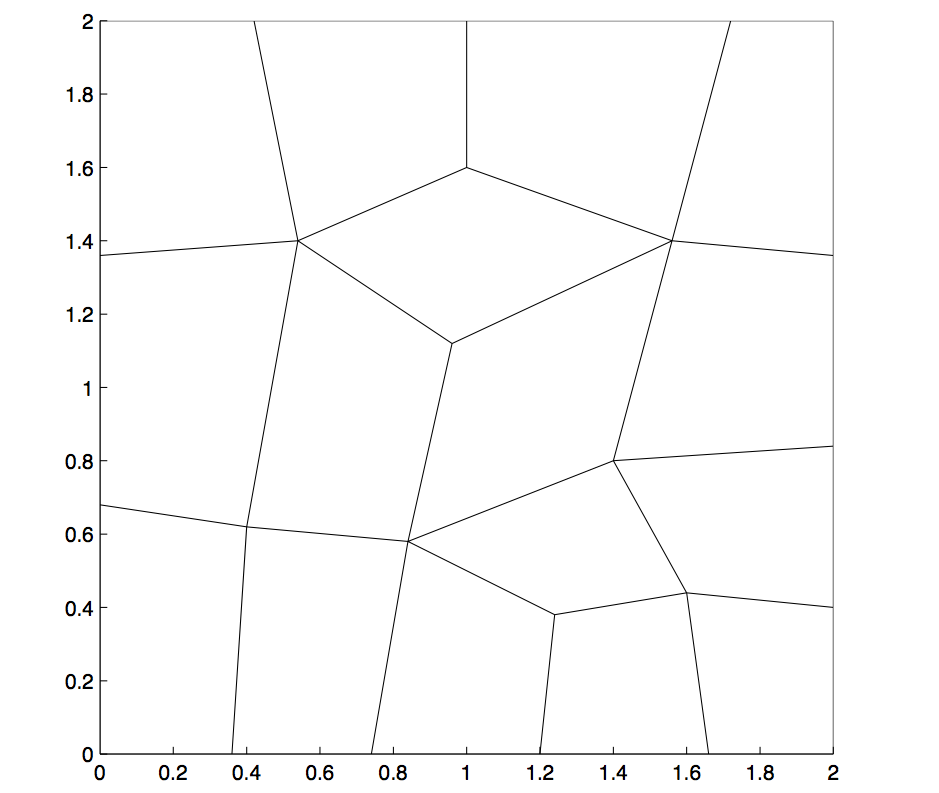}
\caption{Irregular quadrilateral mesh for a square domain.}
\label{fig:fig_tp7}
\end{figure}

Let $a=2 in$, $h=0.04 in$, $\nu=0.3$, $E=40\times 10^6 lb/in^2$;\ \ $f=5\times 10^3 lb/in^2$
for visualization of the surface and $f=5 lb/in^2$ for comparison with the exact solution.
An approximate solution is constructed for the irregular mesh shown in Figure ~\ref{fig:fig_tp7}. 
In the current case global bilinear in-plane parametrisation $\T\Pi^{(bilinear)}$ is considered and it is sufficient to consider the space $\FUN{4}(\T\Pi^{(bilinear)})$. 

\begin{figure}[!pb]
\vspace{-0.2in}
\centering
\includegraphics[clip,totalheight=3.7in]{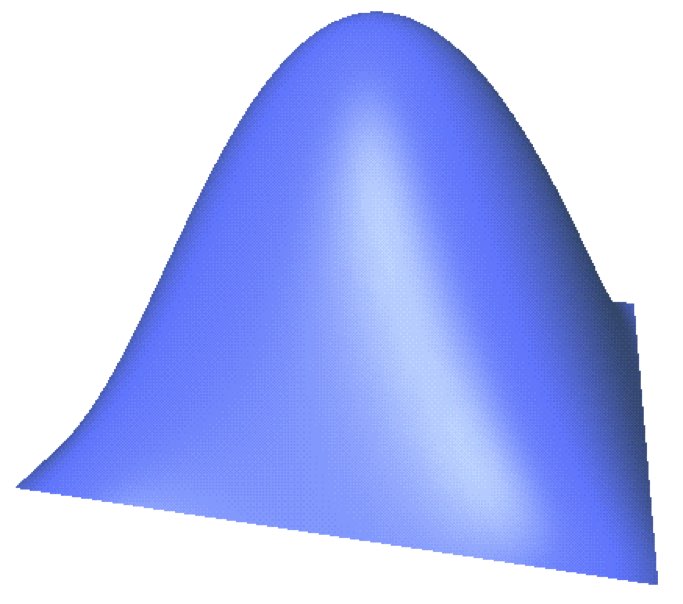}
\caption{The resulting smooth surface (case of the square domain, irregular mesh).}
\label{fig:fig_tp8}
\end{figure}

The resulting smooth surface is shown in Figure ~\ref{fig:fig_tp8}. The total number of the basic free control points is equal to $111$. 
Figure ~\ref{fig:fig_tp9} (see Appendix, Section ~\ref{sect:figures_tp}) presents the level lines for the surface and for its first-order derivatives.
Comparison between the approximate and the exact solution along the segment $X=1 in$ 
is given in Figures ~\ref{fig:fig_tp10} and ~\ref{fig:fig_tp11} (see Appendix, Section ~\ref{sect:figures_tp}). 
It is important to note, that Figure ~\ref{fig:fig_tp11} definitely shows that the approximate solution is not $C^2$-smooth. 

Errors of the approximate solution for the irregular mesh were compared with the errors for the regular $4\times 4$ square mesh. For the regular mesh, the total number of the basic free control points is equal to $144$. (Note, that the number of the basic free control points obtained as the computational result is equal to the number of the basic control points predicted theoretically. Indeed, for the regular mesh, there are $196$ basic control points according to Theorem ~\ref{theorem:mds_dimension_bicubic} and $52$ control points are fixed by the boundary condition). 

Both the irregular and the regular meshes contain $16$ elements and the following table shows that errors along segment $X=1 in$ have the same order for the meshes.  

\vspace{0.05in}
\hspace{-0.23in}
\begin{tabular}{|c||c|c||c|c|}
\hline
Characteristic & 
\multicolumn{2}{c||}
{\begin{tabular}{c}
Maximal absolute value\cr of the error (in inches)
\end{tabular}} &
\multicolumn{2}{c|}
{\begin{tabular}{c}
Error in the center\cr (in percents)
\end{tabular}} \cr
\cline{2-5}
& \begin{tabular}{c}Regular\ \cr\ mesh   \end{tabular} 
& \begin{tabular}{c}Irregular\cr mesh \end{tabular}
& \begin{tabular}{c}\ Regular\ \cr\ mesh   \end{tabular} 
& \begin{tabular}{c}Irregular\cr mesh \end{tabular}
\cr
\hline
$Z$        & $6.72\times 10^{-8}$  & $1.69\times 10^{-7}$  & 
             $8.92\times 10^{-4}\%$ & $2.70\times 10^{-5}\%$ \cr
\hline
$\D{Z}{X}$ & $3.00\times 10^{-15}$ & $1.23\times 10^{-6}$  & 
             $-$                   & $-$                   \cr
\hline
$\D{Z}{Y}$ & $8.91\times 10^{-7}$  & $1.73\times 10^{-6}$  & 
             $-$                   & $-$                   \cr
\hline
$M_{X}$    & $3.29\times 10^{-3}$  & $8.85\times 10^{-3}$  & 
             $0.26\%$               & $0.19\%$               \cr
\hline
$M_{Y}$    & $6.88\times 10^{-3}$  & $5.90\times 10^{-3}$  & 
             $0.26\%$               & $0.10\%$               \cr
\hline
$M_{XY}$   & $3.03\times 10^{-12}$ & $1.62\times 10^{-3}$  & 
             $-$                   & $-$                   \cr
\hline
\end{tabular}

\eject

\FloatBarrier

\FloatBarrier

\eject
\part{Conclusions and  further research.}

\label{part:further_research}
A complete solution for cubic $C^{1}$  boundaries has been constructed. Rigorous poofs show that for $n \ge 5$ there is always a solution. One can also show that the case $n=4$ has a solution provided we respect some conditions on the mesh. All the proofs are constructive and have been implemented for the examples we gave. We also have shown that one can mixes quartic and quintic patches. As a conclusion from a practical point of view, the present work provides a way to solve interpolation/approximation and partial differential problems for arbitrary structures of quadrilateral meshes.
The solution has a linear form, which makes it relatively simple, fast and stable. Construction of the MDS allows operating only with essential degrees of freedom. Application of the approach to the solution of fourth order  partial differential equations results in an approximate solution of a high quality. \\

The following items are s a natural extension of the research, some trivial, others needing further studies:

\begin{itemize}
\item[]

{\it Boundary conditions} Boundary conditions ( BCs) on B\'ezier or B-Splines require special care, since the corresponding basis are non interpolator , cf . Embar et. all
 ~\cite{harari1}, and we have dealt only with the simple homogenous case in our examples. The same remark apply to Robin type of BCs. 
Another topic,periodic BCs are not analysed here, but are essential in treating domains with holes.

\item[]
{\it Extension to non Rationnal B\'ezier patches.} This seems simple for the case of plane bilinear meshes , since the "plane" part is trivially bilinear again, but the case of curved boundaries must be analysed anew.
\item[]
{\it Stability concerns.} Although the geometrical characteristics of a mesh play a very important role in the construction of MDS, they usually define the general type of classification, while the choice of the basic control points in possible cases of ambiguity is made arbitrarily. For example, for $\MDS{4}(\T\Pi^{(bilinear)}$, geometrical {\it "Projections Relation"} defines the number of the basic {\it middle} control points (see Section ~\ref{subsect:example_4_5_middle_bilinear}).  However, the geometry of the elements has no influence on the choice of the basic control point when one of two {\it middle} control points is classified as basic and another one as dependent. Taking into account mesh geometry in cases when an ambiguity is possible and introduction of the tolerance analysis into equations themselves can significantly improve the stability of the solution.

\item[]
{\it Error estimation.}. The results of Part ~\ref{part:illustrative_examples} demonstrate a high quality of the solution.  A rigorous mathematical analysis of the error in the case of an approximate solution of a partial differential equation can be made using the error estimates  methods of IGA, ~\cite{buffa} 

\item[]
{\it Combination of B\'ezier patches with different polynomial degrees.} 
Part ~\ref{part:mds_mixed_4_5} presents the construction of mixed MDS, which combines patches of degree $5$ along the global boundary and patches of degree $4$ in the inner part of the mesh. The approach may be improved so that fewer patches of degree $5$ may be used. It should be possible to consider patches of degree $5$ locally, only in problematic areas of a mesh.

\item[]
{\it Generalisation of MDS.}
Higher order smoothness seems difficult to build on unstructured quadrilateral meshes. But lower orders ( $ n \le 4$)   should be possible provide one uses macro elements  
A study of MDS for lower polynomial degrees of patches and/or higher orders of smoothness is a very important possible domain of the further research. 

\item[]
{\it A study of non-planar initial meshes.} Consideration of non-planar initial meshes gives another direction of a possible research. For example, it may be interesting to start with construction of the weight functions and linearisation of an approximate solution for the Thin Plate problem when the middle surface lies on a cylinder or sphere (see works ~\cite{nielson_sphere}, ~\cite{schumaker_sphere}).
\end{itemize}

Is it possible to eliminate some limitations on the structure of the planar mesh (see Sections ~\ref{sect:mesh_limitations_bilinear} and ~\ref{sect:mesh_limitations_bicubic}) or refine sufficient conditions for the regularity of bicubic in-plane parametrisation (see  Subsection ~\ref{subsect:parametrisation_regularity_bicubic})?

\eject
\addcontentsline{toc}{part}{Bibliography}
\bibliography{bib}
\bibliographystyle{plain}

\eject

\addtocontents{lof}{\vspace{0.45in}}
\addtocontents{lof}{{\bf\large Appendix }}

\setcounter{page}{1}
\addcontentsline{toc}{part}{Appendix }
\appendix{}


\section{Illustrations to approximate solution of the Thin Plate Problem} 
\label{sect:figures_tp}

\vspace*{0.3in}
\begin{figure}[!ph]
\begin{narrow}{-1.0in}{0.0in}
\includegraphics[clip,width=7in]{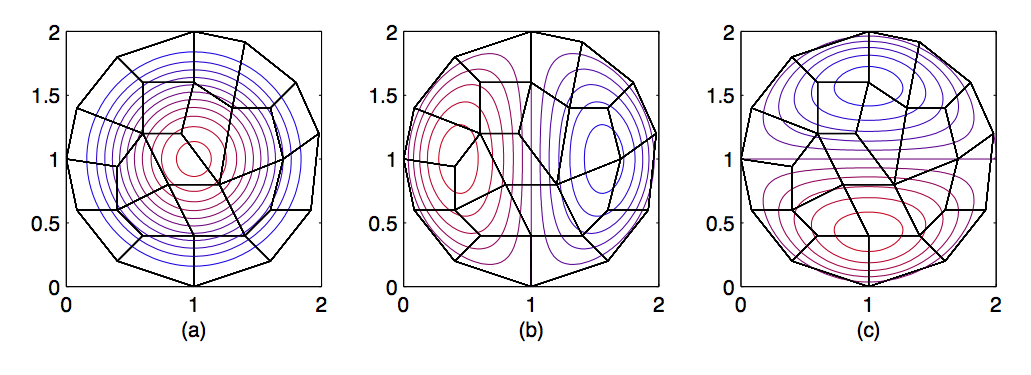}
\end{narrow}
\vspace{-0.3in}
\setcaptionwidth{5.4in}
\caption{
Level lines for the resulting surface and for its first-order derivatives 
(case of the circular domain, irregular mesh)
\ \ \ (a) $Z$ \ \ \ (b) $\D{Z}{X}$ \ \ \ (c) $\D{Z}{Y}$.
}
\label{fig:fig_tp4}
\end{figure}

\begin{figure}[!ph]
\begin{narrow}{-1.0in}{0.0in}
\includegraphics[clip,width=7in]{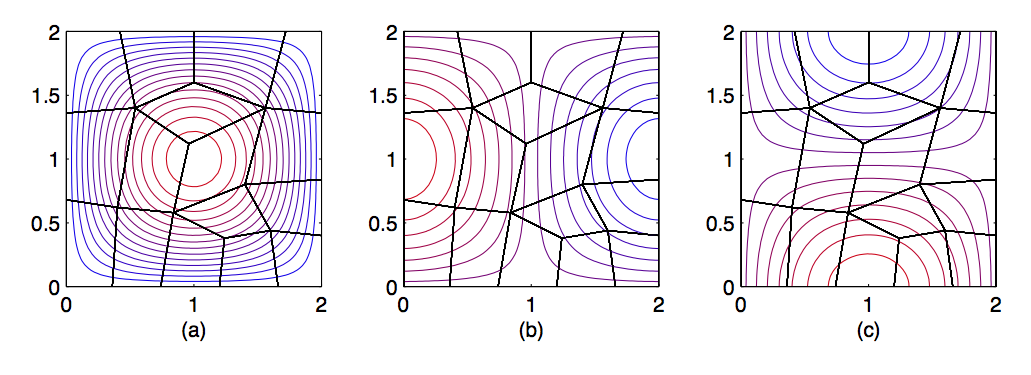}
\end{narrow}
\vspace{-0.3in}
\setcaptionwidth{5.4in}
\caption{
Level lines for the resulting surface and for its first-order derivatives 
(case of the square domain, irregular mesh)
\ \ \ (a) $Z$ \ \ \ (b) $\D{Z}{X}$ \ \ \ (c) $\D{Z}{Y}$.
}
\label{fig:fig_tp9}
\end{figure}

\FloatBarrier

\begin{figure}[!pt]
\vspace{-1.5in}
\hspace{-0.3in}
\includegraphics[clip,totalheight=3.6in,width=5.3in]{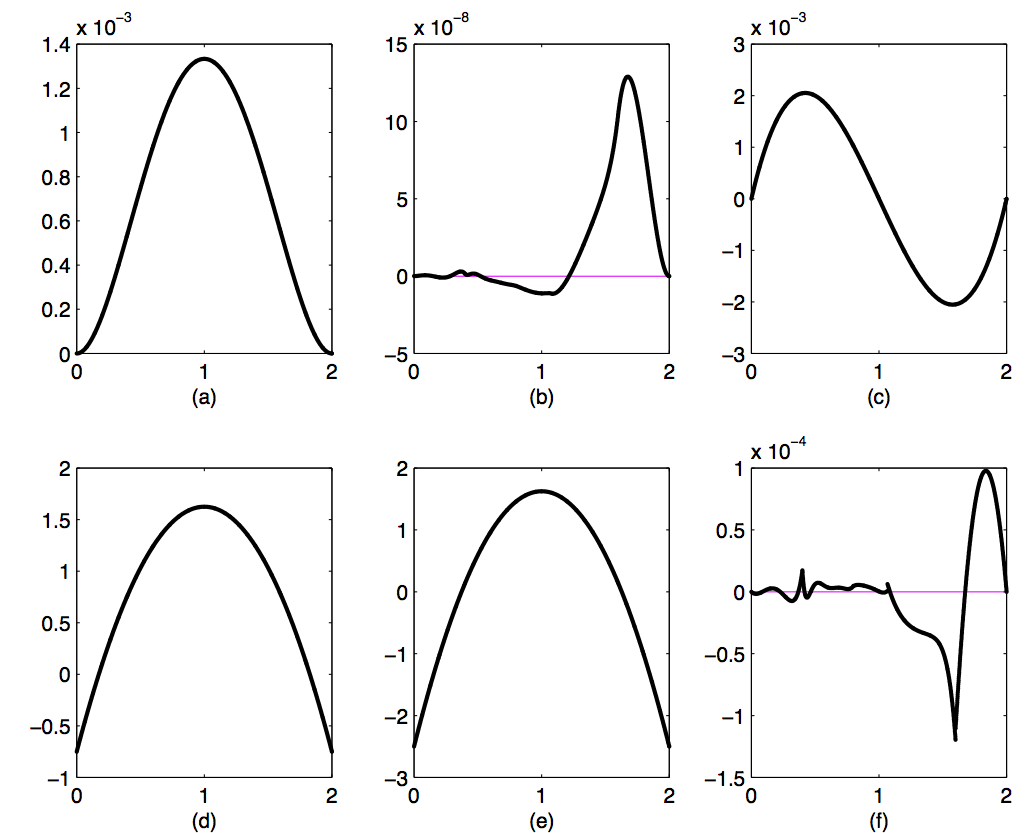}
\vspace{-0.2in}
\setcaptionwidth{5.4in}
\caption{
Values for approximate (bold line) and exact (thin line) solutions, 
their first order derivatives and the bending moments along segment $X=1$ 
(case of the circular domain, irregular mesh)
\ \ \ (a) $Z$ \ \ \ (b) $\D{Z}{X}$ \ \ \ (c) $\D{Z}{Y}$ 
\ \ \ (d) $M_{X}$ \ \ \ (e) $M_{Y}$ \ \ \ (f) $M_{XY}$.
}
\label{fig:fig_tp5}
\end{figure}

\begin{figure}[!pb]
\vspace{-0.5in}
\hspace{-0.3in}
\includegraphics[clip,totalheight=3.6in,width=5.3in]{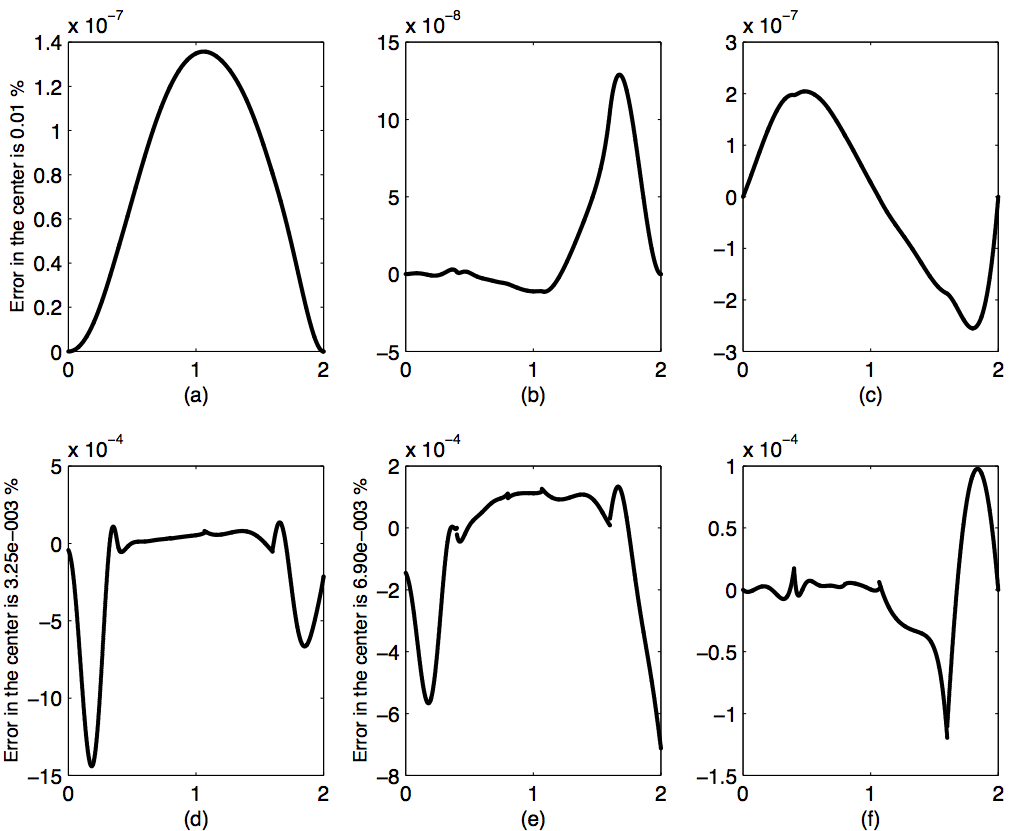}
\vspace{-0.2in}
\setcaptionwidth{5.4in}
\caption{
Difference between the approximate and the exact solutions 
(case of the circular domain, irregular mesh)
\ \ \ (a) $Z$ \ \ \ (b) $\D{Z}{X}$ \ \ \ (c) $\D{Z}{Y}$ 
\ \ \ (d) $M_{X}$ \ \ \ (e) $M_{Y}$ \ \ \ (f) $M_{XY}$.
}
\label{fig:fig_tp6}
\end{figure}

\FloatBarrier

\begin{figure}[!pt]
\vspace{-1.5in}
\hspace{-0.3in}
\includegraphics[clip,totalheight=3.6in,width=5.3in]{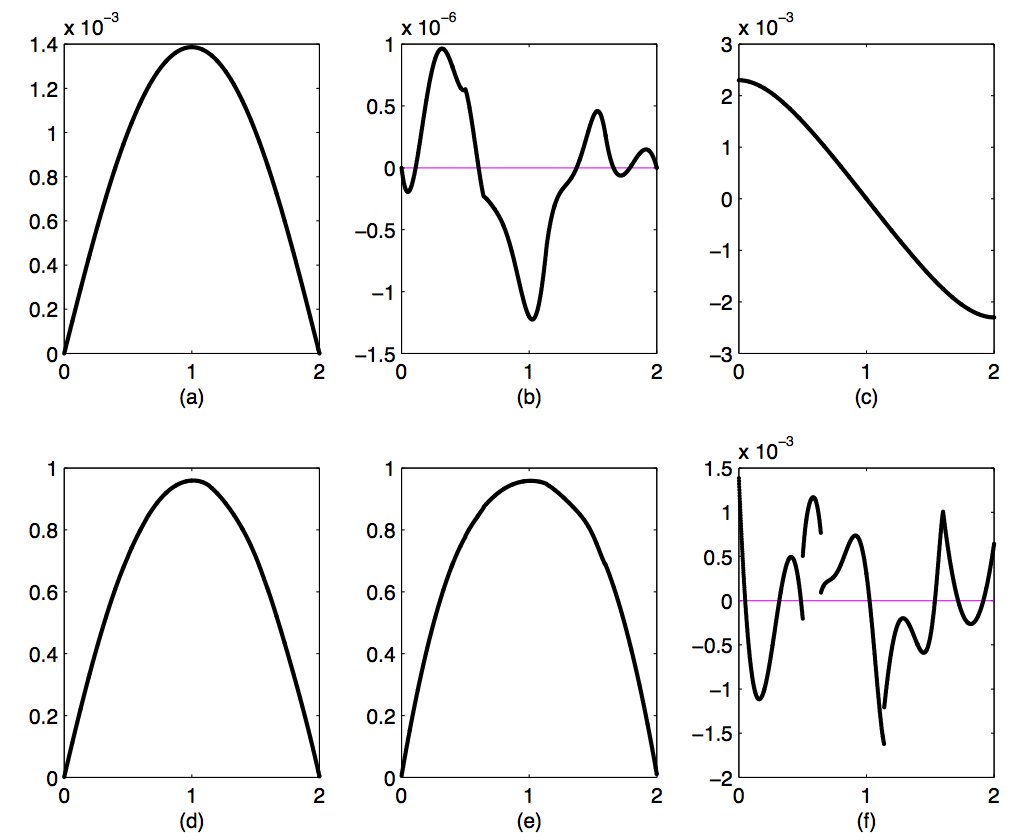}
\vspace{-0.2in}
\setcaptionwidth{5.4in}
\caption{
Values for approximate (bold line) and exact (thin line) solutions, 
their first order derivatives and the bending moments along segment $X=1$ 
(case of the square domain, irregular mesh)
\ \ \ (a) $Z$ \ \ \ (b) $\D{Z}{X}$ \ \ \ (c) $\D{Z}{Y}$ 
\ \ \ (d) $M_{X}$ \ \ \ (e) $M_{Y}$ \ \ \ (f) $M_{XY}$.
}
\label{fig:fig_tp10}
\end{figure}

\begin{figure}[!pb]
\vspace{-0.5in}
\hspace{-0.3in}
\includegraphics[clip,totalheight=3.6in,width=5.3in]{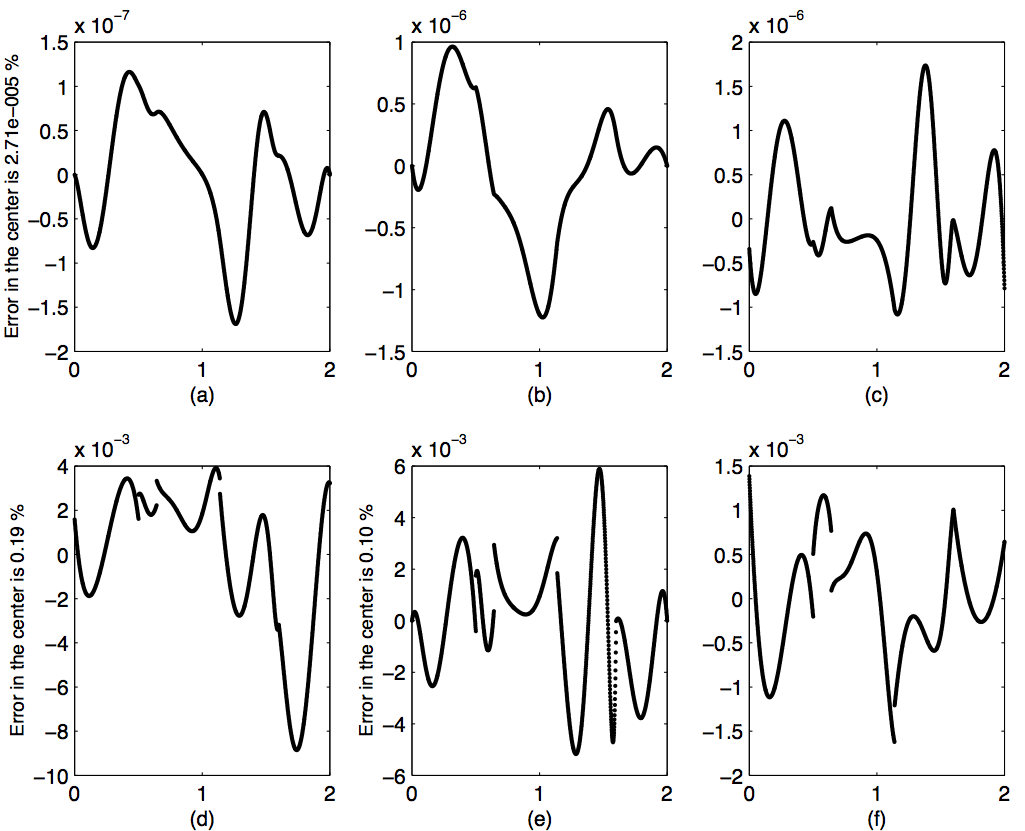}
\vspace{-0.2in}
\setcaptionwidth{5.4in}
\caption{
Difference between approximate and exact solutions 
(case of the square domain, irregular mesh)
\ \ \ (a) $Z$ \ \ \ (b) $\D{Z}{X}$ \ \ \ (c) $\D{Z}{Y}$ 
\ \ \ (d) $M_{X}$ \ \ \ (e) $M_{Y}$ \ \ \ (f) $M_{XY}$.
}
\label{fig:fig_tp11}
\end{figure}

\FloatBarrier

\eject
\section{Algorithms}
\label{sect:algorithms}

\noindent
{\bf Algorithm ~\ref{algorithm:dependency_forest_round}}
In order to clarify the Algorithm, the coefficient of a $D$-type control point $D$ in the {\it "Circular Constraint"} corresponding to vertex $\T V$ will be denoted by $coeff^{D}_{\T V}$.

\vspace{0.15in}

\centerline{\underline{\it Algorithm for construction of the $D$-dependency tree and classification}}
\centerline{\underline{\it of $D$-type and $T$-type control points for a connected component}} 
\centerline{\underline{\it of the $D$-dependency graph which has no dangling half-edge}}

Let connected component $\T{\cal C}$ of the $D$-dependency graph have no dangling half-edge. It means that $\T{\cal C}$ consists of $D$-relevant {\it primary vertices} and {\it full edges} only; let it contain  $|V|$ {\it primary vertices} and $|E|\ge|V|-1$ edges of the initial mesh. 

\vspace{0.1in}
\noindent
$\bullet$
\underline{\bf If ${\bf |E|=|V|-1}$} then $\T{\cal C}$ coincides with its spanning tree and there is no chance of {\it assigning} $D$-type control point to every vertex of the component (the Algorithm {\it fails}).

\vspace{0.1in}
\noindent
$\bullet$ \underline{\bf Otherwise} 

\begin{figure}[!ph]
\centering
\includegraphics[clip,totalheight=2.5in]{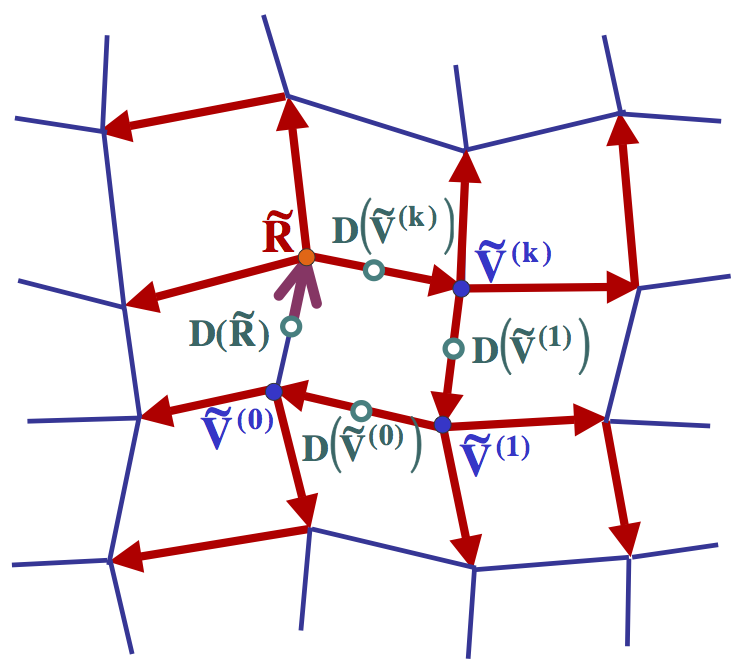}
\caption{An illustration for the construction of a $D$-dependency tree when a component of $D$-dependency graph has no dangling half-edge.}
\label{fig:fig26}
\end{figure}

\begin{narrow}{-0.1in}{0.0in}
\noindent\vspace{-0.2in}
\begin{itemize}
\item[{\bf (1)}] 
Build a \underline{spanning tree} of $\T{\cal C}$. There is at least one edge $\T e$ of the connected component which does not participate in the spanning tree. (Although the choice of the spanning tree and the choice of the edge may affect the stability and even the existence of the solution, the current presentation is restricted to the description of a simple algorithm. A study of more complicated cases is left for a possible further research). 

\item[{\bf (2)}] 
Choose one vertex of $\T e$ to be the \underline{root} of the $D$-dependency tree (vertex $\T R$), denote the second vertex by $\T V^{(0)}$ (see Figure ~\ref{fig:fig26}). {\it Assign} the $D$-type control point of edge $\T e$ to the root vertex $\T R$  (denote the control point by $\T D(\T R)$).

\item[{\bf (3)}] 
Build the \underline{$D$-dependency tree} of $\T{\cal C}$:
\underline{direct} every {\it mesh edge} of the tree from the upper vertex to a lower one (according to the hierarchy of the spanning tree). 

\item[{\bf (4)}] 
\underline{Classify} $D$-type and $T$-type \underline{control points} according to the same principles as in Algorithm ~\ref{algorithm:D_T_classification_degree_4}. The only difference is that now $\T D(\T R)$ can not be considered as a basic control point. Let 
$\T V^{(0)}\rightarrow \T V^{(1)}\rightarrow\ldots\T V^{(k)}\rightarrow\T R$ be the path in the $D$-dependency tree from vertex $\T V^{(0)}$ to the root.
All vertices $\T V^{(0)},\ldots,\T V^{(k)}$ are inner even vertices of the initial mesh, because any boundary vertex uses at most one $D$-type control point and so it has at most one adjacent full edge in the $D$-dependency tree. The following local modifications of Algorithm ~\ref{algorithm:D_T_classification_degree_4} should be made. 
\begin{itemize}
\item[{\bf (a)}]
According to Algorithm ~\ref{algorithm:D_T_classification_degree_4}, one relates to 
$\T D(\T R)$ for the first time when the level of vertex $\T V^{(0)}$ in the $D$-dependency tree (which is traversed down-up) is reached. At this step one should "close" the {\it "Circular Constraint"} of vertex $\T V^{(0)}$ using $D(\T V^{(0)})$.  Both $D(\T R)$ and $D(\T V^{(0)})$ participate in the{\it "Circular Constraint"} with a non-zero coefficients, because both corresponding full edges belong to the $D$-dependency graph. Therefore the {\it "Circular Constraint"} defines dependency of $D(\T V^{(0)})$ on $D(\T R)$ in the following way
\vspace{-0.07in}
\begin{equation}
\begin{array}{l}
D(\T V^{(0)}) = 
-\frac{coeff^{D(\T R)}_{\T V^{(0)}}} {coeff^{D(\T V^{(0)})}_{\T V^{(0)}}}
D(\T R)+f^{ D(\T V^{(0)})}\cr
\end{array}
\end{equation}
\vspace{-0.07in}
where $f^{ D(\T V^{(0)})}$ stands for some linear combination of the basic control points.
\item[{\bf (b)}]
Proceeding to traverse the $D$-dependency tree down-up, one gradually defines the linear dependencies of $D(\T V^{(1)}),\ldots,
D(\T V^{(k)})$ on $D(\T R)$ and some basic control points. At the level of the root vertex one gets
\begin{equation}
D(\T R) = coeff\  D(\T R)+ f^{ D(\T R)}
\label{eq:root_cp_dependency}
\end{equation}
where
\vspace{-0.07in}
\begin{equation}
\begin{array}{l}
coeff =
(-1)^k
\frac{coeff^{D(\T R)}_{\T V^{(0)}}}
     {coeff^{D(\T V^{(0)})}_{\T V^{(0)}}}
\frac{coeff^{D(\T V^{(0)})}_{\T V^{(1)}}}
     {coeff^{D(\T V^{(1)})}_{\T V^{(1)}}}
\ldots
\frac{coeff^{D(\T V^{(k-1)})}_{\T V^{(k)}}}
     {coeff^{D(\T V^{(k)})}_{\T V^{(k)}}}\cr
\end{array}
\end{equation}
\item[{\bf (c)}]
If $coeff\neq 1$ (which holds in a general case) then $\T D(\T R)$ is classified as a dependent control point; its dependency on the basic control points is defined according to Equation ~\ref{eq:root_cp_dependency}. It remains to substitute the expression for $D(\T R)$ into dependency equations for $D(\T V^{(0)}),\ldots, D(\T V^{(k)})$. In this case, the Algorithm {\it succeeds}.

Otherwise it is possible either to try to choose another "free" edge of the spanning tree or even another spanning tree and to run the Algorithm from the beginning or to decide that the Algorithm {\it fails} and to use one of the possibilities listed in the end of Subsection ~\ref{subsect:global_classification_DT_bilinear}.
\end{itemize}
\end{itemize}
\end{narrow}

\vspace{0.4in}
\noindent
{\bf Algorithm ~\ref{algorithm:mixed_mds}}

\centerline{\underline{\it Algorithm for construction of mixed MDS\ $\MDS{4,5}$}} 

\vspace{0.2in}
\noindent
Let a planar mesh have either a polygonal or a smooth global boundary and global in-plane parametrisation $\T\Pi=\T\Pi^{(bilinear)}$ or $\T\Pi=\T\Pi^{(bicubic)}$ is considered. In description of the Algorithm superscripts ${}^{(4)}$ and ${}^{(5)}$ are added in order to specify whether some equation or control point relates to degree $4$ or $5$ of the resulting patch.

\paragraph{"Stage 1"}
Assume (artificially) that there are no $D$-relevant boundary vertices. Classify all $V^{(4)}$, $E^{(4)}$, $D^{(4)}$, $T^{(4)}$-type control points adjacent to the {\it inner vertices} precisely as in the case of global bilinear in-plane parametrisation $\T\Pi^{(bilinear)}$ (see Subsections ~\ref{subsect:local_VE_inner_vertex_bilinear} and Algorithm ~\ref{algorithm:D_T_classification_degree_4}). The classification exists according to Theorem ~\ref{theorem:sufficient_DT_classification_bilinear}.

At the end of this stage all $V^{(4)}$,$E^{(4)}$,$T^{(4)}$-type control points for all {\it inner vertices} and {\it all $D^{(4)}$-type} control points are classified and {\it "Eq(0)"}-type and {\it "Eq(1)"}-type equations at all inner vertices are satisfied.

For an edge with one boundary vertex, control points $(\T L^{(4)}_0,\T C^{(4)}_0,\T R^{(4)}_0)$,\\ $(\T L^{(4)}_1,\T C^{(4)}_1,\T R^{(4)}_1)$ and $\T C^{(4)}_2$ are classified and equations {\it $"Eq^{(4)}(0)"$} and {\it $"Eq^{(4)}(1)"$} are satisfied (see Figure ~\ref{fig:fig36a}).

\paragraph{"Stage 2"}
For every inner edge with two inner vertices, apply a local template for classification of the {\it $middle^{(4)}$} control points, exactly as in the case of global bilinear in-plane parametrisation $\T\Pi^{(bilinear)}$ (see Subsection ~\ref{subsect:local_mds_middle_bilinear}).

At the end of this stage all points from $\GCP{4}$ lying on or adjacent to any {\it edge with two inner vertices} and {\it all $D^{(4)}$-type} control points are classified (see Figure ~\ref{fig:fig36a}). $G^1$-continuity conditions are satisfied for any edge with two inner vertices. For every boundary mesh element, all "degree 4" control points are classified.

\paragraph{"Stage 3"}
For every boundary mesh element, consider $3D$ patch of degree $5$. Use the {\it degree elevation} formulas in order to classify control points $\T P^{(5)}_{i,j}$, $i=0,\ldots,5$, $j=0,1$ and 
$\T P^{(5)}_{0,2}$, $\T P^{(5)}_{5,2}$ and in order re-compute their dependencies on "degree 4" basic control points from $\CP{4,5}$ (see Figure ~\ref{fig:fig36}).

At the end of this stage, {\it all} $G^1$-continuity conditions hold for all edges with {\it two inner vertices}. For an edge with one boundary vertex, control points $(\T L^{(5)}_i,\T C^{(5)}_i,\T R^{(5)}_i)$, $i=0,1$ and $\T C^{(5)}_2$ are classified and equations {\it $"Eq^{(5)}(0)"$} and {\it $"Eq^{(5)}(1)"$} are satisfied (see Technical Lemma ~\ref{tl:degree_combined}, Appendix, Section ~\ref{sect:technical_lemmas})

\paragraph{"Stage 4"}
For an edge with one boundary vertex, locally classify control points $(\T L^{(5)}_i, \T R^{(5)}_i)$ ($i=2,3,4,5$) and $\T C^{(5)}_i$ ($i=3,4,5$) in such a manner that the remaining $G^1$-continuity equations for the edge will be satisfied.

\vspace{0.1in}
\noindent
In the case of a mesh with a polygonal global boundary and global in-plane parametrisation $\T\Pi^{(bilinear)}$
\begin{itemize}
\item[-]
$V$,$E$-type control points adjacent to the boundary vertex ($\T L^{(5)}_5$, $\T C^{(5)}_5$, $\T R^{(5)}_5$ and $C^{(5)}_4$) are classified by a local template described in Subsection ~\ref{subsect:local_VE_boundary_vertex_bilinear}.
\item[-]
$D$,$T$-type control points adjacent to the boundary vertex ($L^{(5)}_4$, $R^{(5)}_4$ and $C^{(5)}_3$) are classified by a local template described in Subsection ~\ref{subsect:local_mds_DT_boundary_bilinear}.
\item[-]
{\it Middle} control points ($\T L^{(5)}_2$, $\T R^{(5)}_2$, $\T L^{(5)}_3$, $\T R^{(5)}_3$) are classified according to a local template described in Subsection ~\ref{subsect:local_mds_middle_bilinear}.
\end{itemize}

\vspace{0.1in}
\noindent
In the case of a mesh with a smooth global boundary and global in-plane parametrisation $\T\Pi^{(bicubic)}$
\begin{itemize}
\item[-]
Control points $\T L^{(5)}_5$, $\T C^{(5)}_5$, $\T R^{(5)}_5$, $i=4,5$ are classified by local template $TB1^{(bicubic)}$, described in Subsection ~\ref{subsect:local_mds_boundary_vertex_bicubic}.
\item[-]
{\it Middle} control points $\T L^{(5)}_3$, $\T C^{(5)}_3$, $\T R^{(5)}_3$, $\T L^{(5)}_2$, $\T R^{(5)}_2$ are classified according to a local template described in Subsection ~\ref{subsect:local_mds_middle_edge_one_boundary_bicubic}.
\end{itemize}

\FloatBarrier

\eject
\section{Technical Lemmas}
\label{sect:technical_lemmas}

\vspace{0.1in}
\begin{technical_lemma}
Let $\T\Pi^{(bilinear)}\in\PAR{1}$ be global bilinear parametrisation and let $\LL$,$\LL'$,$\GG$,$\GG'$,$\RR$,$\RR'$ be vertices of two adjacent mesh elements (see Figure ~\ref{fig:fig9}). Then control points, which belong to  $\GCP{n}(\T\Pi^{(bilinear)})$ and participate in $G^1$-continuity conditions for the common edge of two patches, and first-order differences of the control points have the following representations
\begin{equation}
\begin{array}{l}
\T{C}_j=\GG\left(1-\frac{j}{n}\right)+\frac{j}{n}\GG',\cr
\T{R}_j=\frac{n-1}{n}\left(\GG\frac{n-j}{n}+\GG'\frac{j}{n}\right)+
        \frac{1}{n}\left(\RR\frac{n-j}{n}+\RR'\frac{j}{n}\right),\cr
\T{L}_j=\frac{n-1}{n}\left(\GG\frac{n-j}{n}+\GG'\frac{j}{n}\right)+
        \frac{1}{n}\left(\LL\frac{n-j}{n}+\LL'\frac{j}{n}\right),\cr
\Delta\T{C}_{j}=\frac{1}{n}(\GG'-\GG),\cr
\Delta\T{L}_{j}=-\frac{1}{n}\left((\LL-\GG)\frac{n-j}{n}+(\LL'-\GG')\frac{j}{n}\right),\cr
\Delta\T{R}_{j}=-\frac{1}{n}\left((\RR-\GG)\frac{n-j}{n}+(\RR'-\GG')\frac{j}{n}\right)
\end{array}
\label{eq:degree_elevation_bilinear}
\end{equation}
Here $j=0,\ldots,n$ in expressions for $\T C_j$, $\T R_j$, $\T L_j$, $\Dlt L_j$, $\Dlt R_j$ and $j=0,\ldots,n-1$ for $\Dlt C_j$.
\label{tl:degree_elevation_bilinear}
\end{technical_lemma}

\vspace{0.1in}
\begin{technical_lemma}
The control points of bicubic in-plane parametrisation $\T P(u,v) = \sum_{i,j=0}^3 \T P_{ij} B^3_{ij}(u,v)$ for a boundary mesh element in the case of a planar mesh with a piecewise-cubic $G^1$-smooth global boundary have the following explicit formulas in terms of the planar mesh data (see Figure ~\ref{fig:fig29})
\begin{equation}
\begin{array}{ll}
\T P_{00} = \T A &
\T P_{02} = \frac{1}{3}(\T A + 2\T D)\cr
\T P_{10} = \frac{1}{3}(2\T A + \T B) &
\T P_{12} = \frac{1}{9}(2\T A + \T B + \T C + 2\T D + 3\T E)\cr 
\T P_{20} = \frac{1}{3}(\T A + 2\T B) &
\T P_{22} = \frac{1}{9}(\T A + 2\T B + 2\T C + \T D + 3\T F)\cr
\T P_{30} = \T B &
\T P_{32} = \frac{1}{3}(\T B + 2\T C)\cr
\T P_{01} = \frac{1}{3}(2\T A + \T D) &
\T P_{03} = \T D\cr
\T P_{11} = \frac{1}{9}(4\T A + 2\T B + \T C + 2\T D) &
\T P_{13} = \T E\cr
\T P_{21} = \frac{1}{9}(2\T A + 4\T B + 2\T C + \T D) &
\T P_{23} = \T F\cr
\T P_{31} = \frac{1}{3}(2\T B + \T C) &
\T P_{33} = \T C\cr
\end{array}
\label{eq:control_points_cubic_param}
\end{equation}
\label{tl:control_points_param_bicubic}
\end{technical_lemma}

\begin{figure}[!pb]
\begin{narrow}{-0.3in}{0.0in}
\begin{minipage}[]{0.5\linewidth}
\subfigure[]
{
\includegraphics[clip,totalheight=2.2in]{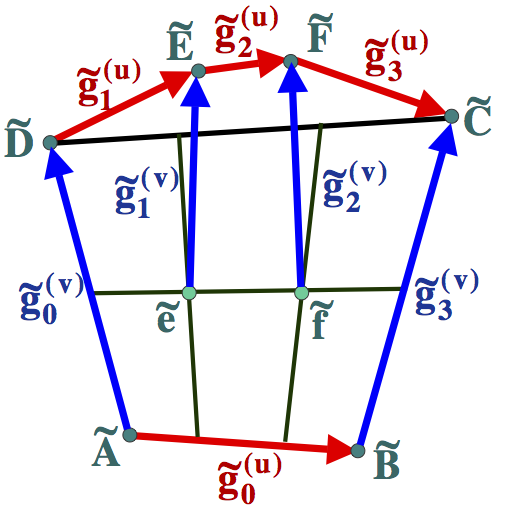}
\label{fig:fig30a}
}
\end{minipage}
\begin{minipage}[]{0.5\linewidth}
\subfigure[]
{
\includegraphics[clip,totalheight=2.2in]{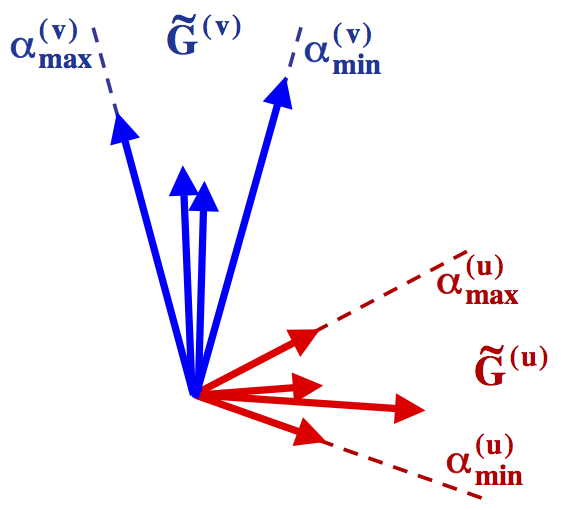}
\label{fig:fig30b}
}
\end{minipage}
\end{narrow}
\vspace{-0.2in}
\caption{An illustration for a sufficient condition for the regularity of bicubic in-plane parametrisation.}
\label{fig:fig30}
\end{figure}

\vspace{0.1in}
\begin{technical_lemma}
Let $\T P(u,v)$ be the bicubic in-plane parametrisation of a boundary mesh element defined according to Technical Lemma ~\ref{tl:control_points_param_bicubic} (see Figure ~\ref{fig:fig29}).
Let
\begin{equation}
\begin{array}{l}
\T e=\frac{1}{2}
\left\{\frac{1}{3}(2\T A+\T B)+\frac{1}{3}(\T C+2\T D)\right\}\cr
\T f=\frac{1}{2}
\left\{\frac{1}{3}(\T A+2\T B)+\frac{1}{3}(2\T C+\T D)\right\}
\end{array}
\end{equation}
and let two families of the planar vectors $\T g^{(u)}_i$ and $\T g^{(v)}_j$ ($i,j=0,\ldots,3$) be defined as follows (see Figure ~\ref{fig:fig30a})
\begin{equation}
\begin{array}{ll}
\T g^{(u)}_0 = \T B-\T A &\ \ \ \ \ \ \     
\T g^{(v)}_0 = \T D-\T A \cr
\T g^{(u)}_1 = \T E-\T D &\ \ \ \ \ \ \
\T g^{(v)}_1 = \T E-\T e \cr
\T g^{(u)}_2 = \T F-\T E &\ \ \ \ \ \ \
\T g^{(v)}_2 = \T F-\T f \cr
\T g^{(u)}_3 = \T C-\T F &\ \ \ \ \ \ \
\T g^{(v)}_3 = \T C-\T B \cr
\end{array}
\end{equation}
Let further $\T G^{(u)}$ and $\T G^{(v)}$ be the minimal infinite sectors which respectively contain all vectors from the first and second families, and let these sectors be bounded by rays corresponding to the angles $\alpha^{(u)}_{min}$, $\alpha^{(u)}_{max}$, $\alpha^{(v)}_{min}$, $\alpha^{(v)}_{max}$ (see Figure ~\ref{fig:fig30b}). Then
\begin{itemize}
\item[\bf (1)]
$\left\{\T g^{(u)}_i\right\}_{i=0}^3$ and 
$\left\{\T g^{(u)}_j\right\}_{j=0}^3$ are respectively generators of $\D{\T P}{u}$ and $\D{\T P}{v}$ in the meaning that for every parametric value $(u,v)\in[0,1]^2$, $\D{\T P}{u}$ and $\D{\T P}{v}$ can be written as a positive linear combination of these vectors. Here a positive linear combination is defined as such a linear combination that all its coefficients are non-negative and at least one coefficient is strictly positive.
\item[\bf (2)]
Relations $\alpha^{(u)}_{min} > -\alpha^{(v)}_{max}$ and
$\alpha^{(u)}_{max} < \alpha^{(v)}_{min}$ provide a sufficient condition for the regularity of parametrisation $\T P(u,v)$.
\end{itemize}
\label{tl:regularity_cubic}
\end{technical_lemma}
{\bf Proof}

\noindent
{\bf (1)} 
The first-order partial derivatives of in-plane parametrisation have the following representations in terms of the control points $\T P_{ij}$, $(i=0,\ldots,2,\ j=0,\ldots,3)$.
\begin{equation}
\begin{array}{l}
\D{\T P}{u}(u,v)=
3\sum_{i=0}^2\sum_{j=0}^3(\T P_{i+1,j}-\T P_{ij})B^2_i(u)B^3_j(v)\cr
\D{\T P}{v}(u,v)=
3\sum_{i=0}^3\sum_{j=0}^2(\T P_{i,j+1}-\T P_{ij})B^3_i(u)B^2_j(v)
\end{array}
\end{equation}
B\'ezier polynomials are always non-negative and for every parametric value $(u,v)\in[0,1]^2]$, at least one of the products $\left\{B^3_i(u)B^2_j(v)\right\}_{\lmtT{i=0,\ldots,3}{j=0,\ldots,2}{}}$ is strictly positive.
It implies that in order to prove that $\left\{g^{(u)}_i\right\}_{i=0}^3$ and 
$\left\{g^{(v)}_j\right\}_{j=0}^3$ serve as generators of $\D{\T P}{u}$ and $\D{\T P}{v}$, it is sufficient to show that every one of the first-order differences of the control points can be represented as a positive linear combination of the vectors from the corresponding family. Explicit formulas for the control points of in-plane parametrisation (see Equation 
~\ref{eq:control_points_cubic_param}) lead to the following expressions.

\noindent
The first-order differences in $u$-direction
\begin{equation}
\hspace{-0.4in}
\begin{array}{ll}
\T P_{10}-\T P_{00} = &\T P_{20}-\T P_{10} = \T P_{30}-\T P_{20} =
\frac{1}{3}(\T B-\T A) = \frac{1}{3}\T g^{(u)}_0\cr

\T P_{11}-\T P_{01} = &\T P_{21}-\T P_{11} = \T P_{31}-\T P_{21} =\cr
&\frac{1}{9}(2(\T B-\T A)+(\T C-\T D))=
\frac{1}{9}(2\T g^{(u)}_0 + \T g^{(u)}_1+ \T g^{(u)}_2+ \T g^{(u)}_3)\cr

\T P_{12}-\T P_{02}=&\frac{1}{9}((\T B-\T A)+(\T C-\T D)+3(\T E-\T D))=\cr &\frac{1}{9}(\T g^{(u)}_0+4\T g^{(u)}_1+\T g^{(u)}_2+\T g^{(u)}_3)\cr
\T P_{22}-\T P_{12}=&\frac{1}{9}((\T B-\T A)+(\T C-\T D)+3(\T F-\T E))=\cr &\frac{1}{9}(\T g^{(u)}_0+\T g^{(u)}_1+4\T g^{(u)}_2+\T g^{(u)}_3)\cr
\T P_{32}-\T P_{22}=&\frac{1}{9}((\T B-\T A)+(\T C-\T D)+3(\T C-\T F))=\cr &\frac{1}{9}(\T g^{(u)}_0+\T g^{(u)}_1+\T g^{(u)}_2+4\T g^{(u)}_3)\cr

\T P_{13}-\T P_{03} = &\T E-\T D = \T g^{(u)}_1\cr
\T P_{23}-\T P_{13} = &\T F-\T E = \T g^{(u)}_2\cr
\T P_{33}-\T P_{23} = &\T C-\T F = \T g^{(u)}_3\cr
\end{array}
\label{eq:dir_u_cubic}
\end{equation}

\noindent

The first-order differences in $v$-direction
\begin{equation}
\hspace{-0.4in}
\begin{array}{ll}
\T P_{01}-\T P_{00} = &\T P_{02}-\T P_{01} = \T P_{03}-\T P_{02}=
\frac{1}{3}(\T D-\T A) = \frac{1}{3} \T g^{(v)}_0\cr

\T P_{11}-\T P_{10} = &\frac{1}{9}(2(\T D-\T A)+(\T C-\T B))=
\frac{1}{9}(2\T g^{(v)}_0+\T g^{(v)}_1)\cr
\T P_{12}-\T P_{11}=&\frac{1}{18}(6(\T E-\T e)+2(\T D-\T A)+(\T C-\T B)=\cr
&\frac{1}{18}(2 g^{(v)}_0+6 g^{(v)}_1+ g^{(v)}_3)\cr
\T P_{13}-\T P_{12} = &\frac{2}{3}(\T E -\T e) =\frac{2}{3}\T g^{(v)}_1\cr

\T P_{21}-\T P_{20} = &\frac{1}{9}((\T D-\T A)+2(\T C-\T B))=
\frac{1}{9}(\T g^{(v)}_0+2\T g^{(v)}_1)\cr
\T P_{22}-\T P_{21}=&\frac{1}{18}(6(\T F-\T f)+(\T D-\T A)+2(\T C-\T B)=\cr
&\frac{1}{18}(g^{(v)}_0+6 g^{(v)}_2+ 2g^{(v)}_3)\cr
\T P_{23}-\T P_{22} = &\frac{2}{3}(\T F -\T f) =\frac{2}{3}\T g^{(v)}_2\cr

\T P_{31}-\T P_{30} = &\T P_{32}-\T P_{31} = \T P_{33}-\T P_{32}=
\frac{1}{3}(\T C-\T B) = \frac{1}{3} \T g^{(v)}_3\cr
\end{array}
\label{eq:dir_v_cubic}
\end{equation}
Relations ~\ref{eq:dir_u_cubic} and ~\ref{eq:dir_v_cubic} imply that for every parametric value of $(u,v)\in[0,1]^2$, the first-order partial derivatives $\D{\T P}{u}(u,v)$ and $\D{\T P}{v}(u,v)$ can be represented as a positive linear combination of the generators. It completes the straightforward proof of the first statement of Technical Lemma ~\ref{tl:regularity_cubic}.

\vspace{0.1in}
\noindent
{\bf(2)}
The proof of the regularity of the bicubic in-plane parametrisation is subdivided into a few parts according to the
different requirements listed in the definition of a regular parametrisation (Definition ~\ref{def:reg_param}).

\paragraph{Bijection.}
One should prove that $\T P(u,v)$ defines a bijection.
On the contrary, let there exist such parametric values $(u_0,v_0)\neq(u_1,v_1)$ that $\T P(u_0,v_0)=\T P(u_1,v_1)$.
Let segment $\sigma(t)$ connect these two values in the parametric domain
\begin{equation}
\sigma(t)=(u_0,v_0)(1-t)+(u_1,v_1)t
\end{equation}
Then $\T P(t)=\T P(\sigma(t))$ is a smooth closed curve in $XY$-plane. Let $f(t)$ be $Z$-coordinate of the vector product of $\T P(t)$ and some constant vector $\T d$
\begin{equation}
f(t)=\<\T d, \T P(t)\>
\end{equation}
The choice of $(u_0,v_0)$ and $(u_1,v_1)$ implies that  $f(0)=f(1)$. According to the Lagrange theorem, there exists such a value $\tau\in[0,1]$ that $f'(\tau)=0$.

In order to show that the assumption is not correct, it is sufficient to choose such a vector $\T d$ that $f'(t)\neq 0$ for every $t\in[0,1]$. In this case, the proof of bijection will be completed. The choice of $\T d$ is based on the formula 
\begin{equation}
f'(t)=\<\T d,\D{\T P}{u}(t)(u_1-u_0)+\D{\T P}{v}(t)(v_1-v_0)\>
\end{equation}
and is made in the following manner
\begin{itemize}
\item[-]
If $(u_1-u_0)(v_1-v_0)\ge 0$ then choose the direction of $\T d$ strictly between $-\alpha^{(v)}_{max}$ and $\alpha^{(u)}_{min}$. In this case 
$\<\T d,\D{\T P}{u}\> >0$, $\<\T d,\D{\T P}{v}\> >0$ for every $t$;\ $f'(t)$ for every $t$ has the same sign as $(u_1-u_0)$ (or as $(v_1-v_0)$ if $(u_1-u_0)=0$) and can not be equal to zero.
\item[-]
If $(u_1-u_0)(v_1-v_0)<0$ then choose the direction of $\T d$ strictly between $\alpha^{(u)}_{max}$ and $\alpha^{(v)}_{min}$. In this case 
$\<\T d,\D{\T P}{u}\> <0$, $\<\T d,\D{\T P}{v}\> >0$ for every $t$ and $f'(t)$ for every $t$ has the same sign as $(v_1-v_0)$.
\end{itemize} 

\paragraph{Smoothness.}
It is clear that $\T P(u,v)$ is $C^1$-smooth as a bicubic B\'ezier polynomial.

\paragraph{Regularity of Jacobian.}
Note, that $det(J^{(\T P)}(u,v))=\<\D{\T P}{u},\D{\T P}{v}\>$.
The regularity of the Jacobian $ J^{(\T P)}$ trivially follows from the facts that for every parametric value $(u,v)\in[0,1]^2$, the partial derivatives $\D{\T P}{u}$ and $\D{\T P}{v}$ can be represented as some positive linear combinations of generators $\left\{\T g^{(u)}_i\right\}_{i=0}^3$
and $\left\{\T g^{(v)}_j\right\}_{j=0}^3$ respectively and that 
$\<\T g^{(u)}_i,\T g^{(v)}_j\> > 0$ for every $i,j=0,\ldots,3$ due to the assumed separation between $\T G^{(u)}$ and $\T G^{(v)}$.

\nopagebreak 
\eop${}_{{\bf Technical Lemma ~\ref{tl:regularity_cubic}}}$

\vspace{0.1in}
\begin{technical_lemma}
Let a planar mesh have piecewise-cubic $G^1$-smooth global boundary and
$\T P(u,v)$ be the restriction of the global in-plane parametrisation 
$\T\Pi^{(bicubic)}$ on such a boundary mesh element, that its upper edge lies on the global boundary (see Figure ~\ref{fig:fig29}).
Then the-first order partial derivatives of $\T P(u,v)$ along the left and the right edges of the element have the following explicit form.

\vspace{0.05in}
\noindent
The partial derivative in the direction along the left edge 
\vspace{-0.1in}
\begin{equation}
\hspace{-0.4in}
\begin{array}{l}
\D{\T P}{v}(0,v) = \T D\!-\!\T A\cr
\end{array}
\end{equation}

\noindent
The partial derivative in the cross direction for the left edge
\vspace{-0.1in}
\begin{equation}
\hspace{-0.4in}
\begin{array}{l}
\D{\T P}{u}(0,v) =
(\T B\!-\!\T A)B^2_0(v)+
\frac{1}{2}
\left((\T B\!-\!\T A)\!+\!(\T C\!-\!\T D)\right)B^2_1(v)+
3(\T E\!-\!\T D)B^2_2(v)\cr
\end{array}
\end{equation}

\noindent
The partial derivative in the direction along the right edge
\vspace{-0.1in}
\begin{equation}
\hspace{-0.4in}
\begin{array}{l}
\D{\T P}{v}(1,v) = \T C\!-\!\T B\cr
\end{array}
\end{equation}

\noindent
The partial derivative in the cross direction for the right edge
\vspace{-0.1in}
\begin{equation}
\hspace{-0.4in}
\begin{array}{l}
\D{\T P}{u}(1,v) =  
(\T B\!-\!\T A)B^2_0(v)+
\frac{1}{2}
\left((\T B\!-\!\T A)\!+\!(\T C\!-\!\T D)\right)B^2_1(v)+
3(\T C\!-\!\T F)B^2_2(v)\cr
\end{array}
\end{equation}
\label{tl:partial_derivatives_in_plane_bicubic}
\end{technical_lemma}

\vspace{0.1in}
\begin{technical_lemma}
Let conventional weight functions $c(v)$, $l(v)$, $r(v)$ for the common edge of two adjacent boundary mesh elements be defined according to global in-plane parametrisation $\T\Pi^{(bicubic)}$. Then the coefficients of the weight functions have the following representations in terms of the initial mesh data (see Figure ~\ref{fig:fig31})
\begin{equation}
\begin{array}{ll}
l_0=\<\RR-\GG,\GG'-\GG\> \cr
l_1=\frac{1}{2}\<(\RR-\GG)+(\RR'-\GG'),\GG'-\GG\> \cr
l_2=3\<\TTR'-\GG',\GG'-\GG\>\cr
r_0=-\<\GG-\LL,\GG'-\GG\>\cr
r_1=-\frac{1}{2}\<(\GG-\LL)+(\GG'-\LL'),\GG'-\GG\>\cr
r_2=-3\<\GG'-\TTL',\GG'-\GG\>\cr
\end{array}
\label{eq:lr_coeff_bicubic}
\end{equation}
\begin{equation}
\begin{array}{ll}
c_0=&\<\GG-\LL,\RR-\GG\>\cr
c_1=&\frac{1}{4}\left(\<\GG-\LL,\RR'-\GG'\>+\<\GG'-\LL',\RR-\GG\>+
\right.\cr&\left.
2\<\GG-\LL,\RR-\GG\>\right)\cr
c_2=&\frac{1}{6}\left(3\<\GG-\LL,\TTR'-\GG'\>+3\<\GG'-\TTL',\RR-\GG\>+
\right.\cr&\left.
\<(\GG-\LL)+(\GG'-\LL'),(\RR-\GG)+(\RR'-\GG')\>\right)\cr
c_3=&\frac{3}{4}\left(\<(\GG-\LL)+(\GG'-\LL'),\TTR'-\GG'\>+
\right.\cr&\left.
\<\GG'-\TTL',(\RR-\GG)+(\RR'-\GG')\>\right)\cr
c_4=&9\<\GG'-\TTL',\TTR'-\GG'\>
\end{array}
\label{eq:c_coeff_bicubic}
\end{equation}
\label{tl:weight_coeff_cubic_explicit}
\end{technical_lemma}

\vspace{0.1in}
\begin{technical_lemma} 
Let conventional weight functions $c(v)$, $l(v)$, $r(v)$ for the common edge of two adjacent boundary mesh elements be defined according to global in-plane parametrisation $\T\Pi^{(bicubic)}$. Then the following relations between the coefficients of the weight functions with respect to B\'ezier and the power bases hold
\begin{equation}
\begin{array}{lcl}
l^{(power)}_0=l_0 &\ \ \ \ \ \ \ \ \ \ \ \ \ &r^{(power)}_0=r_0\cr
l^{(power)}_1=2(l_1-l_0) & & r^{(power)}_1=2(r_1-r_0)\cr
l^{(power)}_2=l_0-2l_1+l_2 & & r^{(power)}_2=r_0-2r_1+r_2
\end{array}
\end{equation}
\begin{equation}
\begin{array}{l}
c^{(power)}_0=c_0\cr
c^{(power)}_1=4(-c_0+c_1)\cr 
c^{(power)}_2=6(c_0-2c_1+c_2)\cr
c^{(power)}_3=4(-c_0+3c_1-3c_2+c_3)\cr
c^{(power)}_4=c_0-4c_1+6c_2-4c_3+c_4
\end{array}
\end{equation}
\label{tl:weight_coeff_cubic_Bezier_power}
\end{technical_lemma}

\vspace{0.1in}
\begin{technical_lemma}
Let a planar mesh have either a polygonal or a smooth global boundary and let global in-plane parametrisation $\T\Pi^{(bilinear)}$ or $\T\Pi^{(bicubic)}$ be considered.  Further, for the inner vertex of an edge with one boundary vertex, let the system of all {\it $"Eq^{(4)}(0)"$-type} and {\it $"Eq^{(4)}(1)"$-type} equations be satisfied by classification of $V^{(4)}$-type, $E^{(4)}$-type, $D^{(4)}$-type and $T^{(4)}$-type control points adjacent to the vertex. In particular, for the edge with one boundary vertex, control points $\T (L^{(4)}_i,\T  C^{(4)}_i,\T R^{(4)}_i)$ for $i=0,1$ and $\T C^{(4)}_2$ are classified (see Figure ~\ref{fig:fig36a}). Then 
\begin{itemize}
\item[{\bf (1)}]
The degree elevation formulas for control points $(\bar L^{(5)}_i, \bar C^{(5)}_i, \bar R^{(5)}_i)$ for $i=0,1$ and $\bar C^{(5)}_2$ (see Figure ~\ref{fig:fig36b}) involve only control points $(\bar L^{(4)}_i, \bar C^{(4)}_i, \bar R^{(4)}_i)$ for $i=0,1$ and $\bar C^{(4)}_2$.
\item[{\bf (2)}] Control points $(L^{(5)}_i, C^{(5)}_i, R^{(5)}_i)$ ($i=0,1$) and $C^{(5)}_2$, computed according to the degree elevation formulas, satisfy equations {\it $"Eq^{(5)}(0)"$} and {\it $"Eq^{(5)}(1)"$} for the edge with one boundary vertex.
\end{itemize}
Here superscripts ${}^{(4)}$ and ${}^{(5)}$ are used in order to specify whether some equation or control point relates to degree $4$ or $5$ of the resulting patch.
\label{tl:degree_combined}
\end{technical_lemma}

\noindent
{\bf Proof}\hspace{4in}

\noindent
{\bf(1)} The first statement of the Technical Lemma immediately follows from the straightforward formulas 
\begin{equation}
\begin{array}{l}
\bar R^{(5)}_0=
1/5(4 \bar R^{(4)}_0 + \bar C^{(4)}_0)\cr
\bar C^{(5)}_0=
\bar C^{(4)}_0\cr
\bar L^{(5)}_0=
1/5(4 \bar L^{(4)}_0 + \bar C^{(4)}_0)\cr
\bar R^{(5)}_1=
1/25(16 \bar R^{(4)}_1+4 \bar R^{(4)}_0+
4 \bar C^{(4)}_1+\bar C^{(4)}_0)\cr
\bar C^{(5)}_1=
1/5(4 \bar C^{(4)}_1 + \bar C^{(4)}_0)\cr
\bar L^{(5)}_1=
1/25(16 \bar L^{(4)}_1+4 \bar L^{(4)}_0+4 \bar C^{(4)}_1+\bar C^{(4)}_0)\cr
\bar C^{(5)}_2=
1/5(3 \bar C^{(4)}_2 + 2 \bar C^{(4)}_1)
\end{array}
\label{eq:degree_elevation_four_five}
\end{equation}
The important implication is that as soon as the control points  $(\T L^{(4)}_i, \T C^{(4)}_i, \T R^{(4)}_i)$ ($i=0,1$) and $\T C^{(4)}_2$ are classified, control points $(\T L^{(5)}_i, \T C^{(5)}_i, \T R^{(5)}_i)$ ($i=0,1$) and $\T C^{(5)}_2$ are also classified and dependencies of the dependent control points from $\GCP{5}$ are fully defined.

\vspace{0.1in}
\noindent
{\bf(2)} It is sufficient to prove the second statement of the Technical Lemma in the case when the conventional weight functions correspond to global bilinear in-plane parametrisation $\T\Pi^{(bicubic)}$. Indeed, Lemma ~\ref{lemma:Eq0_Eq1_bicubic} shows that for an edge with one boundary vertex, the system of equations  {\it "Eq(0)"} and {\it "Eq(1)"} in the case of $\T\Pi^{(bicubic)}$ is equivalent to the system in the case of $\T\Pi^{(bilinear)}$.

Equations {\it $"Eq^{(5)}(0)"$} and {\it $"Eq^{(5)}(1)"$} involve coefficients of the weight functions (which clearly do not depend on chosen degree $4$ or $5$ of MDS) and the first-order differences of the control points 
$(\dL^{(5)}_i, \dC^{(5)}_i, \dR^{(5)}_i)$, $i=0,1$ 
(which can be expressed in terms of $(\dL^{(4)}_i, \dC^{(4)}_i, \dR^{(4)}_i)$, $i=0,1$ according to Equation ~\ref{eq:degree_elevation_four_five}). One can easily verify, that
\begin{equation}
\begin{array}{l}
"Eq^{(5)}(0)" = 4/5 "Eq^{(4)}(0)"\cr
"Eq^{(5)}(0)" = 4/5 \left("Eq^{(4)}(0)" + "Eq^{(4)}(1)"\right)
\end{array}
\label{eq:eq0_eq1_four_five}
\end{equation}
The last two equations clearly imply that statement {\bf (2)} of the Technical Lemma is satisfied.

\nopagebreak 
\eop${}_{{\bf Technical\ Lemma ~\ref{tl:degree_combined}}}$

\FloatBarrier

\eject
\section{Proofs}
\label{sect:proofs}

\vspace{0.15in}
{\bf Proof of Theorem ~\ref{theorem:G1_equiv_Z_equation}}

First, it will be shown that the satisfaction of Equation 
~\ref{eq:G1_Z_equation} with scalar functions given in Equation 
~\ref{eq:weights_general} is \underline{sufficient} in order to guarantee a $G^1$-smooth concatenation  between the adjacent patches. It will be verified, that the scalar functions given in Equation ~\ref{eq:weights_general} satisfies all requirements of Definition ~\ref{def:G1} and may be used as the weight function. 

\underline{\it Equation ~\ref{eq:G1_eq_general_1}.}
It should be verified that Equation ~\ref{eq:G1_eq_general_1} is satisfied for $X$ and $Y$ components of the partial derivatives $\bar L_u$, $\bar R_u$ and $\bar L_v$. For example, for $X$-component, one gets
\begin{equation}
\begin{array}{l}
LX_u(v)l(v)+RX_u(v)r(v)+LX_v(v)c(v)=\cr
LX_u<\tilde R_u,\tilde L_v>-RX_u<\tilde L_u,\tilde L_v>+LX_v<\tilde L_u, \tilde R_u>=\cr
mix
\left(\begin{array}{ccc}\tilde L_u&LX_u\cr \tilde R_u&RX_u\cr \tilde L_v&LX_v\cr\end{array}\right)=
det
\left(\begin{array}{ccc}LX_u&LY_u&LX_u\cr RX_u&RY_u&RX_u\cr LX_v&LY_v&LX_v\cr\end{array}\right)=0.
\end{array}
\end{equation}
The same proof works for $Y$-components of the partial derivatives.

\underline{\it Equation ~\ref{eq:G1_eq_general_2}.}
It is necessary to verify that $l(v)r(v)<0$ for every $v\in[0,1]$.
Equality $l(v)r(v)=0$ is impossible as a trivial consequence of the regularity of parametrisation. Indeed, for example equality to zero of $r(v)$ implies that
\begin{equation}
\hspace{-0.3in}
0=r(v)=-<\tilde L_u,\tilde L_v>=
-det\left(\begin{array}{cc}\D{LX}{u} & \D{LY}{u}\cr \D{LX}{v} & \D{LY}{v}\end{array}\right)=
-det(J^{(\T L)}(1,v))
\end{equation}
that contradicts the regularity of parametrisation for the left patch.

\vspace{0.1in}
\noindent
$l(v)r(v)$ is a continuous function and $l(0)=<\tilde R_u(0),\tilde L_v(0)> > 0$ and $r(0)=-<\tilde L_u(0), \tilde L_v(0)> < 0$ because both $\tilde R_u(0)$ and $\tilde L_u(0)$ clearly point to the right side with respect to the tangent vector of the common boundary $\tilde L_v(0)$. Therefore $l(v)r(v)\le 0$ for every $v\in[0,1]$.

\underline{\it Equation ~\ref{eq:G1_eq_general_3}.}
One should show that $<\bar L_u,\bar L_v>\neq 0$ for every $v\in[0,1]$. It immediately follows
from the fact that $Z$-component of this cross product is equal to
$<\tilde L_u, \tilde L_v>=-r(v)$ which is already proven to be non-zero. 

\vspace{0.15in}
Now it will be shown that if the adjacent patches join the $G^1$-smoothly, then Equation ~\ref{eq:G1_Z_equation} with scalar functions given in Equation ~\ref{eq:weights_general} is \underline{necessarily} satisfied.

Let the adjacent patches join the $G^1$-smoothly. Then, according to Definition ~\ref{def:G1}, there exist such scalar functions $\breve{l}(v)$, $\breve{r}(v)$, $\breve{c}(v)$ that Equations ~\ref{eq:G1_eq_general_1}, ~\ref{eq:G1_eq_general_2}, ~\ref{eq:G1_eq_general_3} are satisfied. In particular, Equation ~\ref{eq:G1_eq_general_1} implies that $XY$-components of the partial derivatives satisfy the equation
\begin{equation}
\T{L}_u\breve{l}(v)+\T{R}_u\breve{r}(v)+\T{L}_v\breve{c}(v)=0
\end{equation}
The result of the cross product of the last expression respectively with $\T L_v$ and $\T R_u$ is clearly equal to zero. Together with the inequality $\<\T R_u,\T L_v\>\neq 0\>$ (which follows from the regularity of in-plane parametrisation), it leads to the conclusion that the following relations between $\breve{l}(v)$, $\breve{r}(v)$, $\breve{c}(v)$ necessarily hold
\begin{equation}
\breve{r}(v)=-\frac{\<\T L_u,\T L_v\>}{\<\T R_u,\T L_v\>}
\breve{l}(v),\ \ \ \ \ \ \ \ \ \ \ \ \ \ \ \ \ \ \ 
\breve{c}(v)=\frac{\<\T L_u,\T R_u\>}{\<\T R_u,\T L_v\>}
\breve{l}(v)
\end{equation}
and Equation ~\ref{eq:G1_eq_general_1} necessarily has the form
\begin{equation}
\bar L_u 
\breve{l}(v)-
\bar R_u 
\frac{\<\T L_u,\T L_v\>}{\<\T R_u,\T L_v\>}\breve{l}(v)+
\bar L_v
\frac{\<\T L_u,\T R_u\>}{\<\T R_u,\T L_v\>}\breve{l}(v)=
0
\label{eq:aux_G1_eq_Z_equation}
\end{equation}

According to Equation ~\ref{eq:G1_eq_general_2},\  $\breve{l}(v)\neq 0$ for any $v\in[0,1]$. Therefore it is possible to divide both parts of Equation ~\ref{eq:aux_G1_eq_Z_equation} by $\frac{\breve{l}(v)}{ \<\T R_u,\T L_v\>}$. Therefore one sees, that Equation ~\ref{eq:G1_Z_equation} (and even Equation ~\ref{eq:G1_eq_general_1}) is satisfied for the scalar functions given in Equation ~\ref{eq:weights_general}.

\nopagebreak
\eop${}_{{\bf Theorem ~\ref{theorem:G1_equiv_Z_equation}}}$

\vspace{0.32in}
\noindent

\vspace{0.32in}
\noindent
{\bf Proof of Lemma ~\ref{lemma:linear_1}}
Proof of the Lemma is straightforward. Let $\bar L=(\T L,L)$ and $\bar R=(\T R,R)$ be the restrictions of $\bar\Psi$ on two adjacent mesh elements. $\bar\Psi$ agrees with degree $m$ global regular in-plane parametrisation $\T\Pi$, therefore $\T L$ and $\T R$ are polynomials of degree $m$. Conventional weight functions $l(v)$, $r(v)$ and $c(v)$ are evidently polynomials of (formal) degrees $2m-1$, $2m-1$ and $2m$ respectively. $\bar\Psi\in\FUN{(n)}$ and so $R$ and $L$ are polynomials of (formal) degree $n$.
Using B\'ezier representation of $Z$-components of the partial derivatives $\bar L_u$, $\bar R_u$, $\bar R_v$ and of the weight functions $l$, $r$, $c$ 

\begin{equation}
\begin{array}{ll}
L_u = n\sum_{j=0}^n \Delta L_j B^n_j & l = \sum_{k=0}^{2m-1} l_kB^{2m-1}_k\cr
R_u = n\sum_{j=0}^n \Delta R_j B^n_j & r = \sum_{k=0}^{2m-1} r_kB^{2m-1}_k\cr
L_v = n\sum_{j=0}^{n-1} \Delta C_j B^{n-1}_j & c = \sum_{k=0}^{2m} c_kB^{2m}_k\cr
\end{array}
\end{equation}
one gets:
\begin{equation}
\hspace{-0.4in}
\begin{array}{l}
L_u(v)l(v)+R_u(v)r(v)+L_v(v)c(v)=\cr
n\!\sum\limits_{s=0}^{n+2m-1}\!\frac{1}{\Cnk{n+2m-1}{s}}\!
\left\{
\sum_{\lmtT{j+k=s}{0\le j\le n}{0\le k\le 2m-1}}
\Cnk{n}{j}\Cnk{2m-1}{k} (l_k\Delta L_j+r_k\Delta R_j)\right. +\cr
\hspace{1.25in}\left.
\sum_{\lmtT{j+k=s}{0\le j\le n-1}{0\le k\le 2m}}
\Cnk{n-1}{j}\Cnk{2m}{k} c_k\Delta C_j
\right\}B^{n+2m-1}_s(v)
\end{array}
\end{equation}

Writing down equality to zero for every coefficient of the resulting B\'ezier polynomial leads to Equation ~\ref{eq:G1_Z_general}.

The number of necessary and sufficient equations is equal to the actual degree of the polynomial $L_u(v)l(v)+R_u(v)r(v)+L_v(v)c(v)$, that is \\
$max\{max\_deg(l,r)+n,\ \ deg(c)+n-1\}$. 

\nopagebreak 
\eop${}_{{\bf Lemma ~\ref{lemma:linear_1}}}$

\vspace{0.32in}
\noindent

\nopagebreak 

\vspace{0.32in}
\noindent
{\bf Proof of Lemma ~\ref{lemma:eq1_geom_edge_j}}
It is sufficient to prove the Lemma for two adjacent patches. 

Let two adjacent patches $\bar L$ and $\bar R$ be defined by a bilinear in-plane parametrisation, the $XY$-components of the first and second order $3D$ partial derivatives, computed at the common vertex (see Subsection ~\ref{subsect:def_deriv_patches_at_common_vertex}), 
have the following representation in terms of the initial mesh data:
\begin{equation}
\hspace{-0.3in}
\begin{array}{lcl}
\bar\epsilon^{(R)}=(\T e^{(R)}, n\Delta R_0)&
\ \ \ \ \ \ \ \ \ &
\bar\tau^{(R)}=(\T t^{(R)}, n^2(\Delta R_1-\Delta R_0))\cr
\bar\epsilon^{(L)}=(\T e^{(L)}, -n\Delta L_0)& &
\bar\tau^{(L)}=(\T t^{(L)}, -n^2(\Delta L_1-\Delta L_0))\cr
\bar\epsilon^{(C)}=(\T e^{(C)}, n\Delta C_0)& &
\end{array}
\end{equation}
Here $\T e^{(R)}$, $\T e^{(R)}$, $\T e^{(R)}$, $\T e^{(R)}$, $\T e^{(R)}$ are directed planar edges and $\bar\tau^{(R)}$, $\bar\tau^{(L)}$ are the twist characteristics of the planar mesh elements (see Subsection ~\ref{subsect:def_vertices_edges_twists_planar}). 

Let {\it "Eq(0)"} be satisfied for the common edge of $\bar L$ and $\bar R$. It is sufficient to show that {\it "Eq(1)"} is equivalent to the following equation
\begin{equation}
tw^{(R)}+tw^{(L)}=coeff^{(C)}\delta^{(C)}
\label{eq:eq_compact_two_patches}
\end{equation}
\vspace{-0.05in}
where
\vspace{-0.05in}
\begin{equation}
\begin{array}{lr}
\hspace{-0.3in}
tw^{(R)} = 
\frac{1}{\<\T e^{(R)}, \T e^{(C)}\>^2} mix\VT{\bar\tau^{(R)}}{\bar\epsilon^{(R)}}{\bar\epsilon^{(C)}}
& 
tw^{(L)} = 
\frac{1}{\<\T e^{(C)}, \T e^{(L)}\>^2} mix\VT{\bar\tau^{(L)}}{\bar\epsilon^{(C)}}{\bar\epsilon^{(L)}}
\end{array}
\end{equation}
\vspace{-0.1in}
\begin{equation}
\hspace{-0.3in}
\ \ coeff^{(C)}=\frac{\<\T e^{(R)},\T e^{(L)}\>}{\<\T e^{(R)},\T e^{(C)}\>\<\T e^{(C)},\T e^{(L)}\>}
\end{equation}
\vspace{-0.05in}
and $\delta^{(C)}$ is $Z$-component of the second-order derivative along the common edge (see Subsection ~\ref{subsect:def_deriv_patches_at_common_vertex}).

Indeed, the following relations between coefficients of the weight functions and geometrical characteristics of the planar mesh hold
\begin{equation}
\begin{array}{ll}
l_0=\ \ \<\T e^{(R)},\T e^{(C)}\>,&  l_1-l_0=\<\T t^{(R)},\T e^{(C)}\>\cr
r_0=-\<\T e^{(C)},\T e^{(L)}\>,&  r_1-r_0=\<\T t^{(L)},\T e^{(C)}\>
\end{array}
\end{equation}
According to these formulas, the first two summands of {\it "Eq(1)"} can be rewritten as follows
\begin{equation}
\hspace{-0.4in}
\begin{array}{l}
n(l_0\Delta L_1+r_0\Delta R_1)+(l_1\Delta L_0+r_1\Delta R_0)=\cr
n(l_0(\Delta L_1-\Delta L_0)+r_0(\Delta R_1-\Delta R_0))+
\<\T t^{(R)},\T e^{(C)}\>\Delta L_0+\<\T t^{(L)},\T e^{(C)}\>\Delta R_0+\cr
(n+1)(l_0\Delta L_0+r_0\Delta R_0)=\cr
\frac{1}{n}\left(\<\T e^{(C)},\T e^{(R)}\>\tau^{(L)}+\<\T t^{(L)},\T e^{(C)}\>\epsilon^{(R)}\right)-
\frac{1}{n}\left(\<\T e^{(C)},\T e^{(L)}\>\tau^{(R)}+\<\T t^{(R)},\T e^{(C)}\>\epsilon^{(L)}\right)+\cr
(n+1)(l_0\Delta L_0+r_0\Delta R_0)=\cr
\frac{1}{n}\left\{mix\VT{\bar\tau^{(L)}}{\bar\epsilon^{(C)}}{\bar\epsilon^{(R)}}+
                 \<\T t^{(L)},\T e^{(R)}\>\epsilon^{(C)}\right\}-
\frac{1}{n}\left\{mix\VT{\bar\tau^{(R)}}{\bar\epsilon^{(C)}}{\bar\epsilon^{(L)}}+\<\T t^{(R)},\T e^{(L)}\>\epsilon^{(C)}\right\}+\cr
(n+1)(l_0\Delta L_0+r_0\Delta R_0)=\cr
\frac{1}{n}\left\{mix\VT{\bar\tau^{(L)}}{\bar\epsilon^{(C)}}{\bar\epsilon^{(R)}}- mix\VT{\bar\tau^{(R)}}{\bar\epsilon^{(C)}}{\bar\epsilon^{(L)}}\right\}+\cr
\left(\<\T t^{(L)},\T e^{(R)}\>-\<\T t^{(R)},\T e^{(L)}\>\right)\Delta C_0+
(n+1)(l_0\Delta L_0+r_0\Delta R_0)
\end{array}
\end{equation}
The assumption that {\it "Eq(0)"} is satisfied implies that $l_0\Delta L_0+r_0\Delta R_0=-c_0\Delta C_0$ and so {\it "Eq(1)"} is
equivalent to the equation
\begin{equation}
\hspace{-0.4in}
\begin{array}{ll}
0=&\frac{1}{n}\left\{mix\VT{\bar\tau^{(L)}}{\bar\epsilon^{(C)}}{\bar\epsilon^{(R)}}-
mix\VT{\bar\tau^{(R)}}{\bar\epsilon^{(C)}}{\bar\epsilon^{(L)}}\right\}+\cr
&\left\{\<\T t^{(L)},\T e^{(R)}\>-\<\T t^{(R)},\T e^{(L)}\>+2c_1-(n+1)c_0\right\}\Delta C_0+(n-1)c_0\Delta C_1
\end{array}
\label{eq:eq_aux_1}
\end{equation}
It can easily be shown, that 
\begin{equation}
\hspace{-0.4in}
\<\T t^{(L)},\T e^{(R)}\>-\<\T t^{(R)},\T e^{(L)}\>+2c_1=2\<\RR-\GG,\LL-\GG\>=2\<\T e^{(R)},\T e^{(L)}\> = 2c_0
\end{equation}
Therefore Equation ~\ref{eq:eq_aux_1} can be further simplified as follows
\begin{equation}
\hspace{-0.35in}
\begin{array}{ll}
0=&\frac{1}{n}\left\{
mix\VT{\bar\tau^{(L)}}{\bar\epsilon^{(C)}}{\bar\epsilon^{(R)}}-	
mix\VT{\bar\tau^{(R)}}{\bar\epsilon^{(C)}}{\bar\epsilon^{(L)}}
\right\}+
(n-1)c_0(\Delta C_1-\Delta C_0)=\cr  &\frac{1}{n}\left\{
mix\VT{\bar\tau^{(L)}}{\bar\epsilon^{(C)}}{\bar\epsilon^{(R)}}-mix\VT{\bar\tau^{(R)}}{\bar\epsilon^{(C)}}{\bar\epsilon^{(L)}}+
\<\T e^{(R)}, \T e^{(L)}\>\delta^{(C)}\right\}
\end{array}
\label{eq:eq_aux_2}
\end{equation}

Now equation {\it "Eq(0)"} is used in order to make the last equation more symmetric. 
The tangent coplanarity in particular means that
\begin{equation}
\begin{array}{ll}
\bar\epsilon^{(R)}=\frac{1}{\<\T e^{(C)},\T e^{(L)}\>}\left(\<\T e^{(R)},\T e^{(L)}\>\bar\epsilon^{(C)}-
                                                \<\T e^{(R)},\T e^{(C)}\>\bar\epsilon^{(L)}\right)\cr
\bar\epsilon^{(L)}=\frac{1}{\<\T e^{(R)},\T e^{(C)}\>}\left(\<\T e^{(R)},\T e^{(L)}\>\bar\epsilon^{(C)}-
                                                \<\T e^{(C)},\T e^{(L)}\>\bar\epsilon^{(R)}\right)\cr
\end{array}
\end{equation}
Substitution of these formulas in the mixed products of Equation ~\ref{eq:eq_aux_2} gives
\begin{equation}
\hspace{-0.35in}
0=\frac{\<\T e^{(R)},\T e^{(C)}\>}{\<\T e^{(C)},\T e^{(L)}\>}
mix\!\VT
{\!\bar\tau^{(L)}\!}
{\!\bar\epsilon^{(L)}\!}
{\!\bar\epsilon^{(C)}\!}+
\!\!
\frac{\<\T e^{(C)},\T e^{(L)}\>}{\<\T e^{(R)},\T e^{(C)}\>}
mix\!\VT
{\!\bar\tau^{(R)}\!}
{\!\bar\epsilon^{(C)}\!}
{\!\bar\epsilon^{(R)}\!}
\!\!
+\<\T e^{(R)}, \T e^{(L)}\>\delta^{(C)}
\end{equation}
In order to complete the proof it only remains to divide all the members of the last equation
by $\<\T e^{(R)},\T e^{(C)}\>\<\T e^{(C)},\T e^{(L)}\>$ which is not equal to zero due to the strict 
convexity of the quadrilateral mesh elements.

\nopagebreak 
\eop${}_{{\bf Lemma ~\ref{lemma:eq1_geom_edge_j}}}$

\vspace{0.32in}
\noindent
{\bf Proof of Lemma ~\ref{lemma:regular_vertex}}
The strict convexity of the mesh elements implies that any inner even vertex $\T V$ has degree $val(V)=4$ at least.

Let $coeff^{(j)}=0$ for every $j=1,\ldots,val(V)$.
In particular, $coeff^{(2)}=0$ and so $\TTe{1}$ and $\TTe{3}$ are colinear and lie on some straight line $\T l^{(1,3)}$;  $coeff^{(3)}=0$ and so $\TTe{2}$ and $\TTe{4}$ are colinear and lie on some straight line $\T l^{(2,4)}$ (see Figure ~\ref{fig:fig16}).
Therefore for $val(V)=4$ the vertex is proven to be regular.

It remains to show that $val(V)$ could not be greater than $4$. Indeed, let $val(V)>4$. Then, $\TTe{2}$ should be colinear to both $\TTe{4}$ and $\TTe{val(V)}$ (because both $coeff^{(3)}$ and $coeff^{(1)}$ are equal to zero). But $\TTe{4}$ and $\TTe{val(V)}$ can not be colinear because $\TTe{val(V)}$ lies strictly between 
$\TTe{4}$ and $\TTe{deg(1)}$ which span an angle less than $\pi$ due to the strict convexity of the mesh elements.


\nopagebreak 
\eop${}_{{\bf Lemma ~\ref{lemma:regular_vertex}}}$

\vspace{0.32in}
\noindent
{\bf Proof of Lemma ~\ref{lemma:C2_sufficient}}

\vspace{0.1in}
\noindent
{\bf (1)}
Compatibility of $\delta^{(j)}$ $(j=1,\ldots,val(V))$ with some second-order functional derivatives in the functional sense means, that there exist such three scalars $Z_{XX}$, $Z_{XY}$, $Z_{YY}$ that the following relation holds for every $j=1,\ldots,val(V)$
\begin{equation}
\begin{array}{ll}
\delta^{(j)}=&(\alpha^{(j)},\beta^{(j)})
\MM{Z_{XX}} {Z_{XY}} {Z_{XY}} {Z_{YY}}
\VV{\alpha^{(j)}} {\beta^{(j)}} =\cr
&\left(\alpha^{(j)}\right)^2 Z_{XX}+
2\alpha^{(j)}\beta^{(j)} Z_{XY}+
\left(\beta^{(j)}\right)^2 Z_{YY}
\end{array}
\label{eq:eq_ff_agreement}
\end{equation} 
where $\alpha^{(j)}$ and $\beta^{(j)}$ are respectively $X$ and $Y$-components of the planar edge $\TTe{j}$.

In order to prove that the {\it "Circular Constraint"} is satisfied one should prove that
\begin{equation}
\hspace{-0.4in}
\begin{array}{l}
0=\sum\limits_{j=1}^{val(V)}
(-1)^j\frac{\<\TTe{j-1},\TTe{j+1}\>}
{\<\TTe{j-1},\TTe{j}\>\<\TTe{j},\TTe{j+1}\>}\delta^{(j)}=\cr
\!\!\sum\limits_{j=1}^{val(V)}\!\! 
(-1)^j\frac{\<\TTe{j-1},\TTe{j+1}\>}
{\<\TTe{j-1},\TTe{j}\>\<\TTe{j},\TTe{j+1}\>}
\left\{\!\left(\alpha^{(j)}\right)^2\!\! Z_{XX}\!\!+
2\alpha^{(j)}\beta^{(j)} Z_{XY}\!\!+
\left(\beta^{(j)}\right)^2\!\! Z_{YY}\!\!\right\}
\end{array}
\end{equation} 
independently of the values of $Z_{XX}$,$Z_{XY}$ and $Z_{YY}$. In other words it should be proven that the coefficient of every one of $Z_{XX}$,$Z_{XY}$ and $Z_{YY}$ is equal to zero. Here it will be shown that the coefficient of $Z_{XY}$ is equal to zero; proof
for the coefficients of $Z_{XX}$ and $Z_{YY}$ can be given in the same way. Coefficient of $Z_{XY}$ is equal to  
\begin{equation}
2\sum_{j=1}^{val(V)}(-1)^j \eta^{(j)}
\end{equation}
where
\begin{equation}
\tiny \eta^{(j)}=
\frac{(\alpha^{(j-1)}\beta^{(j+1)}-\alpha^{(j+1)}\beta^{(j-1)})\alpha^{(j)}\beta^{(j)}}
{(\alpha^{(j-1)}\beta^{(j)}-\alpha^{(j)}\beta^{(j-1)})
 (\alpha^{(j)}\beta^{(j+1)}-\alpha^{(j+1)}\beta^{(j)})}
\end{equation}
It can be shown by a straightforward computation that
\begin{equation}
\eta^{(j)}=1+\xi^{(j-1)}+\xi^{(j)}, {\rm\ \ where\ \ }
\xi^{(j)}=\frac{\beta^{(j)}\alpha^{(j+1)}}
{\alpha^{(j)}\beta^{(j+1)}-\alpha^{(j+1)}\beta^{(j)}}
\label{eq:eq_aux_4}
\end{equation}
The first statement of the Lemma immediately follows from the last equation because
$2\sum_{j=1}^{val(V)}(-1)^j \eta^{(j)}=
2\sum_{j=1}^{val(V)}(-1)^j(1+\xi^{(j-1)}+\xi^{(j)})$ is clearly equal to zero for an even value of $val(V)$.

\vspace{0.1in}
\noindent
{\bf (2)}
For a non-regular $4$-vertex $\T V$ at least one of the coefficients in the circular equation {\it "Circular Constraint"} is not equal to zero.
Let us assume that $coeff^{(4)}\neq 0$ and that $\delta^{(4)}$ is defined as a dependent variable according to the {\it "Circular Constraint"}.
Then for any choice of $\delta^{(1)}$, $\delta^{(2)}$, $\delta^{(3)}$ there exists such a value of $\delta^{(4)}$ that the {\it "Circular Constraint"} is satisfied, and the value of $\delta^{(4)}$ is uniquely defined.

In order to prove the second statement of Lemma ~\ref{lemma:C2_sufficient}, one should show that for every choice of $\delta^{(j)}$ ($j=1,\ldots,4$) which satisfies the {\it "Circular Constraint"}, there exist such scalars $Z_{XX}$,
$Z_{XY}$, $Z_{YY}$ that Equation ~\ref{eq:eq_ff_agreement} holds for every $j=1,\ldots,4$.

Indeed, for every choice of $\delta^{(1)}$, $\delta^{(2)}$, $\delta^{(3)}$ there exist such unique values of $Z_{XX}$, $Z_{XY}$, $Z_{YY}$ that Equation ~\ref{eq:eq_ff_agreement} holds for $j=1,2,3$. In order to prove it, one can write the system of equation for $\delta^{(1)}$, $\delta^{(2)}$, $\delta^{(3)}$
in the matrix form
\begin{equation}
\hspace{-0.4in}
\VT{\delta^{(1)}}{\delta^{(2)}}{\delta^{(3)}}=S\VT{Z_{XX}}{Z_{XY}}{Z_{YY}},
{\rm \ \ where\ \ }
S=\left(\begin{array}{ccc}
\left(\alpha^{(1)}\right)^2 & 2\alpha^{(1)}\beta^{(1)} & \left(\beta^{(1)}\right)^2\cr
\left(\alpha^{(2)}\right)^2 & 2\alpha^{(2)}\beta^{(2)} & \left(\beta^{(2)}\right)^2\cr
\left(\alpha^{(3)}\right)^2 & 2\alpha^{(3)}\beta^{(3)} & \left(\beta^{(3)}\right)^2
\end{array}\right)
\end{equation}
It is easy to show that 
$det(S)=2\<\TTe{1},\TTe{2}\>\<\TTe{2},\TTe{3}\>\<\TTe{3},\TTe{1}\>\neq 0$ because $\<\TTe{1},\TTe{2}\>\neq 0$, $\<\TTe{2},\TTe{3}\>\neq 0$ as the vector
products of consequent edges of the strictly convex quadrilaterals and $\<\TTe{3},\TTe{1}\>\neq 0$ according to the assumption
that $coeff^{(4)}\neq 0$. It means that for any choice of $\delta^{(1)}$, $\delta^{(2)}$, $\delta^{(3)}$ values of $Z_{XX}$, $Z_{XY}$, $Z_{YY}$ are uniquely defined by equation
\begin{equation}
\VT{Z_{XX}}{Z_{XY}}{Z_{YY}}= S^{-1}\VT{\delta^{(1)}}{\delta^{(2)}}{\delta^{(3)}}
\end{equation}

It remains to show that $\delta^{(4)}$ also agrees with the second-order partial derivatives $Z_{XX}$, $Z_{XY}$, $Z_{YY}$.
The agreement holds because, there exists the unique value of $\delta^{(4)}$ so that the {\it "Circular Constraint"} is satisfied, and from the first statement of the Lemma it is already known that\\ 
$\delta^{(4)}=\left(\alpha^{(4)}\right)^2 Z_{XX}+
2\alpha^{(4)}\beta^{(4)} Z_{XY}+\left(\beta^{(4)}\right)^2 Z_{YY}$
satisfies the {\it "Circular Constraint"}. 

\nopagebreak 
\eop${}_{{\bf Lemma ~\ref{lemma:C2_sufficient}}}$

\vspace{0.32in}
\noindent
{\bf Proof of Theorem ~\ref{theorem:middle_equations_bilinear}}
The sufficiently long proof of Theorem  ~\ref{theorem:middle_equations_bilinear} starts with the trivial Auxiliary Lemma ~\ref{aux_lemma:algebraic}. The Auxiliary Lemma just summarizes well-known algebraic facts, which will be useful in the proof. 

\begin{auxillary_lemma}
Let $Av=b$ be a linear system of equations, where $A$ is 
$p\times q$ matrix with rows $A_0,\ldots,A_{p-1}$, $v$ is a vector of variables and $b$ is a vector of free coefficients.
\begin{equation}
A=
\left(\begin{array}{ccccc}
- & - & A_0    & - & - \cr
- & - & A_1    & - & - \cr
  &   & \vdots &   &   \cr
- & - & A_{p-1} & - & - \cr
\end{array}\right)\ \ \ \
v = \left(\begin{array}{c}
v_0\cr v_1\cr \vdots\cr v_{q-1}
\end{array}\right)\ \ \ \
b = \left(\begin{array}{c}
b_0\cr b_1\cr \vdots\cr b_{p-1}
\end{array}\right)
\end{equation}
Then
\begin{itemize}
\item[{\bf (1)}]
The system has a solution (is consistent) if and only if for every non-trivial set of coefficients $\alpha_0,\ldots,\alpha_{p-1}$ so that the linear combination of rows of matrix $A$ is equal to zero $\sum_{s=0}^{p-1}\alpha_s A_s=0$, the corresponding linear combination of the free coefficients is also equal to zero 
$\sum_{s=0}^{p-1}\alpha_s b_s$ = 0. Here a "non-trivial" set means that the set contains at least one non-zero element. 
\item[{\bf (2)}]
The rank of the matrix $A$ is equal to $r$ if and only if for any set of coefficients $\alpha_0,\ldots,\alpha_{p-1}$ so that 
$\sum_{s=0}^{p-1}\alpha_s A_s=0$, the following two conditions hold
\begin{itemize}
\item[{\bf (a)}]
For any subset of $l<p-r$ coefficients $\left\{\alpha_{s_j}\right\}_{j=1}^{l}$, equality of every one of these coefficients to zero $\alpha_{s_j}=0$ 
$(j=1,\ldots,l)$ does not necessarily imply that all other coefficients should also be equal to zero.
\item[{\bf (b)}]
There exists such a subset of $l=p-r$ coefficients $\left\{\alpha_{s_j}\right\}_{j=1}^{l}$ that if $\alpha_{s_j}=0$ 
$(j=1,\ldots,l)$ then all other coefficients are necessarily equal to zero.
\end{itemize}
\item[{\bf (3)}] 
Under the assumption that $q\ge p$ the number of independent (free) variables is equal to $q-r$. Which variables are free and 
which are dependent can be defined, for example, by the Gauss elimination (pivoting) process.
\end{itemize}
\label{aux_lemma:algebraic}
\end{auxillary_lemma}

\vspace{0.1in}
\noindent
Now it is possible to proceed with a formal proof of Theorem
~\ref{theorem:middle_equations_bilinear}. The proof is based on representation of the {\it "Middle"} system of equations in the matrix form and uses a few Auxiliary Lemmas, which provide relations between a pure algebraic analysis of the system of equations and the geometrical essence of theorem.

\vspace{0.1in}
\noindent
{\bf (1)} {\bf Consistency.}
It will be shown that for every inner edge, the {\it "Middle"} system of equations is solvable in terms of the {\it middle} control points. Moreover, it will be shown that the {\it "Middle"} system of equation has a solution in terms of the "side" {\it middle} control points $(L_t,R_t)$ for $t=2,\ldots,n-2$ even if all the "central" {\it middle} control points $\T C_t$ for $t=3,\ldots,n-3$ are classified in advance as basic (free). 

For an inner edge, let the  "central" {\it middle} control points $\T C_t$ for $t=3,\ldots,n-3$ be classified as basic(free). Then solvability of the {\it "Middle"} system in terms of control points $(L_t,R_t)$ $(t=2,\ldots,n-2)$ is equivalent to solvability of the system in terms of differences $\Delta L_t$, $\Delta R_t$ ($t=2,\ldots,n-2$), because every one of the differences contains exactly one non-classified control point $L_t$ or $R_t$.  

The proof of the solvability starts with the representation of the full system of linearized $G^1$-continuity equations {\it "Eq(s)"} $(s=0,\ldots,n+1)$ in the matrix form with respect to the differences $\Delta L_t$, $\Delta R_t$ ($t=0,\ldots,n$). Equation {\it "Eq(s)"} $(s=0,\ldots,n+1)$ can be rewritten as follows 
\begin{equation}
\hspace{-0.35in}
(n+1-s)(l_0\Delta L_s+r_0\Delta R_s)+
s(l_1\Delta L_{s-1}+r_1\Delta R_{s-1})+
sumC_s=0
\end{equation}
where
\begin{equation}
\hspace{-0.4in}
\begin{array}{l}
sumC_s = \frac{\scriptsize (n-s)(n+1-s)}{\scriptsize n}c_0\Delta C_s+
\frac{2s(n+1-s)}{n}c_1\Delta C_{s-1}+
\frac{s(s-1)}{n}c_2\Delta C_{s-2}\cr
\end{array}
\label{eq:def_sumC}
\end{equation}
For the full system of equations {\it "Eq(s)"} ($s=0,\ldots,n+1$), coefficients of $\Delta L_t$ and $\Delta R_t$ ($t=0,\ldots,n$) 
form a $(n+2)\times 2(n+1)$ matrix $A$ and expressions $-sumC_s$ ($s=0,\ldots,n+1$) do not contain non-classified control points and are considered as the right-side coefficients $b_s = -sumC_s$ ($s=0,\ldots,n+1$). One gets the following system, which corresponds to all $n+2$ linearized $G^1$-continuity equations along the given edge 
\begin{equation}
A
\left(\begin{array}{c}
\dL_0\cr
\dR_0\cr
\vdots\cr
\dL_n\cr
\dR_n\cr
\end{array}\right)= 
\left(\begin{array}{c}
b_0\cr
\vdots\cr
b_{n+1}\cr
\end{array}\right)
\end{equation}
In the last system, equations for $s=0, s=1, s=n, s=n+1$ are known to be satisfied and $\Delta L_t$, $\Delta R_t$ for $t=0,1,n-1,n$ do not contain any non-classified control points.

The {\it "Middle"} system of equations with respect to the differences $\Delta L_t$, $\Delta R_t$ ($t=2,\ldots,n-2$) 
has the following representation in the matrix form
\begin{equation}
\hat A
\left(\begin{array}{c}
\dL_2\cr
\dR_2\cr
\vdots\cr
\dL_{n-2}\cr
\dR_{n-2}\cr
\end{array}\right)= 
\left(\begin{array}{c}
\hat b_2\cr
\vdots\cr
\hat b_{n-1}\cr
\end{array}\right)
\label{eq:middle_matrix_bilinear}
\end{equation}
where
\begin{itemize}
\item[-]
Matrix $\hat A$ is submatrix of matrix $A$, composed of rows $s=2,\ldots,n-1$ and columns which correspond to $\Delta L_2,\Delta R_2,\ldots,\Delta L_{n-2},\Delta R_{n-2}$.
\item[-]
Right-side coefficients $\hat b_2,\ldots,\hat b{n-1}$ are equal to the right-side coefficients of the full system plus the result of transfer to the right side of expressions which contain the classified control points only
\begin{equation}
\begin{array}{l}
\hat b_2 = b_2-(2l_1\Delta L_1+2r_1\Delta R_1)\cr
\hat b_s = b_s\ \ \ \ \ {\rm for}\ s=3,\ldots,n-2,\cr
\hat b_{n-1}= b_{n-1}-(2l_0\Delta L_{n-1}+2r_0\Delta R_{n-1})
\end{array}
\label{eq:def_hat_b_bilinear}
\end{equation}
\end{itemize}
 
\noindent
For example in the case of \underline{$n=4$}, matrix $A$ has the form
\begin{equation}
\hspace{-0.3in}
\begin{array}{ccccccccccccccc}
& 
\!\Dlt L_0\! & \!\Dlt R_0\! & \!\Dlt L_1\! & \!\Dlt R_1\! &   
\!\Dlt L_2\! & \!\Dlt R_2\! & \!\Dlt L_3\! & \!\Dlt R_3\! &
\!\Dlt L_4\! & \!\Dlt R_4\!\cr
s\!=\!0& 
5l_0& 5r_0&  0  &  0  &  0  &  0  &  0  &  0  &  0  &  0  \cr 
s\!=\!1& 
l_1 & r_1 & 4l_0& 4r_0&  0  &  0  &  0  &  0  &  0  &  0  \cr 
s\!=\!2&
 0  &  0  & 2l_1& 2r_1& 3l_0& 3r_0&  0  &  0  &  0  &  0  \cr 
s\!=\!3&  
 0  &  0  &  0  &  0  & 3l_1& 3r_1& 2l_0& 2r_0&  0  &  0  \cr 
s\!=\!4& 
 0  &  0  &  0  &  0  &  0  &  0  & 4l_1& 4r_1& l_0 & r_0 \cr 
s\!=\!5& 
 0  &  0  &  0  &  0  &  0  &  0  &  0  &  0  & 5l_1& 5r_1\cr 
\end{array}
\end{equation}
and $\hat A$ is the $2\times 2$ middle submatrix of $A$
\begin{equation}
\begin{array}{ccc}
   &\Dlt L_2 & \Dlt R_2 \cr
s=2&   3l_0  &   3r_0   \cr 
s=3&   3l_1  &   3r_1   \cr 
\end{array}
\end{equation}
In the case of \underline{$n=5$}, $A$ is a $7\times 12$ matrix
and $\hat A$ is the $3\times 4$ middle submatrix 
\begin{equation}
\begin{array}{ccccc}
   &\Dlt L_2 & \Dlt R_2 & \Dlt L_3 & \Dlt R_3 \cr
s=2&   4l_0  &   4r_0   &     0    &     0    \cr 
s=3&   3l_1  &   3r_1   &   3l_0   &    3r_0  \cr 
s=4&     0   &     0    &   4l_1   &    4r_1  \cr 
\end{array}
\label{eq:middle_matrix_n_5}
\end{equation}
In a general case matrix $A$ of the full system and submatrix $\hat A$ of the {\it "Middle"} system have the following structure (here horizontal and vertical lines separate elements which belong to submatrix $\hat A$) 

\vspace{0.1in}
\noindent
The left upper corner of the matrix

\hspace{-0.3in}
\begin{tabular}{ccccc|cccc}
     &$\Delta L_0$&$\Delta R_0$&$\Delta L_1$&$\Delta R_1$
     &$\Delta L_2$&$\Delta R_2$&$\Delta L_3$&$\Delta R_3$\cr
$\!\!\!s\!=\!0$&
$\!\!(n\!+\!1)l_0\!\!$&$\!\!(n\!+\!1)r_0\!\!$&
$   0    $&$   0    $&
$  0     $&$   0    $&$   0    $&$   0    $\cr
$\!\!\!s\!=\!1$&
$  l_1   $&$   r_1  $&$   nl_0 $&$   nr_0 $&
$  0     $&$   0    $&$   0    $&$   0    $\cr
\hline
$\!\!\!s\!=\!2$&
$  0     $&$   0    $&$   2l_1 $&$   2r_1 $&
$\!\!(n\!-\!1)l_0\!\!$&$\!\!(n\!-\!1)r_0\!\!$&
$   0    $&$   0    $\cr
$\!\!\!s\!=\!3$&
$  0     $&$   0    $&$   0    $&$   0    $&
$  3l_1  $&$   3r_1 $&
$\!\!(n\!-\!2)l_0\!\!$&$\!\!(n\!-\!2)r_0\!\!\!$
\end{tabular}

\vspace{0.1in}
\noindent
The middle part of the matrix

\vspace{-0.3in}
\begin{eqnarray}
\hspace{-0.43in}
\begin{array}{lcccc}
     &\Delta L_{t-1}&\Delta R_{t-1}&\Delta L_t    &\Delta R_t  \cr
s\!=\!t\!-\!1&
(n+2-t)l_0    &(n+2-t)r_0    &      0       &      0     \cr
s\!=\!t  &      
      tl_1    &      tr_1    &(n+1-t)l_0    &(n+1-t)r_0  \cr
s\!=\!t\!+\!1&
      0       &      0       &(t+1)l_1      &(t+1)r_1    \cr
\end{array}
\nonumber
\end{eqnarray}

\vspace{0.1in}
\noindent
The right lower corner of the matrix

\hspace{-0.3in}
\begin{tabular}{lcccc|cccc}
&$\!\!\Delta L_{n-3}\!\!$&$\!\!\Delta R_{n-3}\!\!$
&$\!\!\Delta L_{n-2}\!\!$&$\!\!\Delta R_{n-2}\!\!$
&$\!\!\Delta L_{n-1}\!\!$&$\!\!\Delta R_{n-1}\!\!$
&$\!\!\Delta L_n\!\!    $&$\!\!\Delta R_n\!\!\!\!$\cr
$\!\!\!s\!=\!n\!-\!2\!$&
$\!\!\!(n\!-\!2)l_1\!\!$&$\!\!\!(n\!-\!2)r_1\!\!$&
$ 3 l_0  $&$ 3 r_0  $&
$  0     $&$   0    $&$   0    $&$   0    $\cr
$\!\!\!s\!=\!n\!-\!1\!$&
$  0     $&$   0    $&
$\!\!\!(n\!-\!1)l_1\!\!$&$\!\!\!(n\!-\!1)r_1\!\!$&
$  2l_0  $&$   2r_0 $&$   0    $&$   0    $\cr
\hline
$\!\!\!s\!=\!n\!$&
$  0     $&$   0    $&$   0    $&$   0    $&
$  nl_1  $&$   nr_1 $&$   l_0  $&$   r_0  $\cr
$\!\!\!s\!=\!n\!+\!1\!\!
$&$  0     $&$   0    $&$   0    $&$   0    $&
$  0     $&$   0    $&
$\!\!\!\!\!(n\!+\!1)l_1\!\!$&$\!\!\!(n\!+\!1)r_1\!\!\!\!$
\end{tabular}

\vspace{0.15in}
The goal is to prove that the system defined by matrix $\hat A$ and the right-side vector $\hat b$ is consistent. Let some non-trivial linear combination of the rows of matrix $\hat A$ be equal to zero
$\sum_{s=2}^{n-1} \alpha_s \hat A_s = 0$
In order to prove the consistency, one should verify that the corresponding linear combination of the right-side coefficients $\sum_{s=2}^{n-1} \alpha_s \hat b_s$ is also equal to zero. Auxiliary Lemma ~\ref{aux_lemma:necessary_and_sufficient_bilinear} provides the necessary and sufficient conditions for equality to zero of a non-trivial linear combination of the rows of matrix $\hat A$.

\begin{auxillary_lemma}
Let matrix $\hat A$ be the matrix of the {\it "Middle"} system of equations with respect to the differences $\Delta L_t$, $\Delta R_t$ ($t=2,\ldots,n-2$) (see Equation ~\ref{eq:middle_matrix_bilinear}). Then a non-trivial linear combination of its rows $\sum_{s=2}^{n-1} \alpha_s \hat A_s$ is equal to zero if and only if the following two conditions hold
\begin{itemize}
\item[{\bf(1)}]
The {\it "Projections Relation"} is satisfied, or in other words, there exists such a constant $\kappa$ that
\begin{equation} 
\kappa = \frac{l_0}{l_1} = \frac{r_0}{r_1}
\label{eq:def_kappa}
\end{equation}
\item[{\bf(2)}]
Coefficients $\alpha_s$ ($s=2,\ldots,n-1$) are given by the formula
\begin{equation}
\alpha_s = \frac{2}{n(n+1)}\Cnk{n+1}{s}(-1)^s\kappa^{s-2}\alpha_2
\label{eq:def_alpha}
\end{equation}
\end{itemize}
\label{aux_lemma:necessary_and_sufficient_bilinear}
\end{auxillary_lemma}
{\bf Proof} See Appendix, after proof of Theorem ~\ref{theorem:middle_equations_bilinear}.

\vspace{0.15in}
\noindent
Now it remains to prove that if conditions {\bf (1)} and {\bf (2)} of Auxiliary Lemma ~\ref{aux_lemma:necessary_and_sufficient_bilinear} hold then the linear combination of the right-side coefficients 
$\sum_{s=2}^{n-1}\alpha_s \hat b_s$ is equal to zero.
Auxiliary Lemma ~\ref{aux_lemma:right_part_bilinear} presents an important relation between the linear combination, geometrical coefficient of proportionality $\kappa$ defined by the {\it "Projections Relation"} and coefficients of the weight function $c(v)$.

\begin{auxillary_lemma}
For an inner edge, let {\it "Eq(0)"}-type and {\it "Eq(1)"}-type equations be satisfied, matrix $\hat A$ corresponds to the {\it "Middle"} system of equations and let a non-trivial linear combination of the rows of matrix $\hat A$ be equal to zero $\sum_{s=2}^{n-1} \alpha_s \hat A_s=0$. Then the corresponding  linear combination of the right-side coefficients of the {\it "Middle"} system has the following representation in terms of constant $\kappa$ defined by {\it "Projections Relation"} (see Equation ~\ref{eq:def_kappa}) and coefficients of the weight function $c(v)$
\begin{equation}
\sum\limits_{s=2}^{n-1}\alpha_s \hat b_s = 
-\frac{2}{n\kappa^2}\alpha_2
\left\{
\sum\limits_{s=0}^{n-1} (-1)^s\Cnk{n-1}{s}\kappa^s\Dlt C_s
\right\}
\left\{c_0-2\kappa c_1+\kappa^2 c_2\right\}
\label{eq:sumC_kappa_c_relation_bilinear}
\end{equation}
\label{aux_lemma:right_part_bilinear}
\end{auxillary_lemma}
{\bf Proof} See Appendix, after proof of Theorem ~\ref{theorem:middle_equations_bilinear}.

\vspace{0.15in}
\noindent
Auxiliary Lemma ~\ref{lemma:c_kappa_bilinear} contains the geometrical essence of the proof. It derives the relation between coefficients of the weight function $c(v)$ and constant $\kappa$ in the case when the {\it "Projections Relation"} is satisfied.

\begin{auxillary_lemma}
For two adjacent mesh elements with bilinear in-plane parametrisation, let the {\it "Projections Relation"} hold and constant $\kappa$ be defined by Equation ~\ref{eq:def_kappa}.
Then the coefficients of the weight function $c(v)$ satisfy the following equation
\begin{equation}
c_0-2\kappa c_1+\kappa^2 c_2 = 0 
\label{eq:c_kappa_bilinear}
\end{equation}
\label{lemma:c_kappa_bilinear}
\end{auxillary_lemma}
{\bf Proof} See Appendix, after proof of Theorem ~\ref{theorem:middle_equations_bilinear}.

\vspace{0.15in}
\noindent
Auxiliary Lemma ~\ref{lemma:c_kappa_bilinear} clearly allows to complete the proof of the consistency of the {\it "Middle"} system.
Indeed, let a non-trivial linear combination  $\sum_{s=2}^{n-1}\alpha_s \hat A_s$ be equal to zero. Then according to Auxiliary Lemma ~\ref{aux_lemma:necessary_and_sufficient_bilinear} the {\it "Projections Relation"} holds and $\sum_{s=2}^{n-1}\alpha_s \hat b_s$ is equal to zero according to Auxiliary Lemmas 
~\ref{aux_lemma:right_part_bilinear} and 
~\ref{lemma:c_kappa_bilinear}.

\vspace{0.3in}
\noindent
{\bf (2)} {\bf Classification of the control points.}
In the previous part of the proof it was already shown that the {\it "Middle"} system of equations is always solvable in terms of the "side" {\it middle} control points $L_2,\ldots,L_{n-2}$, $R_2,\ldots,R_{n-2}$ and the "central" {\it middle} control points $\T C_3,\ldots,\T C_{n-3}$ can be classified as basic(free) in advance. In addition, Auxiliary Lemma 
~\ref{aux_lemma:necessary_and_sufficient_bilinear} shows that a non-trivial linear combination of rows of matrix $\hat A$ (which corresponds to the {\it "Middle"} system) is equal to zero if and only if the {\it "Projections Relation"} holds and the coefficients of the linear combination are defined by Equation ~\ref{eq:def_alpha}. 

\vspace{0.1in}
It implies that if the \underline{\it "Projections Relation" does not hold}, no non-trivial linear combination of the rows of matrix $\hat A$ is equal to zero and the matrix is of the full row rank $r=n-2$. According to statement {\bf (3)} of Auxiliary Lemma ~\ref{aux_lemma:algebraic}, the number of the basic (free) "side" {\it middle} control points in this case is equal to $2(n-3)-(n-2)=n-4$ and dependencies of the other "side" {\it middle} control points can be defined by the Gauss elimination process applied to the matrix $\hat A$.

\vspace{0.1in}
If the \underline{\it "Projections Relation" holds}, then rank $r$ of the matrix $\hat A$ is equal to $n-3$. Indeed, $r\le n-3$ because there exists a non-trivial zero combination of the rows of $\hat A$ (see Equation ~\ref{eq:def_alpha}). It remains to show that $r\ge n-3$. On the contrary, let $r<n-3$. Then according to statement {\bf (2a)} of Auxiliary Lemma ~\ref{aux_lemma:algebraic}, equality to zero of  $l=1$ coefficient in the zero linear combination of rows 
$\sum_{s=2}^{n-1}\alpha_s \hat A_s=0$
should not necessarily imply that all other coefficients are equal to zero. However, it contradicts the recursive dependency (see Equation ~\ref{eq:def_alpha_recursive}) between coefficients of a non-trivial zero linear combination of the rows. 

The number of basic (free) "side" {\it middle} control points in this case is equal to $n-3$ and from the Gauss elimination process it follows that there is exactly one basic control point in every pair $(\T L_j, \T R_j)$ for $j=2,\ldots,n-2$. Below the Gauss elimination process is illustrated in the case of $n=5$; in case of a general $n$ the resulting matrix has precisely the same structure. For $n=5$, application of two steps of the Gauss elimination process to matrix $\hat A$ of the {\it "Middle"} system of equations (see Equation ~\ref{eq:middle_matrix_n_5}) gives \begin{equation}
\hspace{-0.4in}
\begin{array}{ll}
{\begin{array}{l}
\hat A_2\cr \hat A_3\cr \hat A_4\cr
\end{array}}
&
{\left(
\begin{array}{cccc}
4 l_0 & 4 r_0 &   0   &   0  \cr
3 l_1 & 3 r_1 & 3 l_0 & 3 r_0\cr
  0   &   0   & 4 l_1 & 4 r_1\cr
\end{array}
\right)
\stackrel
{\hat A_2\leftarrow\hat A_2-\frac{3}{4\kappa}\hat A_1}
{\longrightarrow}}
\cr
&
{\left(
\begin{array}{cccc}
4 l_0 & 4 r_0 &   0   &   0  \cr
  0   &   0   & 3 l_0 & 3 r_0\cr
  0   &   0   & 4 l_1 & 4 r_1\cr
\end{array}
\right)
\stackrel
{\hat A_3\leftarrow\hat A_3-\frac{4}{3\kappa}\hat A_2}
{\longrightarrow}
\left(
\begin{array}{cccc}
4 l_0 & 4 r_0 &   0   &   0  \cr
  0   &   0   & 3 l_0 & 3 r_0\cr
  0   &   0   &   0   &   0  \cr
\end{array}
\right)}
\end{array}
\end{equation}
The resulting matrix clearly shows that every second column contains a pivot and so in every pair $(\T L_j,\T R_j)$ $(j=2,3)$ exactly one control point can be classified as a basic control point. In addition, for every pair $(\T L_j,\T R_j)$ $(j=2,3)$, the resulting matrix defines dependency of the dependent control point on the basic one. 

\nopagebreak 
\eop${}_{{\bf Theorem ~\ref{theorem:middle_equations_bilinear}}}$

\vspace{0.32in}
\noindent
{\bf Proof of Auxiliary Lemma ~\ref{aux_lemma:necessary_and_sufficient_bilinear}}
Conditions {\bf (1)} and {\bf (2)} are clearly sufficient in order to guarantee that the linear combination $\sum_{s=2}^{n-1} \alpha_s \hat A_s$ is equal to zero.

It remains to show that the conditions are necessary. Indeed,
equality $\sum_{s=2}^{n-1} \alpha_s \hat A_s = 0$ implies that for columns which correspond to $\Delta L_t$ and $\Delta R_t$ ($t=2,\ldots,n-2$), one gets 
\begin{equation}
\begin{array}{l}
\alpha_t(n+1-t)l_0+\alpha_{t+1}(t+1)l_1=0\cr
\alpha_t(n+1-t)r_0+\alpha_{t+1}(t+1)r_1=0\cr
\end{array}
\label{eq:two_columns_bilinear}
\end{equation}
The existence of the constant $\kappa$ (condition {\bf (1)}) immediately follows from the last pair of equations and the assumption that the linear combination is non-trivial. 
In addition, Equation ~\ref{eq:two_columns_bilinear} defines the recursive dependency between coefficients of the linear combination 
\begin{equation}
\alpha_{t+1}=-\frac{n+1-t}{t+1}\kappa\alpha_t, \ \ \ 
t=2,\ldots,n-2 
\label{eq:def_alpha_recursive}
\end{equation}
A trivial inductive proof leads to the explicit formula for $\alpha_s$ (Equation ~\ref{eq:def_alpha}), which means that condition {\bf (2)} is satisfied.

\nopagebreak
\eop${}_{{\bf Auxiliary Lemma ~\ref{aux_lemma:necessary_and_sufficient_bilinear}}}$

\vspace{0.32in}
\noindent
{\bf Proof of Auxiliary Lemma ~\ref{aux_lemma:right_part_bilinear}}
The proof consists of two parts. The first part shows that the linear combination of the right-side coefficients of the {\it "Middle"} system of equations $\sum_{s=2}^{n-1}\alpha_s \hat b_s$ is equal to a similar linear combination of the right-side coefficients of the full system of equations $\sum_{s=0}^{n+1}\alpha_s b_s$. The second part expresses the linear combination $\sum_{s=0}^{n+1}\alpha_s b_s$ in terms of $\kappa$ and coefficients of the weight function $c(v)$.

Auxiliary Lemma ~\ref{aux_lemma:necessary_and_sufficient_bilinear} implies, that the {\it "Projections Relation"} holds and constant $\kappa$ is correctly defined. According to the definition of $\hat b_2$ (see Equation ~\ref{eq:def_hat_b_bilinear}), definition of $\kappa$  and  the assumption  that equations {\it "Eq(0)"} and {\it "Eq(1)"} are satisfied, one gets 
\begin{equation}
\hspace{-0.3in}
\begin{array}{ll}
\hat b_2 = 
&b_2-2(l_1\Delta L_1+r_1\Delta R_1)
\stackrel{(definition\ of\ \kappa)}{=}\cr
&b_2-2\kappa^{-1}(l_0\Delta L_1+r_0\Delta R_1)
\stackrel{\it "Eq(1)"}{=}\cr
&b_2-\frac{2}{n}\kappa^{-1}
(b_1-\frac{1}{n+1}(l_1\Delta L_0+r_1\Delta R_0))
\stackrel{(definition\ of\ \kappa)}{=}\cr
&b_2-\frac{2}{n}\kappa^{-1}
(b_1-\frac{1}{n+1}\kappa^{-1}(l_0\Delta L_0+r_0\Delta R_0))
\stackrel{\it "Eq(0)"}{=}\cr
&b_2-\frac{2}{n}\kappa^{-1}b_1+\frac{2}{n(n+1)}\kappa^{-2}b_0
\end{array}
\label{eq:hat_b_2}
\end{equation}
In the same manner, it can be shown that
\begin{equation}
\hspace{-0.2in}
\hat b_{n-1}=b_{n-1}-\frac{2}{n}\kappa b_n+
\frac{2}{n(n+1)}\kappa^2 b_{n+1}
\label{eq:hat_b_n_1}
\end{equation}
Let $\alpha_s$ for $s=0,1,n,n+1$ be formally defined according to the recursive relation ~\ref{eq:def_alpha_recursive}
\begin{equation}
\begin{array}{ll}
\alpha_0 = \frac{2}{n(n+1)}\kappa^{-2}\alpha_2, &
\alpha_1 = -\frac{2}{n}\kappa^{-1}\alpha_2\cr
\alpha_n = -\frac{2}{n}\kappa\alpha_{n-1},  &
\alpha_{n+1} = \frac{2}{n(n+1)}\kappa^2\alpha_{n-1}
\end{array}
\label{eq:def_alpha_first_last}
\end{equation}
Then Equations ~\ref{eq:hat_b_2}, ~\ref{eq:hat_b_n_1} allow to conclude that 
\begin{equation}
\sum_{s=2}^{n-1}\alpha_s \hat b_s = \sum_{s=0}^{n+1}\alpha_s b_s
\label{eq:hat_b_s_b_s_bilinear}
\end{equation} 

Now it remains to substitute the explicit formulas for $\alpha_s$ (see Auxiliary Lemma ~\ref{aux_lemma:necessary_and_sufficient_bilinear}, Equation ~\ref{eq:def_alpha}) and $b_s=-sumC_s$ (see Equation ~\ref{eq:def_sumC}) into the right side of Equation ~\ref{eq:hat_b_s_b_s_bilinear}. A straightforward computation of $\sum_{s=0}^{n+1}\alpha_s b_s$ gives
{\small
\begin{equation}
\hspace{-0.4in}
\begin{array}{l}
\sum\limits_{s=0}^{n+1}\alpha_s b_s =  
-\sum\limits_{s=0}^{n+1}\alpha_s sumC_s =\cr
-\frac{1}{n}\sum\limits_{s=0}^{n+1}\alpha_s
\left\{
(n\!-\!s)(n\!+\!1\!-\!s)c_0\Delta C_s
\!+\!2s(n\!+\!1\!-\!s)c_1\Delta C_{s-1}\!+\!
s(s\!-\!1)c_2\Delta C_{s-2}
\right\}\!=\cr
-\frac{1}{n}\sum\limits_{s=0}^{n-1}\Delta C_s
\left\{
\alpha_s c_0(n\!-\!s)(n\!+\!1\!-\!s)\!+\!
2\alpha_{s+1}c_1(s\!+\!1)(n\!-\!s)\!+\!
\alpha_{s+2}c_2(s\!+\!1)(s\!+\!2)
\right\}\!=\!\!\!\!\!\!\!\cr
-\frac{2}{n}{\alpha_2\kappa^2}
\left\{
\sum\limits_{s=0}^{n-1} 
(-1)^s\Cnk{n-1}{s}\kappa^s\Dlt C_s
\right\}
\left\{c_0-2\kappa c_1+\kappa^2 c_2\right\}
\end{array}
\label{eq:sum_computation_bilinear}
\end{equation}}
\noindent
Equation ~\ref{eq:sum_computation_bilinear} completes the proof of Auxiliary Lemma ~\ref{aux_lemma:right_part_bilinear}.

\nopagebreak
\eop${}_{{\bf Auxiliary Lemma ~\ref{aux_lemma:right_part_bilinear}}}$

\vspace{0.32in}
\noindent
{\bf Proof of Auxiliary Lemma ~\ref{lemma:c_kappa_bilinear}}

Let $\varphi,\varphi',\psi,\psi'$ be angles between the common edge of two patches and left and right adjacent edges (see Figure ~\ref{fig:fig27}).
In the case of the bilinear in-plane parametrisation, coefficients of the weight functions $l(v)$ and $r(v)$ have the following representations in terms of the geometrical characteristics of  two adjacent mesh elements (see Equation ~\ref{eq:weights_bilinear}) 
\begin{equation}
\begin{array}{l}
l_0=\<\RR-\GG,\GG'-\GG\>=||\GG'-\GG||\ ||\RR-\GG||sin\psi\cr
l_1=\<\GG-\GG',\RR'-\GG'\>=||\GG'-\GG||\ ||\RR'-\GG'||sin\psi'\cr
r_0=-\<\GG'-\GG,\LL-\GG\>=-||\GG'-\GG||\ ||\LL-\GG||sin\varphi\cr
r_1=-\<\LL'-\GG',\GG-\GG'\>=-||\GG'-\GG||\ ||\LL'-\GG'||sin\varphi'\cr
\end{array}
\end{equation}
Then the {\it "Projections Relation"} implies that
\begin{equation}
\kappa = 
\frac{||\RR-\GG||sin\psi}{||\RR'-\GG'||sin\psi'}=
\frac{||\LL-\GG||sin\varphi}{||\LL'-\GG'||sin\varphi'}
\label{eq:kappa_sin}
\end{equation}

For coefficients of the weight function $c(v)$ the following relations hold
\begin{equation}
\hspace{-0.35in}
\begin{array}{ll}
c_0\!=\!\!&\<\RR-\GG,\LL-\GG\>=||\RR-\GG||\ ||\LL-\GG|| sin(\varphi+\psi)\cr
c_2\!=\!\!&-\<\LL'-\GG',\RR'-\GG'\>=-||\RR'-\GG'||\ ||\LL'-\GG'||
sin(\varphi'+\psi')\cr
c_1\!=\!\!&\frac{1}{2}
\left(\<\RR-\GG,\LL'-\GG'\>-\<\LL-\GG,\RR'-\GG'\>\right)=\cr
&\frac{1}{2}\left(
||\RR-\GG||\ ||\LL'-\GG'|| sin(\varphi'-\psi)-
||\RR'-\GG'||\ ||\LL-\GG|| sin(\varphi-\psi')\right)\!\!\!
\end{array}
\end{equation}
Using the standard trigonometric formulas and Equation ~\ref{eq:kappa_sin}, it is possible to rewrite $2\kappa c_1$ in the form which clearly shows that Auxiliary Lemma ~\ref{lemma:c_kappa_bilinear} is satisfied
\begin{equation}
\begin{array}{ll}
2\kappa c_1=
&\kappa||\RR-\GG||\ ||\LL'-\GG'|| 
(sin\varphi'\ cos\psi-sin\psi\ cos\varphi')-\cr
&\kappa||\RR'-\GG'||\ ||\LL-\GG||
(sin\varphi\ cos\psi'-sin\psi'\ cos\varphi)=\cr
&||\RR-\GG||\ ||\LL-\GG||(sin\varphi\ cos\psi+sin\psi\ cos\varphi)-\cr
&\kappa^2
||\RR'-\GG'||\ ||\LL'-\GG'||
(sin\psi'\ cos\varphi'+sin\varphi'\ cos\psi')=\cr
&c_0+\kappa^2 c_2 
\end{array}
\end{equation}

\nopagebreak 
\eop${}_{{\bf Auxiliary\ Lemma ~\ref{lemma:c_kappa_bilinear}}}$

\vspace{0.32in}
\noindent
{\bf Proof of Lemma ~\ref{lemma:dependency_forest_properties}}

\vspace{0.1in}
\noindent
{\bf (1)} 
It should be shown that the direction of every mesh edge is defined mostly one connected component of the $D$-dependency graph (by mostly one $D$-dependency tree). A mesh edge $\T e$ is directed if it is an edge of the spanning tree of some connected component $\T{\cal C}$ 
of the $D$-dependency graph or if it corresponds to a dangling half-edge of the root vertex of $\T{\cal C}$. In the first case  the direction  of $\T e$ is obviously uniquely defined, because $\T e$ connects two vertices of the same component and can not be influenced by any other component. In the second case $\T e$ clearly could not be an edge of the spanning tree of some other component $\T{\cal C'}$
(otherwise it would connect two vertices of $\T{\cal C'}$). The only situation when component $\T{\cal C'}$ may affect the orientation of $\T e$ is the situation when $\T e$ also corresponds to a dangling half-edge of the root vertex of $\T{\cal C'}$. But it is impossible because in this case roots of $\T{\cal C}$ and $\T{\cal C'}$ should 
belong to the same connected component of the $D$-dependency graph due to the connection by $\T e$.

\vspace{0.1in}
\noindent
{\bf (2)} 
On the contrary, let some vertex $\T V$ use the $D$-type control point of directed edge $\T e = (\T V,\T V')$ and let $\T V$ belong to $D$-dependency tree $\cal T$ while $\T e$ belongs to  $D$-dependency tree $\cal T'$. The only situation when directed edge $\T e$ does not belong to the same $D$-dependency tree as its vertex $\T V$ is the situation when $\T V'$ is the root vertex of $\cal T'$ and $\T e$ corresponds to a dangling half-edge of $\T V'$. Then according to the definition of the dangling half-edge it means that $\T V'$ does not use the $D$-type control point of $\T e$, which leads to a contradiction with the assumption.

\vspace{0.1in}
\noindent
{\bf (3)}
The third statement of Lemma ~\ref{lemma:dependency_forest_properties} is an obvious result of Algorithms ~\ref{algorithm:dependency_graph_construction} and
~\ref{algorithm:dependency_forest_construction} for the construction of the $D$-dependency graph and of the $D$-dependency forest.

\nopagebreak
\eop${}_{{\bf Lemma ~\ref{lemma:dependency_forest_properties}}}$

\vspace{0.32in}
\noindent
{\bf Proof of Theorem ~\ref{theorem:sufficient_DT_classification_bilinear}}

Let $\T{\cal C}$ be a connected component of the
$D$-dependency graph. The only type of {\it primary} vertices which may belong to the component is the inner even vertices, excluding regular $4$-regular vertices, because no other vertex uses any $D$-type control points according to the assumption of theorem. 
Let $\T R$ be the rightmost (geometrically) {\it primary} vertex which belongs to $\T{\cal C}$ or the lowest rightmost vertex in case of ambiguity. It will be shown that $\T R$ has at least one dangling half-edge.

Let $\T r$ be a vertical line which passes through $\T R$. Then in the initial mesh there exist edges $(\T R,\T V^{(1)}),\ldots,(\T R,\T V^{(k)})$ $(k\ge 1)$ which lie strictly on the right from $\T r$ or at the lower half of $\T r$ (see Figure ~\ref{fig:fig24a}).
Indeed, $\T R$ is an inner vertex and so it is surrounded by domain $\T\Omega$. If all mesh edges adjacent to $\T R$ lie on the left of $\T r$ or on the upper half of $\T r$, there should exist a quadrilateral mesh element with inner angle $\ge\pi$, which  contradicts the strict convexity of mesh elements (see Figure ~\ref{fig:fig24b}).


Edges $(\T R,\T V^{(1)}),\ldots,(\T R,\T V^{(k)})$ do not belong to the $D$-dependency graph, otherwise vertices $\T V^{(1)},\ldots,\T V^{(k)}$ would belong to $\T{\cal C}$, which contradicts the fact that $\T R$ is the lowest rightmost primary vertex of $\T{\cal C}$. It means that for every one of edges $(\T R,\T V^{(j)})$ ($j=1,\ldots,k$) at least one of its vertices does not use the $D$-type control point $\T D^{(j)}$ corresponding to the edge. The purpose is to show, that $\T R$ uses the $D$-type control point of at least one of these edges; in this case $\T R$ and the secondary vertex of this edge form a dangling half-edge.

On the contrary, let $\T R$ do not use any $\T D^{(j)}$ ($j=1,\ldots,k$).

Let \underline{$k=1$} and $\T R$ does not use $\T D^{(1)}$. Then in the initial mesh $(\T R, \T W^{(1)})$ and $(\T R, \T W^{(2)})$ - two neighboring edges of $(\T R,\T V^{(1)})$ - should be colinear (see Figure ~\ref{fig:fig24c}). In this case one of the vertices $W^{(1)}$, $W^{(2)}$ lies strictly on the right from $\T r$ or on the lower half of $\T r$, which contradicts the assumption that $k=1$.

Let \underline{$k\ge 3$}. Then angle $\angle\T V^{(1)} \T R \T V^{(k)}<\pi$ (see Figure ~\ref{fig:fig24d}) and for any $j=2,...,k-1$ the $D$-type control point of $(\T R,\T V^{(j)})$ should participate in the {\it "Circular Constraint"} of the vertex $\T R$ with a non-zero coefficient. It means that $\T R$ uses $\T D^{(j)}$ for $j=2,\ldots,k-1$.

The last case is \underline{$k=2$}.
If $\T R$ does not use $D$-type control points of both $(\T R,\T V^{((1)})$ and $(\T R,\T V^{(2)})$, then in the initial mesh there should exist two such pairs of colinear edges that the edges
form a continuous sequence with respect to the order of edges around $\T R$. For example in Figure ~\ref{fig:fig24e},
$\T V^{(1)},\T R,\T W^{(1)}$ and $\T V^{(2)},\T R,\T W^{(2)}$ are triples of the colinear points, and edges $(\T R,\T W^{(2)})$, $(\T R,\T V^{(1)})$, $(\T R,\T V^{(2)})$, $(\T R,\T W^{(1)})$ follow one another in the counter-clockwise order around $\T R$.
If $deg(R) = 4$, then $\T R$ should be a regular $4$-vertex, which contradicts the fact that $\T R$ belongs to the $D$-dependency graph. If $deg(R) \ge 6$ then the arrangement of edges adjacent to $\T R$ contradicts the {\it "Uniform Edge Distribution Condition"} which is assumed to be satisfied.

\nopagebreak 
\eop${}_{{\bf Theorem ~\ref{theorem:sufficient_DT_classification_bilinear}}}$

\vspace{0.32in}
\noindent
{\bf Proof of Theorem ~\ref{theorem:mds_dimension_bilinear}}

Let $|{\cal B}_{V,E\!-\!type}^{(n)}|$,
$|{\cal B}_{D,T\!-\!type}^{(n)}|$,
$|{\cal B}_{middle}^{(n)}|$ denote respectively the number of basic control points of $V$,$E$-type, $D$,$T$-type and of the {\it middle} control points which participate in $G^1$-continuity conditions.
Then
\begin{equation}
|\GMDS{n}| = 
|{\cal B}_{V,E\!-\!type}^{(n)}|+
|{\cal B}_{D,T\!-\!type}^{(n)}|+
|{\cal B}_{middle}^{(n)}|
\end{equation}
Clearly, for any $n$
\begin{equation}
\begin{array}{l}
|{\cal B}_{V,E\!-\!type}^{(n)}|=3|Vert_{non\!-\!corner}|
\end{array}
\label{eq:EV_basic_bilinear}
\end{equation}
In addition, the following useful equality always holds
\begin{equation}
\begin{array}{l}
\sum\limits_{\T V\ non-corner}\!\!\!\!\!\!\!\!\!\! val(V)=
2\left[
|Edge_{inner}|+|Vert_{\twolines{boundary}{non\!-\!corner}}|
\right]
\end{array}
\label{eq:sum_V_expression}
\end{equation}

Let \underline{$n\ge 5$}. Then local templates for classification of $D$,$T$-type control points never intersect and the number of basic $D$,$T$-type control points adjacent to a non-corner vertex $\T V$ is equal to $val(V)+1$ for inner regular $4$-vertices, to $2=val(V)-1$ for boundary vertices and to $val(V)$ for the rest of vertices. Therefore
\begin{equation}
\begin{array}{l}
|{\cal B}_{D,T\!-\!type}^{(n)}| = \!\!
\sum\limits_{\T V\ non-corner}\!\!\!\!\!\!\!\!\!\! val(V)+
|Vert_{\twolines{inner}{4\!-\!regular}}|-|Vert_{\twolines{boundary}{non\!-\!corner}}|
\end{array}
\label{eq:DT_5_basic_bilinear}
\end{equation} 
According to Theorem ~\ref{theorem:middle_equations_bilinear},
\begin{equation}
\begin{array}{l}
|{\cal B}_{middle}^{(n)}| =
(2n-9)|Edge_{inner}|+|Edge_{\twolines{inner,}
                     {"Projections\ Relation"\ holds}}|
\end{array}
\label{eq:middle_5_basic_bilinear}
\end{equation}
Equations ~\ref{eq:EV_basic_bilinear}, ~\ref{eq:DT_5_basic_bilinear}, ~\ref{eq:middle_5_basic_bilinear} and ~\ref{eq:sum_V_expression} allow to conclude that the formula for 
$|\GMDS{n}|$, given in Theorem ~\ref{theorem:mds_dimension_bilinear}, is correct for $n\ge 5$.

Let now \underline{$n=4$}. Then, according to Theorem ~\ref{theorem:middle_equations_bilinear},
\begin{equation}
\begin{array}{l}
|{\cal B}_{middle}^{(4)}| = |Edge_{\twolines{inner,}
                            {"Projections\ Relation"\ holds}}|
\end{array}
\label{eq:middle_4_basic_bilinear}
\end{equation}
The number of basic $D$-type and $T$-type control points can be computed as follows. The total number of $D$-type control points which belong to $\GMDS{4}$ is equal to $|Edge_{inner}|$. One dependent $D$-type control point is {\it assigned} to every $D$-relevant vertex. Therefore
\begin{equation}
\begin{array}{l}
|{\cal B}_{D-type}^{(4)}| = 
|Edge_{inner}|-
|Vert_{\twolines{inner,even,}{not\ 4\!-\!regular}}|-
|Vert_{\twolines{boundary}{D\!-\!relevant}}|
\end{array}
\label{eq:D_4_basic_bilinear}
\end{equation}
The number of basic $T$-type control points is given by
\begin{equation}
\begin{array}{ll}
|{\cal B}_{T-type}^{(4)}| = &
|Vert_{\twolines{inner,even,}{not\ 4\!-\!regular}}|+
|Vert_{\twolines{inner,even,}{4\!-\!regular}}|+\cr
&|Vert_{\twolines{boundary,non-corner,}{not\ D\!-\!relevant}}|+
2|Vert_{\twolines{boundary}{D\!-\!relevant}}|
\end{array}
\label{eq:T_4_basic_bilinear}
\end{equation}
Therefore the number of $D$-type and $T$-type basic control point together do not depend on the presence of the boundary $D$-relevant vertices and can be computed by the formula
\begin{equation}
\begin{array}{l}
|{\cal B}_{D,T-type}^{(4)}| = 
|Edge_{inner}|-
|Vert_{\twolines{inner}{4\!-\!regular}}|-
|Vert_{\twolines{boundary}{non\!-\!corner}}|
\end{array}
\label{eq:DT_4_basic_bilinear}
\end{equation}
From Equations ~\ref{eq:EV_basic_bilinear}, ~\ref{eq:DT_5_basic_bilinear}, ~\ref{eq:middle_5_basic_bilinear} it follows that 
\begin{equation}
\begin{array}{ll}
|\GMDS{4}| = &
3|Vert_{non-corner}|+
|Vert_{\twolines{boundary}{non\!-\!corner}}|+
|Edge_{inner}|+\cr
&|Vert_{\twolines{inner}{4\!-\!regular}}|+
|Edge_{\twolines{inner,}{"Projections\ Relation"\ holds}}|
\end{array}
\label{eq:dimension_gmds_4_bilinear}
\end{equation}
The last expression fits the general formula (Equation ~\ref{eq:dimension_gmds_bilinear}) for $n=4$.

\nopagebreak 
\eop${}_{{\bf Theorem ~\ref{theorem:mds_dimension_bilinear}}}$

\vspace{0.32in}
\noindent
{\bf Proof of Lemma ~\ref{lemma:c_lr_degree_relation}}

Lemma ~\ref{lemma:c_lr_degree_relation} has a very simple formal proof. 

If $max\_deg(l,r)=2$ then the statement of the Lemma is evidently satisfied because $deg(c)\le 4$.

If $max\_deg(l,r)=1$ then $l^{(power)}_2=\<\T\rho^{(power)}_2,\GG'-\GG\>=0$ and
$r^{(power)}_2=-\<\T\lambda^{(power)}_2,\GG'-\GG\>=0$. It means that vectors
$\T\lambda^{(power)}_2$, $\T\rho^{(power)}_2$, $\GG'-\GG$ are parallel and so 
$c^{(power)}_4=\<\lambda^{(power)}_2,\rho^{(power)}_2\>=0$.

If $max\_deg(l,r)=0$ then in the same manner as in the previous case, one sees that vectors $\T\lambda^{(power)}_1$, $\T\rho^{(power)}_1$, $\T\lambda^{(power)}_2$, $\T\rho^{(power)}_2$, $\GG'-\GG$ are parallel and so $ c^{(power)}_4=c^{(power)}_3=0$.

\nopagebreak 
\eop${}_{{\bf Lemma ~\ref{lemma:c_lr_degree_relation}}}$

\vspace{0.32in}
\noindent
{\bf Proof of Lemma ~\ref{lemma:Eq_sumC_linear_system}}

The system of $n+4$ linear equations given in Lemma ~\ref{lemma:G1_cubic_formal} has the following representation in terms of $sumLR_s$ and $sumC_s$ (see Definition ~\ref{def:sumLR_sumC})
\begin{equation}
sumLR_s+sumLR_{s-1}+sumC_s=0
\end{equation}
where $s=0,\ldots,n+3$.

Expressing recursively $sumLR_{s-1}$ from the previous equations for $s=0,\ldots,n+2$ and leaving equation for $s=n+3$ unchanged, one gets an equivalent system 
\begin{equation}
\begin{array}{ll}
(s=0) & sumLR_0 = -sumC_0 \cr
(s=1) & sumLR_1 = -sumLR_0-sumC_1 = sumC_0-sumC_1 \cr
\vdots & \cr
(s=n+2) & sumLR_{n+2} = (-1)^{n+3}\sum_{k=0}^{n+2}(-1)sumC_k \cr
(s=n+3) & sumLR_{n+2} = -sumC_{n+3} \cr
\end{array}
\end{equation}
For $s=0,\ldots,n+2$, the equation has precisely the form of the indexed equation {\it Eq(s)} and it remains to note that the last couple of equations is equivalent to the couple {\it "Eq(n+2)"} and {\it "sumC-equation"}.

\nopagebreak 
\eop${}_{{\bf Lemma ~\ref{lemma:Eq_sumC_linear_system}}}$

\vspace{0.32in}
\noindent
{\bf Proof of Lemma ~\ref{lemma:Eq0_Eq1_bicubic}}

According to Lemma ~\ref{lemma:G1_edge_with_two_inner_vertices_bicubic}, statements of Lemma ~\ref{lemma:Eq0_Eq1_bicubic} clearly hold for an inner vertex which has no adjacent edges with one boundary vertex. Moreover, {\it "Eq(0)"}-type and {\it "Eq(1)"}-type equations remain unchanged for any edge with two inner vertices

It remains to consider equations {"Eq(0)"} and {"Eq(1)"} of an edge with one boundary vertex, applied at the inner vertex of the edge.


In order to distinguish between the case of $\T\Pi^{(bilinear)}$ and the case of $\T\Pi^{(bicubic)}$ global in-plane parametrisations, superscripts ${}^{(bilinear)}$ and ${}^{(bicubic)}$ are added to coefficients of the weight functions and to the indexed equations. It is easy to verify (see Equations ~\ref{eq:weights_bilinear},
~\ref{eq:lr_coeff_bicubic}, ~\ref{eq:c_coeff_bicubic}) that
\begin{equation}
\begin{array}{ll}
l_0^{(bicubic)}=l_0^{bilinear} & l_1^{(bicubic)}=\frac{1}{2}(l_0^{(bilinear)}+l_1^{(bilinear)})\cr
r_0^{(bicubic)}=r_0^{(bilinear)} & r_1^{(bicubic)}=\frac{1}{2}(r_0^{(bilinear)}+r_1^{(bilinear)})\cr
c_0^{(bicubic)}=l_0^{(bilinear)} & 
c_1^{(bicubic)}=\frac{1}{2}(c_0^{(bilinear)}+c_1^{(bilinear)})
\end{array}
\label{eq:bilinear_bicubic_relation}
\end{equation}

\vspace{0.1in}
\noindent
{\bf (1)}
$"Eq(0)^{(bicubic)}"$ involves only zero-indexed coefficients of the weight functions, which are equal for global in-plane parametrisations $\T\Pi^{(bilinear)}$ and $\T\Pi^{(bicubic)}$ .

\vspace{0.1in}
\noindent
{\bf (2)}
From Equation ~\ref{eq:bilinear_bicubic_relation} and formulas 
for {\it "Eq(0)"} and {\it "Eq(1)"} in the case of bilinear and in the case of bicubic in-plane parametrisations (see Equations ~\ref{eq:s_0_bilinear}, ~\ref{eq:s_1_bilinear} and ~\ref{eq:def_Eq_sumC_bicubic}) it follows that
\begin{equation}
\begin{array}{l}
"Eq(0)^{(bicubic)}"="Eq(0)^{(bilinear)}"\cr
"Eq(1)^{(bicubic)}"="Eq(0)^{(bilinear)}"+"Eq(1)^{(bilinear)}"
\end{array}
\end{equation}
Therefore systems of equations {\it "Eq(0)"}, {\it "Eq(1)"} are equivalent for $\T\Pi^{(bilinear)}$ and $\T\Pi{(bicubic)}$ in-plane parametrisations. 

\nopagebreak 
\eop${}_{{\bf Lemma ~\ref{lemma:Eq0_Eq1_bicubic}}}$

\vspace{0.32in}
\noindent
{\bf Proof of Theorem ~\ref{theorem:C_equation_sufficiency} and Theorem ~\ref{theorem:LR_classification} }

\paragraph{General approach to the proof.}

By definition, the {\it "Middle"} system of equations consists of the {\it "sumC-equation"} and the {\it "Restricted Middle"} system of indexed equations. {\it "sumC-equation"} involves the  "central" {\it middle} control points alone, while every one of the indexed equations {\it "Eq(s)"} involves both the "central" and the "side" {\it middle} control points.

One of the basic concepts of the proof is  the assumption  that the "central" {\it middle} control points are classified prior to and independently of classification of the "side" {\it middle} control points. The classification of the "central" {\it middle} control points is assumed to be made is such a manner that {\it "sumC-equation"} (and some additional requirements explained below) are satisfied.
The assumption allows to consider the "side" {\it middle} control points as the only variables of the {\it "Restricted Middle"} system , while the "central" {\it middle} control points participate in the system as components of the right-side coefficients. A consistency analysis of the {\it "Restricted Middle"} system establishes some requirements which should be satisfied by the right-side coefficients of the system. The analysis leads to the necessary and sufficient conditions formulated in Theorem ~\ref{theorem:C_equation_sufficiency}, which define additional requirements to the classification of the "central" {\it middle} control points.

Note, that according to this approach, responsibilities of the "central" and the "side" {\it middle} control points are shared in the following manner. The "central" {\it middle} control points are responsible for the satisfaction of {\it "sumC-equation"} as well as for the consistency of the {\it "Restricted Middle"} system of equations, while the classification of the "side" {\it middle} control points reflects the rank analysis of the {\it "Restricted Middle"} system.

\paragraph{Structure of the {\it "Restricted Middle"} system of equations.}
As it was mentioned above, the proof starts with the consistency and rank analysis of the {\it "Restricted Middle"} system of equations.

Similarly to the case of global bilinear in-plane parametrisation $\T\Pi^{(bilinear)}$, let $A$ be matrix of the full system of linearized $G^1$-continuity equations {\it excluding "sumC-equation"} in terms of differences $\dL_t,\dR_t$ $t=0,\ldots,n$ and let $b_s$ $(s=0,\ldots,n+2)$ be the right-side coefficients of the system. 
In other words,  matrix $A$ and vector $b$ correspond to the system of $n+3$ {\it indexed} equations {\it "Eq(s)"} $s=0,\ldots,n+2$
\begin{equation}
A
\left(\begin{array}{c}
\dL_0\cr
\dR_0\cr
\vdots\cr
\dL_n\cr
\dR_n\cr
\end{array}\right)= 
\left(\begin{array}{c}
b_0\cr
\vdots\cr
b_{n+2}\cr
\end{array}\right)
\label{eq:system_indexed_bicubic}
\end{equation}
According to Lemma ~\ref{lemma:Eq_sumC_linear_system}, $A$ contains coefficients of $\dL_t,\dR_t$ ($t=0,\ldots,n$) in expressions for $sumLR_s$ (see Definition ~\ref{def:sumLR_sumC}) and
\begin{equation}
b_s = (-1)^{s+1}\sum\limits_{k=0}^s (-1)^k sumC_k
\label{eq:def_b_s}
\end{equation}
for $s=0,\ldots,n+2$.

\noindent
Assumptions of Theorems ~\ref{theorem:C_equation_sufficiency} and ~\ref{theorem:LR_classification} mean that 
\begin{itemize}
\item[]
Equations {\it "Eq(s)"} for $s=0, s=1, s=n'+1$ for $n\ge 4$ and equation {\it "Eq(n'+2)"} for $s\ge 5$ are satisfied.
\item[]
All control points that contribute to differences $\dL_t$, $\dR_t$ for $t=0,1,n'-1$ for $n\ge 4$ and for $t=n'$ for $n\ge 5$ are classified.
\item[]
\end{itemize}

\noindent
Let $\hat A$ denote matrix of the {\it "Restricted Middle"} system of equations in terms of differences $\dL_t, \dR_t$ ($t=2,\ldots,n'-2$), $\hat A_s$ $(s=2,\ldots,n')$ denotes a row of the matrix $\hat A$ and let $\hat b_s$ ($s=2,\ldots,n'$) be the right-side coefficients of the system.
Matrix $\hat A$ is a $(n'-1)\times 2(n'-3)$ submatrix of matrix $A$ composed of rows $A_s$ for $s=2,\ldots,n'$ and columns corresponding to $\dL_2,\dR_2,\ldots,\dL_{n'-2},\dR_{n'-2}$. Matrix $\hat A$ has the following form

\begin{equation}
\hspace{-0.45in}
{\small
\begin{array}{lccccccccc}
&\!\!\dL_2 & \!\!\dR_2 & \!\!\dL_3 & \!\!\dR_3 & \ldots 
&\!\!\!\!\dL_{n'\!-\!3}\!\! & \!\!\!\!\dR_{n'\!-\!3}\!\! 
&\!\!\!\!\dL_{n'\!-\!2}\!\! & \!\!\!\!\dR_{n'\!-\!2}\!\!\cr
s\!=\!2    & 
\cnk{n}{2}l_0    & \cnk{n}{2}r_0    &
0                & 0                &
\ldots           &
0                & 0                &
0                & 0                \cr
s\!=\!3    &
2\cnk{n}{2}l_1   &2\cnk{n}{2}r_1    &
\cnk{n}{3}l_0    & \cnk{n}{3}r_0    &
\ldots           &
0                & 0                &
0                & 0                \cr
s\!=\!4    &
\cnk{n}{2}l_2    & \cnk{n}{2}r_2    &
2\cnk{n}{3}l_1   &2\cnk{n}{3}r_1    &
\ldots           &
0                & 0                &
0                & 0                \cr
s\!=\!5    & 
0                & 0                &
\cnk{n}{3}l_2    & \cnk{n}{3}r_2    &
\ldots           &
0                & 0                &
0                & 0                \cr
\vdots & 
\vdots           & \vdots           &
\vdots           & \vdots           &
\ddots           &
\vdots           & \vdots           &
\vdots           & \vdots           \cr
s\!=\!n'\!-\!3\!\! &
0                & 0                &
0                & 0                &
\ldots           &
\cnk{n}{n'-3}l_0 & \cnk{n}{n'-3}r_0 &
0                & 0                \cr
s\!=\!n'\!-\!2\!\! & 
0                & 0                &
0                & 0                &
\ldots           &
2\cnk{n}{n'-3}l_1&2\cnk{n}{n'-3}r_1 &
\cnk{n}{n'-2}l_0 & \cnk{n}{n'-2}r_0 \cr
s\!=\!n'\!-\!1\!\! &
0                & 0                &
0                & 0                &
\ldots           &
\cnk{n}{n'-3}l_2 & \cnk{n}{n'-3}r_2 &
2\cnk{n}{n'-2}l_1&2\cnk{n}{n'-2}r_1 \cr
s\!=\!n'   & 
0                & 0                &
0                & 0                &
\ldots           &
0                & 0                &
\cnk{n}{n'-2}l_2 & \cnk{n}{n'-2}r_2 \cr
\end{array}
}
\end{equation}
Right-side coefficients $\hat b_s$ $(s=2,\ldots,n')$ are given by the next formula

\begin{equation}
\hspace{-0.35in}
{\small
\begin{array}{l}
\hat b_2 = b_2-
\left[\Cnk{n}{0}l_2\dL_0+\Cnk{n}{0}r_2\dR_0+
2\Cnk{n}{1}l_1\dL_1+2\Cnk{n}{1}r_1\dR_1\right]\cr
\hat b_3 = b_3-
\left[\Cnk{n}{1}l_2\dL_1+\Cnk{n}{1}r_2\dR_1\right]\cr
\hat b_s = b_s\ \ {\rm for}\ s=4,\ldots,n'-2\cr
\hat b_{n'-1} = b_{n'-1}-
\left[\Cnk{n}{n'-1}l_0\dL_{n'-1}+
\Cnk{n}{n'-1}r_0\dR_{n'-1}\right]\cr                
\hat b_{n'}\!=\! b_{n'}\!-\!
\left[2\Cnk{n}{n'-1}l_1\dL_{n'\!-\!1}\!+\!
2\Cnk{n}{n'-1}r_1\dR_{n'\!-\!1}\!+\!
\Cnk{n}{n'}l_0\dL_{n'}\!+\!
\Cnk{n}{n'}r_0\dR_{n'}\right]\!\!\!\!\!\!\!
\end{array}
}
\label{eq:def_hat_b_bicubic}
\end{equation}

For example, for $n=4$ and $n=5$, matrix $A$ has the explicit form given below, where the last row and two last columns are relevant only in the case of $n=5$ (the row and column are separated by single lines).

\vspace{0.1in}
\hspace{-0.35in}
{\scriptsize
\begin{tabular}{ccccc||cccc||cc|cc}
     &$\dL_0$&$\dR_0$&$\dL_1$&$\dR_1$&$\dL_2$&$\dR_2$&
      $\dL_3$&$\dR_3$&$\dL_4$&$\dR_4$&$\dL_5$&$\dR_5$$\!\!\!$\cr
$\!\!\!\!s\!=\!0\!\!\!$&
      $ \cnk{n}{0}l_0$&$ \cnk{n}{0}r_0$&
      $    0         $&$    0         $&
      $    0         $&$    0         $&
      $    0         $&$    0         $&
      $    0         $&$    0         $&
      $    0         $&$    0         $$\!\!\!$\cr
$\!\!\!\!s\!=\!1\!\!\!$&
      $\!\!2\cnk{n}{0}l_1\!\!$&$\!\!2\cnk{n}{0}r_1\!\!$&
      $ \cnk{n}{1}l_0$&$ \cnk{n}{1}r_0$&
      $    0         $&$    0         $&
      $    0         $&$    0         $&
      $    0         $&$    0         $&
      $    0         $&$    0         $$\!\!\!$\cr
\hline
\hline
$\!\!\!\!s\!=\!2\!\!\!$&
      $ \cnk{n}{0}l_2$&$ \cnk{n}{0}r_2$&
	 $\!\!2\cnk{n}{1}l_1\!\!$&$\!\!2\cnk{n}{1}r_1\!\!$&
      $ \cnk{n}{2}l_0$&$ \cnk{n}{2}r_0$&
      $    0         $&$    0         $&
      $    0         $&$    0         $&
      $    0         $&$    0         $$\!\!\!$\cr
$\!\!\!\!s\!=\!3\!\!\!$&
      $    0         $&$    0         $&
      $ \cnk{n}{1}l_2$&$ \cnk{n}{1}r_2$&
	 $\!\!2\cnk{n}{2}l_1\!\!$&$\!\!2\cnk{n}{2}r_1\!\!$&
      $ \cnk{n}{3}l_0$&$ \cnk{n}{3}r_0$&
      $    0         $&$    0         $&
      $    0         $&$    0         $$\!\!\!$\cr
$\!\!\!\!s\!=\!4\!\!\!$&
      $    0         $&$    0         $&
      $    0         $&$    0         $&
      $ \cnk{n}{2}l_2$&$ \cnk{n}{2}r_2$&
	 $\!\!2\cnk{n}{3}l_1\!\!$&$\!\!2\cnk{n}{3}r_1\!\!$&
      $ \cnk{n}{4}l_0$&$ \cnk{n}{4}r_0$&
      $    0         $&$    0         $$\!\!\!$\cr
$\!\!\!\!s\!=\!5\!\!\!$&
      $    0         $&$    0         $&
      $    0         $&$    0         $&
      $    0         $&$    0         $&
      $ \cnk{n}{3}l_2$&$ \cnk{n}{3}r_2$&
	 $\!\!2\cnk{n}{4}l_1\!\!$&$\!\!2\cnk{n}{4}r_1\!\!$&
      $ \cnk{n}{5}l_0$&$ \cnk{n}{5}r_0$$\!\!\!$\cr
\hline
\hline
$\!\!\!\!s\!=\!6\!\!\!$&
      $    0         $&$    0         $&
      $    0         $&$    0         $&
      $    0         $&$    0         $&
      $    0         $&$    0         $&
      $ \cnk{n}{4}l_2$&$ \cnk{n}{4}r_2$&
	 $\!\!2\cnk{n}{5}l_1\!\!$&$\!\!2\cnk{n}{5}r_1\!\!$$\!\!\!$\cr
\hline
$\!\!\!\!s\!=\!7\!\!\!$&
      $    0         $&$    0         $&
      $    0         $&$    0         $&
      $    0         $&$    0         $&
      $    0         $&$    0         $&
      $    0         $&$    0         $&
      $ \cnk{n}{5}l_2$&$ \cnk{n}{5}r_2$$\!\!\!$
\end{tabular}
}

\noindent
For both $n=4$ and $n=5$ matrix $\hat A$ is a middle $4\times 4$  submatrix of $A$, composed of rows $s=2,3,4,5$ and columns corresponding to $\dL_2,\dR_2,\dL_3,\dR_3$ ($\hat A$ is separated by double lines).

\paragraph{The possible types of dependencies between coefficients of conventional weight functions $l(v)$ and $r(v)$.}

As shown in the previous paragraph, matrix $\hat A$ contains coefficients of conventional weight functions $l(v)$ and $r(v)$. Different types of dependencies between coefficients of the weight functions lead to the different consistency and rank analysis of the {\it "Restricted Middle"} system of equations. Definition ~\ref{def:cases_classification} presents a refinement of the possible types of dependencies listed in Theorem ~\ref{theorem:LR_classification}. Subdivision into these subcases and auxiliary notations given in the Definition are justified by the series of the following Auxiliary Lemmas.

\begin{definition}
Let conventional weight functions $l(v)$ and $r(v)$ be defined by global in-plane parametrisation $\T\Pi^{(bicubic)}$. Then for an edge with one boundary vertex, the following types of dependencies between coefficients of the weight functions $l(v)$ and $r(v)$ will be considered (here $g^{(ij)}$, $i,j\in\{0,1,2\}$ are given by  Definition ~\ref{def:g_ij})
\begin{itemize}
\item
{\bf "CASE 0"} 
complement of {\bf "CASE 1"} and {\bf "CASE 2"}
\item
{\bf "CASE 1"} 
\begin{eqnarray}
g^{(01)} g^{(12)} g^{(02)}\neq 0,\ \ \ \ \ \ \
\left\{g^{(02)}\right\}^2=4g^{(01)}g^{(12)}
\nonumber
\end{eqnarray}
In this case it is possible to define such a scalar value $\kappa$ that 
\begin{equation}
\kappa = 
-\frac{2 g^{(01)}}{g^{(02)}} = 
-\frac{g^{(02)}}{2 g^{(12)}}
\label{def:kappa}
\end{equation}
The following two subcases are defined
\begin{itemize}
\item[]
{\bf "CASE 1.a"} $\kappa\neq -1$
\item[]
{\bf "CASE 1.b"} $\kappa = -1$
\end{itemize}

\item
{\bf "CASE 2"}
\begin{eqnarray}
g^{(01)}=g^{(12)}=g^{(02)}=0
\nonumber
\end{eqnarray}
In this case it is possible to define such scalars $\kappa^{(ij)}$ $(i,j\in\{0,1,2\})$ that
\begin{equation}
\kappa^{(ij)} = \frac{l_i}{l_j} = \frac{r_i}{r_j}
\label{def:kappa_ij}
\end{equation}
The special notations are used for $\kappa^{(02)}$ and $\kappa^{(12)}$
\begin{equation}
\xi = \kappa^{(02)},\ \ \ \ \eta = \kappa^{(12)}
\label{eq:def_ksi_eta}
\end{equation}
and matrix $G$ is defined as follows
\begin{equation}
G = \left(\begin{array}{cc}
-\xi   & 2\eta\xi \cr
-2\eta & 4\eta^2-\xi
\end{array}\right)
\label{def:G} 
\end{equation}

\begin{itemize}
\item[]
{\bf "CASE 2.1"} $\eta^2\neq\xi$

In this case $G$ has two different eigenvalues $\omega_1$ and $\omega_2$
\begin{itemize}
\item[]
{\bf "CASE 2.1.a"} 
\begin{equation}
\begin{array}{l}
\sqrt{\omega_1}=\eta+\sqrt{\eta^2-\xi}\neq 1\ \ \ and\cr \sqrt{\omega_2}=\eta-\sqrt{\eta^2-\xi}\neq 1
\end{array}
\label{eq:sqrt_omega_1_2}
\end{equation}
\item[]
{\bf "CASE 2.1.b"} 
$\sqrt{\omega_1}=1$ or $\sqrt{\omega_2}=1$
\end{itemize}
\item[]
{\bf "CASE 2.2"} $\eta^2 = \xi$

In this case $G$ has a single eigenvalue $\omega$ 
\begin{itemize}
\item[]
{\bf "CASE 2.2.a"} 
$\sqrt{\omega}=\eta\neq 1$
\item[]
{\bf "CASE 2.2.b"} 
$\sqrt{\omega}=1$
\end{itemize}
\end{itemize}
\end{itemize}
\label{def:cases_classification}
\end{definition}

\paragraph {Explanation regarding the subdivision into the basic cases {\bf "CASE 0"}, {\bf "CASE 1"} and {\bf "CASE 2"}.}

Both proof of the consistency and the rank analysis of {\it "Restricted Middle"} system deal with zero non-trivial linear combinations of rows. Auxiliary Lemma ~\ref{aux_lemma:3_basis_cases} provides necessary and sufficient conditions for the existence of such a linear combination and describes the structure of the set of its coefficients. Subdivision into three basic cases {\bf "CASE 0"}, {\bf "CASE 1"}, {\bf "CASE 2"} directly follows from the classification given in the Auxiliary Lemma. 

\begin{auxillary_lemma}
A non-trivial combination of the rows of matrix $\hat A$ is equal to zero
\begin{equation}
\sum_{s=2}^{n'} \alpha_s\hat A_s = 0
\end{equation}
if and only if either one of the following conditions holds
\begin{itemize}
\item[{\bf(1)}]
\begin{itemize}
\item[]\hspace{-0.2in}
Conditions of {\bf "CASE 1"} are satisfied and 
\item[]\hspace{-0.2in}
Coefficients of the linear combination obey the recursive relation
\begin{equation}
\alpha_{s+1}=\kappa\alpha_s,\ \ \ s=2,\ldots,n'-1
\label{eq:rec_coeff_case_1}
\end{equation}
\end{itemize}
\item[{\bf(2)}]
\begin{itemize}
\item[]\hspace{-0.2in}
Conditions of {\bf "CASE 2"} are satisfied and 
\item[]\hspace{-0.2in}
Coefficients of the linear combination obey the recursive relation
\begin{equation}
\alpha_{s+2}=-\xi\alpha_s-2\eta\alpha_{s+1},\ \ \ s=2,\ldots,n'-2
\label{eq:rec_coeff_case_2}
\end{equation}
\end{itemize}
\end{itemize}
\label{aux_lemma:3_basis_cases}
\end{auxillary_lemma}
{\bf Proof} See Appendix, after proof of Theorem ~\ref{theorem:C_equation_sufficiency} and Theorem ~\ref{theorem:LR_classification}.

\begin{conclusion}
Auxiliary Lemma ~\ref{aux_lemma:3_basis_cases} implies that the rank and consistency of the {\it "Restricted Middle"} system of equations depend in the following way on the subdivision in the basic cases
\begin{itemize}
\item[]{\bf "CASE 0"}
\begin{itemize}
\item[]
Matrix $\hat A$ has the full row rank, $rank(\hat A) = n'-1$.
\item[]
The {\it "Restricted Middle"} system is consistent.
\end{itemize}

\item[]{\bf "CASE 1"}
\begin{itemize}
\item[]
Matrix $\hat A$ has row rank deficiency $1$, $rank(\hat A) = n'-2$.
\item[]
The {\it "Restricted Middle"} system is consistent if and only if for any set of coefficients $\{\alpha_a\}_{s=2}^{n'}$ which satisfy Equation ~\ref{eq:rec_coeff_case_1}, the corresponding linear combination of the right sides is equal to zero $\sum_{s=2}^{n'}\alpha_s\hat b_s=0$.
\end{itemize}

\item[]{\bf "CASE 2"}
\begin{itemize}
\item[]
Matrix $\hat A$ has row rank deficiency $2$, $rank(\hat A) = n'-3$.
\item[]
The {\it "Restricted Middle"} system is consistent if and only if for any set of coefficients $\{\alpha_a\}_{s=2}^{n'}$ which satisfy Equation ~\ref{eq:rec_coeff_case_2}, the corresponding linear combination of the right sides is equal to zero $\sum_{s=2}^{n'}\alpha_s\hat b_s=0$.
\end{itemize}
\end{itemize}
\label{concl:middle_system}
\end{conclusion}
Conclusion ~\ref{concl:middle_system}
allows to study the necessary and sufficient conditions for the consistency of the {\it "Restricted Middle"} system in terms of right-side coefficients, or in other words, in terms of the "central" {\it middle} control points.

Before proceeding with the consistency analysis, let us give some additional explanations concerning subdivision of {\bf "CASE 2"}, which has the most complicated structure.

\paragraph{Explanation regarding the subdivision of {\bf "CASE 2"} into {\bf "CASE 2.1"} and {\bf "CASE 2.2"}.}

Auxiliary Lemma ~\ref{aux_lemma:G_explanation} justifies the definition of matrix $G$ introduced in Equation ~\ref{def:G}. 

\begin{auxillary_lemma}
Let conditions of {\bf "CASE 2"} be satisfied and let matrix $G$ be defined by Equation ~\ref{def:G}. 
Let further $\{\alpha_s\}_{s=2}^{n'}$ be coefficients of a non-trivial zero combination of rows $\sum_{s=2}^{n'}\alpha_s\hat A_s=0$. Then dependency of the coefficients has the following matrix form
\begin{equation}
(\alpha_{s+2},\alpha_{s+3})=(\alpha_s,\alpha_{s+1})G
\label{eq:G_explanation}
\end{equation}
for $s=2,\ldots,n'-3$.
\label{aux_lemma:G_explanation}
\end{auxillary_lemma}
{\bf Proof} 
Auxiliary Lemma ~\ref{aux_lemma:G_explanation}
immediately follows from Equation ~\ref{eq:rec_coeff_case_2} given in Auxiliary Lemma ~\ref{aux_lemma:3_basis_cases}.

\nopagebreak 
\eop${}_{{\bf Auxiliary\ Lemma ~\ref{aux_lemma:G_explanation}}}$

\vspace{0.1in}
Auxiliary Lemma ~\ref{aux_lemma:jordan_form} describes the Jordan form of matrix $G$. The Jordan form appears to play an important role in the consistency analysis and explains the subdivision of {\bf "CASE 2"} into subcases {\bf "CASE 2.1"} and {\bf "CASE 2.2"}.

\begin{auxillary_lemma}
Let the conditions of {\bf "CASE 2"} be satisfied. Then the matrix $G$ (see Equation ~\ref{def:G}) has the following eigenvalues, eigenvectors and Jordan form
\begin{itemize}
\item []
{\hspace{-0.3in}\bf "CASE 2.1"}
\\ 

\vspace{-0.15in}
Matrix $G$ has two different eigenvalues 
\begin{equation}
\omega_{1,2} = -\xi+2\eta(\eta\pm\sqrt{\eta^2-\xi})
\label{eq:lambda_12_def}
\end{equation}
which correspond to the eigenvectors
\begin{equation}
v_{1,2}=
\left(\begin{array}{c}
\xi \cr 
\eta\pm\sqrt{\eta^2-\xi}
\end{array}\right)=
\left(\begin{array}{c}
\xi \cr 
\sqrt{\omega_{1,2}}
\end{array}\right)
\label{eq:v_12_def}
\end{equation}
and $G$ can be represented in the form
\begin{equation}
G=(v_1\ v_2)
\left(\begin{array}{cc}
\omega_1 & 0\cr
0 & \omega_2
\end{array}\right)
(v_1\ v_2)^{-1}
\label{eq:G_structure_2_1}
\end{equation}

\item[]
{\hspace{-0.3in}\bf "CASE 2.2"} 
\\

\vspace{-0.15in}
Matrix $G$ has a single eigenvalue
\begin{equation}
\omega=\eta^2
\label{eq:lambda_def}
\end{equation}
which corresponds to the eigenvector
\begin{equation}
v=
\left(\begin{array}{c}
\eta^2\cr \eta\end{array}\right)=
\left(\begin{array}{c}
\omega\cr \sqrt{\omega}\end{array}\right)
\label{eq:v_def}
\end{equation}
and $G$ can be represented in the form
\begin{equation}
G=(v\ u)
\left(\begin{array}{cc}
\omega & 1\cr
0 & \omega
\end{array}\right)
(v\ u)^{-1}
\label{eq:G_structure_2_2}
\end{equation}
where
\begin{equation}
u=\frac{1}{4}
\left(\begin{array}{c}-1\cr\frac{1}{\eta}\end{array}\right)
\end{equation}
\end{itemize}
\label{aux_lemma:jordan_form}
\end{auxillary_lemma}
Auxiliary Lemma ~\ref{aux_lemma:jordan_form}
has a completely straightforward proof, which is not presented in the current work.

The next two paragraphs, which are directly connected to the consistency analysis of the {\it "Restricted Middle"} system, provide the additional explanations for subdivision into subcases {\bf "CASE 0"},{\bf "CASE 1"}, {\bf "CASE 2.1"} and {\bf "CASE 2.2"}.
According to Conclusion ~\ref{concl:middle_system}, the {\it "Restricted Middle"} system of equations is known to be consistent in {\bf "CASE 0"}, therefore it is sufficient to analyse the consistency in {\bf "CASE 1"} and {\bf "CASE 2"}.

\paragraph{Some special relations the for coefficients of conventional weight function c(v).}

It appears, that the coefficients of the conventional weight function $c(v)$ satisfy some special relations, which involve constant $\kappa$ in  {\bf "CASE 1"} and eigenvalues of the matrix $G$ in {\bf "CASE 2"}. These relations are very useful for the consistency analysis of the {\it "Restricted Middle"} system.

\begin{auxillary_lemma}
Let conventional weight function $c(v)$ be defined by global in-plane parametrisation $\T\Pi^{(bicubic)}$. Then for an edge with one boundary vertex, coefficients of the weight function satisfy the following relations
\begin{itemize}
\item[]{\bf "CASE 1"}
\begin{equation}
\sum_{k=0}^4 \kappa^k \Cnk{4}{k} c_k = 0
\label{eq:coeff_c_relation_1}
\end{equation}
\item[]{\bf "CASE 2.1"}
\begin{equation}
\sum_{k=0}^4\left(-\sqrt{\omega_{1,2}}\right)^k\CNK{4}{k} c_k = 0
\label{eq:coeff_c_relation_2_1}
\end{equation}
\item[]{\bf "CASE 2.2"}
\begin{eqnarray}
\sum_{k=0}^4\left(-\sqrt{\omega}\right)^k\CNK{4}{k} c_k = 0
\label{eq:coeff_c_relation_2_2_coeff}\\
\sum_{k=0}^4\left(-\sqrt{\omega}\right)^k k\CNK{4}{k} c_k = 0
\label{eq:coeff_c_relation_2_2_deriv}
\end{eqnarray}
\end{itemize}
\label{aux_lemma:coeff_c_relations}
\end{auxillary_lemma}
{\bf Proof} See Appendix, after proof of Theorem ~\ref{theorem:C_equation_sufficiency} and Theorem ~\ref{theorem:LR_classification}.

\paragraph{Necessary and sufficient conditions for the consistency of the {\it "Restricted Middle"} system of equations in terms of $b_s$.}

The purpose of the paragraph is to formulate the necessary and sufficient conditions for the consistency of the {\it "Restricted Middle"} system of equations in terms of the right-side coefficients of the full system of the indexed equations (see Equations ~\ref{eq:system_indexed_bicubic} and ~\ref{eq:def_b_s}). 
The next three paragraphs describe the relations between these conditions and the actual degrees of the weight functions and finally explain the rules of classification of the control point $\T C_{n'-2}$.

Auxiliary Lemma ~\ref{aux_lemma:hat_b_s_by_b_s} allows to rewrite a zero linear combination of the right-side coefficients of the {\it "Restricted Middle"} system in terms of the right-side coefficients of the full system of the indexed equations.

\begin{auxillary_lemma}
Let global in-plane parametrisation $\T\Pi^{(bicubic)}$ be considered and for an edge with one boundary vertex, let all non-middle control points be classified and equations {\it "sumC-equation"}, {\it "Eq(s)"} for $s=0,1,n'+1$ and {\it "Eq(n'+2)"} in case of $n\ge 5$ be satisfied. 

Let $\{\alpha_s\}_{s=2}^{n'}$ be coefficients of a non-trivial zero linear combination of rows of the {\it "Restricted Middle"} system of equations $\sum_{s=2}^{n'}\alpha_s \hat A_s = 0$. 
Then the corresponding linear combination of the right-side coefficients of the {\it "Restricted Middle"} system of equations has the following representation in terms of the right-side coefficients of the full system of the indexed equations

\begin{itemize}
\item[]
{\bf "CASE 1"}
\begin{equation}
\sum_{s=2}^{n'}\alpha_s \hat b_s = 
\alpha_0\sum_{s=0}^{n+2}\kappa^s b_s
\label{eq:hat_b_s_by_b_s_1}
\end{equation}
where by the definition 
\begin{equation}
\alpha_0=\kappa^{-2}\alpha_2
\label{eq:def_alpha_0_case_1}
\end{equation}
\item[]
{\bf "CASE 2"}
\begin{equation}
\sum_{s=2}^{n'}\alpha_s \hat b_s = 
(\alpha_0\ \alpha_1)
\sum_{s=0}^{\lceil n'/2\rceil} G^s
\left(\begin{array}{c}
b_{2s}\cr b_{2s+1}\end{array}\right)
\label{eq:hat_b_s_by_b_s_2}
\end{equation}
where by the definition
\begin{equation}
(\alpha_0\ \alpha_1)=(\alpha_2\ \alpha_3)G^{-1}
\end{equation}
and $b_{n'+2}=b_7=0$ in case of $n=4$.
\end{itemize}
\label{aux_lemma:hat_b_s_by_b_s}
\end{auxillary_lemma}
{\bf Proof} See Appendix, after proof of Theorem ~\ref{theorem:C_equation_sufficiency} and Theorem ~\ref{theorem:LR_classification}.

\vspace{0.1in}
Auxiliary Lemma ~\ref{aux_lemma:consisitency_b_s}
provides an elegant necessary and sufficient conditions for the consistency of the {\it "Restricted Middle"} system of equations in terms of the right-side coefficients of the full system of the indexed equations.

\begin{auxillary_lemma}
Let global in-plane parametrisation $\T\Pi^{(bicubic)}$ be considered and for an edge with one boundary vertex, let all non-middle control points be classified and equations {\it "sumC-equation"}, {\it "Eq(s)"} for $s=0,1,n'+1$ and {\it "Eq(n'+2)"} in case of $n\ge 5$ be satisfied. 

Then the {\it "Restricted Middle"} system of equations is consistent if and only if the following conditions hold
\begin{itemize}
\item[]
{\bf "CASE 1"}
\begin{equation}
\sum_{s=0}^{n+2}\kappa^s b_s = 0
\label{eq:consisitency_b_s_1}
\end{equation}
\item[]
{\bf "CASE 2.1"}
\begin{equation}
\sum_{s=0}^{n+2}(-\sqrt{\omega_{1,2}})^s b_s = 0
\label{eq:consisitency_b_s_2_1}
\end{equation}
\item[]
{\bf "CASE 2.2"}
\begin{equation}
\sum_{s=0}^{n+2}(-\sqrt{\omega})^s b_s = 0
\label{eq:consisitency_b_s_2_2_f}
\end{equation}
\begin{equation}
\sum_{s=0}^{n+2}(-\sqrt{\omega})^{s-1} s\ b_s = 0
\label{eq:consisitency_b_s_2_2_s}
\end{equation}
\end{itemize}
\label{aux_lemma:consisitency_b_s}
\end{auxillary_lemma}
{\bf Proof} See Appendix, after proof of Theorem ~\ref{theorem:C_equation_sufficiency} and Theorem ~\ref{theorem:LR_classification}.

\paragraph{Relations between linear combinations of $b_s$ and coefficients of the conventional weight function $c(v)$.}

In the current paragraph, the linear combinations of the right-sides coefficients, described in  Auxiliary Lemma 
~\ref{aux_lemma:consisitency_b_s}, are expanded using the explicit formula for $b_s$  (see Equations ~\ref{eq:def_b_s} and ~\ref{eq:def_sumC}). Straightforward Auxiliary Lemma ~\ref{aux_lemma:phi_formal_expansion} describes relation between the linear combinations and coefficients of the weight function $c(v)$ with respect to the power basis. The Auxiliary Lemma finally justifies the fine subdivision into subcases indexed by letters {\bf (a)} and {\bf (b)} (see Definition ~\ref{def:cases_classification}) and makes the first step towards explanation of the classification rule for the control point $\T C_{n'-2}$ (see Lemma ~\ref{lemma:C_classification}).

\begin{auxillary_lemma}
For a scalar value $\phi$ the following relations hold
\begin{equation}
\hspace{-0.4in}
\sum_{s=0}^{n+2}\phi^s b_s=\left\{
\begin{array}{l}
\begin{array}{l}
"sumC\!-\!equation"\left[-\frac{(-\phi)^{n+3}}{1+\phi}\right]+\cr
\frac{1}{1+\phi}
\left[\sum\limits_{i=0}^{n-1}\phi^i\Cnk{n-1}{i}\Delta C_i\right]
\left[\sum\limits_{j=0}^4 \phi^j \Cnk{4}{j} c_j\right]
\ \ \ {\rm if}\ \  \phi\neq -1
\end{array}\!\!\!\!\!\!\!\!\!\!\!\!\cr
\begin{array}{l}
"sumC\!-\!equation"[-(n+3)] +\cr\cr
\left[\sum\limits_{i=0}^{n-1}(-1)^i i\Cnk{n-1}{i}\Delta C_i\right]
c_4^{(power)}+\cr
"C\!-\!equation"\left[4 c_4^{(power)}+c_3^{(power)}\right]
\ \ \ \ \ \ {\rm if}\ \ \phi = -1
\end{array}\!\!\!\!\!\!\!\!\!\!\!\!
\end{array}
\right.
\end{equation}

\begin{equation}
\hspace{-0.4in}
\sum_{s=0}^{n+2}s\phi^{s-1} b_s=\left\{
\begin{array}{l}
\begin{array}{l}
"sumC\!-\!equation"
\left[
\frac{(-\phi)^{n+2}(n+3)-(-\phi)^{n+3}(n+2)}{(1+\phi)^2}
\right]+\cr
\frac{1}{(1+\phi)^2}
\left[
\sum\limits_{i=0}^{n-1}
\left(\frac{1+\phi}{\phi}i-1\right) \phi^i \Cnk{n-1}{i}\Delta C_i
\right]
\left[\sum\limits_{j=0}^4 \phi^j \Cnk{4}{j} c_j\right]+\cr
\frac{1}{1+\phi}
\left[
\sum\limits_{i=0}^{n-1}
\phi^i \Cnk{n-1}{i}\Delta C_i
\right]
\left[\sum\limits_{j=0}^4 \phi^j j\Cnk{4}{j} c_j\right]
\ \ \ \ \ \ \ \ \ \ {\rm if}\ \ \phi\neq -1
\end{array}\!\!\!\!\!\!\!\!\!\!\!\!\!\!\cr\cr
\begin{array}{l}
"sumC\!-\!equation"\left[\frac{1}{2}(n+2)(n+3)\right] -\cr
\left[
\sum\limits_{i=0}^{n-1}
(-1)^i (i-1)i\Cnk{n-1}{i}\Delta C_i
\right]
\left[\frac{1}{2}c_4^{(power)}\right]-\cr
\left[
\sum\limits_{i=0}^{n-1}
(-1)^i i\Cnk{n-1}{i}\Delta C_i
\right]
\left[4 c_4^{(power)}+c_3^{(power)}\right]-\cr
"C\!-\!equation"
\left[
6 c_4^{(power)}\!\!+\!3 c_3^{(power)}\!\!+\!c_2^{(power)}\right]
\ \ {\rm if}\ \ \phi = -1
\end{array}\!\!\!\!\!\!\!\!\!\!\!\!\!\!
\end{array}
\right.
\end{equation}
\label{aux_lemma:phi_formal_expansion}
\end{auxillary_lemma}
Of course, this technically complicated Auxiliary Lemma becomes much more simple and meaningful if it is applied to $\phi=
\kappa, \omega_{1,2}, \omega$ and if one uses the special relations between the coefficients of the weight function $c(v)$ given in Auxiliary Lemma ~\ref{aux_lemma:coeff_c_relations} and assumes that {\it "sumC-equation"} is satisfied. Further simplification of Auxiliary Lemma ~\ref{aux_lemma:phi_formal_expansion} is based on some properties of the actual degrees of the weight functions, described in the next paragraph.

\paragraph{Actual degrees of the weight functions.}
It appears that the division into the refined cases given in Definition ~\ref{def:cases_classification} is closely connected to the actual degrees of conventional weight functions. 

\begin{auxillary_lemma}
Let the conventional weight functions $l(v)$, $r(v)$ be defined by the global in-plane parametrisation $\T\Pi^{(bicubic)}$. Then for an edge with one boundary vertex, the maximal actual degree of the weight functions obeys the following equality:
\begin{equation}
max\_deg(l,r)=
\left\{
\begin{array}{lll}
2 & {\rm in} &
{\bf "CASE\ 0"},\cr
&&{\bf "CASE\ 1.a"},\cr
&&{\bf "CASE\ 2.1.a"},\cr 
&&{\bf "CASE\ 2.2.a"}\cr
1 & {\rm in} &
{\bf "CASE\ 1.b"},\cr
&&{\bf "CASE\ 2.1.b"}\cr
0 & {\rm in} &
{\bf "CASE\ 2.2.b"}
\end{array}
\right.
\end{equation}
\label{aux_lemma:max_deg_l_r_classification}
\end{auxillary_lemma}
{\bf Proof} See Appendix, after proof of Theorem ~\ref{theorem:C_equation_sufficiency} and Theorem ~\ref{theorem:LR_classification}.

\vspace{0.15in}
\noindent
Auxiliary Lemma ~\ref{aux_lemma:max_deg_l_r_classification} together with Lemma ~\ref{lemma:c_lr_degree_relation} lead to the next Conclusion.

\begin{conclusion}
Let the conventional weight function $c(v)$ be defined by the global in-plane parametrisation $\T\Pi^{(bicubic)}$. Then for an edge with one boundary vertex, the actual degree of the weight function $c(v)$ satisfies the following inequalities
\begin{equation}
deg(c)\le
\left\{
\begin{array}{lll}
3 & {\rm in\ } &
{\bf "CASE\ 1.b"},\cr
&&{\bf "CASE\ 2.1.b"}\cr
2 & {\rm in\ } &
{\bf "CASE\ 2.2.b"}
\end{array}
\right.
\end{equation}
\label{conclusion:actual_deg_c_classification}
\end{conclusion}

\paragraph{Necessary and sufficient conditions for the consistency of the {\it "Middle"} system of equations in terms of 
{\it "C-equation"} and coefficients of conventional weight function $c(v)$ with respect to the power basis.}

Auxiliary Lemmas ~\ref{aux_lemma:consisitency_b_s},
~\ref{aux_lemma:phi_formal_expansion}, 
~\ref{aux_lemma:coeff_c_relations},
~\ref{aux_lemma:max_deg_l_r_classification} and Conclusion
~\ref{conclusion:actual_deg_c_classification} lead to a simple form of necessary and sufficient conditions for the consistency of the {\it "Middle"} system.

\begin{auxillary_lemma}
Consider the  in-plane parametrisation $\T\Pi^{(bicubic)}$ and for an edge with one boundary vertex, let all  non-middle control points be classified and equations {\it "Eq(s)"} for $s=0,1,n'+1$ and {\it "Eq(n'+2)"} in case of $n\ge 5$ be satisfied. Then the {\it "Middle"} system of equations is consistent if and only if:
\begin{itemize}
\item 
\begin{equation}
"C\!-\!equation"c_4^{(power)}=0
\label{eq:middle_consisitency_c_power_4}
\end{equation}
in the case when $max\_deg(l,r)=2$
\item
\begin{equation}
"C\!-\!equation"c_3^{(power)}=0
\label{eq:middle_consisitency_c_power_3}
\end{equation}
in the case when $max\_deg(l,r)=1$
\item
\begin{equation}
"C\!-\!equation"c_2^{(power)}=0
\label{eq:middle_consisitency_c_power_2}
\end{equation}
in the case when $max\_deg(l,r)=0$
\end{itemize}
\label{aux_lemma:middle_consisitency_c_power}
\end{auxillary_lemma}
{\bf Proof} See Appendix, after proof of Theorem ~\ref{theorem:C_equation_sufficiency} and Theorem ~\ref{theorem:LR_classification}.

Results of Auxiliary Lemma ~\ref{aux_lemma:middle_consisitency_c_power}
finally allow to complete the proof of Theorem ~\ref{theorem:C_equation_sufficiency} and 
Theorem ~\ref{theorem:LR_classification}.

\paragraph{Proof of Theorem ~\ref{theorem:C_equation_sufficiency}}

The paragraph presents proof of both statements of Theorem ~\ref{theorem:C_equation_sufficiency}.

\vspace{0.15in}
\noindent
{\bf (1)}
According to Lemma ~\ref{lemma:sumC_expansion}, it is sufficient to show that {\it "sumC-equation"} is satisfied in the case when $deg(c)=4$. From the fact that $deg(l),deg(r)\le 2$ and 
Lemma ~\ref{lemma:c_lr_degree_relation} it follows that in this case $max\_deg(l,r)=2$ and $deg(c)-max\_deg(l,r)=2$. The assumption of Theorem ~\ref{theorem:C_equation_sufficiency} implies that {\it "C-equation"} is satisfied. According to Lemma ~\ref{lemma:sumC_expansion}, it means that {\it "sumC-equation"} is satisfied.

\vspace{0.15in}
\noindent 
{\bf (2)}
The second statement of Theorem ~\ref{theorem:C_equation_sufficiency} immediately follows from  Auxiliary Lemma
~\ref{aux_lemma:middle_consisitency_c_power}. Indeed, let $max\_deg(l,r)=2$. In this case Equation ~\ref{eq:middle_consisitency_c_power_4} provides a necessary and sufficient condition for the consistency of the {\it "Middle"} system.
Equation ~\ref{eq:middle_consisitency_c_power_4} is satisfied if and only if either $deg(c)\le 3$ or $deg(c)=4$ and {\it "C-equation"} is satisfied. In other words, the satisfaction of {\it "C-equation"} in the case when 
$deg(c)-max\_deg(l,r)=2$ is a necessary and sufficient condition for the satisfaction of Equation ~\ref{eq:middle_consisitency_c_power_4} and so for the consistency of the {\it "Middle"} system of equations. Similar analysis of cases when $max\_deg(l,r)=1$ and when $max\_deg(l,r)=0$ completes the proof of the second statement of Theorem ~\ref{theorem:C_equation_sufficiency}.

\paragraph{Proof of Theorem ~\ref{theorem:LR_classification}}
The rank analysis presented in Theorem ~\ref{theorem:LR_classification} directly follows from Conclusion ~\ref{concl:middle_system} and from the fact that classification of the control points $\T C_t$ ($t=3,\ldots,n'-2)$ according to Lemma ~\ref{lemma:C_classification} guarantees the consistency of the {\it "Restricted Middle"} system. Like in the case of global bilinear in-plane parametrisation $\T\Pi^{(bilinear)}$, the correctness of classification of the control points $(L_t,R_t)$ $(t=2,\ldots,n'-2)$ into basic and dependent ones can be easily verified by application of the Gauss elimination process to the matrix of the {\it "Restricted Middle"} system.

\nopagebreak 
\eop${}_{{\bf Theorem ~\ref{theorem:C_equation_sufficiency}\ and\ Theorem ~\ref{theorem:LR_classification}}}$

\vspace{0.32in}
\noindent
{\bf Proof of Auxiliary Lemma ~\ref{aux_lemma:3_basis_cases}}

Let $\alpha_2,\ldots,\alpha_{n'}$ be coefficients of a non-trivial linear combination of rows $\sum_{s=2}^{n'} \alpha_s\hat A_s$. The linear combination is equal to zero if and only if the coefficient of every one of $\dL_s$, $\dR_s$ $(s=2,\ldots,n'-2)$ is equal to zero. In other words, it means that the linear combination is equal to zero if and only if the following system of equations is satisfied
\begin{equation}
\begin{array}{lcl}
(\dL_s) & & \alpha_s l_0+2\alpha_{s+1} l_1+\alpha_{s+2}l_2=0\cr
(\dR_s) & & \alpha_s r_0+2\alpha_{s+1} r_1+\alpha_{s+2}r_2=0
\end{array}
\label{eq:coeff_dL_dR}
\end{equation}
where $s=2,\ldots,n'-2$.

First, it will be shown that the satisfaction of either one of the conditions {\bf (1)} or {\bf (2)} of the Auxiliary Lemma is \underline{necessary} for the satisfaction of the last system of equations. Using the fact that $l_0,r_0\neq 0$  ($l_2,r_2\neq 0$), one can express $\alpha_s$ ($\alpha_{s+2}$) independently of the equation corresponding to $\dL_s$ and the equation corresponding to $\dR_s$
(see Equation  ~\ref{eq:coeff_dL_dR}).
Equality of the expressions for $\alpha_s$ ($\alpha_{s+2}$) leads to the conclusion that the following relations should be satisfied for $s=2,\ldots,n'-2$
\begin{eqnarray}
2\alpha_{s+1}g^{(01)}+\alpha_{s+2} g^{(02)}=0
\label{eq:aux_1}
\\
\alpha_s g^{(02)}+2\alpha_{s+1}g^{(12)}=0
\label{eq:aux_2}
\end{eqnarray}

The following two cases are possible
\begin{description}
\item[]
\underline {\it If $g^{(02)}\neq 0$} then
\begin{itemize}
\item 
$g^{(01)}\neq 0$, $g^{(12)}\neq 0$.

On the contrary, let $g^{(01)}=0$. Then according to Equation ~\ref{eq:aux_1}, $\alpha_s=0$ for $s=4,\ldots,n'$. In addition, from formulas for coefficients of $\dL_2$ and $\dL_3$ it follows that $\alpha_2=0$ and $\alpha_3=0$.
It means that equality $g^{(01)}=0$ leads to a trivial linear combination which contradicts the assumption of Auxiliary Lemma 
~\ref{aux_lemma:3_basis_cases}.
In a similar manner it can be shown that $g^{(12)}\neq 0$.

\item
$\left\{g^{(02)}\right\}^2=4g^{(01)}g^{(12)}$. 

The equality immediately follows from plugging in $s=2$ into Equation ~\ref{eq:aux_1} and $s=3$ into Equation ~\ref{eq:aux_2} for $s=3$.

\item
$\alpha_{s+1}=\kappa\alpha_s$ for $s=2,\ldots,n'-1$, as follows from considering Equation ~\ref{eq:aux_1} for $s=n'-2$ and Equation ~\ref{eq:aux_2} for $s=2,\ldots,n'-2$.
\end{itemize}

It means that \underline{\it condition {\bf (1)}} of Auxiliary Lemma ~\ref{aux_lemma:3_basis_cases} is satisfied.

\item
\underline{\it If $g^{(02)} = 0$} then
\begin{itemize}
\item 
$g^{(01)} = 0$, $g^{(12)} = 0$.

On the contrary, let $g^{(01)}\neq 0$. Then according to Equation ~\ref{eq:aux_1}, $\alpha_s=0$ for $s=3,\ldots,n'-1$. From formulas for coefficients of $\dL_2$ and $\dL_{n'-2}$ it follows that $\alpha_2=0$ and $\alpha_{n'}=0$. It means that one gets a trivial linear combination, which contradicts the assumption
of Auxiliary Lemma ~\ref{aux_lemma:3_basis_cases}.
In a similar manner it can be shown that $g^{(12)}=0$.

\item
$\alpha_{s+2} = -\xi\alpha_s-2\eta\alpha_{s+1}$,
as it immediately follows from definitions of $\xi$, $\eta$ and Equation ~\ref{eq:coeff_dL_dR}. (Equality 
$g^{(01)}=g^{(12)}=g^{(02)}=0$ implies that coefficients of $\dL_s$, $\dR_s$ given in Equation ~\ref{eq:coeff_dL_dR} are described by the proportional expressions.)
\end{itemize}
\end{description}

It means that \underline{\it condition {\bf (2)}} of Auxiliary Lemma ~\ref{aux_lemma:3_basis_cases} is satisfied.

\vspace{0.15in}
\noindent
\underline{Sufficiency} of condition {\bf (1)} in {\bf "CASE 1"} and condition {\bf (2)} in {\bf "CASE 2"} is evident 
since the recursive formulas for coefficients of the linear combination (see Equations ~\ref{eq:rec_coeff_case_1} and ~\ref{eq:rec_coeff_case_2}) clearly implies that Equation ~\ref{eq:coeff_dL_dR} is satisfied.

\nopagebreak 
\eop${}_{{\bf Auxiliary\ Lemma ~\ref{aux_lemma:3_basis_cases}}}$

\vspace{0.32in}
\noindent
{\bf Proof of Auxiliary Lemma ~\ref{aux_lemma:coeff_c_relations}}

\noindent
{\bf "CASE 1"}

\noindent
The proof is based on representation of coefficients of the weight function $c(v)$ given in Equation ~\ref{eq:c_lambda_rho}. This representation implies that 
\begin{equation}
\sum_{k=0}^4 \kappa^k c_k \CNK{4}{k}=
\<\sum_{i=0}^2\kappa^i\LL_i\CNK{2}{i},
  \sum_{j=0}^2\kappa^j\RR_j\CNK{2}{j}\>
\end{equation}
In order to prove that Equation ~\ref{eq:coeff_c_relation_1} holds, it is sufficient to show that two vectors that participate in the last expression are parallel. According to Equation ~\ref{eq:coeff_dL_dR}, coefficients of the weight functions $r(v)$ and $l(v)$ satisfy the following relations
\begin{equation}
\begin{array}{l}
\sum_{i=0}^2 r_i\kappa^i\CNK{2}{i}=0,\cr
\sum_{j=0}^2 l_j\kappa^j\CNK{2}{j}=0
\end{array}
\end{equation}
Using formulas $r_i=-\<\LL_i,\GG'-\GG\>$ and $l_j=\<\RR_j,\GG'-\GG\>$, $i,j=0,1,2$ (see Equation ~\ref{eq:lr_lambda_rho}), one gets 
\begin{equation}
\begin{array}{l}
0=-\<\sum_{i=0}^2\kappa^i\LL_i\CNK{2}{i},\GG'-\GG\>\cr
0=\<\sum_{j=0}^2\kappa^j\RR_j\CNK{2}{j},\GG'-\GG\>
\end{array}
\end{equation}
The last couple of equations means that vectors $\sum_{i=0}^2\kappa^i\LL_i\CNK{2}{i}$ and\\
$\sum_{j=0}^2\kappa^j\RR_j\CNK{2}{j}$ are parallel to $\GG'-\GG$ and to each other.

\vspace{0.2in}
\noindent
{\bf "CASE 2"}

\noindent
It will be shown first that the following two equations hold
\begin{equation}
\eta c_0-2\xi c_1+2\xi^2 c_3-\xi^2\eta c_4 = 0
\label{eq:xi_eta_1}
\end{equation}
\begin{equation}
\frac{1}{\xi^2}(\xi-2\eta^2)c_0+4\frac{1}{\xi}\eta c_1-6 c_2+
4\eta c_3+(\xi-2\eta^2)c_4 = 0
\label{eq:xi_eta_2}
\end{equation}
Like in the previous case, the proof uses expressions for coefficients of the weight function $c(v)$ in terms of the partial derivatives of in-plane parametrisations for the left and the right elements (see Equation ~\ref{eq:c_lambda_rho}).

Let $\alpha_i$ be the angle between $\GG'-\GG$ and $-\LL_j$ and $\beta_j$ be the angle between $\RR_j$ and $\GG'-\GG$ (both angles are measured in the counter clockwise direction). 
Then 
\begin{equation}
\begin{array}{l}
l_j = \<\RR_j,\GG'-\GG\>=|\RR_j||\GG'-\GG|sin\beta_j\cr
r_i = -\<\GG'-\GG,-\LL_i>=-|\LL_i||\GG'-\GG|sin\alpha_i\cr
\<\LL_i,\RR_j\>=|\LL_i||\RR_j|sin(\alpha_i+\beta_j)
\end{array}
\end{equation}

\vspace{0.15in}
\noindent
Proof of \underline{Equation ~\ref{eq:xi_eta_1}} makes use of two different representations of $\<\LL_1,\RR_1\>$. First, $\<\LL_1,\RR_1\>$ may be represented as follows
\begin{equation}
\hspace{-0.2in}
\begin{array}{ll}
\<\LL_1,\RR_1\>=
|\LL_1||\RR_1|(sin\alpha_1 cos\beta_1+cos\alpha_1 sin\beta_1)=\cr
\kappa^{(10)}
\left(|\LL_0||\RR_1|sin\alpha_0 cos\beta_1+
|\LL_1||\RR_0|cos\alpha_1sin\beta_0\right)=\cr
\kappa^{(10)}
\left(\<\LL_0,\RR_1\>\!+\!\<\LL_1,\RR_0\>\!-\!
|\LL_0||\RR_1|cos\alpha_0 sin\beta_1\!-\!
|\LL_1||\RR_0|sin\alpha_1 cos\beta_0\right)\!=\!\!\!\!\!\!\!\cr
\kappa^{(10)}\left(2 c_1-\kappa^{(10)}\<\LL_0,\RR_0\>\right)=
\kappa^{(10)}(2 c_1-\kappa^{(10)} c_0)
\end{array}
\label{eq:lambda_1_rho_1_first}
\end{equation}
In the same manner, it can be shown that 
\begin{equation}
\<\LL_1,\RR_1\>=\kappa^{(12)}(2 c_3-\kappa^{(12)}c_4)
\label{eq:lambda_1_rho_1_second}
\end{equation}
Equations ~\ref{eq:lambda_1_rho_1_second} and ~\ref{eq:lambda_1_rho_1_second} imply that
\begin{equation}
\kappa^{(10)}(2 c_1-\kappa^{(10)} c_0)=
\kappa^{(12)}(2 c_3-\kappa^{(12)} c_4)
\end{equation}
Substitution of $\kappa^{(10)}=\frac{\eta}{\xi}$ and $\kappa^{(12)}=\eta$ in the last equation completes the proof of Equation ~\ref{eq:xi_eta_1}.

\vspace{0.15in}
\noindent
\underline{Equation ~\ref{eq:xi_eta_2}} is proven in a similar manner. Coefficient $c_2$ of the weight function $c(v)$ has the following representation (see Equation ~\ref{eq:c_lambda_rho})
\begin{eqnarray*}
c_2 = \frac{1}{6}
\left(\<\LL_0,\RR_2\>+4\<\LL_1,\RR_1\>+\<\LL_2,\RR_0\>\right)
\end{eqnarray*}
Here
\begin{equation}
\begin{array}{ll}
\<\LL_0,\RR_2\>+\<\LL_2,\RR_0\>=&
|\LL_0||\RR_2|(sin\alpha_0 cos\beta_2+ cos\alpha_0 sin\beta_2)+\cr
&
|\LL_2||\RR_0|(sin\alpha_2 cos\beta_0+ cos\alpha_2 sin\beta_0)=\cr
&\kappa^{(02)}\<\LL_2,\RR_2\>+\kappa^{(20)}\<\LL_0,\RR_0\>=\cr
&\kappa^{(02)}c_4+\kappa^{(20)}c_0
\end{array}
\end{equation}
\begin{equation}
4\<\LL_1,\RR_1\>=2\kappa^{(10)}(2 c_1-\kappa^{(10)}c_0)+
2\kappa^{(12)}(2 c_3-\kappa^{(12)}c_4)
\end{equation}
It means that
\begin{equation}
\hspace{-0.45in}
\begin{array}{ll}
c_2\!=\!\frac{1}{6}\left[
\left(\kappa^{(20)}\!-\!
2\left[\kappa^{(10)}\right]^2\right)c_0\!+\!
4\kappa^{(10)}c_1+4\kappa^{(12)}c_3\!+\!
\left(\kappa^{(02)}\!-\!2\left[\kappa^{(12)}\right]^2\right)c_4
\right]
\end{array}
\end{equation}
Substitution of $\kappa^{(20)}=\frac{1}{\xi}$,
$\kappa^{(10)}=\frac{\eta}{\xi}$,
$\kappa^{(02)}=\xi$, $\kappa^{(12)}=\eta$ in the last formula completes the proof of Equation ~\ref{eq:xi_eta_2}.

\vspace{0.15in}
\noindent
Now it is easy to prove that \underline{Equations ~\ref{eq:coeff_c_relation_2_1},
~\ref{eq:coeff_c_relation_2_2_coeff} and
~\ref{eq:coeff_c_relation_2_2_deriv}} are satisfied.

\vspace{0.1in}
\noindent
{\bf "CASE 2.1"}

\noindent
Note, that $\omega_1\omega_2=\xi^2\neq 0$ and so in order to prove that Equation ~\ref{eq:coeff_c_relation_2_1} is satisfied for $\omega_1$ and $\omega_2$, it is sufficient to prove the equality 
\begin{equation}
\frac{1}{\omega_1}
\sum_{k=0}^4\left(-\sqrt{\omega_1}\right)^k\CNK{4}{k} c_k\pm
\frac{1}{\omega_2}
\sum_{k=0}^4\left(-\sqrt{\omega_2}\right)^k\CNK{4}{k} c_k=0
\label{eq:sum_mixed}
\end{equation}
Using the formulas for $\sqrt{\omega_{1,2}}$ (see Equation ~\ref{eq:sqrt_omega_1_2}), one can rewrite two summands of  Equation ~\ref{eq:sum_mixed} in the following form
\begin{equation}
\begin{array}{l}
\frac{1}{\omega_1}
\sum_{k=0}^4\left(-\sqrt{\omega_1}\right)^k\CNK{4}{k} c_k+
\frac{1}{\omega_2}
\sum_{k=0}^4\left(-\sqrt{\omega_2}\right)^k\CNK{4}{k} c_k=\cr
-2\left\{
\frac{1}{\xi^2}(\xi-2\eta^2)c_0+4\frac{1}{\xi}\eta c_1-6 c_2+
4\eta c_3+(\xi-2\eta^2)c_4
\right\}
\end{array}
\end{equation}
and
\begin{equation}
\begin{array}{l}
\frac{1}{\omega_1}
\sum_{k=0}^4\left(-\sqrt{\omega_1}\right)^k\CNK{4}{k} c_k-
\frac{1}{\omega_2}
\sum_{k=0}^4\left(-\sqrt{\omega_2}\right)^k\CNK{4}{k} c_k=\cr
-\frac{4}{\xi^2}\sqrt{\eta^2-\xi}
\left\{
\eta c_0-2\xi c_1+2\xi^2 c_3-\xi^2\eta c_4
\right\}
\end{array}
\end{equation}
The first expression is proportional to the expression given in Equation ~\ref{eq:xi_eta_2} and the second one is proportional to the expression given in Equation ~\ref{eq:xi_eta_1}, which implies that both of them are equal to zero.

\vspace{0.1in}
\noindent
{\bf "CASE 2.2"}

\noindent
In this case $\omega=\eta^2=\xi\neq 0$, $\sqrt\omega=\eta$ and Equations ~\ref{eq:coeff_c_relation_2_2_coeff} and
~\ref{eq:coeff_c_relation_2_2_deriv} have the following representations in terms of $\eta$
\begin{equation}
\begin{array}{l}
c_0-4\eta c_1+6\eta^2 c_2-4\eta^3 c_3+\eta^4 c_4 = 0\cr
4\eta\left\{-c_1+3\eta c_2-3\eta^2 c_3+\eta^3 c_4\right\}=0
\end{array}
\end{equation}
Using the equality $\eta^2=\xi$, it is easy to verify that the first expression is proportional to Equation ~\ref{eq:xi_eta_2} and the second one is proportional to 
$\eta^2$\{Equation ~\ref{eq:xi_eta_2}\}+
$\frac{1}{\eta}$\{Equation ~\ref{eq:xi_eta_1}\}, which implies that both of the expressions are equal to zero.

\nopagebreak 
\eop${}_{{\bf Auxiliary\ Lemma ~\ref{aux_lemma:coeff_c_relations}}}$

\vspace{0.32in}
\noindent
{\bf Proof of Auxiliary Lemma ~\ref{aux_lemma:hat_b_s_by_b_s}}

Proof of Auxiliary Lemma ~\ref{aux_lemma:hat_b_s_by_b_s} is quite straightforward, therefore only the basic steps of the proof are presented below.

\vspace{0.15in}
\noindent
{\bf "CASE 1"}

\noindent
From the assumption  that {\it "Eq(0)"}, {\it "Eq(1)"}, {\it "Eq($n'+1$)"} and {\it "Eq($n'+2$)"} (in the case of $n\ge 5$) being satisfied, definition of $\kappa$ (see Equation ~\ref{def:kappa})
and formulas for $\hat b_2$, $\hat b_3$, $\hat b{n'-1}$ and $\hat b_{n'}$ (see Equation ~\ref{eq:def_hat_b_bicubic}), it follows that
\begin{equation}
\kappa^2\hat b_2+\kappa^3\hat b_3 = 
b_0+\kappa b_1+\kappa^2 b_2+\kappa^3 b_3
\end{equation}
\begin{equation}
\hat b_{n'-1}+\kappa\hat b_{n'}=
b_{n'-1}+\kappa b_{n'}+\kappa^2 b_{n'+1}+\kappa^3 b_{n'+2}
\end{equation}
For the coefficients $\alpha_s$ ($s=2,\ldots,n'$) of non-trivial zero linear combination of the rows of the {\it "Restricted Middle"} system, the corresponding linear combination of the right-side coefficients can be written as 
\begin{equation}
\begin{array}{ll}
\sum_{s=2}^{n'}\alpha_s \hat b_s =&
\alpha_0 \sum_{s=2}^{n'}\kappa^s \hat b_s =\cr
&\alpha_0 \left\{
\sum_{s=0}^{n'+1}\kappa^s b_s +\kappa^{n'+2} b_{n'+2}\right\}=
\alpha_0 \sum_{s=0}^{n+2}\kappa^s b_s
\end{array}
\label{eq:sum_full_case_1}
\end{equation}
Equation ~\ref{eq:sum_full_case_1} is based on the fact that $\hat b_s = b_s$ for $s=4,\ldots,n'-2$ (see Equation ~\ref{eq:def_hat_b_bicubic}) and on the recursive formula for the coefficients of a non-trivial zero combinations $\alpha_{s+1}=\kappa\alpha_s$ for $s=2,\ldots,n'-1$
(see Equation ~\ref{eq:rec_coeff_case_1}).

\vspace{0.15in}
\noindent
{\bf "CASE 2"}

\noindent
Similarly to the previous case, the main difficulty in the proof is connected to the expression of the first and the last couples of $\hat b_s$ in terms of $b_s$.

The assumption  that equations 
{\it "Eq(s)"} for $s=0,1,n'+1$ and for $s=n'+2$ in case of 
$n\ge 5$ are satisfied leads to the following relations
\begin{equation}
\begin{array}{ll}
\hat b_2 =&
b_2+\frac{4\eta^2-\xi}{\xi^2}b_0-2\frac{\eta}{\xi}b_1\cr
\hat b_3 =&
b_3+\frac{2\eta}{\xi^2}b_0-\frac{1}{\xi}b_1\cr
\hat b_{n'-1}=&
b_{n'-1}-\xi b_{n'+1}+2\xi\eta b_{n'+2}\cr
\hat b_{n'}=&
b_{n'}-2\eta b_{n'+1}+(4\eta^2-\xi)b_{n'+2}
\end{array}
\label{eq:hat_b_ends_case_2}
\end{equation}
The first couple of relations clearly implies that
\begin{equation}
(\alpha_2\ \alpha_3)
\VV{\hat b_2}{\hat b_3}=
(\alpha_0\ \alpha_1)
\left[
\VV{b_0}{b_1}+G\VV{b_2}{b_3}
\right]
\label{eq:first_couple}
\end{equation}
If $n'$ is odd then the last couple of relations (see Equation ~\ref{eq:hat_b_ends_case_2}) can be rewritten in the form
\begin{equation}
\VV{\hat b_{n'-1}}{\hat b_{n'}}=
\VV{b_{n'-1}}{b_{n'}}+G\VV{b_{n'+1}}{b_{n'+2}}
\label{eq:last_couple_odd_n_tag}
\end{equation}
If $n'$ is even then it can be shown that
\begin{equation}
\begin{array}{l}
(\alpha_{n'-2}\ \alpha_{n'-1})
\VV{\hat b_{n'-2}}{\hat b_{n'-1}}+\alpha_{n'}\hat\beta_{n'}=\cr
(\alpha_{n'-2}\ \alpha_{n'-1})\left[
\VV{b_{n'-2}}{b_{n'-1}}+
G\VV{b_{n'}}{b_{n'+1}}\right]
\end{array}
\label{eq:last_couple_even_n_tag}
\end{equation}
Equation ~\ref{eq:last_couple_even_n_tag} makes use of the fact that the satisfaction of {\it "sumC-equation"} leads to equality $b_{n'+2}=0$. Note, that it is the only step of the proof where the satisfaction of {\it "sumC-equation"} is required. It implies that
requirement of the satisfaction of {\it "sumC-equation"} in the formulation of Auxiliary Lemma ~\ref{aux_lemma:hat_b_s_by_b_s}
is necessary  only in case of even $n\ge 6$.

Now Equation ~\ref{eq:hat_b_s_by_b_s_2} follows from Equations ~\ref{eq:first_couple}, ~\ref{eq:last_couple_odd_n_tag}, ~\ref{eq:last_couple_even_n_tag} and relation $(\alpha_{2s}\ \alpha_{2s+1})=(\alpha_0\ \alpha_1)G^s$, $s=2,\ldots,n'-3$ (see Equation ~\ref{eq:G_explanation}).

\nopagebreak 
\eop${}_{{\bf Auxiliary\ Lemma ~\ref{aux_lemma:hat_b_s_by_b_s}}}$

\vspace{0.32in}
\noindent
{\bf Proof of Auxiliary Lemma ~\ref{aux_lemma:consisitency_b_s}}

Proof of Auxiliary Lemma ~\ref{aux_lemma:consisitency_b_s} is based on the first statement of the algebraic Auxiliary Lemma ~\ref{aux_lemma:algebraic}, which provides the necessary and sufficient conditions for the consistency of a system in terms of the linear combinations of rows and right-side coefficients.

\vspace{0.15in}
\noindent
{\bf "CASE 1"}

\noindent
Necessity and sufficiency of Equation ~\ref{eq:consisitency_b_s_1}
immediately follows from Equation 
~\ref{eq:hat_b_s_by_b_s_1}.

\vspace{0.15in}
\noindent
{\bf "CASE 2"}

\noindent
According to Auxiliary Lemma ~\ref{aux_lemma:hat_b_s_by_b_s}, it is clear that the {\it "Restricted Middle"} system is consistent if and only if equation
\begin{equation}
\sum_{s=0}^{\lceil n'/2 \rceil} G^s \VV{b_{2s}}{b_{2s+1}}=
\VV{0}{0}
\label{eq:consistency_G_powers}
\end{equation}
is satisfied.

\vspace{0.1in}
\noindent
{\bf "CASE 2.1"}

\noindent
Proof of the Auxiliary Lemma in {\bf "CASE 2.1"} makes use
of the Jordan form of matrix $G$ given in Equation ~\ref{eq:G_structure_2_1}. Using the formula
\begin{equation}
(v_1\ v_2)^{-1}=
\frac{1}{\xi(\sqrt{\omega_2}-\sqrt{\omega_1})}
\MM{\sqrt{\omega_2}}{-\sqrt{\omega_1\omega_2}}
{-\sqrt{\omega_1}}{\sqrt{\omega_1\omega_2}}
\end{equation}
one sees that Equation ~\ref{eq:consistency_G_powers} is equivalent to the system
\begin{equation}
\hspace{-0.2in}
\begin{array}{ll}
\VV{0}{0}=&
\sum\limits_{s=0}^{\lceil n'/2 \rceil}
\MM{\omega_1^s}{0}{0}{\omega_2^s}
\MM{\sqrt{\omega_2}}{-\sqrt{\omega_1\omega_2}}
{-\sqrt{\omega_1}}{\sqrt{\omega_1\omega_2}}
\VV{b_{2s}}{b_{2s+1}}=\cr
&
\sum\limits_{s=0}^{\lceil n'/2 \rceil}
\VV
{\sqrt{\omega_2}\left(
(-\sqrt{\omega_1})^{2s}b_{2s}+
(-\sqrt{\omega_1})^{2s+1}b_{2s+1}\right)}
{-\sqrt{\omega_1}\left(
(-\sqrt{\omega_2})^{2s}b_{2s}+
(-\sqrt{\omega_2})^{2s+1}b_{2s+1}\right)}
\end{array}
\end{equation}
The last system in turn is clearly equivalent to the following equation, which should be satisfied for $\omega_1$ and $\omega_2$
system of equations  
\begin{equation}
\sum_{s=0}^{2\lceil n'/2 \rceil +1}
\left(-\sqrt{\omega_{1,2}}\right)^s b_s = 0
\end{equation}
It remains to show, that the upper bound of summation can be changed from $2\lceil n'/2 \rceil+1$ to $n+2$. The substitution of the upper bound is valid because $b_7=0$ for $n=4$ and $b_{n+2}=0$ for an even $n>4$ (according to the assumption that {\it "sumC-equation"} is satisfied).

\vspace{0.1in}
\noindent
{\bf "CASE 2.2"}

\noindent
Proof of the Auxiliary Lemma in {\bf "CASE 2.2"} is very similar to the proof in {\bf "CASE 2.1"}. In the current case the Jordan form of the matrix $G$ is given by Equation
~\ref{eq:G_structure_2_2}. 
Using the formula 
\begin{equation}
(v\ u)^{-1}=
\frac{1}{2\sqrt\omega}
\MM
{\frac{1}{4\sqrt\omega}}{\frac{1}{4}}
{-\sqrt\omega}{\omega}
\end{equation}
one sees that Equation ~\ref{eq:consistency_G_powers} is equivalent to the system
\begin{equation}
\hspace{-0.45in}
\begin{array}{ll}
\VV{0}{0}=&
\sum\limits_{s=0}^{\lceil n'/2 \rceil}
\MM
{\omega^s}{s\omega^{s-1}}{0}{\omega^s}
\MM
{\frac{1}{4\sqrt\omega}}{\frac{1}{4}}
{-\sqrt\omega}{\omega}
\VV
{b_{2s}}{b_{2s+1}}=\cr
&
\left(
\begin{array}{l}
\ \ \frac{1}{4\sqrt\omega}
\sum_{s=0}^{2\lceil n'/2 \rceil +1}(-\sqrt\omega)^s b_s-
\frac{1}{2}
\sum_{s=0}^{2\lceil n'/2 \rceil +1}(-\sqrt\omega)^{s-1} s b_s
\cr
-\sqrt\omega
\sum_{s=0}^{2\lceil n'/2 \rceil +1}(-\sqrt\omega)^s b_s
\end{array}
\right)
\end{array}
\end{equation} 
Replacement of the upper bound of summation by $n+2$ (which is valid due to the same reasons as in {\bf "CASE 2.1"}) allows to conclude that  Equations ~\ref{eq:consisitency_b_s_2_2_f} and
~\ref{eq:consisitency_b_s_2_2_s} provide the necessary and sufficient conditions for the consistency of the {\it "Restricted Middle"} system.

\nopagebreak 
\eop${}_{{\bf Auxiliary\ Lemma ~\ref{aux_lemma:consisitency_b_s}}}$

\vspace{0.32in}
\noindent
{\bf Proof of Auxiliary Lemma ~\ref{aux_lemma:max_deg_l_r_classification}}
 
Auxiliary Lemma ~\ref{aux_lemma:max_deg_l_r_classification}
makes use of the following pair of implications

\vspace{0.1in}
\hspace{-0.35in}
\begin{tabular}{ll}
{\it "Impl(1)"}\ \ \ \ \ \ &
$max\_deg(l,r)\le 1$ implies that conditions of either 
{\bf "CASE\ 1"} or\\& {\bf "CASE\ 2"} are satisfied.\\
& \\
{\it "Impl(2)"}\ \ \ \ \ \ &
$max\_deg(l,r) = 0$ implies that conditions of case {\bf "CASE\ 2.2.b"}\\& are satisfied.
\end{tabular}

\vspace{0.15in}
\noindent
A short proof of {\it "Impl(1)"} and {\it "Impl(2)"} is given below.

\noindent
\underline{\it "Impl(1)"}

\noindent
Condition $max\_deg(l,r)\le 1$ means that $l_2^{(power)}=l_0-2l_1+l_2=0$, $r_2^{(power)}=r_0-2r_1+r_2=0$. Therefore $r_0 l_2^{(power)}=l_0 r_2^{(power)}$ and equality $g^{(01)}=\frac{1}{2}g^{(02)}$ holds. In the same manner multiplication of $l_2^{(power)}$ and $r_2^{(power)}$ by $r_2$ and $l_2$ respectively leads to equality $g^{(12)}=\frac{1}{2}g^{(02)}$. It means that the condition $max\_deg(l,r)\le 1$ implies that
\begin{equation}
g^{(01)}=g^{(12)}=\frac{1}{2}g^{(02)}
\label{eq:g_ij_relations_1}
\end{equation}
Equation ~\ref{eq:g_ij_relations_1} holds if either $g^{(01)}=g^{(12)}=g^{(02)}=0$ which corresponds to {\bf "CASE 2"} or $g^{(01)}$,$g^{(12)}$,$g^{(02)}$ simultaneously differ from zero and $\left(g^{(02)}\right)^2=4g^{(01)}g^{(12)}$ which corresponds to {\bf "CASE 1"}.

\vspace{0.1in}
\noindent
\underline{\it "Impl(2)"}

\noindent
Condition $max\_deg(l,r)=0$ means that $l_2^{(power)}=r_2^{(power)}=l_1^{(power)}=r_1^{(power)}=0$ and so 
\begin{equation}
l_0=l_1=l_2,\ \ \ \ r_0=r_1=r_2
\end{equation}
It implies that $g^{(01)}=g^{(12)}=g^{(02)}=0$ and that $\xi=\kappa^{(02)}=1$, $\eta=\kappa^{(12)}=1$, which corresponds to {\bf "CASE 2.2.b"}

\vspace{0.2in}

The remaining part of the proof presents computations of the exact value of $max\_deg(l,r)$ in every one of the possible cases.

\vspace{0.15in}
\noindent
{\underline{\bf "CASE 0"}} 
In this case $max\_deg(l,r)=2$ as it immediately follows from {\it "Impl(1)"}.

\vspace{0.15in}
\noindent
{\underline{\bf "CASE 1"}} 
According to {\it "Impl(2)"}, $max\_deg(l,r)\ge 1$. It remains to show that $max\_deg(l,r)\neq 1$ in {\bf "CASE 1.a"} and $max\_deg(l,r)\neq 2$ in {\bf "CASE 1.b"}.

\begin{itemize}
\item[]
\underline{\bf "CASE 1.a"}
On the contrary, let $max\_deg(l,r)=1$. As it follows from the proof of {\it "Impl(1)"}, in this case $g^{(01)}=g^{(12)}=\frac{1}{2}g^{(02)}$. Conditions of {\bf "CASE 1"} imply that 
$g^{(01)},g^{(12)},g^{(02)}\neq 0$. It means that constant $\kappa$ is correctly defined and by the definition is equal to 
$\kappa=-\frac{2 g^{(01)}}{g^{(02)}}=-1$, which contradicts the conditions of {\bf "CASE 1.a"}. 

\item[]
\underline{\bf "CASE 1.b"}
On the contrary, let  $max\_deg(l,r)=2$. 
From the conditions of {\bf "CASE 1.b"}, it follows that 
$2 g^{(01)}=g^{(02)}$, $2 g^{(12)}=g^{(02)}$. Addition of 
$l_0 r_0$ and $l_2 r_2$ to the both parts of the first and second equalities leads to the conclusion that
\begin{equation}
l_0 r_2^{(power)}=r_0 l_2^{(power)},\ \ \ 
l_2 r_2^{(power)}=r_2 l_2^{(power)}
\label{eq:aux1_case_1_2}
\end{equation}
The last couple of equations means that both $l_2^{(power)}$ and $r_2^{(power)}$ are not equal to zero (otherwise from Equation ~\ref{eq:aux1_case_1_2} it follows that $l_2^{(power)}=r_2^{(power)}=0$ which contradicts the assumption that $max\_deg(l,r)=2$) and that
\begin{equation}
\frac{l_0}{r_0}=
\frac{l_2}{r_2}=\frac{l_2^{(power)}}{r_2^{(power)}}
\label{eq:aux2_case_1_2}
\end{equation}
Equation ~\ref{eq:aux2_case_1_2} implies that $g^{(02)}=0$ which contradicts the conditions of {\bf "CASE 1"}.
\end{itemize}

\noindent
\underline{\bf "CASE 2"}

\begin{itemize}
\item[]
\underline{\bf "CASE 2.1"}
Let conditions of {\bf "CASE 2.1"} be satisfied. From the formula $\sqrt{\omega_{1,2}}=\eta\pm\sqrt{\eta^2-\xi}$ 
(see Equation ~\ref{eq:sqrt_omega_1_2}) it follows that conditions of subcase {\bf "CASE 2.1.a"} are equivalent to inequality
\begin{equation}
1-2\eta+\xi\neq 0
\label{eq:equiv_case_2_1_a}
\end{equation}
and conditions of subcase {\bf "CASE 2.1.b"} are equivalent the equality
\begin{equation}
1-2\eta+\xi = 0
\label{eq:equiv_case_2_1_b}
\end{equation}

Using Equations ~\ref{eq:equiv_case_2_1_a}, ~\ref{eq:equiv_case_2_1_b} and the explicit formulas for $\xi$ and $\eta$ (see Equation ~\ref{eq:def_ksi_eta}), one sees that conditions of subcase {\bf "CASE 2.1.a"} are satisfied if and only if the following inequality holds (in addition to the general conditions of {\bf "CASE 2.1"})
\begin{equation}
\frac{l_2^{(power)}}{l_2}=\frac{r_2^{(power)}}{r_2}\neq 0
\label{eq:max_deg_characteristic_2_1_a}
\end{equation}
and conditions of subcase {\bf "CASE 2.1.b"} are satisfied if and only if the following equality holds (in addition to the general conditions of {\bf "CASE 2.1"}) 
\begin{equation}
\frac{l_2^{(power)}}{l_2}=\frac{r_2^{(power)}}{r_2} = 0
\label{eq:max_deg_characteristic_2_1_b}
\end{equation}
Now it is easy to compute the value of $max\_deg(l,r)$ for both subcases of {\bf "CASE 2.1"}.

\begin{itemize}
\item[]
\underline{\bf "CASE 2.1.a"}
Equation ~\ref{eq:max_deg_characteristic_2_1_a} clearly implies that $max\_deg(l,r)=2$.

\item[]
\underline{\bf "CASE 2.1.b"}
Equation ~\ref{eq:max_deg_characteristic_2_1_b} implies that $max\_deg(l,r)\le 1$. According to {\it "Impl(2)"}, $max\_deg(l,r)\neq 0$. It means that $max\_deg(l,r)=1$.
\end{itemize}

\item[]
\underline{\bf "CASE 2.2"}

\begin{itemize}
\item[]
\underline{\bf "CASE 2.2.a"}
On the contrary, let $max\_deg(l,r)\le 1$.
It implies in particular that 
\begin{equation}
\begin{array}{ll}
0=&\frac{l_2^{(power)}}{l_2}=\frac{l_0-2 l_1+l_2}{l_2}=\cr
&\xi-2\eta+1=\eta^2-2\eta+1=(\eta-1)^2
\end{array}
\end{equation}
The last equality means that $\sqrt\omega=\eta=1$,
which contradicts the conditions of subcase {\bf "CASE 2.2.a"}.

\item[]
\underline{\bf "CASE 2.2.b"}
In this case $\sqrt\omega=\eta=1$ and so $l_1=l_2$, $r_1=r_2$. In addition, $\xi=\eta^2=1$ and so $l_0=l_2$, $r_0=r_2$. It clearly implies that $max\_deg(l,r)=0$.
\end{itemize}
\end{itemize}

\nopagebreak 
\eop${}_{{\bf Auxiliary\ Lemma ~\ref{aux_lemma:max_deg_l_r_classification}}}$

\vspace{0.32in}
\noindent
{\bf Proof of Auxiliary Lemma 
~\ref{aux_lemma:middle_consisitency_c_power}}

The {\it "Middle"} system of equations is composed of the {\it "Restricted Middle"} system and {\it "sumC-equation"}. It means that the consistency of the {\it "Middle"} system is equivalent to the couple of the following conditions
\begin{itemize}
\item
{\it "sumC-equation"} is satisfied.
\item
The necessary and sufficient conditions for the consistency of the {\it "Restricted Middle"} system, given in Auxiliary Lemma ~\ref{aux_lemma:consisitency_b_s}, hold.
\end{itemize}
Now Auxiliary Lemma ~\ref{aux_lemma:middle_consisitency_c_power} can be proven by a simple analysis of all possible subcases. Below the proof in {\bf "CASE 2.1.a"} is given; the rest of the subcases can be analysed in a similar manner.

\vspace{0.1in}
\noindent
{\underline{\bf "CASE 2.1.a"}
According to Auxiliary Lemma ~\ref{aux_lemma:max_deg_l_r_classification}, $max\_deg(l,r)=2$. Auxiliary Lemmas ~\ref{aux_lemma:consisitency_b_s} and ~\ref{aux_lemma:phi_formal_expansion} imply that if
{\it "sumC-equation"} is satisfied then the {\it "Restricted Middle"} system of equations is consistent if and only if 
\begin{equation}
\begin{array}{l}
0=\sum_{s=0}^{n+2}(-\sqrt{\omega_{1,2}})^s b_s=
"sumC\!-\!equation"
\left[-\frac{(\sqrt{\omega_{1,2}})^{n+3}}
{1-\sqrt\omega_{1,2}}\right]+\cr
\frac{1}{1-\sqrt{\omega_{1,2}}}
\left[\sum_{i=0}^{n-1}
(-\sqrt{\omega_{1,2}})^i\Cnk{n-1}{i}\Delta C_i\right]
\left[\sum_{j=0}^4(-\sqrt{\omega_{1,2}})^j \Cnk{4}{j} c_j\right]
\end{array}
\label{eq:example_2_1_a}
\end{equation}
From the relation between coefficients of the weight function $c(v)$ given in Auxiliary Lemma ~\ref{aux_lemma:coeff_c_relations} (see Equation ~\ref{eq:coeff_c_relation_2_1}) it follows that Equation ~\ref{eq:example_2_1_a} holds if and only if {\it "sumC-equation"} is satisfied. It means that the {\it "Middle"} system of equations is consistent if and only if {\it "sumC-equation"} is satisfied. According to Lemma ~\ref{lemma:sumC_expansion}, one sees that the {\it "Middle"} system is consistent if and only if $"C-equation"c_4^{(power)}=0$.

\nopagebreak
\eop${}_{{\bf Auxiliary\ Lemma ~\ref{aux_lemma:middle_consisitency_c_power}}}$

\vspace{0.32in}
\noindent
{\bf Proof of the correctness of Algorithm ~\ref{algorithm:mds_global_bicubic}}

\paragraph{"Stage 1"}
According to Lemma ~\ref{lemma:Eq0_Eq1_bicubic}, at every inner vertex, local templates for classification of $V$,$E$,$D$,$T$-type control points remains unchanged with respect to the case of global bilinear in-plane parametrisation $\T\Pi^{(bilinear)}$. Therefore Algorithms from Section ~\ref{sect:global_mds_bilinear} can be reused.

{\it "Stage 1"} of Algorithm ~\ref{algorithm:mds_global_bicubic} includes {\it "Stage 1"} and {\it "Stage 2"} of algorithm for construction of global MDS in case of global bilinear in-plane parametrisation $\T\Pi^{(bilinear)}$ (Algorithm ~\ref{algorithm:mds_global_deg_ge_5_bilinear}), applied to all inner vertices.

Global classification of $V$,$E$-type control points succeeds for any $n\ge 4$. For $n\ge 5$, application of local templates for classification of $D$,$T$-types control points never leads to a contradictions (see Subsection ~\ref{subsect:global_mds_n_ge_5_bilinear}). For $n=4$, global classification of $D$,$T$-type control points succeeds in the current case according to Theorem ~\ref{theorem:sufficient_DT_classification_bilinear}, because no $D$-relevant boundary vertices are involved.

\paragraph{"Stage 2"}
According to Lemma ~\ref{lemma:unchanged_middle_bicubic}, local templates for classification of the {\it middle} control points for edges with two inner vertices remain unchanged with respect to the case of global bilinear in-plane parametrisation $\T\Pi^{(bilinear)}$. Therefore {\it "Stage 2"} of Algorithm ~\ref{algorithm:mds_global_bicubic} is precisely {\it "Stage 3"} of Algorithm ~\ref{algorithm:mds_global_deg_ge_5_bilinear}, applied to all edges with two inner vertex, which is known to succeed for any $n\ge 4$.

\paragraph{"Stage 3"}
It is the only stage, which requires modifications with respect to construction of MDS in the case of global bilinear in-plane parametrisation $\T\Pi^{(bilinear)}$. Classification of the control points at {\it "Stage 3"} is made locally, traversal order of the edges with one boundary vertex is not important.

For every edge with one boundary vertex, the correctness of classification follows from the fact, that after application of {\it "Step 1"} all initial assumptions of Lemma ~\ref{lemma:C_classification} and then (after classification of the "central" {\it middle} control points) of Theorem ~\ref{theorem:LR_classification} are satisfied.

\nopagebreak
\eop${}_{{\bf Algorithm ~\ref{algorithm:mds_global_bicubic}}}$

\eject

\begin{figure}[!pb]
\centering
\begin{narrow}{-0.3in}{-0.3in}
\begin{minipage}[]{0.33\linewidth}
\subfigure[]
{
\includegraphics[clip,totalheight=2.0in]{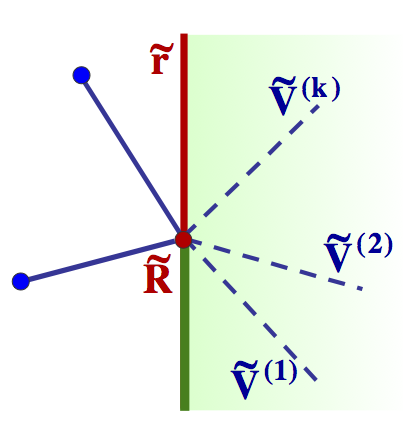}
\label{fig:fig24a}
}
\end{minipage}
\begin{minipage}[]{0.25\linewidth}
\subfigure[]
{
\includegraphics[clip,totalheight=2.0in]{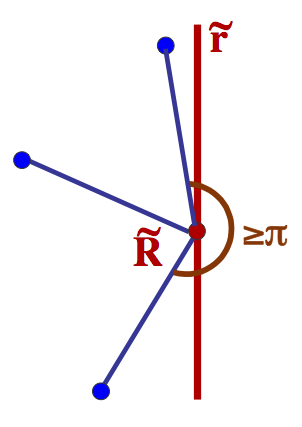}
\label{fig:fig24b}
}
\end{minipage}
\begin{minipage}[]{0.33\linewidth}
\subfigure[]
{
\includegraphics[clip,totalheight=2.0in]{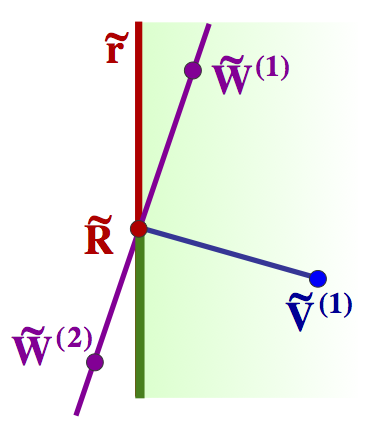}
\label{fig:fig24c}
}
\end{minipage}\\
\centering
\begin{minipage}[]{0.4\linewidth}
\subfigure[]
{
\includegraphics[clip,totalheight=2.0in]{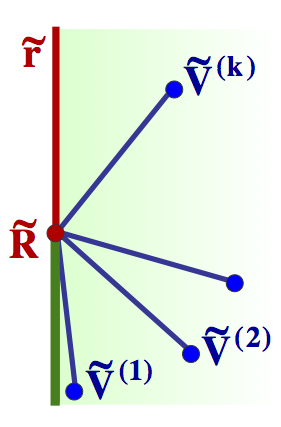}
\label{fig:fig24d}
}
\end{minipage}
\begin{minipage}[]{0.4\linewidth}
\subfigure[]
{
\includegraphics[clip,totalheight=2.0in]{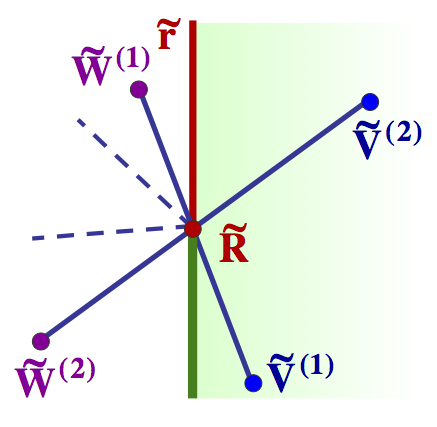}
\label{fig:fig24e}
}
\end{minipage}
\end{narrow}
\vspace{-0.2in}
\caption{An illustrations to proof of Theorem ~\ref{theorem:sufficient_DT_classification_bilinear}.}
\label{fig:fig24}
\end{figure}

\FloatBarrier

\newpage
\section{Computation of Thin Plate energy for bilinear in-plane parametrisation of mesh element}
\label{sect:energy_example_bilinear}

The following example demonstrates the linear form of energy functional in the case of the Thin Plate problem (see Subsection ~\ref{subsect:tp_definition} for the precise definition) and presents explicit computational formulas for the case of bilinear in-plane parametrisation. 

The current Section analyses a single mesh element and therefore it is possible to use slightly different notations than in the rest of the work - $(X,Y,Z)$ instead of $(P_{X},P_{Y},P)$, $J(u,v)$ instead of $J^{(\T P)}(u,v)$ - in order to make formulas more compact and clear. In addition, the superscript of Bernstein polynomials will sometimes be omitted.

Let $\T\Pi$ be a fixed global regular parametrisation,
$\bar\Psi\in\FUN{n}(\T\Pi)$ and $\bar P=(\T P,Z)=\bar\Psi|_{\T p}$ be the restriction of $\bar\Psi$ on some mesh element $\T p\in\T{\cal Q}$.  In-plane parametrisation of the element $\T P=(X,Y)=\T\Pi|_{\T p}$ is fixed.

Energy functional has the following form in terms of $Z$-components of the control points 
\begin{equation}
{\cal E} =\sum_{i,j=0}^n\sum_{k,l=0}^n Z_{i,j}Z_{k,l}a_{ij}^{kl}-\sum_{i,j}^n Z_{ij}b_{ij}. 
\label{eq:energy_quadratic_form_sum}
\end{equation}
where
\begin{equation}
\begin{array}{ll}
a_{ij}^{kl} = & D\int\int_{\tilde p} 
\left(\DT{B_{ij}}{X}+\DT{B_{ij}}{Y}\right)
\left(\DT{B_{kl}}{X}+\DT{B_{kl}}{Y}\right)\cr
              & -(1-\nu)
\left(
\DT{B_{ij}}{X}\DT{B_{kl}}{Y}+
\DT{B_{ij}}{Y}\DT{B_{kl}}{X}
-2\DD{B_{ij}}{X}{Y}\DD{B_{kl}}{X}{Y}
\right) 
dXdY \cr
b_{ij} =      & \int\int_{\tilde p} fB_{ij} dXdY.
\end{array}
\label{eq:energy_matrix_components}
\end{equation}
In  matrix form one gets:
\begin{equation}
{\cal E} = {\cal Z}^T {\cal A} {\cal Z} - {\cal Z}^T {\cal B}
\label{eq:energy_quadratic_form_matrix}
\end{equation}
where
${\cal Z} = (Z_{00},\ldots,Z_{nn})^T$ - vector of $Z$-components of all the  control points of the
patch; symmetric $(n+1)^2\times(n+1)^2$ square matrix ${\cal A}$ is built from components $a_{ij}^{kl}$ and vector ${\cal B} = (b_{00},\ldots,b_{nn})^T$.

Elements of matrix $\cal A$ and vector $\cal B$ can be computed by representing the integrals with respect to variables $(u,v)$ and then by applying any appropriate numerical method. Lemma
~\ref{lemma:partial_derivatives_general} contains some useful general formulas for partial derivatives of $u$, $v$ and Bernstein polynomials with respect to $X$ and $Y$.

\begin{lemma}
Let $\tilde P(u,v)=(X(u,v),Y(u,v))$ be a regular, twice differentiable parametrisation of a planar mesh element and 
\begin{equation}
J(u,v) = \left(\begin{array}{cc} 
\D{X}{u} & \D{Y}{u} \cr
\D{X}{v} & \D{Y}{v} \cr
\end{array}\right)
\end{equation}
is Jacobian corresponding to the parametrisation. Then

\vspace{0.1in}
\noindent
{\bf (1)\ \ }
The first-order and the second-order partial derivatives of $u$ and $v$ with respect to $X$ and $Y$ can be computed according to the following formulas
\begin{equation}
\begin{array}{l}
\MM{\D{u}{X}}{\D{v}{X}}{\D{u}{Y}}{\D{v}{Y}}=J^{-1}(u,v)=
\frac{1}{det(J(u,v))}
\MM{\D{Y}{v}}{-\D{Y}{u}}{-\D{X}{v}}{\D{X}{u}}
\cr
\MM{\DT{u}{X}}{\DT{v}{X}}{\DD{u}{X}{Y}}{\DD{v}{X}{Y}}=\D{J^{-1}(u,v)}{X}=
-J^{-1}(\D{J}{u}\D{u}{X}+\D{J}{v}\D{v}{X})J^{-1}
\cr
\MM{\DD{u}{X}{Y}}{\DD{v}{X}{Y}}{\DT{u}{Y}}{\DT{v}{Y}}=\D{J^{-1}(u,v)}{Y}=
-J^{-1}(\D{J}{u}\D{u}{Y}+\D{J}{v}\D{v}{Y})J^{-1}
\end{array}
\label{eq:partial_uv_xy_general}
\end{equation}

\vspace{0.1in}
\noindent
{\bf (2)\ \ }
The first-order and the second-order partial derivatives of Bernstein polynomial $B(u,v)$ can be computed according to the following formulas
\begin{equation}
\hspace{-0.45in}
\begin{array}{l}
\D{B}{X}=\D{B}{u}\D{u}{X}+\D{B}{v}\D{v}{X}
\cr
\D{B}{Y}=\D{B}{u}\D{u}{Y}+\D{B}{v}\D{v}{Y}
\cr
\DT{B}{X}=\DT{B}{u}\left(\D{u}{X}\right)^2+2\DD{B}{u}{v}\D{u}{X}\D{v}{X}+\DT{B}{v}\left(\D{v}{X}\right)^2+
\D{B}{u}\DT{u}{X}+\D{B}{v}\DT{v}{X}
\cr
\DD{B}{X}{Y}\!=\!\DT{B}{u}\D{u}{X}\D{u}{Y}\!\!+\!
2\DD{B}{u}{v}\!\left(\D{u}{X}\D{v}{Y}\!\!+\!
\D{u}{Y}\D{v}{X}\right)\!\!+\!
\DT{B}{v}\D{v}{X}\D{v}{Y}\!+\!
\D{B}{u}\DD{u}{X}{Y}\!\!+\!\D{B}{v}\DD{v}{X}{Y}\!\!\!\!\!\!\!
\cr
\DT{B}{Y}=\DT{B}{u}\left(\D{u}{Y}\right)^2+2\DD{B}{u}{v}\D{u}{Y}\D{v}{Y}+\DT{B}{v}\left(\D{v}{Y}\right)^2+
\D{B}{u}\DT{u}{Y}+\D{B}{v}\DT{v}{Y}
\end{array}
\label{eq:partial_bernstein_xy_general}
\end{equation}
\label{lemma:partial_derivatives_general}
\end{lemma}

In the case of bilinear in-plane parametrisation $\T\Pi^{(bilinear)}$ (see Section ~\ref{sect:parametrisation_bilinear}) the following formulas take place. For a planar quadrilateral element with vertices $\T{A},\T{B},\T{C},\T{D}$ (see Figure ~\ref{fig:fig8}),  determinant of Jacobian $det(J(u,v))$ may be written as 
\vspace{-0.05in}
\begin{equation}
det(J(u,v))\!=\!\<\T{B}\!-\!\T{A},\T{C}\!-\!\T{D}\>u+
\<\T{C}\!-\!\T{B},\T{D}\!-\!\T{A}\>v+
\<\T{B}\!-\!\T{A},\T{D}\!-\!\T{A}\>
\label{eq:det_jacobian_bilinear}
\end{equation}
\vspace{-0.05in}
$det(J(u,v))$ is a linear function in terms of $u$ and $v$ in a general case and constant in the case of a parallelogram. 
It is important that $det(J(u,v))$ is a linear and not a bilinear function, because $det(J(u,v))$ participates in expressions for partial derivatives and its order is essential for the choice of an appropriate method for numerical or exact integration
when the coefficients of energy matrix are computed. 

Based on Equation ~\ref{eq:det_jacobian_bilinear},
Lemma ~\ref{lemma:partial_derivatives_for_bilinear_parametrisation}
presents the explicit formulas for the first-order and the second-order partial derivatives of $u$ and $v$ with respect to $X$ and $Y$.

\begin{lemma}
Let a convex planar quadrilateral element have bilinear in-plane parametrisation $\tilde P(u,v)=(X(u,v),Y(u,v))$, $J(u,v)$ be Jacobian corresponding to the parametrisation and
$\T{T}=\T t(\T{A},\T{B},\T{C},\T{D})=\T{A}-\T{B}+\T{C}-\T{D}$ be the twist characteristic of the element (see Subsection 
~\ref{subsect:def_vertices_edges_twists_planar}).
Then the first-order and the second-order partial derivatives of $u$ and $v$ with respect to $X$ and $Y$ can be computed according to the following formulas 
\begin{equation}
\hspace{-0.4in}
\MM{\D{u}{X}}{\D{v}{X}}{\D{u}{Y}}{\D{v}{Y}}=
\frac{1}{det(J(u,v))}
{
\footnotesize
\MM{T_{Y} u+(D_{Y}-A_{Y})}{-T_{Y} v-(B_{Y}-A_{Y})}{-T_{X} u-(D_{X}-A_{X})}{T_{X} v+(B_{X}-A_{X})}
}
\end{equation}
\begin{equation}
\hspace{-0.4in}
\left(\begin{array}{cc}
\DT{u}{X}    & \DT{v}{X}\cr
\DD{u}{X}{Y} & \DD{v}{X}{Y}\cr
\DT{u}{Y}    & \DT{v}{Y}
\end{array}\right) =
-\left(\begin{array}{c}
2 \D{u}{X}\D{v}{X}\cr
\D{u}{X}\D{v}{Y}+\D{u}{Y}\D{v}{X}\cr
2 \D{u}{Y}\D{v}{Y}
\end{array}\right) 
(T_{X},T_{Y})
\MM{\D{u}{X}}{\D{v}{X}}{\D{u}{Y}}{\D{v}{Y}}
\end{equation}
\label{lemma:partial_derivatives_for_bilinear_parametrisation}
\end{lemma}
Technical Lemma ~\ref{lemma:partial_derivatives_for_bilinear_parametrisation}
allows concluding that the first-order partial derivatives are the linear rational functions and the second-order partial derivatives are cubic rational functions, where the denominator changes more rapidly when the quadrilateral is less similar to a parallelogram.

\eject
\section{From MDS to solution of the quadratic minimisation problem}
\label{sect:mds_to_algebraic_solution}

Let ${\cal Z}_{all}$ denote $Z$-components of all the  control points of all the  patches, ordered so that
\begin{equation}
{\cal Z}_{all} = (Z^{(1)}_{00}, \ldots, Z^{(1)}_{nn},\ldots, 
Z^{(p_{max})}_{00},\ldots,Z^{(p_{max})}_{nn})^T
\label{Z_all}
\end{equation}
where $Z^{(p)}_{ij}$ denotes a control point belonging to the patch with order number $p=1,
\ldots,p_{max}$. 
Here no dependencies between control points of the different patches are assumed. Let ${\cal A}^{(p)}$ and ${\cal B}^{(p)}$ ($p=1,\ldots,p_{max}$) be respectively the matrix 
and the vector which correspond to the computation of the energy functional for the patch with order number $p$.
Then the global energy functional can be written in the form
\begin{equation}
{\cal E}({\cal Z}_{all}) = \sum_{p=1}^{p_{max}} {\cal E}^{(p)} = 
{\cal Z}_{all}^T {\cal A} {\cal Z}_{all} - {\cal Z}_{all}^T {\cal B}
\end{equation}
where
\begin{equation}
{\cal A} = \left(\begin{array}{cccc}
{\cal A}^{(1)}    & 0                 & \ldots & 0      \cr
0                 & {\cal A}^{(2)}    & \ldots & 0      \cr
\vdots            & \vdots            & \ddots & \vdots \cr
0                 & 0                 & \ldots & {\cal A}^{(p_{max})}
\end{array}\right),\ \ \ \ \ \ \
{\cal B} = \left(\begin{array}{c}
{\cal B}^{(1)}\cr {\cal B}^{(2)}\cr \vdots \cr {\cal B}^{(p_{max})}
\end{array}\right)
\end{equation}
Formal differentiating of the energy functional and setting the result equal to zero yields
\begin{equation}
\D{\cal E}{{\cal Z}_{all}}=2{\cal A}{\cal Z}_{all}-{\cal B}=0
\label{eq:energy_derivative_all}
\end{equation}

Let $\MDS{n}$ be a minimal determining set which fits a chosen "additional" constraints. 
The set of $Z$-coordinates corresponding to in-plane control points from the MDS will be denoted by ${\cal Z}_{basis}$. From the algebraic point of view, ${\cal Z}_{basis}$ describes degrees of freedom of the constrained minimisation problem, where only $G^1$ continuity constraints are applied.
Dependencies of the remaining control points on ${\cal Z}_{basis}$ has a linear form and it is possible to define the dependency matrix ${\cal C}$ so that
\begin{equation}
{\cal Z}_{all} = {\cal C}{\cal Z}_{basis}
\label{eq:Z_all_Z_basis_connection} 
\end{equation}
It leads to the energy functional of the form
\begin{equation}
{\cal E}({\cal Z}_{basis}) = {\cal Z}_{basis}^T {\cal C}^T {\cal A} {\cal C} {\cal Z}_{basis} -
{\cal Z}_{basis}^T{\cal C}^T B
\label{eq:energy_basis}
\end{equation}
Differentiation of the functional gives
\begin{equation}
\D{\cal E}{{\cal Z}_{basis}}=2{\cal C}^T {\cal A} {\cal C} {\cal Z}_{basis}-{\cal C}^T {\cal B}=0
\label{eq:energy_derivative_basis}
\end{equation}
Note, that dependency of a dependent control point on the basic control points is usually defined "step by step". At every stage of the classification process, for a control point which gets dependent status it is sufficient to define how it explicitly depends on the basic control points and control points which are defined as dependent during the previous stages or steps of the classification. It clearly leads to the gradual construction of the final dependency matrix. 

${\cal Z}_{basis}$ is known to fit the considered "additional" constraints. Therefore, it only remains to separate ${\cal Z}_{basis}$ into two subsets ${\cal Z}_{fixed}$ and ${\cal Z}_{free}={\cal Z}_{basis}\setminus{\cal Z}_{fixed}$, where ${\cal Z}_{fixed}$ is $Z$-coordinates of the control points which should be fixed as a result of application of the "additional" constraints. Assuming that 
${\cal Z}_{basis}$ is ordered so that 
\begin{equation}
{\cal Z}_{basis} = \left(\begin{array}{c} {\cal Z}_{free}\cr {\cal Z}_{fixed}\end{array}\right)
\ \ \ {\rm and} \ \ \
{\cal C} = ( {\cal C}_{free}\ \ {\cal C}_{fixed})
\end{equation}
one gets the final linear system of equations with unknowns ${\cal Z}_{free}$
\begin{equation}
\D{\cal E}{{\cal Z}_{free}}=
2{\cal C}^T_{free} {\cal A} {\cal C}_{free} {\cal Z}_{free}+
2{\cal C}^T_{free} {\cal A} {\cal C}_{fixed} {\cal Z}_{fixed}-
{\cal C}^T_{free} {\cal B} = 0
\label{eq:energy_derivative_free}
\end{equation}
Here  ${\cal A}$ is $n_{all}\times n_{all}$ symmetric square matrix,
${\cal B}$ is $n_{all}\times 1$ vector, ${\cal C}_{free}$ is  $n_{all}\times n_{free}$ matrix
and ${\cal C}_{fixed}$ is  $n_{all}\times n_{fixed}$ matrix,  
where $n_{all}$, $n_{free}$, $n_{fixed}$ are numbers of variables in ${\cal Z}_{all}$, ${\cal Z}_{free}$, ${\cal Z}_{fixed}$ respectively. 

\eject
\section{Construction of an interpolating surface:  current approach and interpolation  based on  $3D$ mesh of curves} 
\label{sect:compare_interpolation}

As already mentioned above, from the pure theoretical point of view results of the current work fit the general theory presented in Part ~\ref{part:review_3D_G1} and have a nice and simple geometrical interpretation. 

The current Section concentrates on the more practical aspects of construction of the resulting surface for the (vertex)(tangent plane)-interpolation problem. The (vertex)(tangent plane)-interpolation problem is chosen since it is the "native" problem for the approach based on interpolation of $3D$ mesh of curves and one of the possible applications of our method. The Section highlights the similarities and the differences between the  implementation of the current approach and the interpolation method described in work ~\cite{peters_main}.

\paragraph{Algorithm for construction of a smooth interpolant presented in work ~\cite{peters_main}}
\label{subsect:construction_of_smooth_interpolant}
In the case of a cubic mesh of curves, work ~\cite{peters_main} provides an algorithm for the construction of $G^1$-smooth quadratic interpolant. The algorithm requires that the mesh of curves satisfies the following requirements.
\begin{itemize}
\item[-]
The mesh curves define a unique tangent plane at every mesh vertex.
\item[-]
Sufficient vertex enclosure constraint (see ~\ref{th:suff_vert_enclosure_P}) holds at every even mesh vertex.
\item[-]
At every mesh vertex, tangents of any two sequential curves emanating from the vertex span an angle of less than $\pi$. \end{itemize}
For a quadrilateral mesh, the algorithm defines  $(3,1,1)$-match for every pair of adjacent patches in the general case and $(2,0,0)$-match if the {\it "Tangents Relation"} 
\vspace{-0.05in} 
\begin{equation}
\frac
{\<\bar \epsilon^{(L)},\bar \epsilon^{(R)}\>}
{\<\bar \epsilon^{(R)},\bar \epsilon^{(C)}\>+
\<\bar \epsilon^{(C)},\bar \epsilon^{(L)}\>}=
\frac
{\<\bar \epsilon'^{(L)},\bar \epsilon'^{(R)}\>}
{\<\bar \epsilon'^{(R)},\bar \epsilon'^{(C)}\>+
\<\bar \epsilon'^{(C)},\bar \epsilon'^{(L)}\>}
\label{eq:tangents_relation}
\end{equation}
\vspace{-0.05in}
holds (see Figure ~\ref{fig:fig37}). Here the tangent vectors at two endpoints of the curve are considered in coordinates of the tangent plane they belong to.

Generally, the coefficients of the weight functions are fixed in accordance with the given mesh data. The B\'ezier control points are built locally and linearly, proceeding from the boundary control points to the interior of a patch. Any under-constraint situation is solved by some local heuristic, like averaging or the least-square technique. In some situations,  one additional heuristic constraint has to be imposed : the choice of the twist control points adjacent to an even inner mesh vertex; the choice of two middle "side" control points adjacent to a mesh curve which satisfies the {\it "Tangents Relation"}; the choice of the inner middle control point of a quartic B\'ezier patch.

\paragraph{Comparison with the current approach} 

Although both techniques solve the (vertex)(tangent plane)-interpolation problem by construction of a $G^1$-smooth piecewise B\'ezier surface, they pursue different goals and are {\it applicable in different situations}. The approach presented in work ~\cite{peters_main} serves mainly for the construction of nicely looking smooth surfaces, it may be used as a tool of visualization. The current approach finds the constrained solution of some minimisation problem defined by a given energy functional. There are some similarities and many differences in  these two techniques .

{\it Local/global nature of the solution. Application of the local heuristics.}
Interpolation of $3D$ mesh of curves with a smooth surface uses  local techniques. Control points of the resulting B\'ezier patches are constructed separately for every face of the mesh and depend on the geometry of the mesh. The solution of all the under-constraint situations and the choice of the inner control points is based on some local heuristics. The current approach uses no heuristics, but requires solution of some global system. Construction of a MDS allows to decrease significantly the number of variables and to solve the global system in terms of the basic control points only.

{\it Under-constraint situations and additional degrees of freedom}
There are definite similarities between under-constraint situations of the mesh interpolation algorithm and some special cases studied in the current work. For example, an additional basic twist control point for an inner even vertex (in the current approach, see Subsection ~\ref{subsect:local_DT_inner_vertex_bilinear})  reflects the application of a local heuristic for the choice of the twist control points adjacent to an even inner vertex (in the $3D$ mesh of curves interpolation approach). An additional middle basic control point for an edge, which satisfies the {\it "Projections Relation"} (see Theorem ~\ref{theorem:middle_equations_bilinear}) is the application  of a local heuristic for the choice of the "side" middle control points adjacent to a mesh curve which satisfies the {\it "Tangents Relation"}.

{\it Requirements on the initial data.}
The approach based on interpolation of $3D$ mesh of curves works only if the mesh is admissible; the satisfaction of the vertex enclosure constraint, which involves the first and second-order derivatives, should be verified at every vertex of the mesh. The current approach does not define any requirements on the initial data (excepting minor natural mesh limitations). Moreover curves which participate in $3D$ mesh of curves are fully defined; in particular it means that the first and second-order derivatives of the resulting patches in the boundary directions are initially defined. In the current approach, a boundary curve of a resulting patch has at least one degree of freedom (at least one inner control point of the curve is not fixed by the initial data), which allows to always satisfy the vertex enclosure constraint.

{\it Degrees of the resulting patches.}
According to the current approach, there exist an instance of $\MDS{4}(\T\Pi^{(bilinear)})$ (in the case of a mesh with a polygonal global boundary) and an instance of $\MDS{5}(\T\Pi^{(bicubic)})$ or mixed MDS\ $\MDS{4,5}(\T\Pi^{(bicubic)})$ (in case of a mesh with a smooth global boundary), which fit the (vertex)(tangent-plane)-interpolation condition. 
The approach for interpolation of $3D$ mesh of curves allows to construct the solution of degree $4$ for any admissible mesh; a case of a mesh with a given smooth global boundary is not an exception.

The main reason for the difference is that the approaches use different techniques for definition of the weight functions. As soon as a $3D$ mesh of curves appears to be admissible, the approach of curve mesh interpolation  defines coefficients of the weight functions which participate in the vertex enclosure constraints separately at every mesh vertex and then proceeds with definition of the weight functions independently for every curve. In the current approach, the weight functions of the different edges of the same mesh element are interconnected by an in-plane parametrisation of the element. In the case of a planar mesh with a piecewise-cubic boundary, the bicubic in-plane parametrisation of the boundary mesh elements leads to the high degrees of the weight functions for the edges with one boundary vertex. Of course, for every such edge it might be possible to construct the weight functions separately, precisely in the same manner as in the case of global bilinear in-plane parametrisation. For example, in Figure ~\ref{fig:fig35}, the weight functions of edge $(\GG,\GG')$ are defined by the bilinear parametrisations of virtual mesh elements constructed according to the middle control points of the boundary curve. Although these weight functions and correspondent $G^1$-continuity equations could be analysed similarly to the analysis presented in Section ~\ref{sect:vertex_enclosure}, they lead to a disagreement of in-plane control points corresponding to the weight functions of three inner edges of the boundary mesh element. It implies that in-plane parametrisation is no longer correctly defined for boundary mesh elements, while the definition of in-plane parametrisation a priory plays the principal role for the current approach.

\eject

\begin{figure}[!pt]
\centering
\includegraphics[clip,totalheight=1.8in]{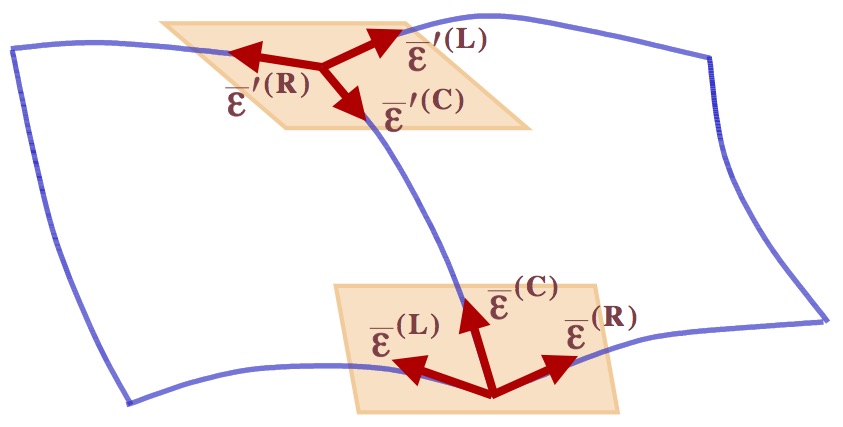}
\vspace{-0.1in}
\caption{An illustration for the {\it "Tangents Relation"}.}
\label{fig:fig37}
\end{figure}

\begin{figure}[!ph]
\centering
\includegraphics[clip,totalheight=2.3in]{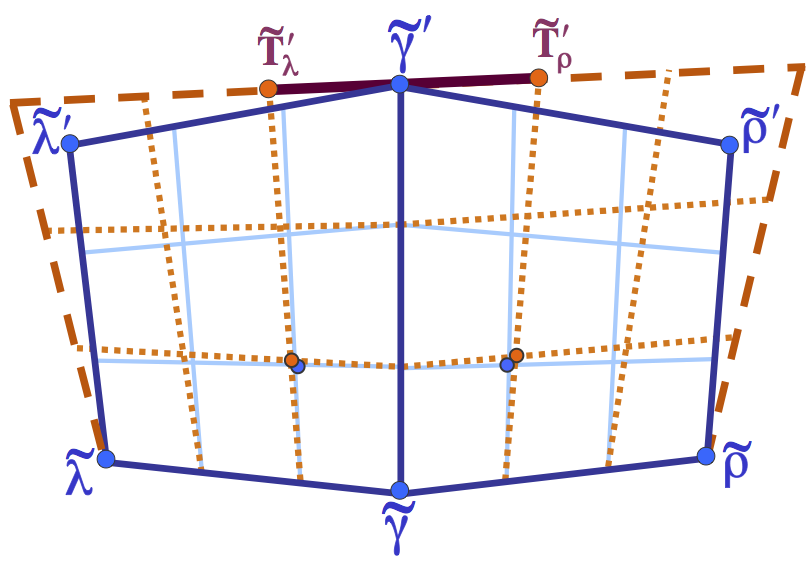}
\vspace{-0.15in}
\caption{Different weight functions correspond to different in-plane parametrisations of the boundary mesh elements.}
\label{fig:fig35}
\end{figure}

\FloatBarrier

\end{document}